\documentclass[11pt]{amsart}
\usepackage[margin=1.1in]{geometry}
\usepackage{amsmath,amssymb,amsthm,mathtools}
\usepackage{enumitem}
\usepackage{booktabs}
\usepackage[dvipsnames]{xcolor}
\usepackage{dsfont}
\usepackage[colorlinks=true,linkcolor=blue!90!black,citecolor=red!90!black,urlcolor=blue!90!black]{hyperref}

\numberwithin{equation}{section}
\allowdisplaybreaks
\theoremstyle{plain}
\newtheorem{theorem}{Theorem}[section]
\newtheorem{proposition}[theorem]{Proposition}
\newtheorem{lemma}[theorem]{Lemma}
\newtheorem{corollary}[theorem]{Corollary}
\newtheorem{mainthm}{Theorem}

\newtheorem{maincor}[mainthm]{Corollary}
\theoremstyle{definition}
\newtheorem{definition}[theorem]{Definition}
\newtheorem{example}[theorem]{Example}
\newtheorem{notation}[theorem]{Notation}
\theoremstyle{remark}
\newtheorem{remark}[theorem]{Remark}

\newcommand{\R}{\mathbb{R}}
\newcommand{\C}{\mathbb{C}}
\newcommand{\Z}{\mathbb{Z}}

\newcommand{\one}{\mathds{1}}
\newcommand{\T}{\mathcal{T}}
\newcommand{\cO}{\mathcal{O}}
\DeclareMathOperator{\GL}{GL}
\DeclareMathOperator{\Hom}{Hom}

\DeclareMathOperator{\Irr}{Irr}
\DeclareMathOperator{\Rep}{Rep}
\DeclareMathOperator{\Nrd}{Nrd}
\DeclareMathOperator{\Trd}{Trd}
\DeclareMathOperator{\St}{St}
\DeclareMathOperator{\Sp}{Sp}
\DeclareMathOperator{\Ind}{Ind}
\DeclareMathOperator{\ind}{ind}
\DeclareMathOperator{\Gr}{Gr}
\DeclareMathOperator{\Lie}{Lie}

\DeclareMathOperator{\JL}{JL}
\DeclareMathOperator{\vol}{vol}
\DeclareMathOperator{\diag}{diag}
\DeclareMathOperator{\Ad}{Ad}
\newcommand{\Qt}{\widetilde{Q}}
\newcommand{\Stil}{\widetilde{S}}

\newcommand{\xione}{\xi^{(1)}}
\newcommand{\xitwo}{\xi^{(2)}}
\newcommand{\mm}{\mathfrak{m}}
\newcommand{\bb}{\mathfrak{b}}

\newcommand{\sss}{\mathfrak{s}}

\newcommand{\aaa}{\mathfrak{a}}
\newcommand{\Ss}{\mathrm{S}}
\newcommand{\Qq}{\mathrm{Q}}
\newcommand{\dI}{d_{\mathrm{I}}}
\newcommand{\dII}{d_{\mathrm{II}}}
\newcommand{\dIII}{d_{\mathrm{III}}}
\newcommand{\OI}{\mathcal{O}_{\mathrm{I}}}
\newcommand{\OII}{\mathcal{O}_{\mathrm{II}}}
\newcommand{\OIII}{\mathcal{O}_{\mathrm{III}}}

\newcommand{\WI}{W_{\mathrm{I}}}
\newcommand{\WII}{W_{\mathrm{II}}}
\newcommand{\WIII}{W_{\mathrm{III}}}
\newcommand{\Kir}{\mathcal{K}}

\begin{document}
\title[On Generalized linear periods for inner forms of $\GL_n$]{On Generalized linear periods for inner forms of $\GL_n$}
\author{Hariom Sharma}
\address{Department of Mathematics, Indian Institute of Technology Bombay, Powai, Mumbai 400076, India}
\email{hariomshrma97@gmail.com, hariom@math.iitb.ac.in}

\begin{abstract}
Let $\mathrm{D}$ be a central division algebra of index $e$ over a non-Archimedean local field $\mathrm{F}$ of characteristic zero. For every integer $n\ge 1$ and every $s\in\R$, we classify the irreducible smooth representations of $\GL_n(\mathrm{D})$ that admit a nonzero linear form transforming under the subgroup $\GL_1(\mathrm{D})\times\GL_{n-1}(\mathrm{D})$ by the character $\nu^{es}\otimes\nu^{-es}$. We further show that such a linear form is unique up to scalars, and we characterize its existence in terms of Langlands parameters.
\end{abstract}

\maketitle
\setcounter{tocdepth}{1}
\tableofcontents

\section{Introduction}\label{sec:intro}

\subsection{Periods and distinction}
Let $G$ be an $\ell$-group in the sense of Bernstein and Zelevinsky \cite[\S 1.1]{BZ76}, that is, a Hausdorff topological group in which the identity admits a basis of neighbourhoods consisting of compact open subgroups. Let $H$ be a closed subgroup of $G$, and let $\chi\colon H\to\C^\times$ be a smooth character. For a smooth representation $(\pi,V)$ of $G$, put
$$\Hom_H(\pi,\chi)=\{\ell\in\Hom_\C(V,\C):\ \ell(\pi(h)v)=\chi(h)\,\ell(v)\ \text{for all }h\in H,\ v\in V\}.
$$
A nonzero element of $\Hom_H(\pi,\chi)$ realizes the one-dimensional representation $\chi$ as a quotient of the restriction $\pi|_H$. We say that $\pi$ is \emph{$(H,\chi)$-distinguished} if $\Hom_H(\pi,\chi)\ne0$, and simply \emph{$H$-distinguished} when $\chi$ is trivial.

The distinction problem asks for a description of the irreducible smooth representations of $G$ that are $(H,\chi)$-distinguished. It comprises the following questions:
\begin{enumerate}[label=\textup{(\roman*)}]
\item which irreducible smooth representations $\pi$ of $G$ satisfy $\Hom_H(\pi,\chi)\ne0$;
\item what the dimension of $\Hom_H(\pi,\chi)$ is;
\item how $(H,\chi)$-distinction is reflected in arithmetic invariants attached to $\pi$, such as its Langlands parameter and its local $L$-functions.
\end{enumerate}
When $H\backslash G$ is a spherical variety, these questions form part of the local side of the relative Langlands program. Their global counterparts concern period integrals of automorphic forms over $H$, whose nonvanishing is often governed by special values or poles of automorphic $L$-functions \cite{FJ93}.

\subsection{The problem}

Let $\mathrm{F}$ be a non-Archimedean local field of characteristic zero, and let $\mathrm{D}$ be a central division algebra of index $e$ over $\mathrm{F}$. For $m\ge1$, we put $G_m=\GL_m(\mathrm{D})$, which is an inner form of $\GL_{em}(\mathrm{F})$, and we let $G_0$ denote the trivial group. Let $\nu=|\cdot|_{\mathrm{F}}\circ\Nrd$, where $\Nrd$ is the reduced norm and $|\cdot|_{\mathrm{F}}$ is the normalized absolute value of $\mathrm{F}$. Consider the subgroup
$$H_{1,n-1}=\{\diag(g_1,g_2):g_1\in G_1,\ g_2\in G_{n-1}\}\simeq\GL_1(\mathrm{D})\times\GL_{n-1}(\mathrm{D})
$$
of $G_n$ and, for $s\in\R$, the character $\chi_s$ of $H_{1,n-1}$ defined by
$$\chi_s(\diag(g_1,g_2))=\nu(g_1)^{es}\,\nu(g_2)^{-es}.$$
The group $H_{1,n-1}$ is the fixed-point subgroup of the involution given by conjugation by $\diag(1,-1,\dots,-1)$. Thus, $(G_n,H_{1,n-1})$ is a symmetric pair, and it is an inner form of the symmetric pair $(\GL_{en}(\mathrm{F}),\GL_e(\mathrm{F})\times\GL_{e(n-1)}(\mathrm{F}))$. In this article, we give a complete solution of the $(H_{1,n-1},\chi_s)$-distinction problem for every index $e$, every integer $n\ge1$, and every $s\in\R$.

When $\mathrm{D}=\mathrm{F}$, the problem is classical. Prasad \cite[Theorem 2]{Pra93} classified the irreducible admissible representations of $\GL_3(\mathrm{F})$ that are $\GL_2(\mathrm{F})$-distinguished, and conjectured the answer for $\GL_n(\mathrm{F})$ in general. Venketasubramanian \cite[Corollary 6.15]{Ven13} completed the classification of $\GL_{n-1}(\mathrm{F})$-distinguished representations, and Smith \cite[Corollary~A.2.2]{Smi17} reformulated these results as a classification of the representations of $\GL_n(\mathrm{F})$ that are distinguished by $\GL_1(\mathrm{F})\times\GL_{n-1}(\mathrm{F})$ with respect to a character. Multiplicity one for twisted linear periods was proved by Chen and Sun \cite[Theorem B]{CS20} for $(\GL_{2n}(\mathrm{F}),\GL_n(\mathrm{F})\times\GL_n(\mathrm{F}))$ and all but countably many characters, and by Lu \cite[Theorem 1.1]{Lu22} for $(\GL_{2p+1}(\mathrm{F}),\GL_p(\mathrm{F})\times\GL_{p+1}(\mathrm{F}))$ and good characters; see also \cite{JR96,AGRS10}.

 For inner forms, Anandavardhanan et al.  \cite[Theorem 2.5]{ALMSY24} proved multiplicity one for linear periods when $|p-q|\le1$, and Lu \cite[Theorem 1.3]{Lu26} studied generalized linear periods over a division algebra. More recently, Li et al. \cite[Theorem 1.1]{LLX26} proved multiplicity one for all twisted linear periods when $e\le2$. Their method relies on the anti-involution $g\mapsto{}^t\bar g$ of $M_n(\mathrm{D})$, which exists only when $e\le2$. For $\mathrm{D}$ quaternion, Dagar, Verma and the author \cite[Conjecture 1.1]{DSV27} classified the $H_{1,n-1}$-distinguished representations of $G_n$ for $n\le4$ and conjectured the answer for all $n$. Dagar and the author \cite[Conjecture 1.1]{DS26} then formulated a conjectural classification of generalized linear periods for $\mathrm{D}$ quaternion. In the present article, we prove these conjectures and establish their analogues for every index $e$. We also prove the multiplicity one theorem for every $e$ without using an anti-involution. 
Problems of this type have recently been considered by Gan et al. in  \cite[Conjecture 1.1]{GGP23}. There they study a twisted variant of their restriction problems for classical groups, in which the trivial representation of the subgroup is replaced by a Weil representation, a member of a one-parameter family of small representations.

\subsection{The main results}

We first fix some notation. Throughout, $\times$ denotes normalized parabolic induction from standard parabolic subgroups; in a product of characters, each power of $\nu$ is a character of the corresponding block. Let $\T$ denote the set of infinite-dimensional irreducible representations of $G_2$ that are $H_{1,1}$-distinguished; this set is described explicitly in Theorem~\ref{thm:B} below. Let $Q_n$ be the unique irreducible quotient of the representation $\nu^{e/2}\times\nu^{e\frac{n+1}2}$ induced from $G_{n-1}\times G_1$, and let $\Qt_n$ be its contragredient; see \S\ref{ss:Qm}. For $\pi\in\Irr(G_m)$, we write $\varphi_\pi$ for its Langlands parameter, that is, the Langlands parameter of the representation $\JL(\pi)$ of $\GL_{em}(\mathrm{F})$ attached to $\pi$ by the Jacquet--Langlands transfer; see Definition~\ref{def:L-param}.

Our main theorem concerns the case $n\ge3$. The answer for $n=2$ takes a different form and is stated separately in Theorem~\ref{thm:B}, while the case $n=1$ is elementary; see  Proposition~\ref{prop:n1}. For $n\ge3$, the distinguished representations form two families of (possibly degenerate) principal series attached to maximal parabolic subgroups, together with one family induced from $\T$.

\begin{mainthm}[The case $n\ge3$]\label{thm:A}
	Let $n\ge3$ and $s\in\R$.
	\begin{enumerate}[label=\textup{(\roman*)}]
		\item \textup{(Classification.)} An irreducible smooth representation $\pi$ of $G_n$ is $(H_{1,n-1},\chi_s)$-distinguished if and only if $\pi$ is isomorphic to one of the following:
		\begin{enumerate}[label=\textup{(\alph*)}]
			\item the representation $\nu^{-es-\frac e2}\times\nu^{es+e\frac{n-1}2}$ induced from $G_{n-1}\times G_1$, if $s\notin\{-\frac n2,0\}$; and its unique irreducible subrepresentation if $s\in\{-\frac n2,0\}$, namely $\nu^{en/2}\Qt_n$ if $s=-\frac n2$ and $\one_n$ if $s=0$;
			\item the representation $\nu^{es-e\frac{n-1}2}\times\nu^{-es+\frac e2}$ induced from $G_1\times G_{n-1}$, if $s\notin\{0,\frac n2\}$; and its unique irreducible subrepresentation if $s\in\{0,\frac n2\}$, namely $\one_n$ if $s=0$ and $\nu^{-en/2}Q_n$ if $s=\frac n2$;
			\item the representation $\nu^{-es}\times\tau$ induced from $G_{n-2}\times G_2$, where $\tau\in\T$ and $\nu^{-es}\times\tau$ is irreducible.
		\end{enumerate}
		\item \textup{(Multiplicity one.)} For every $\pi\in\Irr(G_n)$, we have $\dim_\C\Hom_{H_{1,n-1}}(\pi,\chi_s)\le1$. Hence, this dimension equals $1$ for the representations in \textup{(i)} and $0$ otherwise.
		\item \textup{(Langlands parameters.)}  An irreducible smooth representation $\pi$ of $G_n$ is $(H_{1,n-1},\chi_s)$-distinguished if and only if
		$\varphi_{\pi}\simeq\varphi_{\nu^{-es}}\oplus\varphi_{\tau_0}$,
		where $\nu^{-es}$ is regarded as a character of $G_{n-2}$ and $\varphi_{\tau_0}$ is the $2e$-dimensional Langlands parameter of either the trivial representation $\tau_0=\one_2$, in which case $s=\pm\frac{n-2}2$, or of a representation $\tau_0\in\T$.
	\end{enumerate}
\end{mainthm}

The induced representations in (a) and (b) are irreducible exactly at the non-exceptional values of $s$, and they have length two at the exceptional values; see Lemmas~\ref{lem:dps-irred}--\ref{lem:dps2-exceptional}. The irreducibility condition in (c) is explicit. By Lemma~\ref{lem:reducible-nu-tau}, the representation $\nu^{-es}\times\tau$ is reducible exactly when either $\tau=\St_2$ and $s=\pm\frac n2$, or $\tau=\nu^{ex}\times\nu^{-ex}$ with $x\ne\pm\frac12$ and $\pm x\in\{-s\pm\frac{n-1}2\}$.

To state the result for $n=2$, we need twisted Jacquet modules. Let $N$ be the unipotent radical of the upper triangular parabolic subgroup of $G_2$, and define $\psi(\begin{smallmatrix}1&x\\0&1\end{smallmatrix})=\psi_{\mathrm{F}}(\Trd x)$, where $\psi_{\mathrm{F}}$ is a nontrivial character of $\mathrm{F}$. For $\tau\in\Irr(G_2)$, let $\tau_{N,\psi}$ denote the twisted Jacquet module of $\tau$; it is a representation of the stabilizer of $\psi$ in the Levi subgroup, which is $\mathrm{D}^\times$ embedded diagonally. For $\sigma\in\Irr(\mathrm{D}^\times)$, let $\Sp_2(\sigma)$ and $\St_2(\sigma)$ be the generalized Speh and Steinberg representations; see \S\ref{ss:segments}.

\begin{mainthm}[The case $n=2$]\label{thm:B}
	Let $s\in\R$.
	\begin{enumerate}[label=\textup{(\roman*)}]
		\item \textup{(Independence of $s$.)} For every infinite-dimensional $\tau\in\Irr(G_2)$, we have
		$$
		\dim\Hom_{H_{1,1}}(\tau,\chi_s)=\dim\Hom_{H_{1,1}}(\tau,\one)=\dim\Hom_{\mathrm{D}^\times}(\tau_{N,\psi},\one).$$
		In particular, a representation $\pi\in\Irr(G_2)$ is $(H_{1,1},\chi_s)$-distinguished if and only if either $\pi\simeq\one_2$ and $s=0$, or $\pi$ is infinite-dimensional and $\pi\in\T$, that is, $\Hom_{\mathrm{D}^\times}(\pi_{N,\psi},\one)\ne0$.
		
		\item \textup{(Classification by representation type.)} An infinite-dimensional representation $\tau\in\Irr(G_2)$ lies in $\T$ if and only if one of the following holds:
		\begin{enumerate}[label=\textup{(\alph*)}]
			\item \textup{(Principal series.)} $\tau=\sigma\times\widetilde\sigma$ is irreducible for some $\sigma\in\Irr(\mathrm{D}^\times)$;
			\item \textup{(Twisted Steinberg.)} $\tau=\chi\St_2$ for a character $\chi$ of $\mathrm{D}^\times$ with $\chi^2=\one$;
			\item \textup{(Generalized Speh and Steinberg.)} $\tau\in\{\Sp_2(\sigma),\St_2(\sigma)\}$ for a self-dual $\sigma\in\Irr(\mathrm{D}^\times)$ with $\dim\sigma>1$, and $\tau_{N,\psi}$ admits a nonzero $\mathrm{D}^\times$-invariant linear form;
			\item \textup{(Supercuspidal.)} $\tau$ is supercuspidal, and $\tau_{N,\psi}$ admits a nonzero $\mathrm{D}^\times$-invariant linear form.
		\end{enumerate}
		For each $\sigma$ as in \textup{(c)}, exactly one of $\Sp_2(\sigma)$ and $\St_2(\sigma)$ lies in $\T$.
		
		\item \textup{(Multiplicity one.)} For every $\pi\in\Irr(G_2)$, we have $\dim\Hom_{H_{1,1}}(\pi,\chi_s)\le1$.
		\end{enumerate}
\end{mainthm}

Theorem~B\textup{(i)} is related to a result of Raghuram \cite[Theorem~1.1]{Rag07}. For an infinite-dimensional $\pi\in\Irr(G_2)$ and irreducible representations $\tau_1,\tau_2$ of $\mathrm{D}^\times$, he proved that
$$\dim\Hom_{\mathrm{D}^\times\times\mathrm{D}^\times}(\pi,\tau_1\otimes\tau_2)=\dim\Hom_{\mathrm{D}^\times}(\pi_{N,\psi},\tau_1\otimes\tau_2)$$
when either $\tau_1\otimes\tau_2$ does not intertwine with $\pi_N$, or $\pi$ is irreducibly parabolically induced and $\pi_N$ is semisimple. He expected this identity to hold for all $\pi$ and all $\tau_1\otimes\tau_2$. Taking $\tau_1\otimes\tau_2=\chi_t$, whose restriction to the diagonal copy of $\mathrm{D}^\times$ is trivial, Theorem~B\textup{(i)} proves this identity for every infinite-dimensional $\pi$ and every $t\in\R$, without either assumption. This covers, for instance, $\pi=\St_2$ with $t=1$, where $(\St_2)_N=\chi_1$ but $\St_2$ is not parabolically induced.

Theorem~\ref{thm:B}\textup{(ii)} covers each type of infinite-dimensional $\tau\in\Irr(G_2)$: irreducible principal series, $\St_2(\sigma)$, $\Sp_2(\sigma)$ with $\dim\sigma>1$, and supercuspidals. In \textup{(a)} and \textup{(b)} the criterion of \textup{(i)} becomes explicit in $\sigma$. In \textup{(c)} and \textup{(d)},  it is still stated in terms of $\tau_{N,\psi}$.

\begin{remark}\label{rem:quaternion-dichotomy}
	Let $\mathrm{D}$ be quaternion, so that $e=2$,  and  let $\sigma\in\Irr(\mathrm{D}^\times)$ self-dual with $\dim\sigma>1$. Since $\Sp_2(\sigma)_{N,\psi}\simeq\omega_\sigma\circ\Nrd$ \cite[Theorem~1.1]{NS25}, we have  $\Sp_2(\sigma)\in\T$ if $\omega_\sigma\circ\Nrd=1$ and $\St_2(\sigma)\in\T$ otherwise. See \cite{Rag07,GT10,Pra01} for related discussion of the representations in $\T$.
\end{remark}

\begin{remark}\label{rem:tempered}
	The classification shows that tempered distinguished representations occur only in small rank. For $n=2$, by Theorem~\ref{thm:B}, for every $s\in\R$ the tempered $(H_{1,1},\chi_s)$-distinguished representations are the tempered members of $\T$, namely the irreducible principal series $\sigma\times\widetilde\sigma$ with $\sigma$ unitary, the representations $\chi\St_2$, the representations $\St_2(\sigma)$ lying in $\T$, and the supercuspidal representations in $\T$. For $n=3$, by Theorem~\ref{thm:A}, a tempered representation is $(H_{1,2},\chi_s)$-distinguished if and only if $s=0$ and it is of the form $\one_1\times\tau$ with $\tau\in\T$ tempered. For $n\ge4$, no tempered representation of $G_n$ is $(H_{1,n-1},\chi_s)$-distinguished for any $s\in\R$, since each representation in Theorem~\ref{thm:A} involves a character of $G_m$ with $m\ge2$.
\end{remark}

Specializing Theorem~\ref{thm:A} to $s=0$, we obtain the following untwisted statement. For $\mathrm{D}$ quaternion, it was conjectured in~\cite[Conjecture~1.1]{DSV27}. It has recently been proved by Chen and Huang~\cite[Theorem~A]{CH26}, by a method different from ours.

\begin{maincor}\label{cor:C}
	Let $n\ge3$. An irreducible representation $\pi$ of $G_n$ is $H_{1,n-1}$-distinguished if and only if either $\pi\simeq\one_n$, or $\pi\simeq\one_{n-2}\times\tau$ with $\tau\in\T$ and $\one_{n-2}\times\tau$ irreducible. In that case, $\dim\Hom_{H_{1,n-1}}(\pi,\one)=1$.
\end{maincor}

The family of characters $\chi_s$ may look special. The following remark shows that it covers all characters of the form $\nu^x\otimes\nu^y$.

\begin{remark}[Normalization]\label{rem:normalization}
	The normalization of $\chi_s$ entails no loss of generality. For $x,y\in\R$, let $\nu^x\otimes\nu^y$ denote the character $\diag(g_1,g_2)\mapsto\nu(g_1)^x\nu(g_2)^y$ of $H_{1,n-1}$. For every smooth representation $\pi$ of $G_n$ and every $t\in\R$, we have $$\Hom_{H_{1,n-1}}(\nu^t\pi,\chi_s)=\Hom_{H_{1,n-1}}(\pi,\nu^{es-t}\otimes\nu^{-es-t}).$$ Taking $s=\frac{x-y}{2e}$ and $t=-\frac{x+y}2$, we see that $\pi$ is $(H_{1,n-1},\nu^x\otimes\nu^y)$-distinguished if and only if $\nu^{-\frac{x+y}2}\pi$ is $(H_{1,n-1},\chi_{\frac{x-y}{2e}})$-distinguished. Hence, Theorems~\ref{thm:A} and~\ref{thm:B} determine distinction and multiplicity for all characters $\nu^x\otimes\nu^y$. 
\end{remark}

For $e=1$, Theorem~\ref{thm:A} recovers the results of Prasad \cite{Pra93}, Venketasubramanian \cite{Ven13} and Smith \cite[Corollary~A.2.2]{Smi17} for parts (i) and (iii), and the result of Lu \cite[Theorem 1.3]{Lu22} for part (ii). For $e=2$, parts (i) and (iii) prove the conjecture of \cite[Conjecture 1.1]{DS26}, and part (ii) also follows from \cite[Theorem 1.1]{LLX26}.

\subsection{Method}

The proof starts from an orbit analysis. For a representation $\pi_1\times\pi_2$ induced from the maximal parabolic subgroup $P=P_{(k,n-k)}$, the group $H_{1,n-1}$ has exactly three orbits on $\Gr_k(\mathrm{D}^n)$: two closed orbits and one open orbit; see Section~\ref{sec:orbits}. We work in Offen's framework \cite{Off17} for the symmetric space $G_n/H_{1,n-1}$. For each orbit we determine an explicit representative $x$, its Weyl element and Levi subgroup $L$, the standard parabolic subgroup $Q\subseteq P$ with Levi subgroup $L$, the stabilizers $P_x=Q_x$, and the characters $\delta_{Q_x}\delta_Q^{-1/2}$ in the geometric lemma \cite[Proposition~4.1]{Off17}. This gives the two-sided estimate
$$\dI+\dII\le\dim\Hom_{H_{1,n-1}}(\pi_1\times\pi_2,\chi_s)\le\dI+\dII+\dIII$$
of Theorem~\ref{thm:orbit-bound}. Here $\dI$ and $\dII$ are distinction multiplicities in smaller rank, and $\dIII$ is computed from Jacquet modules. When $e\ge2$, the open orbit contributes through a diagonally embedded copy of $\mathrm{D}^\times$. Its contribution therefore has to be expressed through arbitrary irreducible representations of $\mathrm{D}^\times$, which need not be characters.

These bounds drive an induction on $n$ that starts in rank two and three. In rank two (Section~\ref{sec:rank-two}), the Kirillov model of $\GL_2(\mathrm{D})$ and an extension argument show two things for every $e$: $(H_{1,1},\chi_t)$-distinction of infinite-dimensional representations does not depend on $t$, and $\T$ can be described by twisted Jacquet modules. The case $n=3$ is proved directly for every $e$. For $n\ge4$ (Sections~\ref{sec:n3}--\ref{sec:induction}), every distinguished representation is a quotient of some $\pi'\times\sigma$ with $\pi'\in\Irr(G_{n-1})$ and $\sigma\in\Irr(\mathrm{D}^\times)$. The three orbits then give three cases, and one of them involves a generalized linear period of $\pi'$ for the parameter $s-\frac12$. We analyze the constituents with the Zelevinsky and Langlands classifications of Tadi\'c and Lapid--M\'inguez. The case $n=4$ needs one further argument (Lemma~\ref{lem:n4-special}).

Multiplicity one is proved separately for each of the six kinds of distinguished representations. Each is realized as a quotient of an induced representation whose three orbit contributions add up to one. The characterization by Langlands parameters then follows from the classification (Sections~\ref{sec:multiplicity}--\ref{sec:langlands}). None of these arguments uses an anti-involution of $G_n$ or a duality between $(\pi,s)$ and $(\widetilde\pi,-s)$, since such tools are available only for $e\le2$.

\subsection{Organization}
Section~\ref{sec:prelim} collects notation and known results. Section~\ref{sec:elementary} contains elementary facts on generalized linear periods and the degenerate principal series. Section~\ref{sec:orbits} contains the orbit analysis. Section~\ref{sec:rank-two} treats $n=1$ and $n=2$. Section~\ref{sec:n3} contains general reductions and the case $n=3$. Sections~\ref{sec:aux} and~\ref{sec:induction} prove the classification for $n\ge4$. Section~\ref{sec:multiplicity} proves multiplicity one, and Section~\ref{sec:langlands} the characterization by Langlands parameters. 

\subsection{Acknowledgements} The author would like to thank Dipendra Prasad for many interesting conversations, insightful comments, and constant encouragement throughout the preparation of this article. The author also thanks U. K. Anandavardhanan and Ravi Raghunathan for their helpful suggestions. The author gratefully acknowledges the support of the Indian Institute of Technology Bombay through the Institute Postdoctoral Fellowship.
The author used Claude Opus 5 (Anthropic) and GPT-6.2 Astra (OpenAI) solely for language editing and improving the clarity of the manuscript.

\section{Notation and preliminaries}\label{sec:prelim}

In this section, we fix the notation and collect the facts about the representation theory of $\GL_m(\mathrm{D})$ that will be used later. For $\GL_m(\mathrm{F})$, the standard references are \cite{BZ76,BZ77,Zel80}; for the inner forms, we follow Tadi\'c \cite{Tad90} and Lapid--M\'inguez \cite{LM16}.

\subsection{The groups}\label{ss:groups}

Let $\mathrm{F}$ be a non-Archimedean local field of characteristic zero, with normalized absolute value $|\cdot|_{\mathrm{F}}$. Let $\mathrm{D}$ be a central division algebra of index $e$ over $\mathrm{F}$, so that $\dim_{\mathrm{F}}\mathrm{D}=e^2$, and let $\Nrd$ and $\Trd$ denote the reduced norm and the reduced trace, respectively. For $m\ge 1$, we put
$$
G_m=\GL_m(\mathrm{D}),
$$
and we let $G_0$ be the trivial group. The group $G_m$ is an inner form of $\GL_{em}(\mathrm{F})$, and $G_1=\mathrm{D}^\times$. The reduced norm $\Nrd\colon G_m\to \mathrm{F}^\times$ is a surjective homomorphism whose kernel is the commutator subgroup of $G_m$; this is due to Nakayama--Matsushima for $m=1$ and to Dieudonn\'e for $m\ge2$ (see \cite{Dra83,PR94}). Consequently, every character of $G_m$ is of the form $\mu\circ\Nrd$ for a unique character $\mu$ of $\mathrm{F}^\times$, and we also write $\mu$ for $\mu\circ\Nrd$. We put
$$
\nu(g)=|\Nrd(g)|_{\mathrm{F}}\qquad (g\in G_m),
$$
and, for $x\in\R$, we write $\nu^x$ for the character $g\mapsto \nu(g)^x$ of $G_m$; the integer $m$ will always be clear from the context. The trivial character of $G_m$ is denoted by $\one_m$. The centre of $G_m$ is $\{zI_m:z\in \mathrm{F}^\times\}$, and since $\Nrd(zI_m)=z^{em}$, we have
\begin{equation}\label{eq:nu-centre}
	\nu(zI_m)=|z|_{\mathrm{F}}^{em}\qquad(z\in \mathrm{F}^\times).
\end{equation}

Fix an integer $n\ge1$. We regard $V=\mathrm{D}^n$ as a right $\mathrm{D}$-vector space of column vectors with standard basis $e_1,\dots,e_n$, on which $G_n$ acts by left multiplication.

The following lemma computes the absolute value of the determinant of the natural action of a block-diagonal matrix on an off-diagonal block. It will be used to compute modulus characters of parabolic subgroups.

\begin{lemma}\label{lem:det}
	Let $a,b\ge1$, $A\in G_a$ and $B\in G_b$. The $\mathrm{F}$-linear map $Z\mapsto AZB^{-1}$ of the $\mathrm{F}$-vector space $M_{a\times b}(\mathrm{D})$ has determinant of absolute value $\nu(A)^{eb}\nu(B)^{-ea}$.
\end{lemma}

\begin{proof}
	Let $E/\mathrm{F}$ be a splitting field of $\mathrm{D}$, and fix an isomorphism $M_a(\mathrm{D})\otimes_{\mathrm{F}}E\simeq M_{ea}(E)$. After extension of scalars, the left $M_a(\mathrm{D})$-module $\mathrm{D}^a$ becomes an $M_{ea}(E)$-module of $E$-dimension $e^2a$, hence a direct sum of $e$ copies of the standard module $E^{ea}$. Therefore, left multiplication by $A$ on $\mathrm{D}^a$ has determinant $\Nrd(A)^e$. Since $M_{a\times b}(\mathrm{D})$ is the direct sum of its $b$ columns, each of which is a copy of $\mathrm{D}^a$, left multiplication by $A$ on $M_{a\times b}(\mathrm{D})$ has determinant $\Nrd(A)^{eb}$. Similarly, right multiplication by $B^{-1}$ on the space ${}^b\mathrm{D}$ of row vectors has determinant $\Nrd(B)^{-e}$. Since $M_{a\times b}(\mathrm{D})$ is the direct sum of its $a$ rows, right multiplication by $B^{-1}$ on $M_{a\times b}(\mathrm{D})$ has determinant $\Nrd(B)^{-ea}$. The map $Z\mapsto AZB^{-1}$ is the composition of these two commuting maps, so its determinant is $\Nrd(A)^{eb}\Nrd(B)^{-ea}$. Taking absolute values proves the lemma.
\end{proof}

\subsection{The subgroups \texorpdfstring{$H_{p,q}$}{Hpq} and the characters \texorpdfstring{$\chi_s$}{chi s}}\label{ss:Hpq}

Let $p,q\ge0$ with $p+q=n$, and put
$$
H_{p,q}=\{\diag(g_1,g_2):g_1\in G_p,\ g_2\in G_q\}\simeq G_p\times G_q .
$$
The subgroup $H_{p,q}$ is the group of fixed points of the involution $\Ad(\diag(I_p,-I_q))$ of $G_n$. Equivalently, it is the subgroup of $G_n$ that stabilizes each of the two subspaces $e_1\mathrm{D}\oplus\dots\oplus e_p\mathrm{D}$ and $e_{p+1}\mathrm{D}\oplus\dots\oplus e_n\mathrm{D}$ of $V$. For $s\in\R$, we define the character $\chi_s$ of $H_{p,q}$ by
$$
\chi_s(\diag(g_1,g_2))=\nu(g_1)^{es}\,\nu(g_2)^{-es}.
$$
When $p=0$ or $q=0$, the corresponding factor is absent. Thus, $\chi_s=\nu^{es}$ on $H_{1,0}=G_1$, and $\chi_s=\nu^{-es}$ on $H_{0,m}=G_m$. When we want to specify the pair $(p,q)$, we write $\chi_s^{(p,q)}$.

We now introduce the central notion of this article.

\begin{definition}
	Let $\pi$ be a smooth representation of $G_n$. We say that $\pi$ is \emph{$(H_{p,q},\chi_s)$-distinguished}, or that $\pi$ admits a \emph{generalized linear period} with respect to $(H_{p,q},\chi_s)$, if the space
	$$
	\Hom_{H_{p,q}}(\pi,\chi_s)=\{\ell\colon\pi\to\C\ \text{linear}:\ \ell(\pi(h)v)=\chi_s(h)\ell(v)\ \text{for all } h\in H_{p,q},\,v\in\pi\}
	$$
	is nonzero. When $s=0$, we simply say that $\pi$ is $H_{p,q}$-distinguished.
\end{definition}

Finally, let
$$
w_{p,q}=\begin{pmatrix}0&I_q\\ I_p&0\end{pmatrix}\in G_n,
$$
so that $w_{p,q}\diag(g_1,g_2)w_{p,q}^{-1}=\diag(g_2,g_1)$ for $g_1\in G_p$ and $g_2\in G_q$. In particular, $w_{p,q}H_{p,q}w_{p,q}^{-1}=H_{q,p}$, and $\chi^{(q,p)}_{-s}(w_{p,q}hw_{p,q}^{-1})=\chi^{(p,q)}_s(h)$ for $h\in H_{p,q}$.

\subsection{Representations, parabolic induction and Jacquet modules}\label{ss:reps}

All representations considered in this article are smooth and complex. We denote by $\Rep(G_m)$ the category of smooth representations of $G_m$ of finite length, by $\Irr(G_m)$ the set of isomorphism classes of irreducible smooth representations of $G_m$, and by $\mathcal{C}(G_m)\subseteq\Irr(G_m)$ the subset of supercuspidal ones. Every $\pi\in\Irr(G_m)$ is admissible and has a central character, which we denote by $\omega_\pi$. We write $\widetilde\pi$ for the smooth contragredient of $\pi$. The irreducible subquotients of $\pi\in\Rep(G_m)$, counted with multiplicity, are called the \emph{constituents} of $\pi$. Since $\mathrm{D}^\times/\mathrm{F}^\times$ is compact, every irreducible representation of $G_1=\mathrm{D}^\times$ is finite-dimensional and supercuspidal.

For a composition $\alpha=(n_1,\dots,n_r)$ of $n$, let $P_\alpha=M_\alpha N_\alpha$ be the corresponding standard (block upper triangular) parabolic subgroup of $G_n$, with Levi subgroup $M_\alpha=\{\diag(g_1,\dots,g_r):g_i\in G_{n_i}\}$ and unipotent radical $N_\alpha$. Note that $H_{p,q}=M_{(p,q)}$. The modulus character of $P_\alpha$ is $\delta_{P_\alpha}(x)=|\det(\Ad(x)|_{\Lie N_\alpha})|_{\mathrm{F}}$, where $\Lie N_\alpha$ is regarded as an $\mathrm{F}$-vector space. Using Lemma~\ref{lem:det}, we can compute it explicitly.

\begin{lemma}\label{lem:modulus}
	For $m=\diag(g_1,\dots,g_r)\in M_\alpha$ and $u\in N_\alpha$, we have
	$$
	\delta_{P_\alpha}(mu)=\prod_{i=1}^r\nu(g_i)^{e\left(\sum_{j>i}n_j-\sum_{j<i}n_j\right)}.
	$$
	In particular, $\delta_{P_{(k,n-k)}}(\diag(g,h))=\nu(g)^{e(n-k)}\nu(h)^{-ek}$ for $g\in G_k$ and $h\in G_{n-k}$.
\end{lemma}

\begin{proof}
	Since $N_\alpha$ acts unipotently on $\Lie N_\alpha$, we have $\delta_{P_\alpha}(u)=1$, and we may therefore assume that $u=1$. As an $\mathrm{F}$-vector space, $\Lie N_\alpha=\bigoplus_{i<j}M_{n_i\times n_j}(\mathrm{D})$, and $m$ acts on the $(i,j)$-summand by $Z\mapsto g_iZg_j^{-1}$. By Lemma~\ref{lem:det}, the determinant of this action has absolute value $\nu(g_i)^{en_j}\nu(g_j)^{-en_i}$. Taking the product over all pairs $i<j$, the factor $\nu(g_i)$ appears with exponent $en_j$ for each $j>i$ and with exponent $-en_j$ for each $j<i$. This gives the formula.
\end{proof}

For $\rho_i\in\Rep(G_{n_i})$, we write
$$
\rho_1\times\dots\times\rho_r=\Ind_{P_\alpha}^{G_n}(\rho_1\otimes\dots\otimes\rho_r)
$$
for the normalized parabolic induction. This is the space of locally constant functions $f$ on $G_n$ with values in the space of $\rho_1\otimes\dots\otimes\rho_r$ such that
$$
f(mug)=\delta_{P_\alpha}(m)^{1/2}(\rho_1\otimes\dots\otimes\rho_r)(m)f(g)\qquad(m\in M_\alpha,\ u\in N_\alpha,\ g\in G_n),
$$
on which $G_n$ acts by right translation. For $\pi\in\Rep(G_n)$, the normalized Jacquet module is $r_\alpha(\pi)=\delta_{P_\alpha}^{-1/2}\otimes\pi_{N_\alpha}$, where $\pi_{N_\alpha}$ denotes the space of $N_\alpha$-coinvariants. We use the following standard facts without further comment \cite{BZ77,Tad90}:
\begin{itemize}
	\item the induction and Jacquet functors are exact and preserve finite length;
	\item induction is transitive;
	\item Frobenius reciprocity holds:
	$$
	\Hom_{G_n}(\pi,\rho_1\times\dots\times\rho_r)\simeq\Hom_{M_\alpha}(r_\alpha(\pi),\rho_1\otimes\dots\otimes\rho_r);
	$$
	\item for every character $\chi$ of $G_n$, we have $\chi(\rho_1\times\dots\times\rho_r)\simeq\chi\rho_1\times\dots\times\chi\rho_r$;
	\item induction commutes with taking contragredients:
	\begin{equation}\label{eq:contragredient}
		(\rho_1\times\dots\times\rho_r)^\sim\simeq\widetilde\rho_1\times\dots\times\widetilde\rho_r .
	\end{equation}
\end{itemize}
Let $\mathcal{R}_m$ be the Grothendieck group of $\Rep(G_m)$. Normalized induction makes $\bigoplus_{m\ge0}\mathcal R_m$ a commutative ring \cite{Tad90}. In particular, $\rho_1\times\rho_2$ and $\rho_2\times\rho_1$ have the same constituents. Hence, if $\rho_1\times\rho_2$ is irreducible, then $\rho_1\times\rho_2\simeq\rho_2\times\rho_1$.

We shall often need the Jacquet modules of characters, which are computed as follows.

\begin{lemma}\label{lem:jacquet-char}
	Let $a,b\ge1$ and $x\in\R$. The normalized Jacquet module of the character $\nu^x$ of $G_{a+b}$ with respect to $P_{(a,b)}$ is the character $\nu^{x-eb/2}\otimes\nu^{x+ea/2}$ of $G_a\times G_b$.
\end{lemma}

\begin{proof}
	Since $N_{(a,b)}$ acts trivially on $\nu^x$, the space $(\nu^x)_{N_{(a,b)}}$ is the restriction of $\nu^x$ to $M_{(a,b)}$, namely $\nu^x\otimes\nu^x$. By Lemma~\ref{lem:modulus}, we have $\delta_{P_{(a,b)}}^{-1/2}(\diag(g,h))=\nu(g)^{-eb/2}\nu(h)^{ea/2}$, and the lemma follows.
\end{proof}

By Lemma~\ref{lem:modulus}, for characters $\chi_i$ of $G_{n_i}$, the representation $\chi_1\times\dots\times\chi_r$ is realized on the space of locally constant functions $f\colon G_n\to\C$ such that
\begin{equation}\label{eq:char-induced}
	f(\diag(g_1,\dots,g_r)\,u\,g)=\prod_{i=1}^r\chi_i(g_i)\,\nu(g_i)^{\frac e2\left(\sum_{j>i}n_j-\sum_{j<i}n_j\right)}f(g)\qquad(u\in N_\alpha,\ g\in G_n).
\end{equation}

Finally, we record an elementary observation that will be used repeatedly to compute spaces of homomorphisms into induced representations.

\begin{lemma}\label{lem:socle-hom}
	Let $X\in\Rep(G_m)$ have irreducible socle $\rho$, that is, $\rho$ is the unique irreducible subrepresentation of $X$, and let $\pi\in\Irr(G_m)$. Then $\dim\Hom_{G_m}(\pi,X)=1$ if $\pi\simeq\rho$, and $\Hom_{G_m}(\pi,X)=0$ otherwise.
\end{lemma}

\begin{proof}
	The image of a nonzero homomorphism $\pi\to X$ is an irreducible subrepresentation of $X$, and hence equals $\rho$. Thus, $\Hom_{G_m}(\pi,X)=\Hom_{G_m}(\pi,\rho)$. By Schur's lemma, this space is one-dimensional if $\pi\simeq\rho$ and zero otherwise.
\end{proof}

\subsection{Segments and the classifications of Zelevinsky and Langlands}\label{ss:segments}

We follow \cite[\S3]{LM16} and \cite{Tad90}. Let $\rho\in\mathcal C(G_k)$. By \cite[\S3.1 and Appendix A]{LM16}, there is a positive integer $s_\rho$ with the following property: for $\rho'\in\mathcal C(G_{k'})$, the representation $\rho\times\rho'$ is reducible if and only if $k'=k$ and $\rho'\simeq\nu_\rho^{\pm1}\rho$, where
$$
\nu_\rho=\nu^{s_\rho}.
$$
Consider the trivial character $\one_1$ of $G_1=\mathrm{D}^\times$. Its Jacquet--Langlands transfer is the Steinberg representation of $\GL_e(\mathrm{F})$, and therefore $s_{\one_1}=e$, that is,
$$
\nu_{\one_1}=\nu^e .
$$
Thus, for a character $\mu$ of $\mathrm{D}^\times$, the representation $\one_1\times\mu$ of $G_2$ is reducible if and only if $\mu\simeq\nu^{\pm e}$. Since $s_{\widetilde\rho}=s_\rho$ and $s_{\chi\rho}=s_\rho$ for every character $\chi$, the same holds for every character of $\mathrm{D}^\times$ in place of $\one_1$; that is, $\nu_\chi=\nu^e$ for every character $\chi$ of $\mathrm{D}^\times$.

The reducibility points $\nu_\rho^{\pm1}\rho$ lead to the notion of a segment, which we now recall.

\begin{definition}
	Let $\rho\in\mathcal C(G_k)$, and let $a,b\in\R$ with $b-a\in\Z_{\ge0}$. The \emph{segment} $\Delta=[a,b]_{(\rho)}$ is the sequence $(\nu_\rho^a\rho,\nu_\rho^{a+1}\rho,\dots,\nu_\rho^b\rho)$. Its \emph{length} is $b-a+1$, its \emph{degree} is $k(b-a+1)$, and its \emph{contragredient} is $\widetilde\Delta=[-b,-a]_{(\widetilde\rho)}$. Two segments $\Delta$ and $\Delta'$ are \emph{linked} if their union, as a set, is a segment and neither of them contains the other. If $\Delta$ and $\Delta'$ are linked, we may write $\Delta=[a,b]_{(\rho)}$ and $\Delta'=[a',b']_{(\rho)}$ with the same $\rho$. We then say that $\Delta$ \emph{precedes} $\Delta'$, and write $\Delta\prec\Delta'$, if $a<a'$; in this case, $a<a'\le b+1$ and $b<b'$. Moreover, $\Delta\cup\Delta'$ and $\Delta\cap\Delta'$ denote the segments formed by their union and their intersection, respectively, the latter being possibly empty.
\end{definition}

Each segment determines two irreducible representations, as follows. For a segment $\Delta=[a,b]_{(\rho)}$, the representation $\nu_\rho^a\rho\times\nu_\rho^{a+1}\rho\times\dots\times\nu_\rho^b\rho$ has a unique irreducible subrepresentation $\Ss(\Delta)$ and a unique irreducible quotient $\Qq(\Delta)$ \cite[\S3]{LM16}, \cite{Tad90}. Moreover, $\Ss(\Delta)$ is the unique irreducible quotient of $\nu_\rho^b\rho\times\dots\times\nu_\rho^a\rho$, and we have $\Ss(\Delta)^\sim\simeq\Ss(\widetilde\Delta)$ and $\Qq(\Delta)^\sim\simeq\Qq(\widetilde\Delta)$. We say that an irreducible representation is a \emph{discrete series representation} if it is essentially square-integrable. The representation $\Qq(\Delta)$ is a discrete series representation, and every discrete series representation is isomorphic to $\Qq(\Delta)$ for a unique segment $\Delta$.

\begin{notation}\label{not:line}
	Segments of the form $[a,b]_{(\one_1)}$ will simply be written $[a,b]$, and we put $[a]=[a,a]$. Thus, $[a,b]$ is the sequence of characters $\nu^{ea},\nu^{e(a+1)},\dots,\nu^{eb}$ of $\mathrm{D}^\times$. For $m\ge1$ and $x\in\R$, we put
	$$
	\Delta_m(x)=\left[x-\tfrac{m-1}2,\,x+\tfrac{m-1}2\right],\qquad \Delta_m=\Delta_m(0).
	$$
\end{notation}

In this notation, the characters $\nu^{ex}$ of $G_m$ correspond to the segments $\Delta_m(x)$.

\begin{lemma}\label{lem:characters-segments}
	Let $m\ge1$ and $x\in\R$. The character $\nu^{ex}$ of $G_m$ is isomorphic to $\Ss(\Delta_m(x))$. In particular, $\one_m=\Ss(\Delta_m)$.
\end{lemma}

\begin{proof}
	Apply \eqref{eq:char-induced} to the minimal parabolic subgroup $P_{(1,\dots,1)}$. For $1\le i\le m$, the $i$-th exponent is
	$$
	e\left(x-\tfrac{m-1}2+i-1\right)+\tfrac e2\big((m-i)-(i-1)\big)=ex.
	$$
	Hence, the representation $\nu^{e(x-\frac{m-1}2)}\times\nu^{e(x-\frac{m-3}2)}\times\dots\times\nu^{e(x+\frac{m-1}2)}$ is realized on the space of locally constant functions $f$ on $G_m$ such that $f(\diag(g_1,\dots,g_m)ug)=\prod_i\nu(g_i)^{ex}f(g)$ for all $u\in N_{(1,\dots,1)}$ and $g\in G_m$. The function $g\mapsto\nu(g)^{ex}$ lies in this space, and it spans a subrepresentation isomorphic to $\nu^{ex}$. Therefore, $\nu^{ex}$ is the unique irreducible subrepresentation $\Ss(\Delta_m(x))$.
\end{proof}

Finally, we introduce the Steinberg representations and the generalized Speh and Steinberg representations of $G_2$. We put $\St_m=\Qq(\Delta_m)$, the Steinberg representation of $G_m$, and, for $\sigma\in\Irr(G_1)$,
$$
\Sp_2(\sigma)=\Ss\big([-\tfrac12,\tfrac12]_{(\sigma)}\big),\qquad \St_2(\sigma)=\Qq\big([-\tfrac12,\tfrac12]_{(\sigma)}\big).
$$
These are, respectively, the unique irreducible subrepresentation and the unique irreducible quotient of $\nu_\sigma^{-1/2}\sigma\times\nu_\sigma^{1/2}\sigma$. For $\sigma=\one_1$, they are $\one_2$ and $\St_2$. For a character $\chi$ of $G_2$, we write $\chi\St_2=\chi\otimes\St_2$. By our convention that a character $\mu$ of $\mathrm{F}^\times$ also denotes $\mu\circ\Nrd$, a character $\chi$ of $\mathrm{D}^\times$ determines a character of $G_2$. With this identification, $\chi\St_2=\St_2(\chi)$ and $\chi\one_2=\Sp_2(\chi)$. 

A \emph{multisegment} is a finite multiset of segments. An ordering $(\Delta_1,\dots,\Delta_r)$ of a multisegment $\mm$ is \emph{admissible} if $\Delta_i$ does not precede $\Delta_j$ whenever $i<j$. Admissible orderings exist; for instance, any ordering in which the beginnings of the segments (for a common $\rho$) are non-increasing is admissible. We put
$$
\zeta(\mm)=\Ss(\Delta_1)\times\dots\times\Ss(\Delta_r),\qquad\lambda(\mm)=\Qq(\Delta_1)\times\dots\times\Qq(\Delta_r).
$$
Up to isomorphism, these representations do not depend on the choice of admissible ordering. We write $\widetilde\mm=\{\widetilde\Delta:\Delta\in\mm\}$.

The irreducible representations of $G_m$ are classified by multisegments in two ways, as follows.

\begin{theorem}[Zelevinsky classification; {\cite[Theorem 3.4]{LM16}}, \cite{Tad90}]\label{thm:zel}
	For every multisegment $\mm$, the representation $\zeta(\mm)$ has a unique irreducible subrepresentation $\Ss(\mm)$, and $\Ss(\mm)$ occurs in $\zeta(\mm)$ with multiplicity one. The map $\mm\mapsto\Ss(\mm)$ is a bijection from the set of multisegments of degree $m$ onto $\Irr(G_m)$.
\end{theorem}

\begin{theorem}[Langlands classification; {\cite[Theorem 3.6]{LM16}}, \cite{Tad90}]\label{thm:lang}
	For every multisegment $\mm$, the representation $\lambda(\mm)$ has a unique irreducible quotient $\Qq(\mm)$, and $\Qq(\mm)$ occurs in $\lambda(\mm)$ with multiplicity one. The map $\mm\mapsto\Qq(\mm)$ is a bijection from the set of multisegments of degree $m$ onto $\Irr(G_m)$.
\end{theorem}

We call $\mm$ the \emph{Zelevinsky multisegment} of $\Ss(\mm)$ and the \emph{Langlands multisegment} of $\Qq(\mm)$. By Lemma~\ref{lem:characters-segments}, the Zelevinsky multisegment of the character $\nu^{ex}$ of $G_m$ is $\{\Delta_m(x)\}$, a single segment of length $m$. On the other hand, the Zelevinsky multisegment of $\rho\in\mathcal C(G_m)$ is $\{[0]_{(\rho)}\}$, a single segment of length one. Consequently, a character of $G_m$ with $m\ge2$ is not supercuspidal. We shall also use the isomorphisms $\Ss(\mm)^\sim\simeq\Ss(\widetilde\mm)$ and $\Qq(\mm)^\sim\simeq\Qq(\widetilde\mm)$ \cite[\S3]{LM16}.

An \emph{elementary operation} on a multisegment replaces a pair $\{\Delta,\Delta'\}$ of linked segments by $\{\Delta\cup\Delta',\Delta\cap\Delta'\}$, where the intersection is omitted if it is empty. We write $\bb\le\aaa$ if $\bb$ is obtained from $\aaa$ by a finite number (possibly zero) of elementary operations.

The simplest instance of these classifications is the product of two segment representations, which we recall explicitly.

\begin{lemma}[{\cite[Proposition 3.5 and Theorem 3.6]{LM16}}, {\cite[Proposition 4.3]{Tad90}}]\label{lem:two-segments}
	Let $\Delta_1$ and $\Delta_2$ be segments. Then $\Ss(\Delta_1)\times\Ss(\Delta_2)$ is irreducible if and only if $\Delta_1$ and $\Delta_2$ are not linked, and the same holds for $\Qq(\Delta_1)\times\Qq(\Delta_2)$. Suppose that $\Delta_1$ and $\Delta_2$ are linked, and put $\Delta_3=\Delta_1\cup\Delta_2$ and $\Delta_4=\Delta_1\cap\Delta_2$. Then both products have length two. Moreover, if $\Delta_2\prec\Delta_1$, then:
	\begin{itemize}
		\item $\Ss(\Delta_1)\times\Ss(\Delta_2)$ has unique irreducible subrepresentation $\Ss(\{\Delta_1,\Delta_2\})$ and unique irreducible quotient $\Ss(\Delta_3)\times\Ss(\Delta_4)$;
		\item $\Qq(\Delta_1)\times\Qq(\Delta_2)$ has unique irreducible subrepresentation $\Qq(\Delta_3)\times\Qq(\Delta_4)$ and unique irreducible quotient $\Qq(\{\Delta_1,\Delta_2\})$.
	\end{itemize}
	Here, a factor attached to an empty $\Delta_4$ is omitted, and the representations $\Ss(\Delta_3)\times\Ss(\Delta_4)\simeq\Ss(\{\Delta_3,\Delta_4\})$ and $\Qq(\Delta_3)\times\Qq(\Delta_4)\simeq\Qq(\{\Delta_3,\Delta_4\})$ are irreducible. If $\Delta_1\prec\Delta_2$, then the roles of the subrepresentation and the quotient are interchanged.
\end{lemma}

More generally, the constituents of $\zeta(\aaa)$ and $\lambda(\aaa)$ are controlled by the order $\le$.

\begin{lemma}[{\cite[Proposition 3.5 and Theorem 3.6]{LM16}}, {\cite[Proposition 4.4]{Tad90}}]\label{lem:constituents}
	Let $\aaa$ and $\bb$ be multisegments of the same degree. Then $\Ss(\bb)$ is a constituent of $\zeta(\aaa)$ if and only if $\bb\le\aaa$. Similarly, $\Qq(\bb)$ is a constituent of $\lambda(\aaa)$ if and only if $\bb\le\aaa$.
\end{lemma}

When a product of irreducible representations is irreducible, its multisegment is easy to describe.

\begin{lemma}[{\cite[\S3]{LM16}}]\label{lem:sum-multisegments}
	Let $\pi_i=\Ss(\mm_i)$ for $i=1,2$, and suppose that $\pi_1\times\pi_2$ is irreducible. Then $\pi_1\times\pi_2\simeq\Ss(\mm_1+\mm_2)$. Likewise, if $\pi_i=\Qq(\mm_i)$ for $i=1,2$ and $\pi_1\times\pi_2$ is irreducible, then $\pi_1\times\pi_2\simeq\Qq(\mm_1+\mm_2)$.
\end{lemma}

Finally, representations whose supercuspidal supports lie on different lines never interact.

\begin{lemma}[{\cite{Tad90}}, {\cite[\S2--3]{LM16}}]\label{lem:different-lines}
	Let $\pi_1\in\Irr(G_a)$ and $\pi_2\in\Irr(G_b)$, and suppose that no element of the supercuspidal support of $\pi_1$ is of the form $\nu_\rho^k\rho$ with $\rho$ in the supercuspidal support of $\pi_2$ and $k\in\Z$. Then $\pi_1\times\pi_2$ is irreducible.
\end{lemma}

The following consequences of Theorem~\ref{thm:zel} and Lemma~\ref{lem:two-segments} will be used repeatedly.

\begin{lemma}\label{lem:S-quotient}
	Let $(\Delta_1,\dots,\Delta_r)$ be an admissible ordering of a multisegment $\mm$. Then $\Ss(\mm)$ is the unique irreducible quotient of $\Ss(\Delta_r)\times\dots\times\Ss(\Delta_1)$.
\end{lemma}

\begin{proof}
	For segments $\Delta$ and $\Delta'$, we have $\Delta\prec\Delta'$ if and only if $\widetilde{\Delta'}\prec\widetilde\Delta$. Hence, the ordering $(\widetilde\Delta_r,\dots,\widetilde\Delta_1)$ of $\widetilde\mm$ is admissible. By Theorem~\ref{thm:zel}, the representation $\Ss(\widetilde\mm)\simeq\Ss(\mm)^\sim$ is the unique irreducible subrepresentation of $\Ss(\widetilde\Delta_r)\times\dots\times\Ss(\widetilde\Delta_1)$. By \eqref{eq:contragredient} and the isomorphisms $\Ss(\Delta_i)^\sim\simeq\Ss(\widetilde\Delta_i)$, the latter representation is the contragredient of $\Ss(\Delta_r)\times\dots\times\Ss(\Delta_1)$. Since taking contragredients is an exact involution on $\Rep(G_m)$ that interchanges subrepresentations and quotients, the assertion follows.
\end{proof}

The next lemma allows us to reduce questions about irreducible quotients of induced representations to the case of irreducible inducing data.

\begin{lemma}\label{lem:quotient-constituents}
	Let $X\in\Rep(G_k)$ and $Y\in\Rep(G_{m-k})$, and let $\pi\in\Irr(G_m)$ be a quotient of $X\times Y$. Then $\pi$ is a quotient of $X'\times Y'$ for some constituent $X'$ of $X$ and some constituent $Y'$ of $Y$.
\end{lemma}

\begin{proof}
	We argue by induction on the sum of the lengths of $X$ and $Y$. If $X$ and $Y$ are both irreducible, there is nothing to prove. Suppose that $X$ is reducible, and let $0\to X_1\to X\to X_2\to0$ be an exact sequence with $X_1,X_2\ne0$. Since induction is exact, the sequence $0\to X_1\times Y\to X\times Y\to X_2\times Y\to0$ is exact. A surjection $X\times Y\to\pi$ is either nonzero on $X_1\times Y$, in which case its restriction to $X_1\times Y$ is surjective because $\pi$ is irreducible, or it factors through $X_2\times Y$. In both cases, the assertion follows from the induction hypothesis, since every constituent of $X_1$ or $X_2$ is a constituent of $X$. The same argument applies if $Y$ is reducible.
\end{proof}

If no two segments of a multisegment are linked, the corresponding product is irreducible.

\begin{lemma}\label{lem:unlinked-product}
	If the segments of $\mm=\{\Delta_1,\dots,\Delta_r\}$ are pairwise unlinked, then $\Ss(\Delta_1)\times\dots\times\Ss(\Delta_r)$ is irreducible and isomorphic to $\Ss(\mm)$, for every order of the factors.
\end{lemma}

\begin{proof}
	Since no two segments of $\mm$ are linked, no segment of $\mm$ precedes another, and therefore every ordering of $\mm$ is admissible. Thus, for every order of the factors, the product is $\zeta(\mm)$. By Lemma~\ref{lem:constituents}, the constituents of $\zeta(\mm)$ are among the $\Ss(\bb)$ with $\bb\le\mm$. No elementary operation applies to $\mm$, so the only such $\bb$ is $\mm$ itself. Since $\Ss(\mm)$ occurs in $\zeta(\mm)$ with multiplicity one by Theorem~\ref{thm:zel}, the representation $\zeta(\mm)$ is irreducible and isomorphic to $\Ss(\mm)$.
\end{proof}

We now make Lemma~\ref{lem:two-segments} explicit for two characters of $G_1=\mathrm{D}^\times$ on the line of $\one_1$.

\begin{lemma}\label{lem:St2}
	Let $a\in\R$. The representation $\nu^{ea}\times\nu^{e(a+1)}$ of $G_2$ has length two. Its unique irreducible subrepresentation is the character $\nu^{e(a+\frac12)}$, and its unique irreducible quotient is
	$$
	\nu^{e(a+\frac12)}\St_2=\Ss(\{[a],[a+1]\})=\Qq([a,a+1]).
	$$
	In $\nu^{e(a+1)}\times\nu^{ea}$, the subrepresentation and the quotient are interchanged.
\end{lemma}

\begin{proof}
	We apply Lemma~\ref{lem:two-segments} to $\Delta_1=[a]$ and $\Delta_2=[a+1]$. These segments are linked, with $[a]\prec[a+1]$, $[a]\cup[a+1]=[a,a+1]=\Delta_2(a+\frac12)$ and $[a]\cap[a+1]=\emptyset$. Since $\Delta_1\prec\Delta_2$, the representation $\Ss([a])\times\Ss([a+1])=\nu^{ea}\times\nu^{e(a+1)}$ has unique irreducible subrepresentation $\Ss(\Delta_2(a+\frac12))$, which is the character $\nu^{e(a+\frac12)}$ by Lemma~\ref{lem:characters-segments}, and unique irreducible quotient $\Ss(\{[a],[a+1]\})$. On the other hand, the same representation equals $\Qq([a])\times\Qq([a+1])$, so its unique irreducible quotient is also $\Qq([a,a+1])$. Hence, $\Ss(\{[a],[a+1]\})=\Qq([a,a+1])$. For $a=-\frac12$, this representation is $\St_2=\Qq(\Delta_2)$. Since twisting by $\nu^{ec}$ shifts every segment on the line of $\one_1$ by $c$ and commutes with parabolic induction, twisting by $\nu^{e(a+\frac12)}$ gives the general case. The final assertion follows in the same way from the case $\Delta_2\prec\Delta_1$ of Lemma~\ref{lem:two-segments}.
\end{proof}

For $a\in\R$, we put
\begin{equation}\label{eq:tau-a}
	\tau_a=\nu^{ea}\times\nu^{-ea}\in\Rep(G_2).
\end{equation}
The segments $[a]$ and $[-a]$ are linked if and only if $a=\pm\frac12$. Hence, by Lemma~\ref{lem:two-segments}, $\tau_a$ is irreducible if and only if $a\ne\pm\frac12$, and in that case $\tau_a\simeq\tau_{-a}$.

\subsection{The representations \texorpdfstring{$Q_m$}{Qm} and \texorpdfstring{$\widetilde Q_m$}{Q~m}}\label{ss:Qm}

We now introduce two families of representations that appear as exceptional constituents of degenerate principal series. They are the inner-form analogues of the representations $L_n$ and $Z_n$ of \cite[Definition 3.1]{Ven13}.

\begin{definition}\label{def:Qm}
	For $m\ge2$, we put
	$$
	Q_m=\Ss\big(\{\Delta_{m-1}(\tfrac12),[\tfrac{m+1}2]\}\big),\qquad \Qt_m=\Ss\big(\{[-\tfrac{m+1}2],\Delta_{m-1}(-\tfrac12)\}\big).
	$$
\end{definition}

Explicitly, $\Delta_{m-1}(\frac12)=[-\frac{m-3}2,\frac{m-1}2]$ and $\Delta_{m-1}(-\frac12)=[-\frac{m-1}2,\frac{m-3}2]$. Since $\Delta_{m-1}(\frac12)\prec[\frac{m+1}2]$, Lemma~\ref{lem:S-quotient} shows that $Q_m$ is the unique irreducible quotient of the representation $\nu^{e/2}\times\nu^{e\frac{m+1}2}$ induced from $G_{m-1}\times G_1$. Likewise, $\Qt_m$ is the unique irreducible quotient of the representation $\nu^{-e\frac{m+1}2}\times\nu^{-e/2}$ induced from $G_1\times G_{m-1}$. The contragredient of the multisegment $\{\Delta_{m-1}(\frac12),[\frac{m+1}2]\}$ is $\{\Delta_{m-1}(-\frac12),[-\frac{m+1}2]\}$, and therefore $\Qt_m\simeq Q_m^\sim$. For $t\in\R$, the twist $\nu^{et}Q_m$ is the representation $\Ss$ of the multisegment obtained by adding $t$ to all endpoints, and similarly for $\nu^{et}\Qt_m$. For $m=2$, Lemma~\ref{lem:St2} gives
\begin{equation}\label{eq:Q2}
	Q_2=\nu^e\St_2,\qquad \Qt_2=\nu^{-e}\St_2 .
\end{equation}

We first record that no twist of $Q_m$ or $\Qt_m$ is a character.

\begin{lemma}\label{lem:Qm-not-char}
	Let $m\ge2$ and $t\in\R$. Neither $\nu^{et}Q_m$ nor $\nu^{et}\Qt_m$ is a character.
\end{lemma}

\begin{proof}
	The Zelevinsky multisegments of these representations consist of two segments, whereas the Zelevinsky multisegment of a character consists of a single segment by Lemma~\ref{lem:characters-segments}. Hence, the assertion follows from Theorem~\ref{thm:zel}.
\end{proof}

For $m=2$, the representations $Q_2$ and $\Qt_2$ are twists of each other by \eqref{eq:Q2}. For $m\ge3$, this is no longer the case.

\begin{lemma}\label{lem:Q-vs-Qt}
	Let $m\ge3$. Then $\nu^{et}Q_m\not\simeq\nu^{et'}\Qt_m$ for all $t,t'\in\R$.
\end{lemma}

\begin{proof}
	The Zelevinsky multisegment of a twist of $Q_m$ consists of a segment of length $m-1\ge2$ together with a singleton immediately to the right of its end. The Zelevinsky multisegment of a twist of $\Qt_m$ consists of a segment of length $m-1$ together with a singleton immediately to the left of its beginning. Since $m-1\ge2$, these multisegments are different, and the assertion follows from Theorem~\ref{thm:zel}.
\end{proof}

To realize $Q_m$ and $\Qt_m$ as quotients of standard modules, we need their Langlands multisegments.

\begin{lemma}\label{lem:Qm-Langlands}
	Let $m\ge3$. The Langlands multisegments of $Q_m$ and $\Qt_m$ are given by
	\begin{align*}
		Q_m&=\Qq\Big(\big\{[j]:j=-\tfrac{m-3}2,\dots,\tfrac{m-3}2\big\}\cup\big\{[\tfrac{m-1}2,\tfrac{m+1}2]\big\}\Big),\\
		\Qt_m&=\Qq\Big(\big\{[-\tfrac{m+1}2,-\tfrac{m-1}2]\big\}\cup\big\{[j]:j=-\tfrac{m-3}2,\dots,\tfrac{m-3}2\big\}\Big),
	\end{align*}
	where $j$ runs over $\frac{m-3}2-\Z_{\ge0}$.
\end{lemma}

\begin{proof}
	By the algorithm of M\oe glin and Waldspurger \cite{MW86}, which holds for $G_m$ by \cite{BR07}, we have $\Ss(\aaa)\simeq\Qq(\aaa^\#)$ for a multisegment $\aaa^\#$ computed as follows. Let $d$ be the largest endpoint of a segment of $\aaa$. Choose a shortest segment $\Delta^{(1)}$ of $\aaa$ ending at $d$, then a shortest segment $\Delta^{(2)}$ of $\aaa$ ending at $d-1$ and preceding $\Delta^{(1)}$, and so on, for as long as possible. If this chain has $r$ members, then $[d-r+1,d]$ is a segment of $\aaa^\#$, and the remaining segments of $\aaa^\#$ are those of $(\aaa^-)^\#$, where $\aaa^-$ is obtained from $\aaa$ by removing the last element of each $\Delta^{(i)}$. In particular, for a single segment $\Delta$, the algorithm produces the multisegment consisting of the singletons of $\Delta$.
	
	Let $\aaa=\{\Delta_{m-1}(\frac12),[\frac{m+1}2]\}$. Then $d=\frac{m+1}2$ and $\Delta^{(1)}=[\frac{m+1}2]$. The segment $\Delta^{(2)}=[-\frac{m-3}2,\frac{m-1}2]$ ends at $d-1$ and precedes $\Delta^{(1)}$, and no segment of $\aaa$ ends at $\frac{m-3}2$. Hence, $[\frac{m-1}2,\frac{m+1}2]\in\aaa^\#$ and $\aaa^-=\{[-\frac{m-3}2,\frac{m-3}2]\}$. Since $\aaa^-$ consists of a single segment, $(\aaa^-)^\#$ consists of its singletons. This proves the first formula.
	
	Since $\Qt_m\simeq Q_m^\sim$ and $\Qq(\bb)^\sim\simeq\Qq(\widetilde\bb)$, the second formula follows from the first by replacing each segment with its contragredient. Alternatively, it can be obtained directly from the algorithm. For $\aaa=\{[-\frac{m+1}2],\Delta_{m-1}(-\frac12)\}$, the algorithm successively produces the singletons $[\frac{m-3}2],[\frac{m-5}2],\dots,[-\frac{m-3}2]$, because, as $m\ge3$, the singleton $[-\frac{m+1}2]$ does not end at $d-1$ for any of these steps. After these steps, we have $\aaa^-=\{[-\frac{m+1}2],[-\frac{m-1}2]\}$, which yields the segment $[-\frac{m+1}2,-\frac{m-1}2]$.
\end{proof}

As a consequence, we can realize $Q_m$ and $\Qt_m$ as the unique irreducible quotients of representations induced from $\St_2$ and a trivial character.

\begin{lemma}\label{lem:Qm-quotient}
	Let $m\ge3$. Then $Q_m$ is a quotient of the representation $\nu^{em/2}\St_2\times\one_{m-2}$ induced from $G_2\times G_{m-2}$, and every irreducible quotient of this representation is isomorphic to $Q_m$.
\end{lemma}

\begin{proof}
	By Lemma~\ref{lem:Qm-Langlands}, we have $Q_m=\Qq(\aaa)$, where $\aaa$ consists of $[\frac{m-1}2,\frac{m+1}2]$ and the singletons $[j]$ with $|j|\le\frac{m-3}2$. The ordering $([\frac{m-1}2,\frac{m+1}2],[\frac{m-3}2],[\frac{m-5}2],\dots,[-\frac{m-3}2])$ of $\aaa$ is admissible. Hence, by Theorem~\ref{thm:lang} and Lemma~\ref{lem:St2}, $Q_m$ is the unique irreducible quotient of
	$$
	\lambda(\aaa)=\nu^{em/2}\St_2\times\nu^{e\frac{m-3}2}\times\dots\times\nu^{-e\frac{m-3}2}.
	$$
	The product of the last $m-2$ factors has $\Ss(\Delta_{m-2})=\one_{m-2}$ as a quotient (\S\ref{ss:segments}). By the exactness and transitivity of parabolic induction, $\lambda(\aaa)$ surjects onto $\nu^{em/2}\St_2\times\one_{m-2}$. Therefore, every irreducible quotient of $\nu^{em/2}\St_2\times\one_{m-2}$ is an irreducible quotient of $\lambda(\aaa)$, and hence is isomorphic to $Q_m$. Finally, $\nu^{em/2}\St_2\times\one_{m-2}$ is nonzero and of finite length, so it has an irreducible quotient.
\end{proof}

The analogous statement for $\Qt_m$ is proved in the same way.

\begin{lemma}\label{lem:Qtm-quotient}
	Let $m\ge3$. Then $\Qt_m$ is a quotient of the representation $\one_{m-2}\times\nu^{-em/2}\St_2$ induced from $G_{m-2}\times G_2$, and every irreducible quotient of this representation is isomorphic to $\Qt_m$.
\end{lemma}

\begin{proof}
	By Lemma~\ref{lem:Qm-Langlands}, we have $\Qt_m=\Qq(\aaa)$, where $\aaa$ consists of the singletons $[j]$ with $|j|\le\frac{m-3}2$ and the segment $[-\frac{m+1}2,-\frac{m-1}2]$. The ordering $([\frac{m-3}2],\dots,[-\frac{m-3}2],[-\frac{m+1}2,-\frac{m-1}2])$ of $\aaa$ is admissible. Hence, $\Qt_m$ is the unique irreducible quotient of
	$$
	\lambda(\aaa)=\nu^{e\frac{m-3}2}\times\dots\times\nu^{-e\frac{m-3}2}\times\nu^{-em/2}\St_2,
	$$
	which surjects onto $\one_{m-2}\times\nu^{-em/2}\St_2$. We conclude as in the proof of Lemma~\ref{lem:Qm-quotient}.
\end{proof}

\subsection{Induction of a character and a representation of \texorpdfstring{$\mathrm{D}^\times$}{Dx}}\label{ss:char-times-sigma}

The next two lemmas determine the composition series, and in particular the socles, of the representations that arise from the open orbit in \S\ref{sec:orbits}.

\begin{lemma}\label{lem:char-sigma}
	Let $m\ge2$, $y\in\R$ and $\sigma\in\Irr(G_1)$, and regard $\nu^y$ as a character of $G_{m-1}$. Then $\nu^y\times\sigma$ is reducible if and only if $\sigma\simeq\nu^{y\pm em/2}$. More precisely, the following hold.
	\begin{enumerate}
		\item If $\sigma\simeq\nu^{y+em/2}$, then $\nu^y\times\sigma$ has unique irreducible subrepresentation $\nu^{y+e/2}$ and unique irreducible quotient $\nu^{y-e/2}Q_m$.
		\item If $\sigma\simeq\nu^{y-em/2}$, then $\nu^y\times\sigma$ has unique irreducible subrepresentation $\nu^{y+e/2}\Qt_m$ and unique irreducible quotient $\nu^{y-e/2}$.
		\item Otherwise, $\nu^y\times\sigma$ is irreducible.
	\end{enumerate}
\end{lemma}

\begin{proof}
	By Lemma~\ref{lem:characters-segments}, we have $\nu^y=\Ss(\Delta_{m-1}(y/e))$, where $\Delta_{m-1}(y/e)=[\frac ye-\frac{m-2}2,\frac ye+\frac{m-2}2]$. Moreover, $\sigma=\Ss([0]_{(\sigma)})$. By Lemma~\ref{lem:two-segments}, the representation $\nu^y\times\sigma=\Ss(\Delta_{m-1}(y/e))\times\Ss([0]_{(\sigma)})$ is reducible if and only if the segments $[0]_{(\sigma)}$ and $\Delta_{m-1}(y/e)$ are linked. This forces $\sigma$ to lie on the line of $\one_1$, that is, $\sigma\simeq\nu^{ec}=\Ss([c])$ for some $c\in\R$. Moreover, $[c]$ is linked to $\Delta_{m-1}(y/e)$ exactly when $c$ is adjacent to one of its endpoints, that is, when $c=\frac ye\pm\frac m2$. This proves the first assertion and (3).
	
	Suppose that $c=\frac ye+\frac m2$. Then $\Delta_{m-1}(y/e)\prec[c]$ and $\Delta_{m-1}(y/e)\cup[c]=\Delta_m(\frac ye+\frac12)$. By Lemma~\ref{lem:two-segments} and Lemma~\ref{lem:characters-segments}, the unique irreducible subrepresentation of $\nu^y\times\sigma$ is $\Ss(\Delta_m(\frac ye+\frac12))=\nu^{y+e/2}$. Its unique irreducible quotient is $\Ss(\{\Delta_{m-1}(y/e),[c]\})$. This multisegment is obtained from $\{\Delta_{m-1}(\frac12),[\frac{m+1}2]\}$ by adding $\frac ye-\frac12$ to all endpoints, so the quotient is $\nu^{e(\frac ye-\frac12)}Q_m=\nu^{y-e/2}Q_m$. This proves (1).
	
	Suppose that $c=\frac ye-\frac m2$. Then $[c]\prec\Delta_{m-1}(y/e)$ and $[c]\cup\Delta_{m-1}(y/e)=\Delta_m(\frac ye-\frac12)$. By Lemma~\ref{lem:two-segments}, the unique irreducible subrepresentation of $\nu^y\times\sigma$ is $\Ss(\{[c],\Delta_{m-1}(y/e)\})$. This multisegment is obtained from $\{[-\frac{m+1}2],\Delta_{m-1}(-\frac12)\}$ by adding $\frac ye+\frac12$ to all endpoints, so the subrepresentation is $\nu^{y+e/2}\Qt_m$. The unique irreducible quotient is $\Ss(\Delta_m(\frac ye-\frac12))=\nu^{y-e/2}$. This proves (2).
\end{proof}

Reversing the order of the factors gives the following companion statement.

\begin{lemma}\label{lem:sigma-char}
	Let $m\ge2$, $y\in\R$ and $\sigma\in\Irr(G_1)$, and regard $\nu^y$ as a character of $G_{m-1}$. Then $\sigma\times\nu^y$ is reducible if and only if $\sigma\simeq\nu^{y\pm em/2}$. More precisely, the following hold.
	\begin{enumerate}
		\item If $\sigma\simeq\nu^{y-em/2}$, then $\sigma\times\nu^y$ has unique irreducible subrepresentation $\nu^{y-e/2}$ and unique irreducible quotient $\nu^{y+e/2}\Qt_m$.
		\item If $\sigma\simeq\nu^{y+em/2}$, then $\sigma\times\nu^y$ has unique irreducible subrepresentation $\nu^{y-e/2}Q_m$ and unique irreducible quotient $\nu^{y+e/2}$.
		\item Otherwise, $\sigma\times\nu^y$ is irreducible.
	\end{enumerate}
\end{lemma}

\begin{proof}
	The segments involved are those in the proof of Lemma~\ref{lem:char-sigma}, with the factors in the opposite order. By Lemma~\ref{lem:two-segments}, reversing the order of the factors interchanges the unique irreducible subrepresentation and the unique irreducible quotient, and it does not affect irreducibility.
\end{proof}

In particular, taking $y=\frac e2$ and $\sigma=\nu^{e\frac{m+1}2}$ in Lemma~\ref{lem:char-sigma}(1), we obtain an exact sequence
\begin{equation}\label{eq:Qm-sequence}
	0\longrightarrow\nu^e\longrightarrow\nu^{e/2}\times\nu^{e\frac{m+1}2}\longrightarrow Q_m\longrightarrow0,
\end{equation}
which is the analogue of \cite[(3.2)]{Ven13}. This sequence does not split, since $\nu^{e/2}\times\nu^{e\frac{m+1}2}$ has a unique irreducible subrepresentation.

\subsection{Twists of the Steinberg representation induced with a character}\label{ss:Sn}

In this subsection, we study the representations $\nu^{\pm en/2}\St_2\times\one_{n-2}$ and $\one_{n-2}\times\nu^{\pm en/2}\St_2$, which appear in the analysis of the closed orbits. For $n\ge3$, we put
$$
\sss_n=\big\{\Delta_{n-2},[\tfrac{n-1}2],[\tfrac{n+1}2]\big\},\qquad \widetilde\sss_n=\big\{\Delta_{n-2},[-\tfrac{n-1}2],[-\tfrac{n+1}2]\big\},
$$
and
$$
S_n=\Ss(\sss_n),\qquad \Stil_n=\Ss(\widetilde\sss_n).
$$
Since $\widetilde{\Delta_{n-2}}=\Delta_{n-2}$, the multisegment $\widetilde\sss_n$ is the contragredient of $\sss_n$, and therefore $\Stil_n\simeq S_n^\sim$.

We first determine when a twist of $\St_2$ induced with the trivial character of $G_{n-2}$ is reducible.

\begin{lemma}\label{lem:St2-times-one}
	Let $n\ge3$, and let $\chi$ be a character of $G_2$. Then $\chi\St_2\times\one_{n-2}$ is reducible if and only if $\chi=\nu^{\pm en/2}$.
\end{lemma}
\begin{proof}
	We use the notation of \cite[\S3]{LM16} for the cuspidal line of $\rho=\one_1$, so that $\nu_\rho=\nu^{e}$. Under this normalization, we have $\Ss(\mm)=Z(\mm)$ and $\Qq(\mm)=L(\mm)$, and segments are written in terms of exponents of $\nu^{e}$.
	
	By Lemma~\ref{lem:characters-segments}, we have $\one_{n-2}=\Ss(\Delta)$, where $\Delta=[a,b]$ with $a=-\frac{n-3}2$ and $b=\frac{n-3}2$. If the supercuspidal support of $\chi\St_2$ does not meet the line $\{\nu^{ej}:j\in\frac{n-3}2+\Z\}$, then the product is irreducible by Lemma~\ref{lem:different-lines}. Otherwise, we have $\chi=\nu^{ex}$ with $x\in\frac n2+\Z$, and, by Lemma~\ref{lem:St2}, $\chi\St_2=\Ss(\mm)$, where $\mm=[x-\frac12]+[x+\frac12]$.
	
	By \cite[Corollary~5.14]{LM16}, the representation $\Ss(\Delta)\times\Ss(\mm)$ is irreducible if and only if both $\mathrm{LC}(\Delta,\mm)$ and $\mathrm{RC}(\Delta,\mm)$ hold. Recall from \cite[\S\S5.5--5.6]{LM16} that $\mathrm{LC}(\Delta,\mm)$ is the existence of an injective map
	$$
	f\colon X=\{\Delta'\in\mm:\Delta\prec\Delta'\}\longrightarrow Y=\{\Delta'\in\mm:[a-1,b-1]\prec\Delta'\}
	$$
	such that $f(\Delta')\prec\Delta'$ for all $\Delta'\in X$. Similarly, $\mathrm{RC}(\Delta,\mm)$ is the existence of an injective map
	$$
	g\colon \widetilde X=\{\Delta'\in\mm:\Delta'\prec\Delta\}\longrightarrow \widetilde Y=\{\Delta'\in\mm:\Delta'\prec[a+1,b+1]\}
	$$
	such that $\Delta'\prec g(\Delta')$ for all $\Delta'\in\widetilde X$. For a singleton $[c]$, we have
	$$
	\Delta\prec[c]\iff c=b+1,\qquad [a-1,b-1]\prec[c]\iff c=b,
	$$
	$$
	[c]\prec\Delta\iff c=a-1,\qquad [c]\prec[a+1,b+1]\iff c=a.
	$$
	Hence, $X=\emptyset$ unless $b+1\in\{x-\frac12,x+\frac12\}$. If $x+\frac12=b+1$, then $X=\{[b+1]\}$ and $Y=\{[b]\}$; since $[b]\prec[b+1]$, the condition $\mathrm{LC}(\Delta,\mm)$ holds. If $x-\frac12=b+1$, then $X=\{[b+1]\}$ and $Y=\emptyset$, so $\mathrm{LC}(\Delta,\mm)$ fails. Thus, $\mathrm{LC}(\Delta,\mm)$ fails if and only if $x=b+\frac32=\frac n2$. By the same argument, $\mathrm{RC}(\Delta,\mm)$ fails if and only if $x+\frac12=a-1$, that is, $x=-\frac n2$. Therefore, $\chi\St_2\times\one_{n-2}$ is reducible if and only if $\chi=\nu^{\pm en/2}$.
\end{proof}

At the two reducibility points, the representations in question embed into standard modules attached to $\sss_n$ and $\widetilde\sss_n$, as follows.

\begin{lemma}\label{lem:Sn-embedding}
	Let $n\ge3$. Then $\nu^{en/2}\St_2\times\one_{n-2}$ is a subrepresentation of $\zeta(\sss_n)$, and its constituents are among the representations $\Ss(\bb)$ with $\bb\le\sss_n$.
\end{lemma}

\begin{proof}
	By Lemma~\ref{lem:St2} with $a=\frac{n-1}2$, the representation $\nu^{en/2}\St_2$ is a subrepresentation of $\nu^{e\frac{n+1}2}\times\nu^{e\frac{n-1}2}=\Ss([\frac{n+1}2])\times\Ss([\frac{n-1}2])$. Since parabolic induction is exact, $\nu^{en/2}\St_2\times\one_{n-2}$ is a subrepresentation of $\Ss([\frac{n+1}2])\times\Ss([\frac{n-1}2])\times\Ss(\Delta_{n-2})$. The ordering $([\frac{n+1}2],[\frac{n-1}2],\Delta_{n-2})$ of $\sss_n$ is admissible, so this representation is $\zeta(\sss_n)$. The second assertion follows from Lemma~\ref{lem:constituents}.
\end{proof}

The contragredient situation is handled in the same way.

\begin{lemma}\label{lem:Stn-embedding}
	Let $n\ge3$. Then $\one_{n-2}\times\nu^{-en/2}\St_2$ is a subrepresentation of $\zeta(\widetilde\sss_n)$, and its constituents are among the representations $\Ss(\bb)$ with $\bb\le\widetilde\sss_n$.
\end{lemma}

\begin{proof}
	By Lemma~\ref{lem:St2} with $a=-\frac{n+1}2$, the representation $\nu^{-en/2}\St_2$ is a subrepresentation of $\Ss([-\frac{n-1}2])\times\Ss([-\frac{n+1}2])$. Hence, $\one_{n-2}\times\nu^{-en/2}\St_2$ is a subrepresentation of $\Ss(\Delta_{n-2})\times\Ss([-\frac{n-1}2])\times\Ss([-\frac{n+1}2])$. The ordering $(\Delta_{n-2},[-\frac{n-1}2],[-\frac{n+1}2])$ of $\widetilde\sss_n$ is admissible, so this representation is $\zeta(\widetilde\sss_n)$, and the second assertion follows from Lemma~\ref{lem:constituents}.
\end{proof}

We now list the multisegments that can occur. In $\sss_n$, the only linked pairs are $\{\Delta_{n-2},[\frac{n-1}2]\}$ and $\{[\frac{n-1}2],[\frac{n+1}2]\}$, since $\Delta_{n-2}$ ends at $\frac{n-3}2$. Hence, the multisegments $\bb\le\sss_n$ are
\begin{equation}\label{eq:b-le-sn}
	\sss_n,\quad \big\{[-\tfrac{n-3}2,\tfrac{n-1}2],[\tfrac{n+1}2]\big\},\quad \big\{\Delta_{n-2},[\tfrac{n-1}2,\tfrac{n+1}2]\big\},\quad \big\{[-\tfrac{n-3}2,\tfrac{n+1}2]\big\}.
\end{equation}
Since $[-\frac{n-3}2,\frac{n-1}2]=\Delta_{n-1}(\frac12)$ and $[-\frac{n-3}2,\frac{n+1}2]=\Delta_n(1)$, the corresponding representations are $S_n$, $Q_n$, $\Ss(\{\Delta_{n-2},[\frac{n-1}2,\frac{n+1}2]\})$ and the character $\nu^e$, respectively. Likewise, the multisegments $\bb\le\widetilde\sss_n$ are
$$
\widetilde\sss_n,\quad \big\{[-\tfrac{n+1}2],[-\tfrac{n-1}2,\tfrac{n-3}2]\big\},\quad \big\{\Delta_{n-2},[-\tfrac{n+1}2,-\tfrac{n-1}2]\big\},\quad \big\{[-\tfrac{n+1}2,\tfrac{n-3}2]\big\},
$$
and the corresponding representations are $\Stil_n$, $\Qt_n$, $\Ss(\{\Delta_{n-2},[-\frac{n+1}2,-\frac{n-1}2]\})$ and $\nu^{-e}$, respectively.

\subsection{Two known results}\label{ss:known}

We shall use the following two results of Anandavardhanan, Lu, Matringe, S\'echerre, and Yang. Both hold for the inner forms $\GL_m(\mathrm{D})$, where $\mathrm{D}$ has arbitrary index. The first shows that discrete series representations never contribute to the classification when $n\ge3$.

\begin{theorem}[{\cite[Theorem 3.8]{ALMSY24}}]\label{thm:ds-not-dist}
	Let $n\ge3$, and let $\mu$ be a character of $H_{1,n-1}$. Then no discrete series representation of $G_n$ is $(H_{1,n-1},\mu)$-distinguished.
\end{theorem}

The second result is a multiplicity one theorem for untwisted linear periods when the two blocks have nearly equal size.

\begin{theorem}[{\cite[Theorem 2.5]{ALMSY24}}]\label{thm:ALMSY-mult}
	Let $p,q\ge1$ with $|p-q|\le1$. Then, for every $\pi\in\Irr(G_{p+q})$, we have $\dim\Hom_{H_{p,q}}(\pi,\one)\le1$.
\end{theorem}

\subsection{The Jacquet--Langlands transfer and Langlands parameters}\label{ss:JL}

In this subsection, we deal with representations of two different groups, namely $G_m=\GL_m(\mathrm{D})$ and $\GL_N(\mathrm{F})$. To avoid confusion, representations of $\GL_N(\mathrm{F})$ are denoted by symbols with a subscript $0$ or a prime, such as $\pi_0$ or $\rho'$. The Langlands parameter of a representation $\pi_0$ of $\GL_N(\mathrm{F})$ is denoted by $\varphi_{\pi_0}$, and the Langlands parameter of a representation $\pi$ of $G_m$, defined in Definition~\ref{def:L-param} below, is denoted by $\varphi_\pi$.

 Let $W_{\mathrm{F}}$ be the Weil group of $\mathrm{F}$. The local Langlands correspondence \cite{HT01,Hen00,Sch13} is a bijection $\pi_0\mapsto\varphi_{\pi_0}$ from $\Irr(\GL_N(\mathrm{F}))$ onto the set of isomorphism classes of $N$-dimensional representations of $W_{\mathrm{F}}\times\mathrm{SL}_2(\C)$ that are continuous and semisimple on $W_{\mathrm{F}}$ and algebraic on $\mathrm{SL}_2(\C)$. Such representations are semisimple. On the side of $\GL_N(\mathrm{F})$, we use the same notation for segments as in \S\ref{ss:segments}, with $\nu=|\det|_{\mathrm{F}}$ and $s_{\rho'}=1$. For $r\ge1$, let $\Sp_r$ denote the $r$-dimensional irreducible algebraic representation of $\mathrm{SL}_2(\C)$. If $\pi_0=\Qq(\{\Delta'_1,\dots,\Delta'_k\})\in\Irr(\GL_N(\mathrm{F}))$ with $\Delta'_i=[a_i,b_i]_{(\rho'_i)}$ and $\rho'_i$ supercuspidal, then
\begin{equation}\label{eq:LLC-formula}
	\varphi_{\pi_0}=\bigoplus_{i=1}^k\varphi_{\rho'_i}\otimes|\cdot|^{\frac{a_i+b_i}2}\otimes\Sp_{b_i-a_i+1},
\end{equation}
where $|\cdot|$ denotes the character of $W_{\mathrm{F}}$ corresponding to $|\cdot|_{\mathrm{F}}$ under local class field theory.

The Jacquet--Langlands correspondence $\JL$ is a bijection from the set of discrete series representations of $G_m$ onto the set of discrete series representations of $\GL_{em}(\mathrm{F})$ \cite{DKV84}. If $\rho\in\mathcal C(G_k)$, then $\JL(\rho)=\Qq([-\frac{s_\rho-1}2,\frac{s_\rho-1}2]_{(\rho')})$ for a supercuspidal representation $\rho'$ of $\GL_{ek/s_\rho}(\mathrm{F})$. For a segment $\Delta=[a,b]_{(\rho)}$ of $G_m$, we put
$$
\JL(\Delta)=\Big[\tfrac{(2a-1)s_\rho+1}2,\tfrac{(2b+1)s_\rho-1}2\Big]_{(\rho')}.
$$
This is a segment of $\GL_{em}(\mathrm{F})$ of length $s_\rho(b-a+1)$ and centre $\frac{s_\rho(a+b)}2$. Since $\JL$ is compatible with twisting by characters, we have $\JL(\Qq(\Delta))=\Qq(\JL(\Delta))$. Following \cite[\S2]{ALMSY24}, for a multisegment $\mm=\{\Delta_1,\dots,\Delta_t\}$ of $G_m$, we put $\JL(\mm)=\{\JL(\Delta_1),\dots,\JL(\Delta_t)\}$, which is a multisegment of $\GL_{em}(\mathrm{F})$.

We can now define the Langlands parameter of an arbitrary irreducible representation of $G_m$.

\begin{definition}[{cf.~\cite[Definition 2.1]{ALMSY24}}]\label{def:L-param}
	Let $\pi=\Qq(\mm)\in\Irr(G_m)$. We define the \emph{transfer} of $\pi$ to be the representation
	$$
	\JL(\pi)=\Qq(\JL(\mm))\in\Irr(\GL_{em}(\mathrm{F})),
	$$
	and the \emph{Langlands parameter} of $\pi$ to be
	$$
	\varphi_\pi:=\varphi_{\JL(\pi)}.
	$$
	Thus, $\varphi_\pi$ is an $em$-dimensional representation of $W_{\mathrm{F}}\times\mathrm{SL}_2(\C)$. We say that $\pi$ is \emph{generic} if $\JL(\pi)$ is generic.
\end{definition}

To summarize, $\pi$ is a representation of $G_m$, its transfer $\JL(\pi)$ is a representation of $\GL_{em}(\mathrm{F})$, and $\varphi_\pi$ is, by definition, the Langlands parameter of $\JL(\pi)$. If $\pi$ is generic, then $\pi$ corresponds, up to sign, to $\JL(\pi)$ under the map $\mathrm{LJ}$ of Badulescu, and Definition~\ref{def:L-param} agrees with \cite[Definition 2.1]{ALMSY24}.

For segments on the line of $\one_1$, the parameter is computed as follows.

\begin{lemma}\label{lem:L-segment}
	Let $a,b\in\R$ with $b-a\in\Z_{\ge0}$. Then $\varphi_{\Qq([a,b])}\simeq|\cdot|^{\frac{e(a+b)}2}\otimes\Sp_{e(b-a+1)}$. In particular, $\varphi_{\Qq([a])}\simeq|\cdot|^{ea}\otimes\Sp_e$.
\end{lemma}

\begin{proof}
	Since $s_{\one_1}=e$ and $\JL(\one_1)$ is the Steinberg representation of $\GL_e(\mathrm{F})$, the representation $\rho'$ is the trivial character of $\GL_1(\mathrm{F})$. Hence, $\JL([a,b])=[\frac{(2a-1)e+1}2,\frac{(2b+1)e-1}2]$ is a segment of length $e(b-a+1)$ and centre $\frac{e(a+b)}2$. The assertion now follows from \eqref{eq:LLC-formula}.
\end{proof}

The Langlands parameter is additive with respect to multisegments.

\begin{lemma}\label{lem:L-additive}
	For every multisegment $\mm$, we have $\varphi_{\Qq(\mm)}\simeq\bigoplus_{\Delta\in\mm}\varphi_{\Qq(\Delta)}$.
\end{lemma}

\begin{proof}
	By Definition~\ref{def:L-param}, we have $\JL(\Qq(\mm))=\Qq(\JL(\mm))$ and $\JL(\Qq(\Delta))=\Qq(\JL(\Delta))$ for each $\Delta\in\mm$. The assertion therefore follows from \eqref{eq:LLC-formula}.
\end{proof}

As for $\GL_N(\mathrm{F})$, an irreducible representation of $G_m$ is determined by its Langlands parameter.

\begin{lemma}\label{lem:L-injective}
	The map $\pi\mapsto\varphi_{\pi}$ is injective on $\Irr(G_m)$.
\end{lemma}

\begin{proof}
	The map $\Delta\mapsto\JL(\Delta)$ is injective on segments \cite[\S2]{ALMSY24}, and hence $\mm\mapsto\JL(\mm)$ is injective on multisegments. Moreover, the Langlands classification is bijective (Theorem~\ref{thm:lang}), and so is the local Langlands correspondence for $\GL_{em}(\mathrm{F})$. The assertion follows.
\end{proof}

Finally, we record the Langlands multisegment of a character of $G_m$.

\begin{lemma}\label{lem:char-Langlands}
	Let $m\ge1$ and $x\in\R$. The character $\nu^{ex}$ of $G_m$ is isomorphic to $\Qq(\aaa)$, where $\aaa$ consists of the singletons $[y]$ with $y\in\Delta_m(x)$. Consequently,
	$$
	\varphi_{\nu^{ex}}\simeq\bigoplus_{y\in\Delta_m(x)}|\cdot|^{ey}\otimes\Sp_e.
	$$
\end{lemma}

\begin{proof}
	By \S\ref{ss:segments}, the character $\nu^{ex}=\Ss(\Delta_m(x))$ is the unique irreducible quotient of the product of the characters $\nu^{ey}$ of $G_1$, $y\in\Delta_m(x)$, taken in decreasing order of $y$. Since $\nu^{ey}=\Qq([y])$, this product is $\lambda(\aaa)$ for an admissible ordering of $\aaa$. Hence, $\nu^{ex}\simeq\Qq(\aaa)$ by Theorem~\ref{thm:lang}. The formula for $\varphi_{\nu^{ex}}$ follows from Lemmas~\ref{lem:L-segment} and~\ref{lem:L-additive}.
\end{proof}

\section{Elementary facts about generalized linear periods}\label{sec:elementary}

In this section, we collect some elementary properties of generalized linear periods that will be used throughout the article. Unless stated otherwise, we assume that $n\ge2$, and we put $H=H_{1,n-1}$.

We begin by observing that conjugation by $w_{p,q}$ interchanges the roles of the two blocks, at the cost of replacing $s$ by $-s$.

\begin{lemma}\label{lem:swap}
	Let $p+q=n$ and $s\in\R$, and let $\pi$ be a smooth representation of $G_n$. Then the map $\ell\mapsto\ell\circ\pi(w_{p,q})^{-1}$ is an isomorphism
	$$
	\Hom_{H_{p,q}}\big(\pi,\chi^{(p,q)}_s\big)\xrightarrow{\ \sim\ }\Hom_{H_{q,p}}\big(\pi,\chi^{(q,p)}_{-s}\big).
	$$
	In particular, $\pi$ is $(H_{p,q},\chi_s)$-distinguished if and only if it is $(H_{q,p},\chi_{-s})$-distinguished.
\end{lemma}

\begin{proof}
	Put $w=w_{p,q}$ and $\ell_w=\ell\circ\pi(w)^{-1}$. Let $h'=\diag(g_2,g_1)\in H_{q,p}$ with $g_1\in G_p$ and $g_2\in G_q$. Then $w^{-1}h'w=\diag(g_1,g_2)\in H_{p,q}$, and hence, for every $v\in\pi$,
	$$
	\ell_w(\pi(h')v)=\ell\big(\pi(w^{-1}h'w)\pi(w)^{-1}v\big)=\nu(g_1)^{es}\nu(g_2)^{-es}\ell_w(v)=\chi^{(q,p)}_{-s}(h')\,\ell_w(v).
	$$
	Thus, $\ell_w\in\Hom_{H_{q,p}}(\pi,\chi^{(q,p)}_{-s})$. The same computation shows that $\ell'\mapsto\ell'\circ\pi(w)$ maps $\Hom_{H_{q,p}}(\pi,\chi^{(q,p)}_{-s})$ into $\Hom_{H_{p,q}}(\pi,\chi^{(p,q)}_s)$, and this map is inverse to $\ell\mapsto\ell_w$.
\end{proof}

Since the centre of $G_n$ is contained in $H_{p,q}$, distinction imposes a condition on the central character.

\begin{lemma}\label{lem:central}
	Let $p+q=n$ and $s\in\R$. If $\pi$ admits a central character $\omega_\pi$ and is $(H_{p,q},\chi_s)$-distinguished, then $\omega_\pi=|\cdot|_{\mathrm{F}}^{e^2s(p-q)}$. In particular, every $(H_{1,n-1},\chi_s)$-distinguished representation $\pi$ with a central character satisfies $\omega_\pi=|\cdot|_{\mathrm{F}}^{e^2s(2-n)}$.
\end{lemma}

\begin{proof}
	Let $0\ne\ell\in\Hom_{H_{p,q}}(\pi,\chi_s)$, and let $z\in\mathrm{F}^\times$. Since $zI_n=\diag(zI_p,zI_q)\in H_{p,q}$, we have $\omega_\pi(z)\ell=\ell\circ\pi(zI_n)=\chi_s(zI_n)\ell$. By \eqref{eq:nu-centre},
	$$
	\chi_s(zI_n)=\nu(zI_p)^{es}\,\nu(zI_q)^{-es}=|z|_{\mathrm{F}}^{e^2sp}\,|z|_{\mathrm{F}}^{-e^2sq}=|z|_{\mathrm{F}}^{e^2s(p-q)}.
	$$
	Since $\ell\ne0$, it follows that $\omega_\pi(z)=|z|_{\mathrm{F}}^{e^2s(p-q)}$.
\end{proof}

The next two lemmas describe the behaviour of distinction with respect to constituents and quotients.

\begin{lemma}\label{lem:dist-constituent}
	Let $\xi\in\Rep(G_n)$, $p+q=n$ and $s\in\R$. If $\xi$ is $(H_{p,q},\chi_s)$-distinguished, then so is some constituent of $\xi$.
\end{lemma}

\begin{proof}
	We argue by induction on the length of $\xi$. If $\xi$ is irreducible, there is nothing to prove. If $\xi$ is reducible, choose an exact sequence $0\to\xi'\to\xi\to\xi''\to0$ with $\xi',\xi''\ne0$. The sequence
	$$
	0\to\Hom_{H_{p,q}}(\xi'',\chi_s)\to\Hom_{H_{p,q}}(\xi,\chi_s)\to\Hom_{H_{p,q}}(\xi',\chi_s)
	$$
	is exact. Since the middle term is nonzero, one of $\xi'$ and $\xi''$ is $(H_{p,q},\chi_s)$-distinguished. Since $\xi'$ and $\xi''$ have smaller length and their constituents are constituents of $\xi$, the assertion follows by induction.
\end{proof}

\begin{lemma}\label{lem:dist-quotient}
	Let $\xi$ be a smooth representation of $G_n$, let $p+q=n$ and $s\in\R$, and let $\pi$ be a quotient of $\xi$. Then composition with the quotient map gives an injection $\Hom_{H_{p,q}}(\pi,\chi_s)\hookrightarrow\Hom_{H_{p,q}}(\xi,\chi_s)$. In particular, if $\pi$ is $(H_{p,q},\chi_s)$-distinguished, then so is $\xi$, and $\dim\Hom_{H_{p,q}}(\pi,\chi_s)\le\dim\Hom_{H_{p,q}}(\xi,\chi_s)$.
\end{lemma}

\begin{proof}
	A linear form on $\pi$ whose composition with a surjection $\xi\to\pi$ vanishes is itself zero.
\end{proof}

Finally, we determine which characters of $G_n$ are distinguished.

\begin{lemma}\label{lem:char-dist}
	Let $n\ge2$ and $s\in\R$, and let $\chi$ be a character of $G_n$. Then $\chi$ is $(H_{1,n-1},\chi_s)$-distinguished if and only if $s=0$ and $\chi=\one_n$.
\end{lemma}

\begin{proof}
	Write $\chi=\mu\circ\Nrd$ with $\mu$ a character of $\mathrm{F}^\times$ (\S\ref{ss:groups}). Since $\chi$ is one-dimensional, it is distinguished if and only if $\chi|_{H}=\chi_s$. Evaluating at $\diag(a,I_{n-1})$ with $a\in G_1$ and at $\diag(1,g)$ with $g\in G_{n-1}$, and using the surjectivity of $\Nrd$ on $G_1$ and on $G_{n-1}$, we see that this condition holds if and only if $\mu=|\cdot|_{\mathrm{F}}^{es}$ and $\mu=|\cdot|_{\mathrm{F}}^{-es}$. Hence, $|\cdot|_{\mathrm{F}}^{2es}=1$, so $s=0$ and $\mu=1$, that is, $\chi=\one_n$. Conversely, $\one_n$ is $(H_{1,n-1},\chi_0)$-distinguished, since $\chi_0=\one$.
\end{proof}

\subsection{The degenerate principal series}\label{ss:dps}

In this subsection, we study the two degenerate principal series that carry the obvious $(H,\chi_s)$-invariant linear forms. For $n\ge2$ and $s\in\R$, we introduce the following real numbers and segments on the line of $\one_1$ (Notation~\ref{not:line}):
$$
L_s=-s-\tfrac{n-1}2,\qquad R_s=-s+\tfrac{n-1}2,\qquad A_s=\Delta_{n-2}(-s)=[L_s+1,R_s-1],
$$
$$
A_s^-=[L_s,R_s-1]=\Delta_{n-1}(-s-\tfrac12),\qquad A_s^+=[L_s+1,R_s]=\Delta_{n-1}(-s+\tfrac12),
$$
where $A_s$ is empty if $n=2$. By Lemma~\ref{lem:characters-segments}, we have $\nu^{-es}=\Ss(A_s)$ as a character of $G_{n-2}$ (for $n\ge3$), $\nu^{-es\mp e/2}=\Ss(A_s^\mp)$ as characters of $G_{n-1}$, and $\nu^{-eL_s}=\Ss([-L_s])$ and $\nu^{-eR_s}=\Ss([-R_s])$ as characters of $G_1$. For $n\ge3$, the only singleton segments linked to $A_s$ are $[L_s]$ and $[R_s]$, and we have $A_s^-=[L_s]\cup A_s$ and $A_s^+=A_s\cup[R_s]$. We put
\begin{align*}
	\xione_s&=\nu^{-es-\frac e2}\times\nu^{es+e\frac{n-1}2}=\Ss(A_s^-)\times\Ss([-L_s])&&\text{(induced from }G_{n-1}\times G_1),\\
	\xitwo_s&=\nu^{es-e\frac{n-1}2}\times\nu^{-es+\frac e2}=\Ss([-R_s])\times\Ss(A_s^+)&&\text{(induced from }G_1\times G_{n-1}),
\end{align*}
and we write $\xi^{(1)\prime}_s=\Ss([-L_s])\times\Ss(A_s^-)$ and $\xi^{(2)\prime}_s=\Ss(A_s^+)\times\Ss([-R_s])$ for the same products with the factors in the opposite order.

Both representations admit explicit nonzero $(H,\chi_s)$-equivariant linear forms, as we now explain.

\begin{example}\label{ex:dps-functionals}
	By \eqref{eq:char-induced}, the representation $\xione_s$ is realized on the space of locally constant functions $f$ on $G_n$ such that
	$$
	f(\diag(g,a)\,u\,x)=\nu(g)^{-es}\nu(a)^{es}f(x)\qquad(g\in G_{n-1},\ a\in G_1,\ u\in N_{(n-1,1)},\ x\in G_n),
	$$
	since $-es-\frac e2+\frac e2=-es$ and $es+e\frac{n-1}2-e\frac{n-1}2=es$. Likewise, $\xitwo_s$ is realized on the space of locally constant functions $f$ such that
	$$
	f(h\,u\,x)=\chi_s(h)f(x)\qquad(h\in H=M_{(1,n-1)},\ u\in N_{(1,n-1)},\ x\in G_n).
	$$
	In both cases, the normalized induction is the unnormalized induction of the displayed character. This character is trivial exactly when $s=0$, and in that case the constant functions form a subrepresentation isomorphic to $\one_n$. Put $w=w_{1,n-1}$. We claim that the linear forms
	$$
	\ell^{(1)}(f)=f(w)\ \text{ on }\xione_s,\qquad \ell^{(2)}(f)=f(I_n)\ \text{ on }\xitwo_s
	$$
	are nonzero elements of $\Hom_H(\xi^{(i)}_s,\chi_s)$. They are nonzero, since an induced representation contains functions that do not vanish at a given point. For $\ell^{(2)}$, the equivariance is immediate from the displayed transformation rule. For $\ell^{(1)}$, let $h=\diag(a,g)\in H$ with $a\in G_1$ and $g\in G_{n-1}$. Then $whw^{-1}=\diag(g,a)$, and therefore
	$$
	\ell^{(1)}(hf)=f(wh)=f(\diag(g,a)w)=\nu(g)^{-es}\nu(a)^{es}f(w)=\chi_s(h)\ell^{(1)}(f).
	$$
\end{example}

We now determine when these representations are irreducible.

\begin{lemma}\label{lem:dps-irred}
	Let $n\ge2$ and $s\in\R$. Then $\xione_s$ is irreducible if and only if $s\notin\{-\frac n2,0\}$. In that case, $\xione_s\simeq\xi^{(1)\prime}_s$.
\end{lemma}

\begin{proof}
	By Lemma~\ref{lem:two-segments}, $\xione_s$ is irreducible if and only if the segments $A_s^-=[L_s,R_s-1]$ and $[-L_s]$ are not linked. These segments are linked if and only if $-L_s\in\{L_s-1,R_s\}$. Now, $-L_s=R_s$ means that $s=0$, and $-L_s=L_s-1$ means that $L_s=\frac12$, that is, $s=-\frac n2$. If $\xione_s$ is irreducible, then $\xione_s\simeq\xi^{(1)\prime}_s$ by the commutativity of the Grothendieck ring (\S\ref{ss:reps}).
\end{proof}

\begin{lemma}\label{lem:dps2-irred}
	Let $n\ge2$ and $s\in\R$. Then $\xitwo_s$ is irreducible if and only if $s\notin\{0,\frac n2\}$. In that case, $\xitwo_s\simeq\xi^{(2)\prime}_s$.
\end{lemma}

\begin{proof}
	By Lemma~\ref{lem:two-segments}, $\xitwo_s$ is irreducible if and only if the segments $A_s^+=[L_s+1,R_s]$ and $[-R_s]$ are not linked. These segments are linked if and only if $-R_s\in\{L_s,R_s+1\}$. Now, $-R_s=L_s$ means that $s=0$, and $-R_s=R_s+1$ means that $s=\frac n2$. The last assertion follows as in the proof of Lemma~\ref{lem:dps-irred}.
\end{proof}

At the exceptional values of $s$, the structure of these representations is as follows.

\begin{lemma}\label{lem:dps1-exceptional}
	Let $n\ge2$. Each of the following representations has length two, with the indicated unique irreducible subrepresentation and unique irreducible quotient:
	$$
	\one_n\hookrightarrow\xione_0\twoheadrightarrow\nu^{-e}Q_n,\qquad
	\nu^{\frac{en}2}\Qt_n\hookrightarrow\xione_{-n/2}\twoheadrightarrow\nu^{\frac{e(n-2)}2},
	$$
	$$
	\nu^{-e}Q_n\hookrightarrow\xi^{(1)\prime}_0\twoheadrightarrow\one_n,\qquad
	\nu^{\frac{e(n-2)}2}\hookrightarrow\xi^{(1)\prime}_{-n/2}\twoheadrightarrow\nu^{\frac{en}2}\Qt_n.
	$$
\end{lemma}

\begin{proof}
	We have $\xione_s=\nu^y\times\sigma$, where $y=-es-\frac e2$ and $\nu^y$ is regarded as a character of $G_{n-1}$, and $\sigma=\nu^{es+e\frac{n-1}2}$. Now, $y+\frac{en}2=-es+e\frac{n-1}2$, which equals the exponent of $\sigma$ exactly when $s=0$. Similarly, $y-\frac{en}2=-es-e\frac{n+1}2$, which equals the exponent of $\sigma$ exactly when $s=-\frac n2$. We apply Lemma~\ref{lem:char-sigma} with $m=n$. For $s=0$, the subrepresentation is $\nu^{y+e/2}=\one_n$ and the quotient is $\nu^{y-e/2}Q_n=\nu^{-e}Q_n$. For $s=-\frac n2$, we have $y=e\frac{n-1}2$; the subrepresentation is $\nu^{y+e/2}\Qt_n=\nu^{en/2}\Qt_n$ and the quotient is $\nu^{y-e/2}=\nu^{e(n-2)/2}$. This proves the first line. Applying Lemma~\ref{lem:sigma-char} to $\xi^{(1)\prime}_s=\sigma\times\nu^y$ gives the second line.
\end{proof}

\begin{lemma}\label{lem:dps2-exceptional}
	Let $n\ge2$. Each of the following representations has length two, with the indicated unique irreducible subrepresentation and unique irreducible quotient:
	$$
	\one_n\hookrightarrow\xitwo_0\twoheadrightarrow\nu^{e}\Qt_n,\qquad
	\nu^{-\frac{en}2}Q_n\hookrightarrow\xitwo_{n/2}\twoheadrightarrow\nu^{\frac{e(2-n)}2},
	$$
	$$
	\nu^{e}\Qt_n\hookrightarrow\xi^{(2)\prime}_0\twoheadrightarrow\one_n,\qquad
	\nu^{\frac{e(2-n)}2}\hookrightarrow\xi^{(2)\prime}_{n/2}\twoheadrightarrow\nu^{-\frac{en}2}Q_n.
	$$
\end{lemma}

\begin{proof}
	We have $\xitwo_s=\sigma\times\nu^y$, where $\sigma=\nu^{es-e\frac{n-1}2}$, $y=-es+\frac e2$, and $\nu^y$ is regarded as a character of $G_{n-1}$. The exponent of $\sigma$ equals $y-\frac{en}2$ exactly when $s=0$, and it equals $y+\frac{en}2$ exactly when $s=\frac n2$. Hence, Lemma~\ref{lem:sigma-char} with $m=n$ gives the first line, and Lemma~\ref{lem:char-sigma} applied to $\xi^{(2)\prime}_s=\nu^y\times\sigma$ gives the second line.
\end{proof}

For later reference, we single out the subrepresentations at the exceptional values.

\begin{remark}\label{rem:exceptional-sub}
	By Lemmas~\ref{lem:dps1-exceptional} and~\ref{lem:dps2-exceptional}, at $s=0$ the unique irreducible subrepresentations of $\xione_0$ and $\xitwo_0$ are both equal to $\one_n$. At $s=-\frac n2$, the unique irreducible subrepresentation of $\xione_s$ is $\nu^{en/2}\Qt_n=\nu^{-es}\Qt_n$, and at $s=\frac n2$, that of $\xitwo_s$ is $\nu^{-en/2}Q_n=\nu^{-es}Q_n$. By Lemma~\ref{lem:Qm-not-char}, the last two representations are not characters. For $n=2$, \eqref{eq:Q2} shows that $\nu^{e}\Qt_2=\nu^{-e}Q_2=\St_2$.
\end{remark}

We now show that these two subrepresentations are distinguished exactly at the corresponding exceptional value of $s$.

\begin{lemma}\label{lem:Qtilde-dist}
	Let $n\ge2$ and $s\in\R$. Then $\nu^{-es}\Qt_n$ is $(H_{1,n-1},\chi_s)$-distinguished if and only if $s=-\frac n2$.
\end{lemma}

\begin{proof}
	By \eqref{eq:Qm-sequence} and \eqref{eq:nu-centre}, we have
	$$
	\omega_{Q_n}(z)=\nu(zI_{n-1})^{e/2}\,\nu(z)^{e\frac{n+1}2}=|z|_{\mathrm{F}}^{e^2n}\qquad(z\in\mathrm{F}^\times).
	$$
	Hence, $\omega_{\Qt_n}=|\cdot|_{\mathrm{F}}^{-e^2n}$ and $\omega_{\nu^{-es}\Qt_n}=|\cdot|_{\mathrm{F}}^{-e^2n(s+1)}$. By Lemma~\ref{lem:central}, distinction forces $-n(s+1)=s(2-n)$, that is, $s=-\frac n2$.
	
	Conversely, let $s=-\frac n2$. By Example~\ref{ex:dps-functionals}, there is a nonzero $\ell\in\Hom_H(\xione_s,\chi_s)$. By Lemma~\ref{lem:dps1-exceptional}, there is an exact sequence
	$$
	0\to\nu^{en/2}\Qt_n\to\xione_s\to\nu^{e(n-2)/2}\to0.
	$$
	If $\ell$ vanished on $\nu^{en/2}\Qt_n$, then it would factor through the character $\nu^{e(n-2)/2}$, which would therefore be $(H,\chi_s)$-distinguished. This contradicts Lemma~\ref{lem:char-dist}, since $s\ne0$. Hence, the restriction of $\ell$ to $\nu^{en/2}\Qt_n=\nu^{-es}\Qt_n$ is nonzero.
\end{proof}

\begin{lemma}\label{lem:Q-dist}
	Let $n\ge2$ and $s\in\R$. Then $\nu^{-es}Q_n$ is $(H_{1,n-1},\chi_s)$-distinguished if and only if $s=\frac n2$.
\end{lemma}

\begin{proof}
	By the computation in the proof of Lemma~\ref{lem:Qtilde-dist}, we have $\omega_{\nu^{-es}Q_n}=|\cdot|_{\mathrm{F}}^{e^2n(1-s)}$. Hence, Lemma~\ref{lem:central} forces $n(1-s)=s(2-n)$, that is, $s=\frac n2$. Conversely, for $s=\frac n2$, we use the linear form $\ell^{(2)}$ of Example~\ref{ex:dps-functionals}, the exact sequence $0\to\nu^{-en/2}Q_n\to\xitwo_{n/2}\to\nu^{e(2-n)/2}\to0$ of Lemma~\ref{lem:dps2-exceptional}, and Lemma~\ref{lem:char-dist}, exactly as in the proof of Lemma~\ref{lem:Qtilde-dist}.
\end{proof}

\section{The orbit structure and Mackey theory for maximal parabolic induction}\label{sec:orbits}

In this section, we consider a representation $\pi_1\times\pi_2$ induced from a maximal parabolic subgroup of $G_n$. We bound the space $\Hom_{H_{1,n-1}}(\pi_1\times\pi_2,\chi_s)$ in terms of three explicit contributions, one for each $H_{1,n-1}$-orbit on the corresponding Grassmannian. Since this geometric input governs everything that follows, we give full details.

We first describe the orbits on the Grassmannian. Then, in \S\S\ref{ss:offen}--\ref{ss:contributions}, we compute for each orbit the data that enter Offen's geometric lemma \cite[Proposition 4.1]{Off17}, namely:
\begin{itemize}
	\item the representative $x$;
	\item the Levi subgroup $L$ and the parabolic subgroup $Q=LV\subseteq P$ attached to $x$;
	\item the stabilizers $P_x=Q_x$;
	\item the character $\delta_{Q_x}\delta_Q^{-1/2}\chi^{\eta^{-1}}$.
\end{itemize}

Throughout this section, we fix $n\ge2$ and an integer $k$ with $1\le k\le n-1$, and we put
$$
H=H_{1,n-1},\qquad P=P_{(k,n-k)}=MN,\qquad M=M_{(k,n-k)}\simeq G_k\times G_{n-k}.
$$
The group $P$ is the stabilizer in $G_n$ of the subspace $e_1\mathrm{D}\oplus\dots\oplus e_k\mathrm{D}$ of $V=\mathrm{D}^n$. Hence, $G_n/P$ is identified with the Grassmannian of $k$-dimensional $\mathrm{D}$-subspaces of $V$, and the double cosets $P\backslash G_n/H$ correspond to the $H$-orbits on this Grassmannian.

We use the convention of \S\ref{ss:Hpq} for the characters $\chi_c$ of the degenerate groups $H_{1,0}=G_1$ and $H_{0,m}=G_m$, and we omit factors $G_0$. Moreover, the Jacquet functors $r_{(0,m)}$ and $r_{(m,0)}$ are the identity functors. With these conventions, all statements below also hold for $k=1$ and $k=n-1$, where some blocks are empty. 

\subsection{The Grassmannian and the invariants of a subspace}\label{ss:grassmannian}

We begin by describing $P\backslash G_n$ geometrically. Let $V=\mathrm{D}^n$ with its standard basis $e_1,\dots,e_n$ (\S\ref{ss:groups}), and put
$$
V_+=e_1\mathrm{D},\qquad V_-=e_2\mathrm{D}\oplus\dots\oplus e_n\mathrm{D},\qquad V_j=e_1\mathrm{D}\oplus\dots\oplus e_j\mathrm{D}\quad(0\le j\le n).
$$
Then $V=V_+\oplus V_-$, and
$$
H=\{g\in G_n:\ gV_+=V_+,\ gV_-=V_-\},\qquad P=\{g\in G_n:\ gV_k=V_k\}.
$$
An element $h=\diag(a,g)\in H$, with $a\in G_1$ and $g\in G_{n-1}$, acts on $V_+$ by left multiplication by $a$ on the first coordinate, and on $V_-\simeq \mathrm{D}^{n-1}$ by $g$. In this way, we identify $H$ with $\GL_{\mathrm{D}}(V_+)\times\GL_{\mathrm{D}}(V_-)$. Let $\Gr_k(V)$ denote the set of $k$-dimensional right $\mathrm{D}$-subspaces of $V$.

The following lemma identifies $P\backslash G_n$ with $\Gr_k(V)$.

\begin{lemma}\label{lem:grass}
	The group $G_n$ acts transitively on $\Gr_k(V)$, and the stabilizer of $V_k$ is $P$. The map
	$$
	P\backslash G_n\longrightarrow\Gr_k(V),\qquad Pg\longmapsto g^{-1}V_k,
	$$
	is a bijection. It transforms the action of $H$ on $P\backslash G_n$ by right translation into the action $(W,h)\mapsto h^{-1}W$ of $H$ on $\Gr_k(V)$. Consequently, the $(P,H)$-double cosets in $G_n$ correspond bijectively to the $H$-orbits on $\Gr_k(V)$. Endowed with the topology transported from $P\backslash G_n$, the set $\Gr_k(V)$ is a compact Hausdorff space.
\end{lemma}

\begin{proof}
	Every $k$-dimensional subspace has a basis that can be completed to a basis of $V$, and this gives transitivity. The stabilizer of $V_k$ is $P$ by definition. The map is well defined, since $(pg)^{-1}V_k=g^{-1}p^{-1}V_k=g^{-1}V_k$ for $p\in P$. It is surjective by transitivity, and it is injective because $g^{-1}V_k=g'^{-1}V_k$ implies $g'g^{-1}\in P$. Moreover, $Pgh\mapsto h^{-1}g^{-1}V_k$ for $h\in H$, which proves the assertion about the action of $H$. Finally, $P$ is closed in $G_n$, so $P\backslash G_n$ is Hausdorff, and $P\backslash G_n$ is compact since $G_n=PK$ for a maximal compact subgroup $K$ of $G_n$. The last assertion follows.
\end{proof}

To classify the $H$-orbits, we attach two numerical invariants to a subspace.

\begin{definition}\label{def:invariants}
	For $W\in\Gr_k(V)$, we put
	$$
	r(W)=\dim_{\mathrm{D}}(W\cap V_+),\qquad t(W)=\dim_{\mathrm{D}}(W\cap V_-),
	$$
	and we let $\pi_+\colon V\to V_+$ denote the projection along $V_-$.
\end{definition}

Since $H$ preserves $V_+$ and $V_-$, the pair $(r(W),t(W))$ depends only on the $H$-orbit of $W$. Only three values of this pair are possible.

\begin{lemma}\label{lem:invariants}
	Let $W\in\Gr_k(V)$. Then $r(W)\le1$ and $t(W)\in\{k-1,k\}$. Moreover, if $t(W)=k$, then $r(W)=0$. Consequently,
	$$
	(r(W),t(W))\in\{(0,k),\ (1,k-1),\ (0,k-1)\}.
	$$
\end{lemma}

\begin{proof}
	First, $r(W)\le\dim V_+=1$. The kernel of $\pi_+|_W$ is $W\cap V_-$, and its image lies in $V_+$, which has dimension one. Hence, $k-t(W)\le1$. If $t(W)=k$, then $W\subseteq V_-$, and therefore $W\cap V_+\subseteq V_-\cap V_+=0$.
\end{proof}

Accordingly, we divide the subspaces into three types.

\begin{definition}\label{def:types}
	We say that $W\in\Gr_k(V)$ is of
	\begin{itemize}
		\item \emph{type I} if $(r(W),t(W))=(0,k)$, that is, if $W\subseteq V_-$;
		\item \emph{type II} if $(r(W),t(W))=(1,k-1)$, that is, if $V_+\subseteq W$;
		\item \emph{type III} if $(r(W),t(W))=(0,k-1)$, that is, if $W\not\subseteq V_-$ and $V_+\not\subseteq W$.
	\end{itemize}
\end{definition}

The equivalences in Definition~\ref{def:types} hold because $t(W)=k$ means $W\subseteq V_-$, and $r(W)=1$ means $V_+\subseteq W$. By Lemma~\ref{lem:invariants}, $r(W)=1$ forces $t(W)=k-1$.

We now give normal forms for subspaces of each type.

\begin{lemma}\label{lem:normal-forms}
	Let $W\in\Gr_k(V)$.
	\begin{enumerate}
		\item If $W$ is of type I, then $W$ is a $k$-dimensional subspace of $V_-$.
		\item If $W$ is of type II, then $W=V_+\oplus(W\cap V_-)$ with $W\cap V_-\in\Gr_{k-1}(V_-)$. Conversely, $V_+\oplus W'$ is of type II for every $W'\in\Gr_{k-1}(V_-)$.
		\item If $W$ is of type III, put $W'=W\cap V_-$. Then $\dim W'=k-1$, and there is $v\in V_-$, unique modulo $W'$, such that $e_1+v\in W$. Moreover, $v\notin W'$ and
		$$
		W=W'\oplus(e_1+v)\mathrm{D}.
		$$
		Conversely, for every $W'\in\Gr_{k-1}(V_-)$ and every $v\in V_-\smallsetminus W'$, the subspace $W'\oplus(e_1+v)\mathrm{D}$ is of type III.
	\end{enumerate}
\end{lemma}

\begin{proof}
	(1) This is the definition.
	
	(2) If $V_+\subseteq W$, then $W=V_+\oplus(W\cap V_-)$, since every $w\in W$ can be written as $w=\pi_+(w)+(w-\pi_+(w))$ with $\pi_+(w)\in V_+\subseteq W$. The converse is clear.
	
	(3) Since $t(W)=k-1$, the image of $\pi_+|_W$ has dimension one, so it equals $V_+$. Hence, there is $w\in W$ with $\pi_+(w)=e_1$, that is, $w=e_1+v$ with $v\in V_-$. Counting dimensions, we obtain $W=W'\oplus w\mathrm{D}$. If $v$ belonged to $W'$, then $e_1=w-v$ would lie in $W$, contradicting $r(W)=0$. If also $e_1+v'\in W$ with $v'\in V_-$, then $v-v'\in W\cap V_-=W'$, which proves uniqueness modulo $W'$.
	
	Conversely, let $W=W'\oplus(e_1+v)\mathrm{D}$ with $v\in V_-\smallsetminus W'$. If $w'+(e_1+v)c\in V_+$ with $w'\in W'$ and $c\in\mathrm{D}$, then $w'+vc=0$. Hence, $vc\in W'$, and since $v\notin W'$, we get $c=0$ and then $w'=0$. Thus, $r(W)=0$. If $w'+(e_1+v)c\in V_-$, then $c=0$. Thus, $W\cap V_-=W'$ and $t(W)=k-1$.
\end{proof}

We can now describe the $H$-orbits on $\Gr_k(V)$.

\begin{proposition}\label{prop:three-orbits}
	The group $H$ has exactly three orbits on $\Gr_k(V)$, namely the sets
	$$
	\OI=\{W:W\subseteq V_-\},\qquad\OII=\{W:V_+\subseteq W\},\qquad\OIII=\{W:W\not\subseteq V_-,\ V_+\not\subseteq W\}
	$$
	of subspaces of types I, II and III, respectively. Representatives of these orbits are
	$$
	\WI=e_2\mathrm{D}\oplus\dots\oplus e_{k+1}\mathrm{D},\qquad\WII=V_k=e_1\mathrm{D}\oplus\dots\oplus e_k\mathrm{D},\qquad\WIII=(e_1+e_{k+1})\mathrm{D}\oplus e_2\mathrm{D}\oplus\dots\oplus e_k\mathrm{D}.
	$$
\end{proposition}

\begin{proof}
	By Lemma~\ref{lem:invariants}, every subspace is of exactly one type, and the type is $H$-invariant. It remains to show that each type forms a single orbit and that $W_*$ is of type $*$. Recall that $\GL_{\mathrm{D}}(V_-)$ acts transitively on $\Gr_j(V_-)$ for every $0\le j\le n-1$.
	
	\emph{Type I.} We have $\OI=\Gr_k(V_-)$ and $k\le n-1=\dim V_-$. Hence, the subgroup $\{1\}\times\GL_{\mathrm{D}}(V_-)$ of $H$ acts transitively on $\OI$. Clearly, $\WI\subseteq V_-$.
	
	\emph{Type II.} By Lemma~\ref{lem:normal-forms}(2), the map $W\mapsto W\cap V_-$ is an $H$-equivariant bijection from $\OII$ onto $\Gr_{k-1}(V_-)$, on which $\GL_{\mathrm{D}}(V_-)$ acts transitively. Clearly, $V_+\subseteq\WII$.
	
	\emph{Type III.} Put $W'_0=e_2\mathrm{D}\oplus\dots\oplus e_k\mathrm{D}$ and $v_0=e_{k+1}$, which makes sense since $k+1\le n$. Then $\WIII=W'_0\oplus(e_1+v_0)\mathrm{D}$ with $v_0\notin W'_0$, so $\WIII$ is of type III by Lemma~\ref{lem:normal-forms}(3). Now let $W=W'\oplus(e_1+v)\mathrm{D}$ be of type III, as in Lemma~\ref{lem:normal-forms}(3). Since $v\notin W'$, a basis of $W'$ together with $v$ is linearly independent and can be completed to a basis of $V_-$. Hence, there is $g\in\GL_{\mathrm{D}}(V_-)$ with $gW'=W'_0$ and $gv=v_0$. Then $h=\diag(1,g)\in H$ satisfies $h(e_1+v)=e_1+v_0$, and therefore $hW=\WIII$.
\end{proof}

Next, we determine the topology of the three orbits.

\begin{lemma}\label{lem:orbit-topology}
	The orbits $\OI$ and $\OII$ are compact, hence closed, and $\OIII$ is open in $\Gr_k(V)$. Moreover, $\OI\simeq\Gr_k(V_-)$ and $\OII\simeq\Gr_{k-1}(V_-)$ as $H$-spaces, and $\OI\cap\OII=\emptyset$.
\end{lemma}

\begin{proof}
	Equip $\Gr_k(V_-)$ and $\Gr_{k-1}(V_-)$ with the topologies given by Lemma~\ref{lem:grass}, applied to $G_{n-1}=\GL_{\mathrm{D}}(V_-)$. The maps
	$$
	\Gr_k(V_-)\to\Gr_k(V),\quad W\mapsto W,\qquad\text{and}\qquad\Gr_{k-1}(V_-)\to\Gr_k(V),\quad W'\mapsto V_+\oplus W',
	$$
	are continuous, since they are induced by the continuous homomorphism $g\mapsto\diag(1,g)$ from $G_{n-1}$ to $G_n$. For instance, $\diag(1,g)\WII=V_+\oplus g(e_2\mathrm{D}\oplus\dots\oplus e_k\mathrm{D})$. Their images are $\OI$ and $\OII$, respectively, and their sources are compact. Since $\Gr_k(V)$ is Hausdorff, these continuous injections are homeomorphisms onto their images, and $\OI$ and $\OII$ are closed. Consequently, $\OIII$, the complement of $\OI\cup\OII$, is open. Finally, the two closed orbits are disjoint because they consist of subspaces of different types.
\end{proof}

We illustrate Proposition~\ref{prop:three-orbits} in small rank.

\begin{example}\label{ex:orbits-small}
	(1) Let $n=2$ and $k=1$. Then $\Gr_1(\mathrm{D}^2)$ is the projective line over $\mathrm{D}$, and
	$$
	\OI=\{e_2\mathrm{D}\},\qquad \OII=\{e_1\mathrm{D}\},\qquad \OIII=\{(e_1+e_2x)\mathrm{D}:x\in \mathrm{D}^\times\}.
	$$
	The element $h=\diag(a,b)\in H_{1,1}$ maps $(e_1+e_2x)\mathrm{D}$ to $(e_1+e_2bxa^{-1})\mathrm{D}$. Thus, $\OIII\simeq \mathrm{D}^\times$, with $H_{1,1}=\mathrm{D}^\times\times \mathrm{D}^\times$ acting by $x\mapsto bxa^{-1}$, and the stabilizer of $x=1$ is the diagonal subgroup $\{\diag(a,a):a\in\mathrm{D}^\times\}$. These are the three orbits of the diagonal Levi subgroup $H_{1,1}$ of $G_2$ on $P\backslash G_2$, namely the two points $P$ and $Pw$ and the open orbit, which underlie the theory of Section~\ref{sec:rank-two}.
	
	(2) Let $n=3$. For $k=1$, the orbit $\OI$ is the set of lines in $V_-$ (a projective line over $\mathrm{D}$), $\OII=\{V_+\}$ is a point, and $\OIII=\{(e_1+v)\mathrm{D}:0\ne v\in V_-\}$. For $k=2$, the orbit $\OI=\{V_-\}$ is a point, $\OII$ is the set of planes containing $V_+$ (a projective line over $\mathrm{D}$, via $W\mapsto W\cap V_-$), and $\OIII$ consists of the planes $W'\oplus(e_1+v)\mathrm{D}$ with $W'$ a line in $V_-$ and $v\in V_-\smallsetminus W'$.
\end{example}

\subsection{The symmetric space and the geometric lemma of Offen}\label{ss:offen}

We now pass from the Grassmannian to the framework of Offen \cite{Off17}, in which the Mackey theory of an induced representation restricted to a symmetric subgroup is formulated. Put
$$
\varepsilon=\diag(1,-I_{n-1}),\qquad\theta=\Ad(\varepsilon)\colon g\mapsto\varepsilon g\varepsilon^{-1}.
$$
Then $\theta$ is an involution of $G_n$, and $H=G_n^\theta$ (\S\ref{ss:Hpq}). Following \cite[\S2.1]{Off17}, we let
$$
X=\{g\in G_n:\ g\theta(g)=I_n\},\qquad g\cdot x=gx\theta(g)^{-1}\quad(g\in G_n,\ x\in X),
$$
and, for a subgroup $S$ of $G_n$ and $x\in X$, we denote by $S_x$ the stabilizer of $x$ in $S$. The map $g\mapsto g\cdot I_n=g\theta(g)^{-1}$ identifies $G_n/H$ with the orbit $G_n\cdot I_n\subseteq X$. Since $\theta(g)=\varepsilon g\varepsilon^{-1}$ and $\varepsilon^2=I_n$, we have, for $x\in X$ and $g\in G_n$,
\begin{equation}\label{eq:x-epsilon}
	(x\varepsilon)^2=x\theta(x)=I_n,\qquad (g\cdot x)\varepsilon=g\,(x\varepsilon)\,g^{-1}.
\end{equation}
Hence, $x\varepsilon$ is a $\mathrm{D}$-linear involution of $V$, and $S_x$ is the centralizer of $x\varepsilon$ in $S$, that is, the subgroup of $S$ preserving the two eigenspaces of $x\varepsilon$. If $x=\eta\cdot I_n$ with $\eta\in G_n$, then $x\varepsilon=\eta\varepsilon\eta^{-1}$ has $(+1)$-eigenspace $\eta V_+$ and $(-1)$-eigenspace $\eta V_-$, and therefore
\begin{equation}\label{eq:Px}
	P_x=P\cap\eta H\eta^{-1}=\{p\in P:\ p(\eta V_+)=\eta V_+,\ p(\eta V_-)=\eta V_-\}.
\end{equation}
The map $P\eta H\mapsto P\cdot(\eta\cdot I_n)$ is a bijection from $P\backslash G_n/H$ onto the set of $P$-orbits in $G_n\cdot I_n$. By Lemma~\ref{lem:grass}, these sets are also in bijection with the set of $H$-orbits on $\Gr_k(V)$, via the map sending $P\eta H$ to the $H$-orbit of $\eta^{-1}V_k$. Thus, the $P$-orbit of $x=\eta\cdot I_n$ is determined by the type of $\eta^{-1}V_k$ (Definition~\ref{def:types}). Equivalently, it is determined by the relative position of $V_k$ with respect to the decomposition $V=\eta V_+\oplus\eta V_-$.

\smallskip
\emph{Weyl group data.} Let $T=\{\diag(t_1,\dots,t_n):t_i\in \mathrm{F}^\times\}$, which is a maximal $\mathrm{F}$-split torus of $G_n$, and let $P_0$ be the subgroup of upper triangular matrices in $G_n$. Then $M_0=C_{G_n}(T)$ is the group of diagonal matrices with entries in $\mathrm{D}^\times$, and $P_0=M_0U_0$ is a minimal parabolic subgroup containing $T$; the standard parabolic subgroups are the subgroups $P_\alpha$ of \S\ref{ss:reps}. The Weyl group $W=N_{G_n}(T)/M_0$ is the symmetric group $S_n$, realized by permutation matrices. Since $\varepsilon\in M_0$, the involution $\theta$ preserves $T$, $M_0$, $P_0$ and every standard parabolic subgroup, and it acts trivially on $W$. Consequently, in the notation of \cite[\S\S2.3--3]{Off17}, the Weyl element $\tau$ with $\theta(P_0)=\tau P_0\tau^{-1}$ is trivial, $\theta'=\theta$, and the twisted involutions in $W$ are exactly the involutions of $S_n$.

Let $\mathcal O$ be a $P$-orbit in $X$. Following \cite[\S3]{Off17}, we attach to $\mathcal O$ the following data:
\begin{itemize}
	\item the Weyl element $w=\iota_M(\mathcal O)$, which is characterized by the conditions that $PxP=PwP$ for $x\in\mathcal O$ and that $w$ has minimal length in $W_MwW_M$;
	\item the standard Levi subgroup
	$$
	L=M\cap w\theta(M)w^{-1}=M\cap wMw^{-1};
	$$
	\item the standard parabolic subgroup $Q=LU_Q$ of $G_n$ with Levi subgroup $L$ and unipotent radical $U_Q$. Since $L\subseteq M$, we have $Q\subseteq P$.
\end{itemize}
By \cite[Lemma 3.2]{Off17}, the representative $x$ of $\mathcal O$ can be chosen in $Lw$, and then $P_x=L_x\ltimes Z_x$, where $Z_x$ is the unipotent radical of $P_x$. By \cite[Lemma 3.3]{Off17}, the projection of $Z_x$ to $M$ is $M\cap U_Q$. Finally, as in the proof of \cite[Proposition 4.1]{Off17}, we have $P_x=Q_x$. In \S\ref{ss:representatives}, we choose such representatives explicitly, and in \S\ref{ss:stabilizers}, we verify all these properties directly.

\subsection{Representatives in \texorpdfstring{$G_n$}{Gn}}\label{ss:representatives}

We now choose explicit representatives of the three double cosets. We define $\eta_{\mathrm I},\eta_{\mathrm{II}},\eta_{\mathrm{III}}\in G_n$ by their action on the standard basis:
\begin{align*}
	\eta_{\mathrm{II}}&=I_n;\\
	\eta_{\mathrm I}e_1&=e_{k+1},\quad \eta_{\mathrm I}e_i=e_{i-1}\ (2\le i\le k+1),\quad \eta_{\mathrm I}e_i=e_i\ (i\ge k+2);\\
	\eta_{\mathrm{III}}e_1&=e_k+e_{k+1},\quad \eta_{\mathrm{III}}e_i=e_{i-1}\ (2\le i\le k),\quad \eta_{\mathrm{III}}e_{k+1}=e_k-e_{k+1},\quad \eta_{\mathrm{III}}e_i=e_i\ (i\ge k+2).
\end{align*}
Since $2\in \mathrm{F}^\times$, the vectors $e_k+e_{k+1}$ and $e_k-e_{k+1}$ span $e_k\mathrm{D}\oplus e_{k+1}\mathrm{D}$, and hence $\eta_{\mathrm{III}}$ is invertible. For $*\in\{\mathrm I,\mathrm{II},\mathrm{III}\}$, we put
$$
x_*=\eta_*\cdot I_n=\eta_*\theta(\eta_*)^{-1}\in G_n\cdot I_n.
$$

The following lemma shows that these elements represent the three orbits of Proposition~\ref{prop:three-orbits}.

\begin{lemma}\label{lem:representatives}
	For $*\in\{\mathrm I,\mathrm{II},\mathrm{III}\}$, we have $\eta_*^{-1}V_k=W_*$. Moreover, the following hold.
	\begin{enumerate}
		\item $\eta_{\mathrm I}V_+=e_{k+1}\mathrm{D}$ and $\eta_{\mathrm I}V_-=V_k\oplus e_{k+2}\mathrm{D}\oplus\dots\oplus e_n\mathrm{D}$.
		\item $\eta_{\mathrm{II}}V_+=V_+$ and $\eta_{\mathrm{II}}V_-=V_-$.
		\item $\eta_{\mathrm{III}}V_+=(e_k+e_{k+1})\mathrm{D}$ and $\eta_{\mathrm{III}}V_-=V_{k-1}\oplus(e_k-e_{k+1})\mathrm{D}\oplus e_{k+2}\mathrm{D}\oplus\dots\oplus e_n\mathrm{D}$.
	\end{enumerate}
	Consequently, $G_n=P\eta_{\mathrm I}H\sqcup P\eta_{\mathrm{II}}H\sqcup P\eta_{\mathrm{III}}H$, the double cosets $P\eta_{\mathrm I}H$ and $P\eta_{\mathrm{II}}H$ are closed in $G_n$, and $P\eta_{\mathrm{III}}H$ is open in $G_n$. Moreover, $P\cdot x_{\mathrm I}$, $P\cdot x_{\mathrm{II}}$ and $P\cdot x_{\mathrm{III}}$ are the three $P$-orbits in $G_n\cdot I_n$.
\end{lemma}

\begin{proof}
	For $\eta_{\mathrm I}$, we have $e_j=\eta_{\mathrm I}e_{j+1}$ for $1\le j\le k$, and therefore $\eta_{\mathrm I}^{-1}V_k=e_2\mathrm{D}\oplus\dots\oplus e_{k+1}\mathrm{D}=\WI$. For $\eta_{\mathrm{III}}$, we have $e_j=\eta_{\mathrm{III}}e_{j+1}$ for $1\le j\le k-1$ and $e_k=\eta_{\mathrm{III}}\big((e_1+e_{k+1})\tfrac12\big)$. Therefore, $\eta_{\mathrm{III}}^{-1}V_k=e_2\mathrm{D}\oplus\dots\oplus e_k\mathrm{D}\oplus(e_1+e_{k+1})\mathrm{D}=\WIII$. The case of $\eta_{\mathrm{II}}=I_n$ is trivial, and the descriptions of $\eta_*V_\pm$ follow directly from the definitions.
	
	Under the bijection of Lemma~\ref{lem:grass}, the image of $P\eta_*H$ in $P\backslash G_n$ corresponds to the $H$-orbit of $\eta_*^{-1}V_k=W_*$. Hence, the decomposition of $G_n$ follows from Proposition~\ref{prop:three-orbits}. Since the projection $G_n\to P\backslash G_n$ is continuous and open, and the double cosets are unions of fibres, a double coset is closed (respectively open) in $G_n$ if and only if its image is closed (respectively open) in $P\backslash G_n$. Thus, the topological assertions follow from Lemma~\ref{lem:orbit-topology}. The last assertion follows from the bijection between $P\backslash G_n/H$ and the set of $P$-orbits in $G_n\cdot I_n$ described in \S\ref{ss:offen}.
\end{proof}

We next compute the involutions $x_*\varepsilon$ explicitly. By \eqref{eq:x-epsilon} and Lemma~\ref{lem:representatives}, we have
\begin{equation}\label{eq:x-epsilon-explicit}
	x_{\mathrm I}\varepsilon=\diag(-I_k,1,-I_{n-k-1}),\qquad x_{\mathrm{II}}\varepsilon=\varepsilon,\qquad x_{\mathrm{III}}\varepsilon=\diag\Big(-I_{k-1},\begin{pmatrix}0&1\\1&0\end{pmatrix},-I_{n-k-1}\Big),
\end{equation}
in blocks of sizes $(k,1,n-k-1)$, $(1,n-1)$ and $(k-1,2,n-k-1)$, respectively. Indeed, $x_*\varepsilon=\eta_*\varepsilon\eta_*^{-1}$ is the involution with $(+1)$-eigenspace $\eta_*V_+$ and $(-1)$-eigenspace $\eta_*V_-$, and $(e_k+e_{k+1})\mathrm{D}$ and $(e_k-e_{k+1})\mathrm{D}$ are the $(+1)$- and $(-1)$-eigenspaces of $\left(\begin{smallmatrix}0&1\\1&0\end{smallmatrix}\right)$ on $e_k\mathrm{D}\oplus e_{k+1}\mathrm{D}$. In particular, multiplying on the right by $\varepsilon$, we obtain
$$
x_{\mathrm I}=\diag(-1,I_{k-1},-1,I_{n-k-1}),\qquad x_{\mathrm{II}}=I_n.
$$
Let $s_k\in W$ denote the transposition of $k$ and $k+1$, realized by the corresponding permutation matrix.

With these representatives, the data of Offen's geometric lemma can be read off directly.

\begin{lemma}\label{lem:offen-data}
	The data of \cite[\S3]{Off17} attached to the three orbits are as follows.
	\begin{enumerate}
		\item For $\OI$ and $\OII$, we have $w=1$, $L=M$, $Q=P$ and $U_Q=N$. Moreover, $x_{\mathrm I},x_{\mathrm{II}}\in M_0\subseteq Lw$.
		\item For $\OIII$, we have $w=s_k$, $L=M_{(k-1,1,1,n-k-1)}$ and $Q=P_{(k-1,1,1,n-k-1)}$, and $U_Q=N_Q$ is the unipotent radical of $Q$. Moreover, $x_{\mathrm{III}}\in M_0s_k\subseteq Lw$.
	\end{enumerate}
	In all three cases, $M\cap Q=L\ltimes(M\cap U_Q)$, where $M\cap U_Q=\{1\}$ in case \textup{(1)} and $M\cap U_Q=N_{(k-1,1)}\times N_{(1,n-k-1)}$ in case \textup{(2)}.
\end{lemma}

\begin{proof}
	(1) By \eqref{eq:x-epsilon-explicit}, the elements $x_{\mathrm I}=(x_{\mathrm I}\varepsilon)\varepsilon$ and $x_{\mathrm{II}}=I_n$ are diagonal with entries $\pm1$. Hence, they lie in $M_0\subseteq P$, and $Px_*\theta(P)=Px_*P=P$ for $*\in\{\mathrm I,\mathrm{II}\}$. Therefore, $\iota_M(P\cdot x_*)=1$, $L=M$, $Q=P$ and $U_Q=N$. In this case, $M\cap U_Q=M\cap N=\{1\}$.
	
	(2) The element $x_{\mathrm{III}}=(x_{\mathrm{III}}\varepsilon)\varepsilon$ is the product of the permutation matrix of $s_k$ with a diagonal matrix with entries $\pm1$. Hence, $x_{\mathrm{III}}\in M_0s_k$ and $Px_{\mathrm{III}}P=Ps_kP$. The transposition $s_k$ is the only simple reflection not in $W_M=S_k\times S_{n-k}$. Thus, the double coset $W_Ms_kW_M$ does not contain $1$, and $s_k$, which has length one, has minimal length in it. Consequently, $\iota_M(P\cdot x_{\mathrm{III}})=s_k$. Conjugation by $s_k$ interchanges $e_k$ and $e_{k+1}$, so $s_kMs_k$ is the stabilizer of the decomposition
	$$
	V=(V_{k-1}\oplus e_{k+1}\mathrm{D})\oplus(e_k\mathrm{D}\oplus e_{k+2}\mathrm{D}\oplus\dots\oplus e_n\mathrm{D}).
	$$
	An element of $L=M\cap s_kMs_k$ preserves the four pairwise intersections $V_{k-1}$, $e_k\mathrm{D}$, $e_{k+1}\mathrm{D}$ and $e_{k+2}\mathrm{D}\oplus\dots\oplus e_n\mathrm{D}$, whose direct sum is $V$. Hence, $L$ consists of the matrices $\diag(g_1,a,b,g_3)$ with $g_1\in G_{k-1}$, $a,b\in G_1$ and $g_3\in G_{n-k-1}$, that is, $L=M_{(k-1,1,1,n-k-1)}$, and $Q$ is the standard parabolic subgroup with Levi subgroup $L$. Since $M_0\subseteq L$, we have $x_{\mathrm{III}}\in M_0s_k\subseteq Ls_k=Lw$. Finally, $M\cap Q=P_{(k-1,1)}\times P_{(1,n-k-1)}$, so $M\cap U_Q$ consists of the elements $(m_1,m_2)\in M$ with $m_1\in N_{(k-1,1)}$ and $m_2\in N_{(1,n-k-1)}$.
\end{proof}

\subsection{The stabilizers}\label{ss:stabilizers}

We now determine the stabilizers $P_{x_*}=P\cap\eta_*H\eta_*^{-1}$ explicitly, using \eqref{eq:Px} and Lemma~\ref{lem:representatives}. Let $x=\eta\cdot I_n$, and put
$$
H_\eta=H\cap\eta^{-1}P\eta=\eta^{-1}P_x\eta.
$$
Following \cite{Off17}, for a character $\chi$ of $H$, we write $\chi^{\eta^{-1}}(p)=\chi(\eta^{-1}p\eta)$ for $p\in P_x=\eta H_\eta\eta^{-1}$. If $p\in P_x$, then $\eta^{-1}p\eta\in H$ preserves $V_+$ and $V_-$, and its restrictions to $V_+$ and $V_-$ are conjugate, via $\eta$, to the restrictions of $p$ to $\eta V_+$ and $\eta V_-$, respectively. Since the reduced norm of an automorphism of a finite-dimensional $\mathrm{D}$-vector space does not depend on the choice of a basis, we obtain
\begin{equation}\label{eq:chi-y}
	\chi_s^{\eta^{-1}}(p)=\nu\big(p|_{\eta V_+}\big)^{es}\,\nu\big(p|_{\eta V_-}\big)^{-es}\qquad(p\in P_x).
\end{equation}

We treat the three orbits in turn, beginning with the closed orbit $\OII$, whose representative is $x_{\mathrm{II}}=I_n$.

\begin{lemma}[The orbit $\OII$]\label{lem:stab-II}
	In blocks of sizes $(1,k-1,n-k)$, we have
	$$
	P_{x_{\mathrm{II}}}=P\cap H=\left\{\begin{pmatrix}a&0&0\\0&g_1&X\\0&0&g_3\end{pmatrix}:\ a\in G_1,\ g_1\in G_{k-1},\ g_3\in G_{n-k},\ X\in M_{(k-1)\times(n-k)}(\mathrm{D})\right\}.
	$$
	Thus, $P_{x_{\mathrm{II}}}=L_x\ltimes Z_x$, where $L_x=M_{x_{\mathrm{II}}}=\{\diag(a,g_1,g_3)\}=H_{1,k-1}\times G_{n-k}$ and $Z_x$ is the subgroup of elements with $a=1$, $g_1=1$ and $g_3=1$, which is contained in $N$. In particular, $P_{x_{\mathrm{II}}}=Q_{x_{\mathrm{II}}}$, and the projection of $Z_x$ to $M$ is trivial, that is, it equals $M\cap U_Q$.
\end{lemma}

\begin{proof}
	An element of $H$ has the form $\diag(a,B)$ with $a\in G_1$ and $B\in G_{n-1}$. It lies in $P$ if and only if it stabilizes $V_k=V_+\oplus(e_2\mathrm{D}\oplus\dots\oplus e_k\mathrm{D})$, that is, if and only if $B$ stabilizes $e_2\mathrm{D}\oplus\dots\oplus e_k\mathrm{D}$. This means that $B$ is block upper triangular with blocks of sizes $(k-1,n-k)$, which gives the displayed form. Since $Q=P$ by Lemma~\ref{lem:offen-data}(1), we have $P_x=Q_x$.
\end{proof}

The other closed orbit $\OI$ is treated similarly.

\begin{lemma}[The orbit $\OI$]\label{lem:stab-I}
	In blocks of sizes $(k,1,n-k-1)$, we have
	$$
	P_{x_{\mathrm I}}=\left\{\begin{pmatrix}g_1&0&Y\\0&c&0\\0&0&g_3\end{pmatrix}:\ g_1\in G_k,\ c\in G_1,\ g_3\in G_{n-k-1},\ Y\in M_{k\times(n-k-1)}(\mathrm{D})\right\}.
	$$
	Thus, $P_{x_{\mathrm I}}=L_x\ltimes Z_x$, where $L_x=M_{x_{\mathrm I}}=\{\diag(g_1,c,g_3)\}=G_k\times H_{1,n-k-1}$, with $H_{1,n-k-1}=\{\diag(c,g_3)\}\subseteq G_{n-k}$, and $Z_x$ is the subgroup of elements with $g_1=1$, $c=1$ and $g_3=1$, which is contained in $N$. In particular, $P_{x_{\mathrm I}}=Q_{x_{\mathrm I}}$, and the projection of $Z_x$ to $M$ is trivial.
\end{lemma}

\begin{proof}
	Write $p\in P$ as $\left(\begin{smallmatrix}A&B\\0&C\end{smallmatrix}\right)$ in blocks of sizes $(k,n-k)$. By \eqref{eq:Px} and Lemma~\ref{lem:representatives}(1), $p\in P_{x_{\mathrm I}}$ if and only if $p(e_{k+1}\mathrm{D})=e_{k+1}\mathrm{D}$ and $p$ stabilizes $V_k\oplus e_{k+2}\mathrm{D}\oplus\dots\oplus e_n\mathrm{D}$. Equivalently, by \eqref{eq:x-epsilon-explicit}, $p$ commutes with $\diag(-I_k,1,-I_{n-k-1})$. The vector $pe_{k+1}$ is the $(k+1)$-st column of $p$, which is formed by the first columns of $B$ and of $C$. It lies in $e_{k+1}\mathrm{D}$ if and only if the first column of $B$ vanishes and the first column of $C$ is $(c,0,\dots,0)^t$ for some $c\in G_1$. Since $pV_k\subseteq V_k$ automatically, the second condition says that $pe_j$ has vanishing $(k+1)$-st coordinate for $j\ge k+2$, that is, that the first row of $C$ vanishes outside its first entry. Hence, $C=\diag(c,g_3)$ and $B=(0\mid Y)$, which gives the displayed form.
\end{proof}

Finally, we treat the open orbit $\OIII$, which is the only case in which $Q\ne P$.

\begin{lemma}[The orbit $\OIII$]\label{lem:stab-III}
	In blocks of sizes $(k-1,1,1,n-k-1)$, we have
	$$
	P_{x_{\mathrm{III}}}=\left\{\begin{pmatrix}g_1&Z_{12}&-Z_{12}&Z_{14}\\0&c&0&Z_{24}\\0&0&c&-Z_{24}\\0&0&0&g_3\end{pmatrix}:\ \begin{array}{l}g_1\in G_{k-1},\ c\in G_1,\ g_3\in G_{n-k-1},\\ Z_{12}\in M_{(k-1)\times1}(\mathrm{D}),\ Z_{14}\in M_{(k-1)\times(n-k-1)}(\mathrm{D}),\\ Z_{24}\in M_{1\times(n-k-1)}(\mathrm{D})\end{array}\right\}.
	$$
	Thus, $P_{x_{\mathrm{III}}}=L_x\ltimes Z_x$, where $L_x=L\cap P_{x_{\mathrm{III}}}=\{\diag(g_1,c,c,g_3)\}\simeq G_{k-1}\times\Delta G_1\times G_{n-k-1}$ and $Z_x$ is the subgroup of elements with $g_1=1$, $c=1$ and $g_3=1$. Moreover, the following hold.
	\begin{enumerate}
		\item $P_{x_{\mathrm{III}}}\subseteq Q=P_{(k-1,1,1,n-k-1)}$, so that $P_{x_{\mathrm{III}}}=Q_{x_{\mathrm{III}}}$.
		\item The image of $p\in P_{x_{\mathrm{III}}}$ in $M=G_k\times G_{n-k}$ is
		$$
		\left(\begin{pmatrix}g_1&Z_{12}\\0&c\end{pmatrix},\ \begin{pmatrix}c&-Z_{24}\\0&g_3\end{pmatrix}\right).
		$$
		Consequently, the projection of $Z_x$ to $M$ is $N_{(k-1,1)}\times N_{(1,n-k-1)}=M\cap U_Q$, and the image of $P_{x_{\mathrm{III}}}$ in $M$ is $L'\ltimes(M\cap U_Q)$, where $L'=\{(\diag(g_1,c),\diag(c,g_3))\}\simeq G_{k-1}\times\Delta G_1\times G_{n-k-1}$.
		\item The kernel of the map $P_{x_{\mathrm{III}}}\to M$ is $P_{x_{\mathrm{III}}}\cap N$, which consists of the elements with $g_1=1$, $c=1$, $g_3=1$ and $Z_{12}=0$, $Z_{24}=0$.
	\end{enumerate}
\end{lemma}

\begin{proof}
	Put $\ell_+=(e_k+e_{k+1})\mathrm{D}=\eta_{\mathrm{III}}V_+$ and $W_-=\eta_{\mathrm{III}}V_-$. By Lemma~\ref{lem:representatives}(3), $W_-$ is the set of vectors whose $k$-th and $(k+1)$-st coordinates are opposite. Let $p\in P_{x_{\mathrm{III}}}$. Since $p$ stabilizes $V_k$ and $W_-$, it stabilizes $V_k\cap W_-=V_{k-1}$. Since $p$ stabilizes $V_k$ and $\ell_+$, it stabilizes $V_k+\ell_+=V_{k+1}$. Hence, $p\in Q$, which proves (1). We write $p$ in blocks of sizes $(k-1,1,1,n-k-1)$ as
	$$
	p=\begin{pmatrix}g_1&A_{12}&A_{13}&A_{14}\\0&a_{22}&a_{23}&a_{24}\\0&0&a_{33}&a_{34}\\0&0&0&g_3\end{pmatrix}.
	$$
	We have $p(e_k+e_{k+1})=(A_{12}+A_{13},\,a_{22}+a_{23},\,a_{33},\,0)^t$, which lies in $\ell_+$ if and only if $A_{13}=-A_{12}$ and $a_{22}+a_{23}=a_{33}$. Next, $p(e_k-e_{k+1})=(A_{12}-A_{13},\,a_{22}-a_{23},\,-a_{33},\,0)^t$ lies in $W_-$ if and only if $a_{22}-a_{23}=a_{33}$. For $j\le k-1$, we have $pe_j\in V_{k-1}\subseteq W_-$ automatically, and for $j\ge k+2$, the vector $pe_j$ lies in $W_-$ if and only if $a_{24}=-a_{34}$. Since $2\in\mathrm{F}^\times$, the two equations for $a_{22}$ and $a_{23}$ give $a_{23}=0$ and $a_{22}=a_{33}=:c$. Writing $Z_{12}=A_{12}$, $Z_{14}=A_{14}$ and $Z_{24}=a_{24}$, we obtain the displayed form. Conversely, every matrix of this form satisfies these conditions, and hence lies in $P_{x_{\mathrm{III}}}$.
	
	Every $p\in P_{x_{\mathrm{III}}}$ can be written as $p=\diag(g_1,c,c,g_3)\cdot z$ with $z=\diag(g_1,c,c,g_3)^{-1}p\in Z_x$. Since $Z_x$ is normal in $P_{x_{\mathrm{III}}}$ and meets $\{\diag(g_1,c,c,g_3)\}$ trivially, this gives the semidirect product decomposition, and $L_x=L\cap P_{x_{\mathrm{III}}}$ by inspection. The image of $p$ in $M$ is given by its diagonal blocks of sizes $k$ and $n-k$, which yields the formula in (2). Since $Z_{12}$ and $Z_{24}$ vary independently, the remaining assertions of (2) follow, and (3) follows from (2).
\end{proof}

In summary, the three lemmas verify, for our representatives, the properties used in \cite[Lemmas 3.2, 3.3 and the proof of Proposition 4.1]{Off17}. Namely, $P_x=L_x\ltimes Z_x$ with $L_x=L\cap P_x$, $P_x=Q_x$, and the projection of $Z_x$ to $M$ equals $M\cap U_Q$.

\subsection{The characters attached to the orbits}\label{ss:orbit-characters}

For a closed subgroup $S$ of $G_n$, we put $\delta_S(y)=|\det(\Ad(y)|_{\Lie S})|_{\mathrm{F}}$ for $y\in S$; this is the modulus function used in \cite{Off17}. Suppose that $S=L\ltimes U$, where $U$ is unipotent and $L$ is isomorphic to a product of groups $G_j$, embedded block-diagonally, possibly with some blocks repeated. Then
$$
\delta_S(lu)=|\det(\Ad(l)|_{\Lie U})|_{\mathrm{F}}\qquad(l\in L,\ u\in U).
$$
Indeed, $\Ad(u)$ is unipotent, and $\Ad(l)$ acts on $\Lie L$ by conjugation on each factor $M_j(\mathrm{D})$. By Lemma~\ref{lem:det}, $|\det(\Ad(g)|_{M_j(\mathrm{D})})|_{\mathrm{F}}=\nu(g)^{ej}\nu(g)^{-ej}=1$ for $g\in G_j$. In \cite[Proposition 4.1]{Off17}, the character of $L_x$ that governs the contribution of the orbit of $x=\eta\cdot I_n$ is
$$
\delta_{Q_x}\,\delta_Q^{-1/2}\,\chi^{\eta^{-1}}.
$$
We now compute this character for the three orbits, beginning with the closed orbit $\OII$.

\begin{lemma}\label{lem:char-II}
	For $m=\diag(a,g_1,g_3)\in L_{x_{\mathrm{II}}}$, written in blocks of sizes $(1,k-1,n-k)$, we have
	\begin{align*}
		\delta_{Q_{x_{\mathrm{II}}}}(m)&=\nu(g_1)^{e(n-k)}\nu(g_3)^{-e(k-1)},\\
		\delta_Q(m)=\delta_P(m)&=\nu(a)^{e(n-k)}\nu(g_1)^{e(n-k)}\nu(g_3)^{-ek},\\
		\chi_s^{\eta_{\mathrm{II}}^{-1}}(m)&=\nu(a)^{es}\nu(g_1)^{-es}\nu(g_3)^{-es}.
	\end{align*}
	Consequently,
	$$
	\delta_{Q_{x_{\mathrm{II}}}}\,\delta_Q^{-1/2}\,\chi_s^{\eta_{\mathrm{II}}^{-1}}(m)=\chi^{(1,k-1)}_{s-\frac{n-k}2}\big(\diag(a,g_1)\big)\cdot\nu(g_3)^{\frac{e(2-k)}2-es}.
	$$
\end{lemma}

\begin{proof}
	By Lemma~\ref{lem:stab-II}, we have $Q_{x_{\mathrm{II}}}=P_{x_{\mathrm{II}}}=L_x\ltimes Z_x$, and $m$ acts on $\Lie Z_x=M_{(k-1)\times(n-k)}(\mathrm{D})$ by $X\mapsto g_1Xg_3^{-1}$. Hence, Lemma~\ref{lem:det} gives the formula for $\delta_{Q_{x_{\mathrm{II}}}}$. The formula for $\delta_Q=\delta_P$ follows from Lemma~\ref{lem:modulus}, applied with $\diag(a,g_1)\in G_k$ and $g_3\in G_{n-k}$. Since $\eta_{\mathrm{II}}=I_n$, we have $\chi_s^{\eta_{\mathrm{II}}^{-1}}=\chi_s$. Multiplying, we find that the exponent of $\nu(a)$ is $es-\frac{e(n-k)}2=e(s-\frac{n-k}2)$, the exponent of $\nu(g_1)$ is $-es+e(n-k)-\frac{e(n-k)}2=-e(s-\frac{n-k}2)$, and the exponent of $\nu(g_3)$ is $-es-e(k-1)+\frac{ek}2=\frac{e(2-k)}2-es$.
\end{proof}

The computation for the other closed orbit $\OI$ is similar.

\begin{lemma}\label{lem:char-I}
	For $m=\diag(g_1,c,g_3)\in L_{x_{\mathrm I}}$, written in blocks of sizes $(k,1,n-k-1)$, we have
	\begin{align*}
		\delta_{Q_{x_{\mathrm I}}}(m)&=\nu(g_1)^{e(n-k-1)}\nu(g_3)^{-ek},\\
		\delta_Q(m)=\delta_P(m)&=\nu(g_1)^{e(n-k)}\nu(c)^{-ek}\nu(g_3)^{-ek},\\
		\chi_s^{\eta_{\mathrm I}^{-1}}(m)&=\nu(c)^{es}\nu(g_1)^{-es}\nu(g_3)^{-es}.
	\end{align*}
	Consequently,
	$$
	\delta_{Q_{x_{\mathrm I}}}\,\delta_Q^{-1/2}\,\chi_s^{\eta_{\mathrm I}^{-1}}(m)=\nu(g_1)^{\frac{e(n-k-2)}2-es}\cdot\chi^{(1,n-k-1)}_{s+\frac k2}\big(\diag(c,g_3)\big).
	$$
\end{lemma}

\begin{proof}
	By Lemma~\ref{lem:stab-I}, $m$ acts on $\Lie Z_x=M_{k\times(n-k-1)}(\mathrm{D})$ by $Y\mapsto g_1Yg_3^{-1}$, and Lemma~\ref{lem:det} gives the formula for $\delta_{Q_{x_{\mathrm I}}}$. The formula for $\delta_Q=\delta_P$ follows from Lemma~\ref{lem:modulus}, applied with $g_1\in G_k$ and $\diag(c,g_3)\in G_{n-k}$. By Lemma~\ref{lem:representatives}(1), $m$ acts on $\eta_{\mathrm I}V_+=e_{k+1}\mathrm{D}$ by $c$, and on $\eta_{\mathrm I}V_-$, in the basis $e_1,\dots,e_k,e_{k+2},\dots,e_n$, by $\diag(g_1,g_3)$. Hence, \eqref{eq:chi-y} gives the formula for $\chi_s^{\eta_{\mathrm I}^{-1}}$. Multiplying, we find that the exponent of $\nu(g_1)$ is $-es+e(n-k-1)-\frac{e(n-k)}2=\frac{e(n-k-2)}2-es$, the exponent of $\nu(c)$ is $es+\frac{ek}2=e(s+\frac k2)$, and the exponent of $\nu(g_3)$ is $-es-ek+\frac{ek}2=-e(s+\frac k2)$.
\end{proof}

For the open orbit $\OIII$, the contributions of the two copies of $c$ cancel.

\begin{lemma}\label{lem:char-III}
	For $l=\diag(g_1,c,c,g_3)\in L_{x_{\mathrm{III}}}$, written in blocks of sizes $(k-1,1,1,n-k-1)$, we have
	\begin{align*}
		\delta_{Q_{x_{\mathrm{III}}}}(l)&=\nu(g_1)^{e(n-k)}\nu(c)^{e(n-2k)}\nu(g_3)^{-ek},\\
		\delta_Q(l)&=\nu(g_1)^{e(n-k+1)}\nu(c)^{2e(n-2k)}\nu(g_3)^{-e(k+1)},\\
		\chi_s^{\eta_{\mathrm{III}}^{-1}}(l)&=\nu(g_1)^{-es}\nu(g_3)^{-es}.
	\end{align*}
	Consequently, the powers of $\nu(c)$ cancel, and
	$$
	\delta_{Q_{x_{\mathrm{III}}}}\,\delta_Q^{-1/2}\,\chi_s^{\eta_{\mathrm{III}}^{-1}}(l)=\nu(g_1)^{\frac{e(n-k-1)}2-es}\,\nu(g_3)^{\frac{e(1-k)}2-es}.
	$$
\end{lemma}

\begin{proof}
	By Lemma~\ref{lem:stab-III}, we have $Q_{x_{\mathrm{III}}}=P_{x_{\mathrm{III}}}=L_x\ltimes Z_x$, and $\Lie Z_x$ is identified with the $\mathrm{F}$-vector space of triples $(Z_{12},Z_{14},Z_{24})$, on which $l$ acts by
	$$
	(Z_{12},Z_{14},Z_{24})\longmapsto(g_1Z_{12}c^{-1},\ g_1Z_{14}g_3^{-1},\ cZ_{24}g_3^{-1}).
	$$
	By Lemma~\ref{lem:det}, the determinants of these three maps have absolute values
	$$
	\nu(g_1)^{e}\nu(c)^{-e(k-1)},\qquad\nu(g_1)^{e(n-k-1)}\nu(g_3)^{-e(k-1)},\qquad\nu(c)^{e(n-k-1)}\nu(g_3)^{-e},
	$$
	respectively, and their product is $\delta_{Q_{x_{\mathrm{III}}}}(l)$. Next, Lemma~\ref{lem:modulus} for the composition $(k-1,1,1,n-k-1)$ gives the exponent $e(n-k+1)$ for $g_1$, the exponent $e\big((n-k)-(k-1)\big)+e\big((n-k-1)-k\big)=2e(n-2k)$ for the two copies of $c$, and the exponent $-e(k+1)$ for $g_3$. By Lemma~\ref{lem:representatives}(3), $l$ acts on $\eta_{\mathrm{III}}V_+=(e_k+e_{k+1})\mathrm{D}$ by $c$, and on $\eta_{\mathrm{III}}V_-$, in the basis $e_1,\dots,e_{k-1},e_k-e_{k+1},e_{k+2},\dots,e_n$, by $\diag(g_1,c,g_3)$. Hence, by \eqref{eq:chi-y},
	$$
	\chi_s^{\eta_{\mathrm{III}}^{-1}}(l)=\nu(c)^{es}\big(\nu(g_1)\nu(c)\nu(g_3)\big)^{-es}=\nu(g_1)^{-es}\nu(g_3)^{-es}.
	$$
	Multiplying, we find that the exponent of $\nu(g_1)$ is $e(n-k)-\frac{e(n-k+1)}2-es=\frac{e(n-k-1)}2-es$, the exponent of $\nu(c)$ is $e(n-2k)-e(n-2k)=0$, and the exponent of $\nu(g_3)$ is $-ek+\frac{e(k+1)}2-es=\frac{e(1-k)}2-es$.
\end{proof}

\subsection{Mackey theory and the geometric lemma}\label{ss:mackey}

We now combine the orbit computations with Mackey theory. For a closed subgroup $K$ of an $\ell$-group $S$ and a smooth representation $(\rho,E)$ of $K$, let $\ind_K^S(\rho)$ denote the unnormalized compact induction. This is the space of locally constant functions $\phi\colon S\to E$ such that $\phi(ks)=\rho(k)\phi(s)$ for $k\in K$ and $s\in S$, and whose support is compact modulo $K$, with $S$ acting by right translation.

We shall use Frobenius reciprocity for compact induction in the following form.

\begin{lemma}[Frobenius reciprocity for compact induction; {\cite[Proposition 2.29]{BZ76}}]\label{lem:frob-compact}
	Let $S$ be a unimodular $\ell$-group, let $K$ be a closed subgroup of $S$ that is the group of $\mathrm{F}$-points of a linear algebraic group, let $\rho$ be a smooth representation of $K$, and let $\chi$ be a character of $S$. Then
	$$
	\Hom_S\big(\ind_K^S(\rho),\chi\big)\simeq\Hom_K\big(\rho,\ \chi|_K\cdot\delta_K\big).
	$$
\end{lemma}

For instance, let $B$ be the Borel subgroup of $\GL_2(\mathrm{F})$ and $\rho=\delta_B$. Then the lemma says that the trivial character is a quotient of the unnormalized induction of $\delta_B$, which is the normalized induction $|\cdot|^{1/2}\times|\cdot|^{-1/2}$. This example fixes the normalization of $\delta_K$.

Let $\sigma\in\Rep(M)$ and $I=\Ind_P^{G_n}(\sigma)$, and regard $\sigma\delta_P^{1/2}$ as a representation of $P$ trivial on $N$. For $\eta\in G_n$ and $x=\eta\cdot I_n$, let $(\sigma\delta_P^{1/2}|_{P_x})^\eta$ denote the representation of $H_\eta=\eta^{-1}P_x\eta$ on the space of $\sigma$ given by $h\mapsto(\sigma\delta_P^{1/2})(\eta h\eta^{-1})$. The restriction of $I$ to $H$ is then described as follows.

\begin{lemma}[Mackey filtration; {\cite[\S4, (4)]{Off17}}, {\cite[Theorem 5.2]{BZ77}}]\label{lem:mackey}
	Let $I^{\mathrm{op}}\subseteq I$ be the subspace of functions vanishing on the closed set $P\eta_{\mathrm I}H\sqcup P\eta_{\mathrm{II}}H$. Then $I^{\mathrm{op}}$ is $H$-stable, and there is an exact sequence of representations of $H$
	$$
	0\longrightarrow I^{\mathrm{op}}\longrightarrow I|_H\longrightarrow I_{\mathrm I}\oplus I_{\mathrm{II}}\longrightarrow0,
	$$
	where $I_*$ denotes the space of restrictions of elements of $I$ to $P\eta_*H$. Moreover, for each $*$, the map $f\mapsto\big(h\mapsto f(\eta_*h)\big)$ induces isomorphisms of representations of $H$
	$$
	I_{\mathrm I}\simeq\ind_{H_{\eta_{\mathrm I}}}^H\big((\sigma\delta_P^{1/2}|_{P_{x_{\mathrm I}}})^{\eta_{\mathrm I}}\big),\quad I_{\mathrm{II}}\simeq\ind_{H_{\eta_{\mathrm{II}}}}^H\big((\sigma\delta_P^{1/2}|_{P_{x_{\mathrm{II}}}})^{\eta_{\mathrm{II}}}\big),\quad I^{\mathrm{op}}\simeq\ind_{H_{\eta_{\mathrm{III}}}}^H\big((\sigma\delta_P^{1/2}|_{P_{x_{\mathrm{III}}}})^{\eta_{\mathrm{III}}}\big).
	$$
\end{lemma}

\begin{proof}
	This is the filtration of $I|_H$ by the $H$-orbits on $P\backslash G_n$ described in \cite[\S4, (4)]{Off17}, for the order $(P\eta_{\mathrm{III}}H,P\eta_{\mathrm{I}}H,P\eta_{\mathrm{II}}H)$ of the double cosets. For this order, the union of the first $i$ double cosets is open for every $i$, by Lemma~\ref{lem:representatives}; see also \cite[\S\S1.8, 2.23]{BZ76} and \cite[\S5]{BZ77}. For the convenience of the reader, we recall the argument.
	
	The elements of $I$ are the locally constant sections of an $\ell$-sheaf on the compact space $P\backslash G_n$. The set $\mathcal Z=P\backslash(P\eta_{\mathrm I}H\sqcup P\eta_{\mathrm{II}}H)$ is closed, its complement $\mathcal U=P\backslash P\eta_{\mathrm{III}}H$ is open, and $\mathcal Z$ is the disjoint union of the two closed sets $P\backslash P\eta_{\mathrm I}H$ and $P\backslash P\eta_{\mathrm{II}}H$. By \cite[Proposition 1.8]{BZ76}, restriction of sections gives an exact sequence $0\to I^{\mathrm{op}}\to I\to I_{\mathcal Z}\to0$. Here, $I^{\mathrm{op}}$ is the space of sections with compact support in $\mathcal U$, that is, the elements of $I$ vanishing on $\mathcal Z$, and $I_{\mathcal Z}=I_{\mathrm I}\oplus I_{\mathrm{II}}$.
	
	Each orbit $P\backslash P\eta_*H$ is locally closed and homeomorphic to $H_{\eta_*}\backslash H$ via $h\mapsto P\eta_*h$ \cite[\S1.5]{BZ76}. If $f\in I$, then the function $\phi(h)=f(\eta_*h)$ satisfies
	$$
	\phi(h_0h)=f\big((\eta_*h_0\eta_*^{-1})\eta_*h\big)=(\sigma\delta_P^{1/2})(\eta_*h_0\eta_*^{-1})\phi(h)\qquad(h_0\in H_{\eta_*}).
	$$
	For the closed orbits, every section is compactly supported modulo $H_{\eta_*}$, because the orbit is compact. For the open orbit, the sections with compact support in $\mathcal U$ correspond to the functions with compact support modulo $H_{\eta_{\mathrm{III}}}$. Conversely, every such $\phi$ arises from a section, by the extension property of locally constant sections on closed and open subsets, respectively \cite[\S1]{BZ76}.
\end{proof}

The contribution of each orbit is computed by the following result of Offen.

\begin{proposition}[{\cite[Proposition 4.1]{Off17}}]\label{prop:offen41}
	Let $*\in\{\mathrm I,\mathrm{II},\mathrm{III}\}$, put $\eta=\eta_*$ and $x=x_*$, and let $L$ and $Q=LU_Q$ be as in Lemma~\ref{lem:offen-data}. Then, for every character $\chi$ of $H$,
	$$
	\Hom_H\Big(\ind_{H_\eta}^H\big((\sigma\delta_P^{1/2}|_{P_x})^\eta\big),\chi\Big)\simeq\Hom_{L_x}\big(r_{L,M}(\sigma),\ \delta_{Q_x}\,\delta_Q^{-1/2}\,\chi^{\eta^{-1}}\big),
	$$
	where $r_{L,M}(\sigma)$ is the normalized Jacquet module of $\sigma$ with respect to the standard parabolic subgroup $M\cap Q=L\ltimes(M\cap U_Q)$ of $M$.
\end{proposition}

\begin{proof}
	This is \cite[Proposition 4.1]{Off17}; by Lemma~\ref{lem:offen-data}, our representatives satisfy $x\in Lw$, as required there. For the convenience of the reader, we recall the argument in our situation, using Lemmas~\ref{lem:stab-II}--\ref{lem:stab-III}.
	
	The groups $H$ and $H_\eta$ are groups of $\mathrm{F}$-points of linear algebraic groups, and $H$ is unimodular. By Lemma~\ref{lem:frob-compact}, the left-hand side is isomorphic to $\Hom_{H_\eta}\big((\sigma\delta_P^{1/2})^\eta,\chi\,\delta_{H_\eta}\big)$. Conjugation by $\eta$ is an isomorphism $H_\eta\to P_x$, $h\mapsto p=\eta h\eta^{-1}$, under which
	$$
	(\sigma\delta_P^{1/2})^\eta(h)=\sigma(p)\delta_P^{1/2}(p),\qquad\chi(h)=\chi^{\eta^{-1}}(p),\qquad\delta_{H_\eta}(h)=\delta_{P_x}(p).
	$$
	Hence, the left-hand side is isomorphic to
	$$
	\Hom_{P_x}\big(\sigma,\ \delta_{P_x}\,\delta_P^{-1/2}\,\chi^{\eta^{-1}}\big),
	$$
	where $P_x$ acts on the space of $\sigma$ through its image in $M$. By Lemmas~\ref{lem:stab-II}--\ref{lem:stab-III}, we have $P_x=Q_x=L_x\ltimes Z_x$, the kernel of $P_x\to M$ is contained in $N$ and therefore acts trivially on $\sigma$, and the projection of $Z_x$ to $M$ is $M\cap U_Q$. The character $\delta_{P_x}\delta_P^{-1/2}\chi^{\eta^{-1}}$ is trivial on the unipotent group $Z_x$. Hence, every linear form in the last space factors through the coinvariants $\sigma_{M\cap U_Q}$, and the space is isomorphic to
	$$
	\Hom_{L_x}\big(\sigma_{M\cap U_Q},\ \delta_{Q_x}\,\delta_P^{-1/2}\,\chi^{\eta^{-1}}\big).
	$$
	Finally, we have $\sigma_{M\cap U_Q}=\delta_{M\cap Q}^{1/2}\,r_{L,M}(\sigma)$. Moreover, since $\Lie U_Q=\Lie N\oplus\Lie(M\cap U_Q)$ as $L$-modules, we have $\delta_Q|_L=\delta_P|_L\,\delta_{M\cap Q}|_L$. Therefore, $\delta_P^{-1/2}\delta_{M\cap Q}^{-1/2}=\delta_Q^{-1/2}$ on $L$, and the assertion follows.
\end{proof}

When the inducing representation is a tensor product, we shall use the following elementary lemma to compute the resulting spaces of homomorphisms.

\begin{lemma}\label{lem:tensor-hom}
	Let $A_1$ and $A_2$ be groups, and for $i=1,2$, let $V_i$ be a representation of $A_i$ and $\lambda_i$ a character of $A_i$. Suppose that $W_1=\Hom_{A_1}(V_1,\lambda_1)$ is finite-dimensional. Then
	$$
	\Hom_{A_1\times A_2}(V_1\otimes V_2,\lambda_1\otimes\lambda_2)\simeq W_1\otimes\Hom_{A_2}(V_2,\lambda_2).
	$$
	In particular, the left-hand side is nonzero if and only if both spaces $\Hom_{A_i}(V_i,\lambda_i)$ are nonzero, and its dimension is $\dim W_1\cdot\dim\Hom_{A_2}(V_2,\lambda_2)$, with the convention $0\cdot\infty=0$.
\end{lemma}

\begin{proof}
	Consider the map sending $\ell_1\otimes\ell_2$ to the linear form $v_1\otimes v_2\mapsto\ell_1(v_1)\ell_2(v_2)$. Choose a basis $\ell_1,\dots,\ell_r$ of $W_1$ and vectors $v^{(1)},\dots,v^{(r)}\in V_1$ with $\ell_i(v^{(j)})=\delta_{ij}$, which is possible since the $\ell_i$ are linearly independent. If $\sum_i\ell_i\otimes\mu_i$ maps to zero, then evaluating at $v^{(j)}\otimes v_2$ gives $\mu_j(v_2)=0$ for all $j$ and $v_2$. Hence, the map is injective.
	
	Conversely, let $\ell$ be an element of the left-hand side. For every $v_2\in V_2$, the linear form $\ell(\cdot\otimes v_2)$ is $(A_1\times1)$-equivariant, and hence lies in $W_1$. Therefore, $\ell(v_1\otimes v_2)=\sum_i\ell_i(v_1)\mu_i(v_2)$ with $\mu_i(v_2)=\ell(v^{(i)}\otimes v_2)$. Equivariance of $\ell$ under $1\times A_2$ shows that $\mu_i\in\Hom_{A_2}(V_2,\lambda_2)$, so the map is surjective.
\end{proof}

\subsection{The contributions of the three orbits}\label{ss:contributions}

We now apply the preceding results to a representation induced from the maximal parabolic subgroup $P$. Let $\pi_1\in\Rep(G_k)$ and $\pi_2\in\Rep(G_{n-k})$, and let $\sigma=\pi_1\otimes\pi_2$, so that $I=\pi_1\times\pi_2$. Let
$$
L_{n,k}=\{\diag(g_1,a,a,g_3):g_1\in G_{k-1},\ a\in G_1,\ g_3\in G_{n-k-1}\}\simeq G_{k-1}\times\Delta G_1\times G_{n-k-1}.
$$
This is the group $L_{x_{\mathrm{III}}}$ of Lemma~\ref{lem:stab-III}, and we regard it as a subgroup of $M_{(k-1,1)}\times M_{(1,n-k-1)}$ via the embedding $\diag(g_1,a,a,g_3)\mapsto(\diag(g_1,a),\diag(a,g_3))$. We define the character $\kappa_s$ of $L_{n,k}$ by
$$
\kappa_s\big(\diag(g_1,a,a,g_3)\big)=\nu(g_1)^{\frac{e(n-k-1)}2-es}\,\nu(g_3)^{\frac{e(1-k)}2-es}.
$$

We first compute the contribution of the closed orbit $\OI$.

\begin{proposition}[Contribution of $\OI$]\label{prop:contrib-I}
	We have
	$$
	\Hom_H(I_{\mathrm I},\chi_s)\simeq\Hom_{G_k\times H_{1,n-k-1}}\Big(\pi_1\otimes\pi_2,\ \nu^{\frac{e(n-k-2)}2-es}\otimes\chi_{s+\frac k2}\Big).
	$$
\end{proposition}

\begin{proof}
	By Lemma~\ref{lem:offen-data}, we have $L=M$ and $Q=P$, so that $r_{M,M}(\sigma)=\sigma$. Hence, by Lemma~\ref{lem:mackey} and Proposition~\ref{prop:offen41}, the left-hand side is isomorphic to $\Hom_{L_{x_{\mathrm I}}}\big(\sigma,\delta_{Q_x}\delta_Q^{-1/2}\chi_s^{\eta_{\mathrm I}^{-1}}\big)$. By Lemma~\ref{lem:stab-I}, the group $L_{x_{\mathrm I}}=G_k\times H_{1,n-k-1}$ acts on $\sigma$ through its inclusion in $M$, and Lemma~\ref{lem:char-I} identifies the character $\delta_{Q_x}\delta_Q^{-1/2}\chi_s^{\eta_{\mathrm I}^{-1}}$ with $\nu^{\frac{e(n-k-2)}2-es}\otimes\chi_{s+\frac k2}$.
\end{proof}

The other closed orbit $\OII$ is treated in the same way.

\begin{proposition}[Contribution of $\OII$]\label{prop:contrib-II}
	We have
	$$
	\Hom_H(I_{\mathrm{II}},\chi_s)\simeq\Hom_{H_{1,k-1}\times G_{n-k}}\Big(\pi_1\otimes\pi_2,\ \chi_{s-\frac{n-k}2}\otimes\nu^{\frac{e(2-k)}2-es}\Big).
	$$
\end{proposition}

\begin{proof}
	The argument is the same as in the proof of Proposition~\ref{prop:contrib-I}. By Lemma~\ref{lem:offen-data}, we again have $L=M$ and $Q=P$. By Lemma~\ref{lem:stab-II}, the group $L_{x_{\mathrm{II}}}=H_{1,k-1}\times G_{n-k}$ acts on $\sigma$ through its inclusion in $M$, and Lemma~\ref{lem:char-II} identifies the character.
\end{proof}

For the open orbit $\OIII$, a Jacquet module appears.

\begin{proposition}[Contribution of $\OIII$]\label{prop:contrib-III}
	We have
	$$
	\Hom_H(I^{\mathrm{op}},\chi_s)\simeq\Hom_{L_{n,k}}\big(r_{(k-1,1)}(\pi_1)\otimes r_{(1,n-k-1)}(\pi_2),\ \kappa_s\big).
	$$
\end{proposition}

\begin{proof}
	By Lemma~\ref{lem:mackey} and Proposition~\ref{prop:offen41}, the left-hand side is isomorphic to $\Hom_{L_{x_{\mathrm{III}}}}\big(r_{L,M}(\sigma),\delta_{Q_x}\delta_Q^{-1/2}\chi_s^{\eta_{\mathrm{III}}^{-1}}\big)$, where $L=M_{(k-1,1,1,n-k-1)}$. Since $M\cap Q=P_{(k-1,1)}\times P_{(1,n-k-1)}$ by Lemma~\ref{lem:offen-data}, we have
	$$
	r_{L,M}(\pi_1\otimes\pi_2)=r_{(k-1,1)}(\pi_1)\otimes r_{(1,n-k-1)}(\pi_2).
	$$
	By Lemma~\ref{lem:stab-III}, the group $L_{x_{\mathrm{III}}}=L_{n,k}$ acts on this space through the embedding $\diag(g_1,c,c,g_3)\mapsto(\diag(g_1,c),\diag(c,g_3))$. Finally, by Lemma~\ref{lem:char-III}, the character is $\kappa_s$.
\end{proof}

For $\pi_1\in\Rep(G_k)$ and $\pi_2\in\Rep(G_{n-k})$, we now put
\begin{align*}
	\dI&=\dim\Hom_{G_k}\big(\pi_1,\nu^{\frac{e(n-k-2)}2-es}\big)\cdot\dim\Hom_{H_{1,n-k-1}}\big(\pi_2,\chi_{s+\frac k2}\big),\\
	\dII&=\dim\Hom_{H_{1,k-1}}\big(\pi_1,\chi_{s-\frac{n-k}2}\big)\cdot\dim\Hom_{G_{n-k}}\big(\pi_2,\nu^{\frac{e(2-k)}2-es}\big),\\
	\dIII&=\dim\Hom_{L_{n,k}}\big(r_{(k-1,1)}(\pi_1)\otimes r_{(1,n-k-1)}(\pi_2),\kappa_s\big),
\end{align*}
with the convention $0\cdot\infty=0$. Since $\pi_1$ and $\pi_2$ have finite length, the spaces $\Hom_{G_k}(\pi_1,\nu^c)$ and $\Hom_{G_{n-k}}(\pi_2,\nu^c)$ are finite-dimensional for every $c\in\R$; indeed, their dimensions are bounded by the multiplicity of $\nu^c$ as a constituent. Hence, by Lemma~\ref{lem:tensor-hom} (applied with the roles of the two factors interchanged in the case of $\dII$), the numbers $\dI$ and $\dII$ are the dimensions of the spaces in Propositions~\ref{prop:contrib-I} and~\ref{prop:contrib-II}, respectively.

Combining the three contributions, we obtain the main result of this section.

\begin{theorem}\label{thm:orbit-bound}
	Let $\pi_1\in\Rep(G_k)$, $\pi_2\in\Rep(G_{n-k})$ and $s\in\R$. Then
	$$
	\dI+\dII\ \le\ \dim\Hom_H(\pi_1\times\pi_2,\chi_s)\ \le\ \dI+\dII+\dIII.
	$$
	In particular, if $\Hom_H(\pi_1\times\pi_2,\chi_s)\ne0$, then one of $\dI$, $\dII$ and $\dIII$ is nonzero. Moreover, if $\dI\ne0$ or $\dII\ne0$, then $\pi_1\times\pi_2$ is $(H,\chi_s)$-distinguished.
\end{theorem}

\begin{proof}
	Applying the left exact functor $\Hom_H(-,\chi_s)$ to the exact sequence of Lemma~\ref{lem:mackey}, we obtain an exact sequence
	$$
	0\to\Hom_H(I_{\mathrm I},\chi_s)\oplus\Hom_H(I_{\mathrm{II}},\chi_s)\to\Hom_H(\pi_1\times\pi_2,\chi_s)\to\Hom_H(I^{\mathrm{op}},\chi_s).
	$$
	By Propositions~\ref{prop:contrib-I}, \ref{prop:contrib-II} and~\ref{prop:contrib-III}, the dimensions of the outer terms are $\dI+\dII$ and $\dIII$, respectively. The assertions follow.
\end{proof}

We now spell out the three conditions in the form in which they will be used. We say that condition (I), (II) or (III) holds for $(\pi_1,\pi_2,s)$ if $\dI\ne0$, $\dII\ne0$ or $\dIII\ne0$, respectively. By Lemma~\ref{lem:tensor-hom}, conditions (I) and (II) mean the following:
\begin{itemize}
	\item[(I)] $\Hom_{G_k}\big(\pi_1,\nu^{\frac{e(n-k-2)}2-es}\big)\ne0$ and $\Hom_{H_{1,n-k-1}}\big(\pi_2,\chi_{s+\frac k2}\big)\ne0$;
	\item[(II)] $\Hom_{H_{1,k-1}}\big(\pi_1,\chi_{s-\frac{n-k}2}\big)\ne0$ and $\Hom_{G_{n-k}}\big(\pi_2,\nu^{\frac{e(2-k)}2-es}\big)\ne0$.
\end{itemize}
For $k=n-1$, the second part of (I) says that $\nu^{e(s+\frac{n-1}2)}$ is a quotient of $\pi_2\in\Rep(G_1)$. For $k=1$, the first part of (II) says that $\nu^{e(s-\frac{n-1}2)}$ is a quotient of $\pi_1\in\Rep(G_1)$. Finally, if $\pi_1$ is irreducible, then $\Hom_{G_k}(\pi_1,\nu^c)\ne0$ if and only if $\pi_1\simeq\nu^c$.

\subsection{The open orbit and representations of \texorpdfstring{$\mathrm{D}^\times$}{Dx}}\label{ss:open-orbit}

We now analyze $\dIII$. In $L_{n,k}$, the two middle factors $G_1=\mathrm{D}^\times$ are identified. Since irreducible representations of $\mathrm{D}^\times$ need not be characters when $e\ge2$, the resulting condition must be expressed in terms of arbitrary irreducible representations of $\mathrm{D}^\times$. We put
$$
A=G_{k-1},\qquad B=G_1=\mathrm{D}^\times,\qquad C=G_{n-k-1},
$$
so that $L_{n,k}=A\times\Delta B\times C$.

We first record two elementary facts about representations of $\mathrm{D}^\times$, which rely on the compactness of $\mathrm{D}^\times/\mathrm{F}^\times$.

\begin{lemma}\label{lem:Dx-semisimple}
	Let $E$ be a finite-dimensional smooth representation of $\mathrm{D}^\times$ on which the centre $\mathrm{F}^\times$ acts by a character. Then $E$ is semisimple.
\end{lemma}

\begin{proof}
	Let $\omega$ be the character by which $\mathrm{F}^\times$ acts on $E$. The character $|\omega|$ of $\mathrm{F}^\times$ is positive, hence trivial on the compact group $\mathcal O_{\mathrm{F}}^\times$, and therefore $|\omega|=|\cdot|_{\mathrm{F}}^{t}$ for some $t\in\R$. Choose $x\in\R$ with $ex=-t$, so that $|\omega(z)|\,|z|_{\mathrm{F}}^{ex}=1$ for all $z\in \mathrm{F}^\times$. Twisting by the character $\nu^x$ of $\mathrm{D}^\times$ does not change the lattice of subrepresentations. Hence, by \eqref{eq:nu-centre}, we may assume that $\omega$ is unitary. Choose any inner product $\langle\,,\rangle_0$ on $E$. For $v,w\in E$, the function $g\mapsto\langle gv,gw\rangle_0$ is locally constant and invariant under $\mathrm{F}^\times$, and therefore defines a function on the compact group $\mathrm{D}^\times/\mathrm{F}^\times$. Averaging it over $\mathrm{D}^\times/\mathrm{F}^\times$ gives a $\mathrm{D}^\times$-invariant inner product on $E$. The orthogonal complement of a subrepresentation with respect to this inner product is a complementary subrepresentation.
\end{proof}

As a consequence, every constituent of a finite-length representation of $\mathrm{D}^\times$ occurs as a quotient.

\begin{lemma}\label{lem:Dx-quotient}
	Every constituent of a finite-length representation of $\mathrm{D}^\times$ is a quotient of that representation.
\end{lemma}

\begin{proof}
	We first show that every extension $0\to\tau\to E\to\sigma\to0$ of non-isomorphic irreducible representations of $\mathrm{D}^\times$ splits. The representation $E$ is finite-dimensional.
	
	Suppose first that $\omega_\tau\ne\omega_\sigma$, and choose $z\in \mathrm{F}^\times$ with $\omega_\tau(z)\ne\omega_\sigma(z)$. Since $z$ is central, the generalized $\omega_\sigma(z)$-eigenspace of $E(z)$ is $\mathrm{D}^\times$-stable, and it is a complement to $\tau$.
	
	Suppose now that $\omega_\tau=\omega_\sigma=\omega$. For $z\in \mathrm{F}^\times$, the operator $E(z)-\omega(z)$ maps $E$ into $\tau$ and vanishes on $\tau$. Hence, it induces a $\mathrm{D}^\times$-homomorphism $\sigma\to\tau$, which is zero by Schur's lemma. Thus, $\mathrm{F}^\times$ acts on $E$ by $\omega$, and $E$ is semisimple by Lemma~\ref{lem:Dx-semisimple}.
	
	Now let $X$ be a finite-length representation of $\mathrm{D}^\times$ with a constituent $\tau$. We argue by induction on the length of $X$. Choose an exact sequence $0\to X'\to X\to\sigma\to0$ with $\sigma$ irreducible. If $\sigma\simeq\tau$, there is nothing to prove. Otherwise, $\tau$ is a constituent of $X'$, so by the induction hypothesis there is a surjection $X'\to\tau$, with kernel $K$, say. The exact sequence $0\to\tau\to X/K\to\sigma\to0$ splits by the first part of the proof. Hence, the composition $X\to X/K\to\tau$ is surjective.
\end{proof}

We now isolate the linear algebra underlying the contribution of the open orbit. Let $J_1\in\Rep(A\times B)$ and $J_2\in\Rep(B\times C)$, and let $\chi_A$ and $\chi_C$ be characters of $A$ and $C$, respectively. Let $U$ be the largest quotient of $J_1$ on which $A$ acts by $\chi_A$, and let $T$ be the largest quotient of $J_2$ on which $C$ acts by $\chi_C$. The kernel of $J_1\to U$, which is spanned by the vectors $(a,1)v-\chi_A(a)v$ with $a\in A$ and $v\in J_1$, is $(A\times B)$-stable, so $B$ acts on $U$. Likewise, $B$ acts on $T$.

\begin{proposition}\label{prop:U-T}
	The representations $U$ and $T$ of $B$ are finite-dimensional, and there are canonical isomorphisms
	$$
	\Hom_{A\times\Delta B\times C}\big(J_1\otimes J_2,\ \chi_A\otimes\one\otimes\chi_C\big)\simeq\Hom_{\Delta B}(U\otimes T,\C)\simeq\Hom_B(U,\widetilde T).
	$$
\end{proposition}

\begin{proof}
	Since $A$ acts on $U$ by scalars, the $B$-subrepresentations of $U$ are $(A\times B)$-subrepresentations. Since $U$ is a quotient of $J_1$, which has finite length, it follows that $U$ has finite length as a representation of $B$. Hence, $U$ is finite-dimensional, since every irreducible representation of $B$ is finite-dimensional. The same argument applies to $T$. By equivariance under $A$ and $C$, every linear form in the left-hand side vanishes on $\ker(J_1\to U)\otimes J_2$ and on $J_1\otimes\ker(J_2\to T)$. Hence, it descends to a $\Delta B$-invariant bilinear form $\beta$ on $U\times T$, and conversely, every such form pulls back to an element of the left-hand side. Finally, since $T$ is finite-dimensional, the map $\beta\mapsto(u\mapsto\beta(u,\cdot))$ identifies these forms with $\Hom_B(U,\widetilde T)$.
\end{proof}

The nonvanishing of this space can be detected by irreducible representations of $\mathrm{D}^\times$, as follows.

\begin{lemma}\label{lem:U-T-criterion}
	The space of Proposition~\ref{prop:U-T} is nonzero if and only if there is $\sigma\in\Irr(\mathrm{D}^\times)$ such that
	$$
	\Hom_{A\times B}(J_1,\chi_A\otimes\sigma)\ne0\qquad\text{and}\qquad\Hom_{B\times C}(J_2,\widetilde\sigma\otimes\chi_C)\ne0.
	$$
\end{lemma}

\begin{proof}
	By construction, we have $\Hom_{A\times B}(J_1,\chi_A\otimes\sigma)=\Hom_B(U,\sigma)$ and $\Hom_{B\times C}(J_2,\widetilde\sigma\otimes\chi_C)=\Hom_B(T,\widetilde\sigma)$. Suppose that $0\ne f\in\Hom_B(U,\widetilde T)$, and let $\sigma$ be an irreducible quotient of the finite-dimensional representation $f(U)$. Then $\sigma$ is a quotient of $U$ and a constituent of $\widetilde T$. Hence, $\widetilde\sigma$ is a constituent of $T$, and therefore a quotient of $T$ by Lemma~\ref{lem:Dx-quotient}. Conversely, surjections $p\colon U\to\sigma$ and $q\colon T\to\widetilde\sigma$ give the nonzero $\Delta B$-invariant bilinear form $\beta(u,t)=\langle p(u),q(t)\rangle$, where $\langle\,,\rangle$ is the canonical pairing between $\sigma$ and $\widetilde\sigma$.
\end{proof}

For the rest of this section, we put
$$
x_s=\frac{e(n-k-1)}2-es,\qquad y_s=\frac{e(1-k)}2-es;
$$
these are the exponents of $\kappa_s$, and they depend on $n$, $k$ and $s$. Since $\kappa_s$ is trivial on $\Delta B$, we may apply Proposition~\ref{prop:U-T} and Lemma~\ref{lem:U-T-criterion} with $J_1=r_{(k-1,1)}(\pi_1)$, $J_2=r_{(1,n-k-1)}(\pi_2)$, $\chi_A=\nu^{x_s}$ and $\chi_C=\nu^{y_s}$. Combined with Frobenius reciprocity, this gives the following description of condition (III).

\begin{proposition}\label{prop:open-condition}
	Let $\pi_1\in\Rep(G_k)$ and $\pi_2\in\Rep(G_{n-k})$. Then $\dIII=\dim\Hom_B(U,\widetilde T)$, and condition \textup{(III)} holds if and only if there is $\sigma\in\Irr(\mathrm{D}^\times)$ such that
	$$
	\Hom_{G_k}\big(\pi_1,\ \nu^{x_s}\times\sigma\big)\ne0\qquad\text{and}\qquad\Hom_{G_{n-k}}\big(\pi_2,\ \widetilde\sigma\times\nu^{y_s}\big)\ne0,
	$$
	where $\nu^{x_s}\times\sigma$ is induced from $G_{k-1}\times G_1$ and $\widetilde\sigma\times\nu^{y_s}$ is induced from $G_1\times G_{n-k-1}$. If $\pi_1$ and $\pi_2$ are irreducible, this means that $\pi_1\hookrightarrow\nu^{x_s}\times\sigma$ and $\pi_2\hookrightarrow\widetilde\sigma\times\nu^{y_s}$.
\end{proposition}

\begin{proof}
	The first assertion follows from Proposition~\ref{prop:contrib-III} and Proposition~\ref{prop:U-T}. By Lemma~\ref{lem:U-T-criterion}, condition (III) holds if and only if there is $\sigma\in\Irr(\mathrm{D}^\times)$ with $\Hom_{A\times B}(r_{(k-1,1)}(\pi_1),\nu^{x_s}\otimes\sigma)\ne0$ and $\Hom_{B\times C}(r_{(1,n-k-1)}(\pi_2),\widetilde\sigma\otimes\nu^{y_s})\ne0$. By Frobenius reciprocity, these spaces are isomorphic to $\Hom_{G_k}(\pi_1,\nu^{x_s}\times\sigma)$ and $\Hom_{G_{n-k}}(\pi_2,\widetilde\sigma\times\nu^{y_s})$, respectively. Finally, a nonzero homomorphism from an irreducible representation is injective.
\end{proof}

For $e\ge2$, this is where the inner form differs essentially from the split case: the representation $\sigma\in\Irr(\mathrm{D}^\times)$ is arbitrary and need not be a character.

When one of the two blocks has size one, $\dIII$ can be expressed by a single Hom space.

\begin{corollary}\label{cor:dIII-k1}
	Let $k=1$ and $\pi_1\in\Irr(G_1)$. Then $\dIII=\dim\Hom_{G_{n-1}}\big(\pi_2,\ \widetilde\pi_1\times\nu^{-es}\big)$, where $\widetilde\pi_1\times\nu^{-es}$ is induced from $G_1\times G_{n-2}$.
\end{corollary}

\begin{proof}
	Here $A=G_0$, so that $U=J_1=\pi_1$, and $y_s=-es$. Since $\pi_1$ and $T$ are finite-dimensional, we have $\Hom_B(\pi_1,\widetilde T)\simeq\Hom_B(T,\widetilde\pi_1)$. By construction of $T$, the latter space equals
	$$
	\Hom_{B\times C}(J_2,\widetilde\pi_1\otimes\nu^{-es})=\Hom_{G_1\times G_{n-2}}\big(r_{(1,n-2)}(\pi_2),\widetilde\pi_1\otimes\nu^{-es}\big),
	$$
	which is isomorphic to $\Hom_{G_{n-1}}(\pi_2,\widetilde\pi_1\times\nu^{-es})$ by Frobenius reciprocity.
\end{proof}

\begin{corollary}\label{cor:dIII-kn1}
	Let $k=n-1$ and $\pi_2\in\Irr(G_1)$. Then $\dIII=\dim\Hom_{G_{n-1}}\big(\pi_1,\ \nu^{-es}\times\widetilde\pi_2\big)$, where $\nu^{-es}\times\widetilde\pi_2$ is induced from $G_{n-2}\times G_1$.
\end{corollary}

\begin{proof}
	Here $C=G_0$, so that $T=J_2=\pi_2$, and $x_s=-es$. Thus, $\dIII=\dim\Hom_B(U,\widetilde\pi_2)=\dim\Hom_{A\times B}(J_1,\nu^{-es}\otimes\widetilde\pi_2)$, and Frobenius reciprocity gives the assertion.
\end{proof}

\subsection{Irreducible inducing data}\label{ss:irreducible-data}

To make Proposition~\ref{prop:open-condition} explicit for irreducible $\pi_1$ and $\pi_2$, we use Lemmas~\ref{lem:char-sigma} and~\ref{lem:sigma-char}. We first introduce notation for the possible irreducible subrepresentations that occur.

\begin{definition}\label{def:O-sets}
	Let $\sigma\in\Irr(G_1)$. For $k\ge2$, let $\mathcal O_1(s,k;\sigma)$ be the set consisting of the unique irreducible subrepresentation of the representation $\nu^{x_s}\times\sigma$ induced from $G_{k-1}\times G_1$, and put $\mathcal O_1(s,1;\sigma)=\{\sigma\}$. For $m=n-k\ge2$, let $\mathcal O_2(s,m;\sigma)$ be the set consisting of the unique irreducible subrepresentation of the representation $\widetilde\sigma\times\nu^{y_s}$ induced from $G_1\times G_{m-1}$, and put $\mathcal O_2(s,1;\sigma)=\{\widetilde\sigma\}$.
\end{definition}

These sets can be described explicitly. We apply Lemma~\ref{lem:char-sigma} with $m=k$ and $y=x_s$, using the identities
$$x_s+\tfrac{ek}2=\tfrac{e(n-1)}2-es,\qquad x_s-\tfrac{ek}2=\tfrac{e(n-2k-1)}2-es,\qquad x_s+\tfrac e2=\tfrac{e(n-k)}2-es.$$
We obtain, for $k\ge2$,
\begin{equation}\label{eq:O1}
	\mathcal O_1(s,k;\sigma)=\begin{cases}\{\nu^{\frac{e(n-k)}2-es}\}&\text{if }\sigma\simeq\nu^{\frac{e(n-1)}2-es},\\[2pt] \{\nu^{\frac{e(n-k)}2-es}\Qt_k\}&\text{if }\sigma\simeq\nu^{\frac{e(n-2k-1)}2-es},\\[2pt] \{\nu^{x_s}\times\sigma\}\ \text{(irreducible)}&\text{otherwise.}\end{cases}
\end{equation}
Similarly, we apply Lemma~\ref{lem:sigma-char} with $m=n-k$ and $y=y_s$. We have $\widetilde\sigma\simeq\nu^{y_s-\frac{e(n-k)}2}$ if and only if $\sigma\simeq\nu^{es+\frac{e(n-1)}2}$, and $\widetilde\sigma\simeq\nu^{y_s+\frac{e(n-k)}2}$ if and only if $\sigma\simeq\nu^{es-\frac{e(n-2k+1)}2}$. Moreover, $y_s-\frac e2=-\frac{ek}2-es$. Hence, for $n-k\ge2$,
\begin{equation}\label{eq:O2}
	\mathcal O_2(s,n-k;\sigma)=\begin{cases}\{\nu^{-\frac{ek}2-es}\}&\text{if }\sigma\simeq\nu^{es+\frac{e(n-1)}2},\\[2pt] \{\nu^{-\frac{ek}2-es}Q_{n-k}\}&\text{if }\sigma\simeq\nu^{es-\frac{e(n-2k+1)}2},\\[2pt] \{\widetilde\sigma\times\nu^{y_s}\}\ \text{(irreducible)}&\text{otherwise.}\end{cases}
\end{equation}

In particular, the sets $\mathcal O_1$ and $\mathcal O_2$ contain very few characters.

\begin{remark}\label{rem:O-characters}
	For $k\ge2$, the only character that can lie in a set $\mathcal O_1(s,k;\sigma)$ is $\nu^{\frac{e(n-k)}2-es}$. Similarly, for $n-k\ge2$, the only character that can lie in a set $\mathcal O_2(s,n-k;\sigma)$ is $\nu^{-\frac{ek}2-es}$. Indeed, twists of $\Qt_k$ and of $Q_{n-k}$ are not characters by Lemma~\ref{lem:Qm-not-char}. Moreover, an irreducible representation $\nu^{x}\times\sigma$ or $\widetilde\sigma\times\nu^y$ of $G_m$, with $m\ge2$, is not a character. To see this, note that by Lemma~\ref{lem:sum-multisegments}, its Zelevinsky multisegment is the sum of the multisegment of the character $\nu^x$ (respectively $\nu^y$) and that of $\sigma$ (respectively $\widetilde\sigma$). It therefore consists of two segments, whereas the Zelevinsky multisegment of a character consists of a single segment. Similarly, such an irreducible induced representation is never a twist of $Q_m$, $\Qt_m$ or $\St_2$. Indeed, the two segments of its Zelevinsky multisegment are not linked by Lemma~\ref{lem:two-segments}, whereas the two segments of the multisegments of $Q_m$, $\Qt_m$ and $\St_2$ are linked.
\end{remark}

Combining Theorem~\ref{thm:orbit-bound} with the above description, we obtain the following trichotomy for irreducible inducing data.

\begin{proposition}\label{prop:three-alternatives}
	Let $\pi_1\in\Irr(G_k)$ and $\pi_2\in\Irr(G_{n-k})$, and suppose that $\pi_1\times\pi_2$ is $(H,\chi_s)$-distinguished. Then at least one of the following holds:
	\begin{itemize}
		\item[(I)] $\pi_1\simeq\nu^{\frac{e(n-k-2)}2-es}$, and $\pi_2$ is $(H_{1,n-k-1},\chi_{s+\frac k2})$-distinguished;
		\item[(II)] $\pi_1$ is $(H_{1,k-1},\chi_{s-\frac{n-k}2})$-distinguished, and $\pi_2\simeq\nu^{\frac{e(2-k)}2-es}$;
		\item[(III)] there is $\sigma\in\Irr(\mathrm{D}^\times)$ such that $\pi_1\in\mathcal O_1(s,k;\sigma)$ and $\pi_2\in\mathcal O_2(s,n-k;\sigma)$.
	\end{itemize}
	Conversely, each of the conditions \textup{(I)} and \textup{(II)} implies that $\pi_1\times\pi_2$ is $(H,\chi_s)$-distinguished.
\end{proposition}

\begin{proof}
	By Theorem~\ref{thm:orbit-bound}, one of the conditions (I), (II) and (III) of \S\ref{ss:contributions} holds. Since $\pi_1$ and $\pi_2$ are irreducible, conditions (I) and (II) of \S\ref{ss:contributions} are exactly the first two alternatives. Suppose that condition (III) holds. Then Proposition~\ref{prop:open-condition} gives $\sigma\in\Irr(\mathrm{D}^\times)$ with $\pi_1\hookrightarrow\nu^{x_s}\times\sigma$ and $\pi_2\hookrightarrow\widetilde\sigma\times\nu^{y_s}$, where the factors $G_0$ are omitted when $k=1$ or $n-k=1$. Since $\pi_1$ and $\pi_2$ are irreducible and these induced representations have unique irreducible subrepresentations by Lemmas~\ref{lem:char-sigma} and~\ref{lem:sigma-char}, we obtain $\pi_1\in\mathcal O_1(s,k;\sigma)$ and $\pi_2\in\mathcal O_2(s,n-k;\sigma)$. The converse is the last assertion of Theorem~\ref{thm:orbit-bound}.
\end{proof}

We illustrate Proposition~\ref{prop:three-alternatives} in the smallest case.

\begin{example}\label{ex:n2-orbits}
	Let $n=2$ and $k=1$, and let $\sigma_1,\sigma_2\in\Irr(G_1)$. Then the three conditions read as follows:
	\begin{itemize}
		\item[(I)] $\sigma_1\simeq\nu^{-e/2-es}$ and $\sigma_2\simeq\nu^{e/2+es}$;
		\item[(II)] $\sigma_1\simeq\nu^{es-e/2}$ and $\sigma_2\simeq\nu^{e/2-es}$;
		\item[(III)] $\sigma_2\simeq\widetilde\sigma_1$, by Corollary~\ref{cor:dIII-k1}, where the factor $G_0$ is omitted.
	\end{itemize}
	In each of the three cases, we have $\sigma_2\simeq\widetilde\sigma_1$. Thus, if $\sigma_1\times\sigma_2$ is irreducible and $(H_{1,1},\chi_s)$-distinguished, then $\sigma_2\simeq\widetilde\sigma_1$. This is consistent with the classification of \S\ref{ss:rank2-classification}.
\end{example}

\section{The cases \texorpdfstring{$n=1$}{n=1} and \texorpdfstring{$n=2$}{n=2}}\label{sec:rank-two}

In this section, we treat the two smallest cases. The case $n=1$ is elementary, whereas the case $n=2$ requires a detailed analysis of twisted Jacquet modules.

\subsection{The case \texorpdfstring{$n=1$}{n=1}}\label{ss:n1}

For $n=1$, we have $H_{1,0}=G_1=\mathrm{D}^\times$ and $\chi_s=\nu^{es}$. In this case, distinction simply means isomorphism with a character.

\begin{proposition}\label{prop:n1}
	Let $s\in\R$ and $\pi\in\Irr(\mathrm{D}^\times)$. Then $\pi$ is $(H_{1,0},\chi_s)$-distinguished if and only if $\pi\simeq\nu^{es}$, and in that case $\dim\Hom_{\mathrm{D}^\times}(\pi,\nu^{es})=1$. Equivalently, $\pi$ is $(H_{1,0},\chi_s)$-distinguished if and only if $\varphi_{\pi}\simeq|\cdot|^{es}\otimes\Sp_e$.
\end{proposition}

\begin{proof}
	A nonzero $\mathrm{D}^\times$-homomorphism from the irreducible representation $\pi$ to the character $\nu^{es}$ is surjective, and hence an isomorphism. Moreover, the space of such homomorphisms is one-dimensional by Schur's lemma. By Lemma~\ref{lem:L-segment}, we have $\varphi_{\nu^{es}}=\varphi_{\Qq([s])}=|\cdot|^{es}\otimes\Sp_e$. Therefore, the last assertion follows from Lemma~\ref{lem:L-injective}.
\end{proof}

\subsection{Notation for \texorpdfstring{$G_2$}{G2}}\label{ss:G2-notation}

In the rest of this section, we put $G=G_2$, $P=P_{(1,1)}=MN$ and $M=H_{1,1}=\mathrm{D}^\times\times \mathrm{D}^\times$. For $a,b\in \mathrm{D}^\times$ and $x\in \mathrm{D}$, we write
$$m(a,b)=\begin{pmatrix}a&0\\0&b\end{pmatrix},\qquad n(x)=\begin{pmatrix}1&x\\0&1\end{pmatrix},\qquad w=\begin{pmatrix}0&1\\1&0\end{pmatrix},$$
so that $N=\{n(x):x\in\mathrm{D}\}\simeq \mathrm{D}$ and $m(a,b)^{-1}n(x)m(a,b)=n(a^{-1}xb)$. Let $\Delta=\{m(c,c):c\in \mathrm{D}^\times\}\simeq \mathrm{D}^\times$. We fix a Haar measure $dx$ on $\mathrm{D}$ and a Haar measure $d^\times a$ on $\mathrm{D}^\times$. By Lemma~\ref{lem:det}, applied to the map $x\mapsto ax(b^{-1})^{-1}$, we have $d(axb)=\nu(a)^e\nu(b)^e\,dx$ for $a,b\in \mathrm{D}^\times$. Moreover, $\delta_P(m(a,b))=\nu(a)^e\nu(b)^{-e}$ by Lemma~\ref{lem:modulus}.

Fix a nontrivial character $\psi_{\mathrm{F}}$ of $\mathrm{F}$, and put $\psi(x)=\psi_{\mathrm{F}}(\Trd x)$ and $\psi_y(x)=\psi(yx)$ for $x,y\in \mathrm{D}$. Since $\Trd(xy)=\Trd(yx)$, we have $\psi_y(x)=\psi(xy)$ and $\psi(cxc^{-1})=\psi(x)$ for $c\in\mathrm{D}^\times$. We shall use the following standard facts of local duality \cite[Ch.~II]{Wei74}. A \emph{lattice} is a compact open subgroup of $(\mathrm{D},+)$. For a lattice $L$, we put $L^\perp=\{y\in \mathrm{D}:\psi_y|_L=1\}$.
\begin{itemize}
	\item[(D1)] The map $y\mapsto\psi_y$ is a topological isomorphism of $\mathrm{D}$ onto its Pontryagin dual.
	\item[(D2)] If $L$ is a lattice, then $L^\perp$ is a lattice and $(L^\perp)^\perp=L$. If $L\subseteq L'$, then $L'^\perp\subseteq L^\perp$. If $L_0\subseteq L$ are lattices, then $y\mapsto\psi_y|_L$ identifies $L_0^\perp/L^\perp$ with the character group of $L/L_0$.
	\item[(D3)] For every lattice $L$ and every $y\in\mathrm{D}$, we have $\int_L\psi(yx)\,dx=\vol(L)\,\one_{L^\perp}(y)$.
	\item[(D4)] Every lattice is contained in a strictly larger one. The powers $\mathfrak p_{\mathrm{D}}^r$ of the maximal ideal of the maximal order of $\mathrm{D}$ form a basis of neighbourhoods of $0$ consisting of lattices. Moreover, the subgroups $U_r=1+\mathfrak p_{\mathrm{D}}^r$, $r\ge1$, form a basis of neighbourhoods of $1$ in $\mathrm{D}^\times$ consisting of compact open subgroups. Note that $U_r=1+L_r^\perp$ for $L_r=(\mathfrak p_{\mathrm{D}}^r)^\perp$, by (D2).
\end{itemize}
For $a,b\in\mathrm{D}^\times$, we have $\psi(a^{-1}xb)=\psi(ba^{-1}x)=\psi_{ba^{-1}}(x)$. By (D1), this equals $\psi(x)$ for all $x\in\mathrm{D}$ if and only if $a=b$. Hence, $\Delta$ is the stabilizer of $\psi$ in $M$.

Finally, we recall the definition of (twisted) Jacquet modules. Let $(\pi,V)$ be a smooth representation of a group containing $N$, and let $\chi$ be a character of $N\simeq \mathrm{D}$. We put
$$
V(N,\chi)=\operatorname{span}\{\pi(n)v-\chi(n)v:n\in N,\ v\in V\},\qquad V_{N,\chi}=V/V(N,\chi),
$$
and we let $p_\chi\colon V\to V_{N,\chi}$ denote the projection. For $\chi=1$, we write $V(N)$, $V_N$ and $p_1$. If $V$ is a representation of $P$, then $\pi_N=V_N$ is the (unnormalized) Jacquet module, which is a representation of $M$. Since $\Delta$ stabilizes $\psi$, the space $\pi_{N,\psi}=V_{N,\psi}$ is a representation of $\Delta\simeq \mathrm{D}^\times$, called the \emph{twisted Jacquet module} of $\pi$. By the definition of the normalized Jacquet module (\S\ref{ss:reps}), we have $\pi_N=\delta_P^{1/2}\otimes r_{(1,1)}(\pi)$.

\subsection{Coinvariants for \texorpdfstring{$N$}{N}}\label{ss:coinvariants}

In this subsection, we collect some elementary facts about (twisted) coinvariants for the group $N\simeq\mathrm{D}$, which we use repeatedly in the analysis of $G_2$. Let $(\pi,V)$ be a smooth representation of $N\simeq \mathrm{D}$. For a lattice $L$ and a (continuous) character $\chi$ of $\mathrm{D}$, we put
$$E_{L,\chi}v=\vol(L)^{-1}\int_L\chi(x)^{-1}\pi(n(x))v\,dx\qquad(v\in V).$$
Since $v$ is fixed by $n(L_1)$ for some lattice $L_1$, and $\chi$ is trivial on some lattice, the integrand is locally constant, and the integral is a finite sum.

We first show that the vanishing of $E_{L,\chi}v$ is preserved when $L$ is enlarged.

\begin{lemma}\label{lem:coinv-enlarge}
	Let $\chi$ be a character of $\mathrm{D}$, let $v\in V$, and let $L\subseteq L'$ be lattices. If $E_{L,\chi}v=0$, then $E_{L',\chi}v=0$.
\end{lemma}

\begin{proof}
	Since $L$ is open and $L'$ is compact, $L'$ is a finite disjoint union of cosets $x_i+L$. For each $i$, the substitution $x\mapsto x+x_i$ gives
	$$\int_{x_i+L}\chi(x)^{-1}\pi(n(x))v\,dx=\chi(x_i)^{-1}\pi(n(x_i))\,\vol(L)\,E_{L,\chi}v=0.$$
	Summing over $i$, we obtain $E_{L',\chi}v=0$.
\end{proof}

The operators $E_{L,\chi}$ detect the subspace $V(N,\chi)$, as follows.

\begin{lemma}\label{lem:coinv-kernel}
	Let $\chi$ be a character of $\mathrm{D}$ and $v\in V$. Then $v\in V(N,\chi)$ if and only if $E_{L,\chi}v=0$ for some lattice $L$.
\end{lemma}

\begin{proof}
	Suppose that $E_{L,\chi}v=0$. Then
	$$v=v-E_{L,\chi}v=\vol(L)^{-1}\int_L\chi(x)^{-1}\big(\chi(x)v-\pi(n(x))v\big)\,dx,$$
	which is a finite linear combination of elements of $V(N,\chi)$. Conversely, let $v=\sum_i\big(\pi(n(x_i))u_i-\chi(x_i)u_i\big)$, and choose a lattice $L$ containing all the $x_i$. Since $x_i+L=L$, the substitution $x\mapsto x-x_i$ gives $E_{L,\chi}\pi(n(x_i))u_i=\chi(x_i)E_{L,\chi}u_i$ for each $i$. Hence, $E_{L,\chi}v=0$.
\end{proof}

Next, a vector fixed by $n(L_0)$ decomposes according to the characters of $L/L_0$.

\begin{lemma}\label{lem:coinv-decomp}
	Let $v\in V$ be fixed by $n(L_0)$ for a lattice $L_0$, and let $L\supseteq L_0$ be a lattice. Then
	$$
	v=\sum_{y\in L_0^\perp/L^\perp}E_{L,\psi_y}v.
	$$
\end{lemma}

\begin{proof}
	The operator $E_{L,\psi_y}$ depends only on the class of $y$ modulo $L^\perp$, so the sum is well defined. By (D2), the functions $x\mapsto\psi(yx)$, with $y\in L_0^\perp/L^\perp$, are exactly the characters of $L/L_0$. Hence, by the orthogonality relations, $\sum_y\psi(-yx)=[L:L_0]\,\one_{L_0}(x)$ for $x\in L$. Since $\vol(L)=[L:L_0]\vol(L_0)$ and $v$ is fixed by $n(L_0)$, we obtain
	$$\sum_yE_{L,\psi_y}v=\frac{[L:L_0]}{\vol(L)}\int_{L_0}\pi(n(x))v\,dx=v.$$
\end{proof}

As a consequence of Lemma~\ref{lem:coinv-kernel}, the (twisted) coinvariant functors are exact.

\begin{lemma}\label{lem:coinv-exact}
	For every character $\chi$ of $\mathrm{D}$, the functor $V\mapsto V_{N,\chi}$ is exact on the category of smooth representations of $N$. Hence, it is exact on smooth representations of $P$ (for $\chi=1$) and on smooth representations of $\Delta N$ (for $\chi=\psi$).
\end{lemma}

\begin{proof}
	Right exactness is formal. For left exactness, let $U\subseteq V$ be a subrepresentation and let $u\in U\cap V(N,\chi)$. By Lemma~\ref{lem:coinv-kernel} applied in $V$, there is a lattice $L$ with $E_{L,\chi}u=0$. Since $E_{L,\chi}u$ is computed in $U$, the same lemma applied in $U$ shows that $u\in U(N,\chi)$. Hence, $U_{N,\chi}\to V_{N,\chi}$ is injective.
\end{proof}

The following lemma shows that a vector is determined by its images in all twisted Jacquet modules.

\begin{lemma}\label{lem:coinv-vanish}
	If $v\in V$ satisfies $p_{\psi_y}(v)=0$ for every $y\in \mathrm{D}$, then $v=0$.
\end{lemma}

\begin{proof}
	Let $L_0$ be a lattice such that $v$ is fixed by $n(L_0)$. For each $y\in L_0^\perp$, Lemmas~\ref{lem:coinv-kernel} and~\ref{lem:coinv-enlarge} give a lattice $L_y\supseteq L_0$ with $E_{L_y,\psi_y}v=0$. If $y'\in y+L_y^\perp$, then $\psi_{y'}=\psi_y$ on $L_y$, and therefore $E_{L_y,\psi_{y'}}v=0$. Since $L_0^\perp$ is compact and the sets $y+L_y^\perp$ are open, there are $y_1,\dots,y_r\in L_0^\perp$ with $L_0^\perp\subseteq\bigcup_{i=1}^r(y_i+L_{y_i}^\perp)$. Let $L=L_{y_1}+\dots+L_{y_r}$, which is a lattice containing $L_0$. By Lemma~\ref{lem:coinv-enlarge}, we have $E_{L,\psi_y}v=0$ for all $y\in L_0^\perp$. Therefore, Lemma~\ref{lem:coinv-decomp} gives $v=0$.
\end{proof}

Finally, we compute the coinvariants of the space of compactly supported functions on $\mathrm{D}$.

\begin{lemma}\label{lem:coinv-Cc}
	Let $W$ be a finite-dimensional vector space, and let $N$ act on $C_c^\infty(\mathrm{D},W)$ by $(n(y)\phi)(x)=\phi(x+y)$. Then
	$$C_c^\infty(\mathrm{D},W)(N)=\Big\{\phi:\int_{\mathrm{D}}\phi(x)\,dx=0\Big\},$$
	and the map $\phi\mapsto\int_{\mathrm{D}}\phi(x)\,dx$ induces an isomorphism $C_c^\infty(\mathrm{D},W)_N\simeq W$.
\end{lemma}

\begin{proof}
	Since the Haar measure is translation invariant, the integral vanishes on $C_c^\infty(\mathrm{D},W)(N)$. Conversely, suppose that $\int_{\mathrm{D}}\phi=0$, and choose a lattice $L$ and $x_0\in \mathrm{D}$ with $\operatorname{supp}\phi\subseteq x_0+L$. Then $(E_{L,1}\phi)(x)=\vol(L)^{-1}\int_{x+L}\phi$. If $x\in x_0+L$, then $x+L=x_0+L$, and this equals $\vol(L)^{-1}\int_{\mathrm{D}}\phi=0$; otherwise, $x+L$ does not meet the support of $\phi$, and it vanishes. Thus, $E_{L,1}\phi=0$, and $\phi\in C_c^\infty(\mathrm{D},W)(N)$ by Lemma~\ref{lem:coinv-kernel}. Finally, the map is surjective, since $\int_{\mathrm{D}}\vol(L)^{-1}\one_L\,w\,dx=w$ for every $w\in W$.
\end{proof}

\subsection{The Kirillov isomorphism}\label{ss:kirillov}

In this subsection, we describe the kernel $V(N)$ of the Jacquet module of a smooth representation of $P$ in terms of its twisted Jacquet module, in the spirit of the Kirillov model. Let $(\pi,V)$ be a smooth representation of $P$, and let $\sigma$ denote the representation of $\mathrm{D}^\times\simeq\Delta$ on $V_{N,\psi}$, so that $p_\psi(\pi(m(c,c))v)=\sigma(c)p_\psi(v)$ for $c\in\mathrm{D}^\times$ and $v\in V$. For $v\in V$, we define the \emph{Kirillov function}
$$\Kir_v\colon \mathrm{D}^\times\to V_{N,\psi},\qquad\Kir_v(a)=p_\psi\big(\pi(m(a,1))v\big).$$
Let $K_\sigma=C_c^\infty(\mathrm{D}^\times,V_{N,\psi})$, with the action of $P$ given by
\begin{equation}\label{eq:kirillov-action}
	(n(x)\Phi)(a)=\psi(ax)\Phi(a),\qquad(m(b,c)\Phi)(a)=\sigma(c)\Phi(c^{-1}ab).
\end{equation}
A direct computation, using $\psi(c^{-1}abx)=\psi(abxc^{-1})$, shows that this defines a representation of $P$.

The map $v\mapsto\Kir_v$ intertwines the action of $P$ on $V$ with \eqref{eq:kirillov-action}.

\begin{lemma}\label{lem:kirillov-equiv}
	For $v\in V$, $x\in \mathrm{D}$ and $a,b,c\in \mathrm{D}^\times$, we have
	$$\Kir_{\pi(n(x))v}(a)=\psi(ax)\Kir_v(a)\qquad\text{and}\qquad\Kir_{\pi(m(b,c))v}(a)=\sigma(c)\Kir_v(c^{-1}ab).$$
\end{lemma}

\begin{proof}
	We have $m(a,1)n(x)=n(ax)m(a,1)$ and $p_\psi\circ\pi(n(ax))=\psi(ax)p_\psi$, which gives the first identity. The second identity follows from $m(a,1)m(b,c)=m(c,c)m(c^{-1}ab,1)$.
\end{proof}

The values of $\Kir_v$ detect the twisted Jacquet modules with respect to all nontrivial characters of $N$.

\begin{lemma}\label{lem:kirillov-zero}
	For $y\in \mathrm{D}^\times$, the operator $\pi(m(y,1))$ maps $V(N,\psi_y)$ onto $V(N,\psi)$. In particular, for $v\in V$, we have $p_{\psi_y}(v)=0$ if and only if $\Kir_v(y)=0$.
\end{lemma}

\begin{proof}
	For $u\in V$ and $x\in\mathrm{D}$, put $u'=\pi(m(y,1))u$. Then
	$$\pi(m(y,1))\big(\pi(n(x))u-\psi_y(x)u\big)=\pi(n(yx))u'-\psi(yx)u'.$$
	Since $x\mapsto yx$ and $u\mapsto u'$ are bijective, the first assertion follows. Hence, $p_{\psi_y}(v)=0$ if and only if $\pi(m(y,1))v\in V(N,\psi)$, that is, if and only if $\Kir_v(y)=0$.
\end{proof}

We can now state the main result of this subsection.

\begin{theorem}[Kirillov isomorphism]\label{thm:kirillov}
	For every smooth representation $(\pi,V)$ of $P$, the map $v\mapsto\Kir_v$ restricts to an isomorphism of $P$-modules $V(N)\xrightarrow{\sim}K_\sigma$. Consequently, there is an exact sequence of $P$-modules
	$$0\to K_\sigma\to V\to V_N\to0.$$
\end{theorem}

\begin{proof}
	We first show that $\Kir_v\in K_\sigma$ for every $v\in V(N)$. Let $U\subseteq\cO_{\mathrm{D}}^\times$ be an open subgroup with $\pi(m(u,1))v=v$ for all $u\in U$. Then $\Kir_v(au)=\Kir_v(a)$ for $a\in\mathrm{D}^\times$ and $u\in U$, so $\Kir_v$ is locally constant. By Lemma~\ref{lem:coinv-kernel}, there is a lattice $L_1$ with $E_{L_1,1}v=0$, and Lemma~\ref{lem:kirillov-equiv} together with (D3) gives
	$$0=\Kir_{E_{L_1,1}v}(a)=\vol(L_1)^{-1}\int_{L_1}\psi(ax)\,dx\ \Kir_v(a)=\one_{L_1^\perp}(a)\Kir_v(a).$$
	So $\Kir_v$ vanishes on $L_1^\perp\cap\mathrm{D}^\times$. On the other hand, if $v$ is fixed by $n(L_2)$ for a lattice $L_2$, then $\Kir_v(a)=\psi(ax)\Kir_v(a)$ for all $x\in L_2$, so $\Kir_v$ vanishes outside $L_2^\perp$. The support of $\Kir_v$ therefore lies in the compact subset $L_2^\perp\smallsetminus L_1^\perp$ of $\mathrm{D}^\times$. Hence $\Kir_v\in K_\sigma$, and the map $v\mapsto\Kir_v$ is $P$-equivariant by Lemma~\ref{lem:kirillov-equiv}.
	
	The map is injective. If $v\in V(N)$ and $\Kir_v=0$, then $p_{\psi_y}(v)=0$ for all $y\ne0$ by Lemma~\ref{lem:kirillov-zero}. Also $p_1(v)=0$ because $v\in V(N)$, so $v=0$ by Lemma~\ref{lem:coinv-vanish}.
	
	For surjectivity, we first produce enough elements in the image. Let $\bar w\in V_{N,\psi}$, and choose $v\in V$ with $p_\psi(v)=\bar w$. As above, $\Kir_v$ is locally constant, and $\Kir_v(1)=\bar w$. By (D4) there is $r_0\ge1$ such that $\Kir_v=\bar w$ on $U_r=1+\mathfrak p_{\mathrm{D}}^r$ for every $r\ge r_0$. Fix such an $r$, and put $L=L_r$, so that $L^\perp=\mathfrak p_{\mathrm{D}}^r$. Since $1\notin\mathfrak p_{\mathrm{D}}^r$, we have $\psi|_L\ne1$. Put $v'=E_{L,\psi}v$. By Lemma~\ref{lem:kirillov-equiv} and (D3),
	\[
	\Kir_{v'}(a)=\vol(L)^{-1}\int_L\psi((a-1)x)\,dx\ \Kir_v(a)=\one_{\mathfrak p_{\mathrm{D}}^r}(a-1)\Kir_v(a)=\one_{U_r}(a)\,\bar w.
	\]
	Moreover, (D3) gives $p_1(v')=\vol(L)^{-1}\int_L\psi(x)^{-1}\,dx\ p_1(v)=0$, because $-1\notin L^\perp$, so $v'\in V(N)$. The image therefore contains $\one_{U_r}\bar w$ for every $\bar w\in V_{N,\psi}$ and every $r\ge r_0(\bar w)$. The image is $M$-stable, so applying $m(b,1)$ shows that it also contains the function $a\mapsto\one_{U_rb^{-1}}(a)\,\bar w$ for every $b\in\mathrm{D}^\times$.
	
	Now let $\Phi\in K_\sigma$. Since $\Phi$ is locally constant and compactly supported, its values span a finite-dimensional subspace of $V_{N,\psi}$. Choose a basis of this subspace, and take $r$ at least as large as the integers $r_0$ attached to its elements. After enlarging $r$ if necessary, $\Phi$ is invariant under left translation by $U_r$. Then $\Phi=\sum_i\one_{U_rb_i}\Phi(b_i)$ is a finite sum of functions of the form above, so $\Phi$ lies in the image. This proves the isomorphism. The exact sequence follows, since $V(N)$ is the kernel of the projection $V\to V_N$.
\end{proof}

Next, we compute the $M$-equivariant linear forms on $K_\sigma$.

\begin{proposition}\label{prop:kirillov-frobenius}
	Let $\sigma$ be a smooth representation of $\mathrm{D}^\times$ on a space $V_\sigma$, let $K_\sigma=C_c^\infty(\mathrm{D}^\times,V_\sigma)$ with the action \eqref{eq:kirillov-action} of $M$, and let $\chi=\chi_1\otimes\chi_2$ be a character of $M$. For $\lambda\in\Hom_{\mathrm{D}^\times}(\sigma,\chi_1\chi_2)$, we put
	\[
	L_\lambda(\Phi)=\int_{\mathrm{D}^\times}\chi_1(a)^{-1}\lambda(\Phi(a))\,d^\times a .
	\]
	Then $\lambda\mapsto L_\lambda$ is an isomorphism $\Hom_{\mathrm{D}^\times}(\sigma,\chi_1\chi_2)\xrightarrow{\sim}\Hom_M(K_\sigma,\chi)$.
\end{proposition}

\begin{proof}
	We first check that $L_\lambda$ is $M$-equivariant. We have $\lambda\circ\sigma(c)=\chi_1(c)\chi_2(c)\lambda$. Since $\mathrm{D}^\times$ is unimodular, the substitution $a=ca'b^{-1}$ preserves $d^\times a$. Hence
	\[
	L_\lambda(m(b,c)\Phi)=\int\chi_1(ca'b^{-1})^{-1}\chi_1(c)\chi_2(c)\lambda(\Phi(a'))\,d^\times a'=\chi_1(b)\chi_2(c)L_\lambda(\Phi),
	\]
	so $L_\lambda\in\Hom_M(K_\sigma,\chi)$.
	
	Fix $r$ so large that $\chi_1$ is trivial on $U=U_r=1+\mathfrak p_{\mathrm{D}}^r$, and note that $U$ is normal in $\mathrm{D}^\times$. Since $\chi_1$ is trivial on $U$, we have $\lambda(w)=\vol(U)^{-1}L_\lambda(\one_U\otimes w)$ for every $w\in V_\sigma$. So $L_\lambda=0$ forces $\lambda=0$, and the map $\lambda\mapsto L_\lambda$ is injective.
	
	For surjectivity, let $L\in\Hom_M(K_\sigma,\chi)$. For $w\in V_\sigma$, consider the linear form $T_w(f)=L(f\otimes w)$ on $C_c^\infty(\mathrm{D}^\times)$. Write $(R_bf)(a)=f(ab)$. Since $m(b,1)(f\otimes w)=(R_bf)\otimes w$, we have $T_w(R_bf)=\chi_1(b)T_w(f)$. Put $S_w(f)=T_w(f\chi_1)$. Since $(R_bf)\chi_1=\chi_1(b)^{-1}R_b(f\chi_1)$, we get
	\[
	S_w(R_bf)=\chi_1(b)^{-1}T_w\big(R_b(f\chi_1)\big)=S_w(f).
	\]
	So $S_w$ is a right-invariant linear form on $C_c^\infty(\mathrm{D}^\times)$. It is therefore a multiple $\lambda(w)$ of the Haar integral \cite[\S1.18]{BZ76}, that is,
	\[
	T_w(f)=\lambda(w)\int f(a)\chi_1(a)^{-1}\,d^\times a,
	\]
	where $\lambda(w)\in\C$ depends linearly on $w$. We show that $\lambda$ is $\mathrm{D}^\times$-equivariant. Because $U$ is normal, $m(c,c)(\one_U\otimes w)=\one_U\otimes\sigma(c)w$. Applying $L$ to this identity and using $L(\one_U\otimes w)=\vol(U)\lambda(w)$ gives $\lambda(\sigma(c)w)=\chi_1(c)\chi_2(c)\lambda(w)$. Thus $\lambda\in\Hom_{\mathrm{D}^\times}(\sigma,\chi_1\chi_2)$, and $L$ agrees with $L_\lambda$ on every function $f\otimes w$. Since these functions span $K_\sigma$, it follows that $L=L_\lambda$.
\end{proof}

Combining the two previous results, we obtain the following description of $M$-equivariant linear forms on $\pi(N)$.

\begin{corollary}\label{cor:kirillov-hom}
	Let $\pi$ be a smooth representation of $G_2$, and let $\chi=\chi_1\otimes\chi_2$ be a character of $M$. Then
	$$\Hom_M(\pi(N),\chi)\simeq\Hom_{\mathrm{D}^\times}(\pi_{N,\psi},\chi_1\chi_2).$$
\end{corollary}

\begin{proof}
	By Theorem~\ref{thm:kirillov}, applied to the restriction of $\pi$ to $P$, we have $\pi(N)\simeq K_\sigma$ as representations of $P$, and hence of $M$, with $\sigma=\pi_{N,\psi}$. The assertion now follows from Proposition~\ref{prop:kirillov-frobenius}.
\end{proof}

\subsection{Independence of the character}\label{ss:independence}

In this subsection, we show that, for infinite-dimensional irreducible representations of $G_2$, the space of $(H_{1,1},\chi_t)$-equivariant linear forms does not depend on $t$. We first need a splitting result for extensions of finite-dimensional representations of $M$.

\begin{lemma}\label{lem:ext-splitting}
	Let $V$ be a finite-dimensional smooth representation of $M$, and let $\chi$ be a character of $M$ that is not a constituent of $V$. Then $\Hom_M(V,\chi)=0$, and every exact sequence $0\to\chi\to E\to V\to0$ of smooth representations of $M$ splits.
\end{lemma}

\begin{proof}
	The first assertion is clear, since the image of a nonzero homomorphism $V\to\chi$ would be a constituent of $V$ isomorphic to $\chi$. For the second assertion, we argue by induction on the length of $V$. Let $Z_M=\mathrm{F}^\times\times \mathrm{F}^\times$ denote the centre of $M$.
	
	Suppose first that $V=\rho$ is irreducible, so that $\rho\not\simeq\chi$. If $\omega_\rho\ne\chi|_{Z_M}$, choose $z\in Z_M$ with $\omega_\rho(z)\ne\chi(z)$. Then the generalized $\omega_\rho(z)$-eigenspace of $E(z)$ is $M$-stable, and it is a complement of $\chi$ in $E$. If $\omega_\rho=\chi|_{Z_M}$, then, for $z\in Z_M$, the endomorphism $E(z)-\chi(z)$ vanishes on $\chi$ and maps $E$ into $\chi$. Hence, it induces an $M$-homomorphism $\rho\to\chi$, which is zero since $\rho\not\simeq\chi$. Therefore, $Z_M$ acts trivially on $E\otimes\chi^{-1}$, and $E\otimes\chi^{-1}$ is a finite-dimensional smooth representation of the compact group $M/Z_M$. Consequently, $E\otimes\chi^{-1}$ is semisimple, and so is $E$.
	
	In general, choose an exact sequence $0\to V'\to V\to\rho\to0$ with $\rho$ irreducible, and let $E'\subseteq E$ be the preimage of $V'$. Since $\chi$ is not a constituent of $V'$, the induction hypothesis gives a decomposition $E'=\chi\oplus C'$ of representations of $M$. The quotient $E/C'$ is an extension of $\rho$ by $E'/C'\simeq\chi$, which splits by the first part of the proof. Hence, the composition $E\to E/C'\to\chi$ is a retraction of the inclusion $\chi\hookrightarrow E$.
\end{proof}

For $t\in\R$, the character $\chi_t=\nu^{et}\otimes\nu^{-et}$ of $M=H_{1,1}$ is trivial on $\Delta$. For $\tau\in\Irr(G_2)$, let $E(\tau)$ denote the set of $t\in\R$ such that $\chi_t$ is a constituent of the finite-dimensional representation $\tau_N$ of $M$. We can now prove the main result of this subsection.

\begin{theorem}\label{thm:independence}
	Let $\tau\in\Irr(G_2)$ be infinite-dimensional. Then, for every $t\in\R$,
	$$\dim\Hom_{H_{1,1}}(\tau,\chi_t)=\dim\Hom_{\mathrm{D}^\times}(\tau_{N,\psi},\one)=\dim\Hom_{H_{1,1}}(\tau,\one)\le1.$$
	In particular, $\tau$ is $(H_{1,1},\chi_t)$-distinguished if and only if it is $H_{1,1}$-distinguished.
\end{theorem}

\begin{proof}
	We first show that $\Hom_M(\tau,\chi_t)\simeq\Hom_{\mathrm{D}^\times}(\tau_{N,\psi},\one)$ whenever $t\notin E(\tau)$. Applying $\Hom_M(-,\chi_t)$ to the exact sequence $0\to\tau(N)\to\tau\to\tau_N\to0$ gives an exact sequence
	$$0\to\Hom_M(\tau_N,\chi_t)\to\Hom_M(\tau,\chi_t)\to\Hom_M(\tau(N),\chi_t),$$
	and its first term vanishes by Lemma~\ref{lem:ext-splitting}. Now let $f\in\Hom_M(\tau(N),\chi_t)$. The push-out of $0\to\tau(N)\to\tau\to\tau_N\to0$ along $f$ is an extension of $\tau_N$ by $\chi_t$, and it splits by Lemma~\ref{lem:ext-splitting}. Hence $f$ extends to an element of $\Hom_M(\tau,\chi_t)$, so restriction is an isomorphism $\Hom_M(\tau,\chi_t)\simeq\Hom_M(\tau(N),\chi_t)$. Since $\nu^{et}\nu^{-et}=\one$, Corollary~\ref{cor:kirillov-hom} identifies the latter space with $\Hom_{\mathrm{D}^\times}(\tau_{N,\psi},\one)$.
	
	Next, we show that $E(\tau)$ never contains both $t$ and $-t$, and that $0\notin E(\tau)$. If $\tau$ is supercuspidal, then $\tau_N=0$ and $E(\tau)=\emptyset$. Otherwise $r_{(1,1)}(\tau)\ne0$. By Frobenius reciprocity, $\tau$ embeds into $\sigma_1\times\sigma_2$ for some irreducible quotient $\sigma_1\otimes\sigma_2$ of $r_{(1,1)}(\tau)$, where $\sigma_1,\sigma_2\in\Irr(G_1)$. By the geometric lemma \cite[Theorem~5.2]{BZ77}, \cite{Tad90}, the constituents of $r_{(1,1)}(\sigma_1\times\sigma_2)$ are $\sigma_1\otimes\sigma_2$ and $\sigma_2\otimes\sigma_1$, so by exactness every constituent of $r_{(1,1)}(\tau)$ is one of these. Since $\tau_N=\delta_P^{1/2}\otimes r_{(1,1)}(\tau)$ and $\delta_P^{1/2}=\nu^{e/2}\otimes\nu^{-e/2}$, the condition $t\in E(\tau)$ forces
	$$\{\sigma_1,\sigma_2\}=\{\nu^{e(t-\frac12)},\nu^{e(\frac12-t)}\}.$$
	If both $t$ and $-t$ lay in $E(\tau)$, then $\{t-\frac12,\frac12-t\}=\{-t-\frac12,\frac12+t\}$, which forces $t=0$. It therefore suffices to rule out $0\in E(\tau)$. Suppose that $0\in E(\tau)$. Then $\tau$ is an infinite-dimensional constituent of $\nu^{-e/2}\times\nu^{e/2}$, so $\tau=\St_2$ by Lemma~\ref{lem:St2}. The other constituent is $\one_2$, and $r_{(1,1)}(\one_2)=\nu^{-e/2}\otimes\nu^{e/2}$ by Lemma~\ref{lem:jacquet-char}. Hence $r_{(1,1)}(\St_2)=\nu^{e/2}\otimes\nu^{-e/2}$ and $(\St_2)_N=\nu^e\otimes\nu^{-e}=\chi_1$. So $E(\St_2)=\{1\}$, which contradicts $0\in E(\tau)$.
	
	We now combine these facts. By Lemma~\ref{lem:swap} with $p=q=1$, we have $\dim\Hom_M(\tau,\chi_t)=\dim\Hom_M(\tau,\chi_{-t})$. We have shown that $0\notin E(\tau)$ and that at least one of $t$ and $-t$ lies outside $E(\tau)$. By the first paragraph, the spaces $\Hom_M(\tau,\chi_t)$, $\Hom_M(\tau,\chi_{-t})$ and $\Hom_M(\tau,\chi_0)$ therefore all have the same dimension as $\Hom_{\mathrm{D}^\times}(\tau_{N,\psi},\one)$. Finally, $\dim\Hom_{H_{1,1}}(\tau,\one)\le1$ by Theorem~\ref{thm:ALMSY-mult}.
\end{proof}

\subsection{Twisted Jacquet modules}\label{ss:twisted-jacquet}
In this subsection, we compute the twisted Jacquet modules of the non-supercuspidal irreducible representations of $G_2$, beginning with those of representations parabolically induced from $P$, which are given by \cite[Theorem 2.1]{PR00}. For more details, see \cite{Pra99}.

\begin{lemma}\label{lem:tw-jacquet-ps}
	Let $\sigma_1,\sigma_2\in\Irr(\mathrm{D}^\times)$ act on spaces $W_1$ and $W_2$, respectively, and put $W=W_1\otimes W_2$. Then $(\sigma_1\times\sigma_2)_{N,\psi}\simeq(W,\sigma_1\otimes\sigma_2|_\Delta)$ as representations of $\mathrm{D}^\times$, where $c\in\mathrm{D}^\times$ acts on $W$ by $\sigma_1(c)\otimes\sigma_2(c)$.
\end{lemma}

We next recall the reducibility of the representations $\sigma_1\times\sigma_2$. For $\sigma\in\Irr(\mathrm{D}^\times)$, the following facts are due to Tadi\'c \cite{Tad90}. They also follow from Lemma~\ref{lem:two-segments}, applied to the segments $[\pm\frac12]_{(\sigma)}$, exactly as in the proof of Lemma~\ref{lem:St2}.
\begin{itemize}
	\item[(T1)] For $\sigma_1,\sigma_2\in\Irr(\mathrm{D}^\times)$, the representation $\sigma_1\times\sigma_2$ is reducible if and only if $\sigma_2\simeq\nu_{\sigma_1}^{\pm1}\sigma_1$.
	\item[(T2)] There are exact sequences
	$$
	0\to\St_2(\sigma)\to\nu_\sigma^{1/2}\sigma\times\nu_\sigma^{-1/2}\sigma\to\Sp_2(\sigma)\to0,\qquad0\to\Sp_2(\sigma)\to\nu_\sigma^{-1/2}\sigma\times\nu_\sigma^{1/2}\sigma\to\St_2(\sigma)\to0,
	$$
	and $\St_2(\sigma)\not\simeq\Sp_2(\sigma)$.
	\item[(T3)] If $\sigma=\chi$ is a character, then $\nu_\chi=\nu^e$, $\Sp_2(\chi)=\chi\cdot\one_2$ is a character of $G_2$, and $\St_2(\chi)=\chi\St_2$.
\end{itemize}

Combining (T2) with Lemma~\ref{lem:tw-jacquet-ps}, we obtain the twisted Jacquet modules of $\St_2(\sigma)$ and $\Sp_2(\sigma)$.

\begin{proposition}\label{prop:tw-jacquet-St-Sp}
	Let $\sigma\in\Irr(\mathrm{D}^\times)$. Then there is an isomorphism of representations of $\mathrm{D}^\times$
	$$(\sigma\otimes\sigma)|_\Delta\simeq\St_2(\sigma)_{N,\psi}\oplus\Sp_2(\sigma)_{N,\psi},
	$$
	under which $\St_2(\sigma)_{N,\psi}$ corresponds to the image of the injection $\St_2(\sigma)_{N,\psi}\hookrightarrow(\nu_\sigma^{1/2}\sigma\times\nu_\sigma^{-1/2}\sigma)_{N,\psi}=(\sigma\otimes\sigma)|_\Delta$ induced by (T2). If $\sigma$ is a character, then $\Sp_2(\sigma)_{N,\psi}=0$ and $\St_2(\sigma)_{N,\psi}\simeq\sigma^2$.
\end{proposition}

\begin{proof}
	We apply the exact functor $(\cdot)_{N,\psi}$ (Lemma~\ref{lem:coinv-exact}) to the first exact sequence of (T2). By Lemma~\ref{lem:tw-jacquet-ps}, the middle term becomes $\sigma\otimes\sigma$, on which $c\in\mathrm{D}^\times$ acts by $\nu_\sigma^{1/2}(c)\sigma(c)\otimes\nu_\sigma^{-1/2}(c)\sigma(c)=\sigma(c)\otimes\sigma(c)$. The centre $\mathrm{F}^\times$ acts on $\sigma\otimes\sigma$ by the character $\omega_\sigma^2$. Hence, the resulting exact sequence splits by Lemma~\ref{lem:Dx-semisimple}. If $\sigma$ is a character, then $\Sp_2(\sigma)$ is a character of $G_2$ by (T3). Therefore, $N$ acts trivially on $\Sp_2(\sigma)$, and $\Sp_2(\sigma)_{N,\psi}=0$. Consequently, $\St_2(\sigma)_{N,\psi}\simeq(\sigma\otimes\sigma)|_\Delta=\sigma^2$.
\end{proof}

To exploit Proposition~\ref{prop:tw-jacquet-St-Sp}, we need the invariant bilinear forms on tensor products of irreducible representations of $\mathrm{D}^\times$.

\begin{lemma}\label{lem:invariant-forms}
	Let $\sigma_1,\sigma_2\in\Irr(\mathrm{D}^\times)$, with $\mathrm{D}^\times$ acting diagonally on $\sigma_1\otimes\sigma_2$. Then $\dim\Hom_{\mathrm{D}^\times}(\sigma_1\otimes\sigma_2,\one)$ is $1$ if $\sigma_2\simeq\widetilde\sigma_1$, and $0$ otherwise. If $\sigma\simeq\widetilde\sigma$, then a nonzero element of $\Hom_{\mathrm{D}^\times}(\sigma\otimes\sigma,\one)$ is a nondegenerate invariant bilinear form $B_\sigma$ on $\sigma$, which is unique up to scalars.
\end{lemma}

\begin{proof}
	Since $\sigma_1$ is finite-dimensional, a linear form on $\sigma_1\otimes\sigma_2$ is the same as a linear map $\sigma_2\to\sigma_1^*=\widetilde\sigma_1$, and invariance of the form corresponds to equivariance of the map. Hence, the first assertion follows from Schur's lemma. If $\sigma\simeq\widetilde\sigma$, then the map $\sigma\to\widetilde\sigma$ corresponding to a nonzero invariant form is a nonzero homomorphism between irreducible representations, and hence an isomorphism. Thus, $B_\sigma$ is nondegenerate.
\end{proof}

The finite-dimensional irreducible representations of $G_2$ are easily described; see \cite[Proposition 3.2.1]{Rag99}.

\begin{lemma}\label{lem:G2-finite-dim}
	Let $\pi\in\Irr(G_2)$ be finite-dimensional. Then $\pi$ is a character $\chi\circ\Nrd$ of $G_2$. In particular, $N$ acts trivially on $\pi$, $\pi_N=\pi$, and $\pi_{N,\psi}=0$.
\end{lemma}

Finally, we list the irreducible representations of $G_2$.

\begin{lemma}\label{lem:G2-kinds}
	Every $\pi\in\Irr(G_2)$ is isomorphic to a representation of one of the following kinds:
	\begin{enumerate}
		\item a supercuspidal representation;
		\item an irreducible induced representation $\sigma_1\times\sigma_2$ with $\sigma_1,\sigma_2\in\Irr(\mathrm{D}^\times)$;
		\item $\St_2(\sigma)$ or $\Sp_2(\sigma)$ for some $\sigma\in\Irr(\mathrm{D}^\times)$.
	\end{enumerate}
\end{lemma}

\begin{proof}
	Suppose that $\pi$ is not supercuspidal. Then $r_{(1,1)}(\pi)\ne0$ has an irreducible quotient $\sigma_1\otimes\sigma_2$, and Frobenius reciprocity gives an embedding $\pi\hookrightarrow\sigma_1\times\sigma_2$. If $\sigma_1\times\sigma_2$ is irreducible, we are done. Otherwise, $\sigma_2\simeq\nu_{\sigma_1}^{\pm1}\sigma_1$ by (T1). Put $\sigma=\nu_{\sigma_1}^{\pm1/2}\sigma_1$, so that $\nu_\sigma=\nu_{\sigma_1}$. Then $\sigma_1\times\sigma_2=\nu_\sigma^{\mp1/2}\sigma\times\nu_\sigma^{\pm1/2}\sigma$, whose constituents are $\St_2(\sigma)$ and $\Sp_2(\sigma)$ by (T2).
\end{proof}

\subsection{The classification for \texorpdfstring{$n=2$}{n=2}}\label{ss:rank2-classification}

We can now classify the $H_{1,1}$-distinguished irreducible representations of $G_2$. For every infinite-dimensional $\pi\in\Irr(G_2)$, Theorem~\ref{thm:independence} gives
\begin{equation}\label{eq:shalika}
	\dim\Hom_{H_{1,1}}(\pi,\chi_t)=\dim\Hom_{H_{1,1}}(\pi,\one)=\dim\Hom_{\mathrm{D}^\times}(\pi_{N,\psi},\one)\qquad(t\in\R).
\end{equation}
Combining \eqref{eq:shalika} with the computation of twisted Jacquet modules in \S\ref{ss:twisted-jacquet}, we obtain the following result.

\begin{theorem}\label{thm:rank2-classification}
	An irreducible representation $\pi$ of $G_2$ is $H_{1,1}$-distinguished if and only if it is isomorphic to one of the following:
	\begin{enumerate}[label=\textup{(\alph*)}]
		\item the trivial representation $\one_2$;
		\item $\chi\St_2$, where $\chi$ is a character of $\mathrm{D}^\times$ of order at most two;
		\item $\Sp_2(\sigma)$, where $\sigma\in\Irr(\mathrm{D}^\times)$ is self-dual with $\dim\sigma>1$, and the twisted Jacquet module $\Sp_2(\sigma)_{N,\psi}$ admits a nonzero $\mathrm{D}^\times$-invariant linear form;
		\item $\St_2(\sigma)$, where $\sigma\in\Irr(\mathrm{D}^\times)$ is self-dual with $\dim\sigma>1$, and the twisted Jacquet module $\St_2(\sigma)_{N,\psi}$ admits a nonzero $\mathrm{D}^\times$-invariant linear form;
		\item $\sigma\times\widetilde\sigma$, where $\sigma\in\Irr(\mathrm{D}^\times)$ satisfies $\widetilde\sigma\not\simeq\nu_\sigma^{\pm1}\sigma$, that is, $\sigma\times\widetilde\sigma$ is irreducible;
		\item a supercuspidal representation whose twisted Jacquet module admits a nonzero $\mathrm{D}^\times$-invariant linear form.
	\end{enumerate}
	The representations in \textup{(b)--(f)} are exactly the infinite-dimensional ones in this list. Moreover, let $\sigma\in\Irr(\mathrm{D}^\times)$ with $\dim\sigma>1$. If $\sigma$ is self-dual, then exactly one of $\Sp_2(\sigma)$ and $\St_2(\sigma)$ is $H_{1,1}$-distinguished; otherwise, neither of them is.
\end{theorem}

\begin{proof}
	We go through the kinds of representations listed in Lemma~\ref{lem:G2-kinds}, starting with the finite-dimensional ones. If $\pi$ is finite-dimensional, then $\pi$ is a character by Lemma~\ref{lem:G2-finite-dim}, and by Lemma~\ref{lem:char-dist} with $s=0$, it is distinguished if and only if $\pi\simeq\one_2$. This gives type~(a). From now on, $\pi$ is infinite-dimensional, so \eqref{eq:shalika} applies.
	
	A supercuspidal representation is infinite-dimensional, since a character of $G_2$ has a nonzero Jacquet module. So, for supercuspidal $\pi$, \eqref{eq:shalika} shows that $\pi$ is distinguished if and only if $\Hom_{\mathrm{D}^\times}(\pi_{N,\psi},\one)\ne0$, that is, if and only if the trivial character of $\mathrm{D}^\times$ is a quotient of $\pi_{N,\psi}$. This gives type~(f).
	
	Next, let $\pi=\sigma_1\times\sigma_2$ be irreducible. By Lemma~\ref{lem:tw-jacquet-ps}, $\pi_{N,\psi}\simeq(\sigma_1\otimes\sigma_2)|_\Delta\ne0$, so $\pi$ is infinite-dimensional by Lemma~\ref{lem:G2-finite-dim}. By \eqref{eq:shalika} and Lemma~\ref{lem:invariant-forms}, $\pi$ is distinguished if and only if $\sigma_2\simeq\widetilde\sigma_1$. By (T1), $\sigma\times\widetilde\sigma$ is irreducible if and only if $\widetilde\sigma\not\simeq\nu_\sigma^{\pm1}\sigma$, which gives type~(e).
	
	Now let $\chi$ be a character of $\mathrm{D}^\times$. By (T3), $\Sp_2(\chi)$ is a character, and these have already been treated. By Proposition~\ref{prop:tw-jacquet-St-Sp}, $(\chi\St_2)_{N,\psi}\simeq\chi^2\ne0$. So $\chi\St_2$ is infinite-dimensional, and by \eqref{eq:shalika}, it is distinguished if and only if $\chi^2=\one_1$. This gives type~(b).
	
	Finally, let $\sigma\in\Irr(\mathrm{D}^\times)$ with $\dim\sigma>1$. By (T2) and the geometric lemma, the Jacquet modules of $\St_2(\sigma)$ and $\Sp_2(\sigma)$ are irreducible representations of $M$ of dimension $(\dim\sigma)^2>1$. Hence neither representation is a character, and both are infinite-dimensional by Lemma~\ref{lem:G2-finite-dim}. By \eqref{eq:shalika}, Proposition~\ref{prop:tw-jacquet-St-Sp} and Lemma~\ref{lem:invariant-forms},
	$$\dim\Hom_{H_{1,1}}(\St_2(\sigma),\one)+\dim\Hom_{H_{1,1}}(\Sp_2(\sigma),\one)=\dim\Hom_{\mathrm{D}^\times}(\sigma\otimes\sigma,\one),$$
	and the right-hand side is $1$ if $\sigma\simeq\widetilde\sigma$, and $0$ otherwise. So if $\sigma$ is not self-dual, neither representation is distinguished. If $\sigma$ is self-dual, exactly one of them is distinguished, and by \eqref{eq:shalika}, it is the one whose twisted Jacquet module has the trivial character as a quotient. This gives types~(c) and~(d).
	
	By Lemma~\ref{lem:G2-kinds}, every $\pi\in\Irr(G_2)$ has now been examined. Along the way, we have seen that the representation of type~(a) is one-dimensional, and that the representations of types~(b)--(f) are infinite-dimensional.
\end{proof}

We give a name to the set of representations in \textup{(b)--(f)}.

\begin{definition}\label{def:T}
	We denote by $\T$ the set of infinite-dimensional $H_{1,1}$-distinguished representations in $\Irr(G_2)$. By Theorem~\ref{thm:rank2-classification}, $\T$ consists of the representations of types \textup{(b)--(f)}.
\end{definition}

By Theorem~\ref{thm:independence}, membership in $\T$ detects $(H_{1,1},\chi_t)$-distinction for every $t$.

\begin{corollary}\label{cor:T-chi-t}
	Let $\tau\in\Irr(G_2)$ be infinite-dimensional, and let $t\in\R$. Then $\tau$ is $(H_{1,1},\chi_t)$-distinguished if and only if $\tau\in\T$.
\end{corollary}

\begin{proof}
	This follows from Theorem~\ref{thm:independence} and Definition~\ref{def:T}.
\end{proof}

We record two families of elements of $\T$ that will be used repeatedly.

\begin{corollary}\label{cor:T-members}
	We have $\St_2\in\T$, and $\tau_a=\nu^{ea}\times\nu^{-ea}\in\T$ for every $a\in\R\smallsetminus\{\pm\frac12\}$.
\end{corollary}

\begin{proof}
	The representation $\St_2$ is of type (b) with $\chi=\one_1$. For $a\ne\pm\frac12$, the representation $\tau_a=\sigma\times\widetilde\sigma$, with $\sigma=\nu^{ea}$, is irreducible (\S\ref{ss:segments}), and hence of type (e).
\end{proof}

\subsection{Multiplicity one for \texorpdfstring{$n=2$}{n=2}}\label{ss:rank2-mult}

We now prove multiplicity one by treating each type of Theorem~\ref{thm:rank2-classification} separately. For infinite-dimensional $\pi$, it suffices by \eqref{eq:shalika} to consider the trivial character. For types \textup{(b)--(e)} of that theorem, multiplicity one follows directly from the computation of twisted Jacquet modules. Theorem~\ref{thm:ALMSY-mult} is needed only for supercuspidal representations.

\begin{proposition}\label{prop:rank2-mult}
	Let $t\in\R$ and $\pi\in\Irr(G_2)$. Then $\dim\Hom_{H_{1,1}}(\pi,\chi_t)\le1$. More precisely, the following hold.
	\begin{enumerate}[label=\textup{(\alph*)}]
		\item If $\pi$ is a character, then $\dim\Hom_{H_{1,1}}(\pi,\chi_t)\le1$, with equality if and only if $\pi=\one_2$ and $t=0$.
		\item If $\pi=\chi\St_2$ with $\chi^2=\one_1$, then $\dim\Hom_{H_{1,1}}(\pi,\chi_t)=\dim\Hom_{\mathrm{D}^\times}(\chi^2,\one)=1$.
		\item If $\sigma$ is self-dual with $\dim\sigma>1$, then
		$$\dim\Hom_{H_{1,1}}(\St_2(\sigma),\chi_t)+\dim\Hom_{H_{1,1}}(\Sp_2(\sigma),\chi_t)=\dim\Hom_{\mathrm{D}^\times}(\sigma\otimes\sigma,\one)=1.$$
		In particular, each of these dimensions is at most one.
		\item If $\pi=\sigma\times\widetilde\sigma$ is irreducible, then $\dim\Hom_{H_{1,1}}(\pi,\chi_t)=\dim\Hom_{\mathrm{D}^\times}(\sigma\otimes\widetilde\sigma,\one)=1$.
		\item If $\pi$ is supercuspidal, then $\dim\Hom_{H_{1,1}}(\pi,\chi_t)=\dim\Hom_{H_{1,1}}(\pi,\one)\le1$.
	\end{enumerate}
\end{proposition}

\begin{proof}
	If $\pi$ is a character, then $\pi$ is one-dimensional, so $\Hom_{H_{1,1}}(\pi,\chi_t)$ has dimension at most one, and Lemma~\ref{lem:char-dist} determines when it is nonzero. This proves~(a). Assertions (b)--(d) follow from \eqref{eq:shalika}, Proposition~\ref{prop:tw-jacquet-St-Sp}, Lemma~\ref{lem:tw-jacquet-ps} and Lemma~\ref{lem:invariant-forms}, as in the proof of Theorem~\ref{thm:rank2-classification}. Assertion~(e) follows from \eqref{eq:shalika} and Theorem~\ref{thm:ALMSY-mult} with $p=q=1$. Finally, suppose that $\pi$ is not of any of the types (a)--(e). Then $\pi$ is infinite-dimensional and not $H_{1,1}$-distinguished by Theorem~\ref{thm:rank2-classification}, so $\Hom_{H_{1,1}}(\pi,\chi_t)=0$ by \eqref{eq:shalika}.
\end{proof}

\subsection{Proof of Theorem~\ref{thm:B}}\label{ss:rank2-main}

We now combine the results of this section. The following theorem summarizes the case $n=2$.

\begin{theorem}\label{thm:main-n2}
	Let $s\in\R$ and $\pi\in\Irr(G_2)$.
	\begin{enumerate}
		\item The representation $\pi$ is $(H_{1,1},\chi_s)$-distinguished if and only if either $\pi\simeq\one_2$ and $s=0$, or $\pi\in\T$.
		\item We have $\dim\Hom_{H_{1,1}}(\pi,\chi_s)\le1$.
	\end{enumerate}
\end{theorem}

\begin{proof}
	(1) If $\pi$ is finite-dimensional, then $\pi$ is a character by Lemma~\ref{lem:G2-finite-dim}, and the assertion follows from Lemma~\ref{lem:char-dist}. If $\pi$ is infinite-dimensional, then the assertion follows from Corollary~\ref{cor:T-chi-t}.
	
	(2) This is Proposition~\ref{prop:rank2-mult}.
\end{proof}

Theorem~\ref{thm:B} is now proved. Indeed, part (i) of Theorem~\ref{thm:B} is Theorem~\ref{thm:independence}, part (ii) is Theorem~\ref{thm:main-n2}(1) together with Theorem~\ref{thm:rank2-classification}, and part (iii) is Theorem~\ref{thm:main-n2}(2).

We conclude this section by comparing the case $n=2$ with the case $n\ge3$.

\begin{remark}\label{rem:n2-vs-n3}
	The answer for $n=2$ has a different shape from that of Theorem~\ref{thm:A}. For $n=2$, the distinguished representations are all the elements of $\T$, independently of $s$, together with $\one_2$ when $s=0$. In contrast, for $n\ge3$, the set of distinguished representations depends on $s$, both through the exceptional values $0$ and $\pm\frac n2$ and through the character $\nu^{-es}$ of $G_{n-2}$.
\end{remark}

\section{General reductions and the case \texorpdfstring{$n=3$}{n=3}}\label{sec:n3}

Throughout this section, unless stated otherwise, we assume that $n\ge3$ and put $H=H_{1,n-1}$. We keep the notation $L_s$, $R_s$, $A_s$, $A_s^\pm$, $\xione_s$ and $\xitwo_s$ of \S\ref{ss:dps}.

\subsection{The listed representations}\label{ss:listed}

We first give a name to the representations that appear in Theorem~\ref{thm:A}(i).

\begin{definition}\label{def:listed}
	Let $n\ge3$ and $s\in\R$. We say that $\pi\in\Irr(G_n)$ is \emph{listed for $(n,s)$} if $\pi$ is isomorphic to one of the following representations:
	\begin{enumerate}[label=\textup{(\alph*)}]
		\item $\xione_s$ if $s\notin\{-\frac n2,0\}$; $\one_n$ if $s=0$; and $\nu^{en/2}\Qt_n$ if $s=-\frac n2$;
		\item $\xitwo_s$ if $s\notin\{0,\frac n2\}$; $\one_n$ if $s=0$; and $\nu^{-en/2}Q_n$ if $s=\frac n2$;
		\item $\nu^{-es}\times\tau$, induced from $G_{n-2}\times G_2$, where $\tau\in\T$ and $\nu^{-es}\times\tau$ is irreducible.
	\end{enumerate}
\end{definition}

By Lemmas~\ref{lem:dps-irred}, \ref{lem:dps2-irred}, \ref{lem:dps1-exceptional} and~\ref{lem:dps2-exceptional} (see Remark~\ref{rem:exceptional-sub}), the representations in (a) are $\xione_s$ when $\xione_s$ is irreducible, and its unique irreducible subrepresentation otherwise. Similarly, the representations in (b) are $\xitwo_s$ when $\xitwo_s$ is irreducible, and its unique irreducible subrepresentation otherwise. Thus, the listed representations are exactly those of Theorem~\ref{thm:A}(i).

We now prove the easy direction of Theorem~\ref{thm:A}(i).

\begin{proposition}\label{prop:listed-dist}
	Let $n\ge3$ and $s\in\R$. Then every representation listed for $(n,s)$ is $(H_{1,n-1},\chi_s)$-distinguished.
\end{proposition}

\begin{proof}
	Suppose first that $s\notin\{-\frac n2,0\}$ (respectively, $s\notin\{0,\frac n2\}$). Then the representation $\xione_s$ (respectively, $\xitwo_s$) is irreducible by Lemma~\ref{lem:dps-irred} (respectively, Lemma~\ref{lem:dps2-irred}), and Example~\ref{ex:dps-functionals} gives a nonzero element of $\Hom_H(\xi^{(i)}_s,\chi_s)$. For $s=0$, the trivial representation $\one_n$ is $(H,\chi_0)$-distinguished, since $\chi_0$ is trivial. Moreover, the representations $\nu^{en/2}\Qt_n$ for $s=-\frac n2$ and $\nu^{-en/2}Q_n$ for $s=\frac n2$ are distinguished by Lemmas~\ref{lem:Qtilde-dist} and~\ref{lem:Q-dist}, respectively.
	
	Finally, let $\pi=\nu^{-es}\times\tau$ be as in (c). Since $\tau\in\T$ is infinite-dimensional, it is $(H_{1,1},\chi_{s+\frac{n-2}2})$-distinguished by Corollary~\ref{cor:T-chi-t}. Put $k=n-2$, so that $\nu^{-es}=\nu^{\frac{e(n-k-2)}2-es}$ and $\chi_{s+\frac{n-2}2}=\chi_{s+\frac k2}$. Hence, condition (I) of Proposition~\ref{prop:three-alternatives} holds for $\pi_1=\nu^{-es}\in\Irr(G_{n-2})$ and $\pi_2=\tau\in\Irr(G_2)$. Therefore, $\pi=\pi_1\times\pi_2$ is distinguished by the converse part of that proposition.
\end{proof}

\subsection{Reduction to induction from \texorpdfstring{$P_{(n-1,1)}$}{P(n-1,1)}}\label{ss:reduction}

In this subsection, we show that every distinguished irreducible representation of $G_n$, $n\ge3$, is a quotient of a representation induced from the maximal parabolic subgroup $P_{(n-1,1)}$. We first restrict the supercuspidal support of such representations.

\begin{lemma}\label{lem:sc-support}
	Let $n\ge3$ and $s\in\R$, and let $\pi\in\Irr(G_n)$ be $(H,\chi_s)$-distinguished. Then the supercuspidal support of $\pi$ consists of representations of $G_1$, except for at most one representation of $G_2$.
\end{lemma}

\begin{proof}
	By Theorem~\ref{thm:ds-not-dist}, $\pi$ is not supercuspidal. Hence, there are $r\ge2$ and $\rho_i\in\mathcal C(G_{n_i})$, with $n_1+\dots+n_r=n$, such that $\pi$ is a quotient of $\Pi=\rho_1\times\dots\times\rho_r$. Indeed, $\widetilde\pi$ embeds into a product of supercuspidal representations, and we take contragredients using \eqref{eq:contragredient}. Supercuspidal representations of groups $G_{n_i}$ and $G_{n_j}$ with $n_i\ne n_j$ lie on different cuspidal lines. Hence, by Lemma~\ref{lem:different-lines}, two adjacent factors $\rho_i$ and $\rho_j$ with $n_i\ne n_j$ commute. Therefore, we may assume that $n_1\ge n_2\ge\dots\ge n_r$. By Lemma~\ref{lem:dist-quotient}, $\Pi$ is distinguished. We apply Theorem~\ref{thm:orbit-bound} with $k=n_1$, $\pi_1=\rho_1$ and $\pi_2=\Sigma=\rho_2\times\dots\times\rho_r$.
	
	If $n_1=1$, there is nothing to prove. Suppose that $n_1\ge2$. Condition (I) would give a nonzero homomorphism from $\rho_1$ to a character of $G_{n_1}$. This is impossible, since a character of $G_{n_1}$ with $n_1\ge2$ is not supercuspidal. Condition (III) fails, since $r_{(n_1-1,1)}(\rho_1)=0$. Hence, condition (II) holds, that is, $\rho_1$ is $(H_{1,n_1-1},\chi_{s-\frac{n-n_1}2})$-distinguished and $\Hom_{G_{n-n_1}}(\Sigma,\nu^{\frac{e(2-n_1)}2-es})\ne0$. Since $\rho_1$ is a discrete series representation, Theorem~\ref{thm:ds-not-dist} gives $n_1\le2$, and therefore $n_1=2$. In this case, the second part of (II) says that $\Sigma$ has the character $\nu^{-es}$ of $G_{n-2}$ as a quotient. The supercuspidal support of a character of $G_{n-2}$ consists of characters of $G_1$ (Lemma~\ref{lem:characters-segments}), and every constituent of $\Sigma$ has supercuspidal support $\{\rho_2,\dots,\rho_r\}$. Hence, $n_i=1$ for all $i\ge2$.
\end{proof}

As a consequence, every distinguished representation is a quotient of a representation induced from $G_{n-1}\times G_1$.

\begin{lemma}\label{lem:quotient-pi-sigma}
	Let $n\ge3$ and $s\in\R$, and let $\pi\in\Irr(G_n)$ be $(H,\chi_s)$-distinguished. Then there exist $\pi'\in\Irr(G_{n-1})$ and $\sigma\in\Irr(G_1)$ such that $\pi$ is a quotient of $\pi'\times\sigma$.
\end{lemma}

\begin{proof}
	Write $\pi$ as a quotient of $\rho_1\times\dots\times\rho_r$ with $\rho_i\in\mathcal C(G_{n_i})$. By Lemma~\ref{lem:sc-support}, we have $n_i\in\{1,2\}$ for all $i$, and at most one $n_i$ equals $2$. Since $n\ge3$, some $n_i$ equals $1$. If the last factor has rank two, we commute it with the neighbouring factor of rank one (Lemma~\ref{lem:different-lines}). Thus, we may assume that $n_r=1$. Then $\pi$ is a quotient of $\Pi'\times\rho_r$, where $\Pi'=\rho_1\times\dots\times\rho_{r-1}$. Taking contragredients, we obtain an embedding $\widetilde\pi\hookrightarrow\widetilde{\Pi'}\times\widetilde\rho_r$ by \eqref{eq:contragredient}. By Frobenius reciprocity, $r_{(n-1,1)}(\widetilde\pi)\ne0$. Since $r_{(n-1,1)}(\widetilde\pi)$ has finite length, it has an irreducible quotient $\tau'\otimes\tau$ with $\tau'\in\Irr(G_{n-1})$ and $\tau\in\Irr(G_1)$. By Frobenius reciprocity again, there is a nonzero homomorphism $\widetilde\pi\to\tau'\times\tau$, which is injective since $\widetilde\pi$ is irreducible. Taking contragredients, we see that $\pi$ is a quotient of $\pi'\times\sigma$, where $\pi'=\widetilde{\tau'}$ and $\sigma=\widetilde\tau$.
\end{proof}

Combining Lemma~\ref{lem:quotient-pi-sigma} with Proposition~\ref{prop:three-alternatives} for $k=n-1$, we obtain the following trichotomy.

\begin{proposition}\label{prop:three-cases}
	Let $n\ge3$ and $s\in\R$, and let $\pi\in\Irr(G_n)$ be $(H,\chi_s)$-distinguished. Let $\pi'\in\Irr(G_{n-1})$ and $\sigma\in\Irr(G_1)$ be such that $\pi$ is a quotient of $\pi'\times\sigma$. Then at least one of the following holds:
	\begin{itemize}
		\item[(I)] $\pi'\simeq\nu^{-es-\frac e2}$ and $\sigma\simeq\nu^{e(s+\frac{n-1}2)}$; in this case, $\pi'\times\sigma=\xione_s$;
		\item[(II)] $\pi'$ is $(H_{1,n-2},\chi_{s-\frac12})$-distinguished, and $\sigma\simeq\nu^{\frac{e(3-n)}2-es}=\nu^{e(L_s+1)}$;
		\item[(III)] there is $\rho\in\Irr(G_1)$ such that $\sigma\simeq\widetilde\rho$ and
		$$
		\pi'\simeq\begin{cases}\nu^{\frac e2-es}=\Ss(A_s^+)&\text{if }\rho\simeq\nu^{eR_s},\\ \nu^{\frac e2-es}\Qt_{n-1}=\Ss(\{[L_s],A_s\})&\text{if }\rho\simeq\nu^{eL_s},\\ \nu^{-es}\times\rho\ \text{(irreducible)}&\text{otherwise}.\end{cases}$$
	\end{itemize}
\end{proposition}

\begin{proof}
	By Lemma~\ref{lem:dist-quotient}, $\pi'\times\sigma$ is distinguished. We apply Proposition~\ref{prop:three-alternatives} with $k=n-1$, $\pi_1=\pi'$ and $\pi_2=\sigma$. Condition (I) reads $\pi'\simeq\nu^{\frac{e(n-k-2)}2-es}=\nu^{-es-\frac e2}$ and $\sigma\simeq\chi_{s+\frac{n-1}2}=\nu^{e(s+\frac{n-1}2)}$ on $H_{1,0}=G_1$. Condition (II) reads as stated, since $\frac{e(2-k)}2-es=\frac{e(3-n)}2-es$ and $L_s+1=\frac{3-n}2-s$. In condition (III), we have $\mathcal O_2(s,1;\rho)=\{\widetilde\rho\}$. Since $k=n-1\ge2$, we may apply \eqref{eq:O1} with $k=n-1$; then $x_s=-es$, $\frac{e(n-1)}2-es=eR_s$ and $\frac{e(n-2k-1)}2-es=eL_s$, and \eqref{eq:O1} gives the three possibilities. Finally, we check that $\nu^{\frac e2-es}\Qt_{n-1}=\Ss(\{[L_s],A_s\})$. Twisting $\Qt_{n-1}=\Ss(\{[-\frac n2],\Delta_{n-2}(-\frac12)\})$ by $\nu^{e(\frac12-s)}$ adds $\frac12-s$ to all endpoints. This turns the singleton $[-\frac n2]$ into $[-\frac n2+\frac12-s]=[L_s]$, and the segment $\Delta_{n-2}(-\frac12)$ into $\Delta_{n-2}(-s)=A_s$.
\end{proof}

\subsection{Reducibility of \texorpdfstring{$\nu^{-es}\times\tau$}{nu times tau}}\label{ss:reducibility}

We now determine exactly when the representations in Definition~\ref{def:listed}(c) are reducible. Recall that $\tau_a=\nu^{ea}\times\nu^{-ea}$ for $a\in\R$, and that $A_s=[L_s+1,R_s-1]$.

\begin{lemma}\label{lem:reducible-nu-tau}
	Let $n\ge3$, $s\in\R$ and $\tau\in\T$. If the representation $\nu^{-es}\times\tau$, induced from $G_{n-2}\times G_2$, is reducible, then one of the following holds:
	\begin{enumerate}
		\item $\tau=\St_2$ and $s=\pm\frac n2$;
		\item $\tau=\tau_x=\nu^{ex}\times\nu^{-ex}$ for some $x\in\{L_s,R_s\}$ with $x\ne\pm\frac12$.
	\end{enumerate}
	Conversely, in each of these cases, $\nu^{-es}\times\tau$ is reducible.
\end{lemma}

\begin{proof}
	The supercuspidal support of the character $\nu^{-es}=\Ss(A_s)$ of $G_{n-2}$ consists of the characters $\nu^{ej}$ of $G_1$ with $j\in A_s$. By Lemma~\ref{lem:different-lines}, if $\nu^{-es}\times\tau$ is reducible, then the supercuspidal support of $\tau$ contains a character $\nu^{ea}$ of $G_1$ with $a\in\R$ (indeed, with $a\in-s+\frac{n-3}2+\Z$). In particular, $\tau$ is not supercuspidal, and its supercuspidal support is contained in $\Irr(G_1)$. Let $\rho\in\Irr(G_1)$ denote the other element of this support. Since $\tau$ is $H_{1,1}$-distinguished, its central character $\omega_\tau$ is trivial by Lemma~\ref{lem:central}.
	
	Suppose first that $\nu^{ea}\times\rho$ is irreducible. Then $\tau=\nu^{ea}\times\rho$, and by Theorem~\ref{thm:rank2-classification}(e), we have $\rho\simeq\nu^{-ea}$. Thus, $\tau=\tau_a$ with $a\ne\pm\frac12$. The representation $\nu^{-es}\times\tau_a$ has the same constituents as $\zeta(\{A_s,[a],[-a]\})$. Hence, by Lemma~\ref{lem:constituents} and Theorem~\ref{thm:zel}, it is reducible if and only if two of the segments $A_s$, $[a]$ and $[-a]$ are linked. Since $a\ne\pm\frac12$, the segments $[a]$ and $[-a]$ are not linked. Therefore, $[a]$ or $[-a]$ is linked to $A_s$, that is, $a$ or $-a$ lies in $\{L_s,R_s\}$. Since $\tau_a\simeq\tau_{-a}$, this is (2).
	
	Suppose now that $\nu^{ea}\times\rho$ is reducible. Then $\rho\simeq\nu^{e(a\pm1)}$ by (T1), and the infinite-dimensional constituent of $\nu^{ea}\times\rho$ is $\nu^{eb}\St_2$ for some $b\in\R$, by Lemma~\ref{lem:St2}. Its central character is $z\mapsto|z|_{\mathrm{F}}^{2e^2b}$, which is trivial only if $b=0$. Hence, $\tau=\St_2$. Finally, $\nu^{-es}\times\St_2=\nu^{-es}(\one_{n-2}\times\nu^{es}\St_2)$, and this representation has the same constituents as $\nu^{-es}(\nu^{es}\St_2\times\one_{n-2})$. Therefore, by Lemma~\ref{lem:St2-times-one}, it is reducible only if $s=\pm\frac n2$. This is (1).
	
	Conversely, suppose that we are in case (2). Then $[x]$ is linked to $A_s$. Let $\bb$ be the multisegment obtained from $\{A_s,[x],[-x]\}$ by merging $[x]$ with $A_s$. By Lemma~\ref{lem:constituents}, $\Ss(\{A_s,[x],[-x]\})$ and $\Ss(\bb)$ are two different constituents of $\nu^{-es}\times\tau_x$, which is therefore reducible. Suppose now that we are in case (1) with $s=\frac n2$. Then $\nu^{-es}\times\St_2$ and $\nu^{-en/2}(\nu^{en/2}\St_2\times\one_{n-2})$ have the same constituents. By Lemma~\ref{lem:Sn-embedding}, the representation $\nu^{en/2}\St_2\times\one_{n-2}$ is a nonzero subrepresentation of $\zeta(\sss_n)$. It therefore contains the unique irreducible subrepresentation $S_n$ of $\zeta(\sss_n)$ (Theorem~\ref{thm:zel}). Moreover, $Q_n$ is a quotient of it by Lemma~\ref{lem:Qm-quotient}. Since $S_n\not\simeq Q_n$ by Theorem~\ref{thm:zel}, this representation is reducible. The case $s=-\frac n2$ follows by taking contragredients, since $\St_2$ is self-dual.
\end{proof}

In particular, the irreducibility condition in Theorem~\ref{thm:A}(i)(c) can be made completely explicit.

\begin{remark}\label{rem:irreducibility-explicit}
	By Lemma~\ref{lem:reducible-nu-tau}, for $\tau\in\T$, the representation $\nu^{-es}\times\tau$ is irreducible except exactly in the following cases: $\tau=\St_2$ and $s=\pm\frac n2$; or $\tau=\tau_x$ for some $x\in\{-s-\frac{n-1}2,\,-s+\frac{n-1}2\}$ with $x\ne\pm\frac12$. In particular, $\nu^{-es}\times\tau$ is irreducible for every $s$ in each of the following cases:
	\begin{itemize}
		\item $\tau$ is supercuspidal;
		\item $\tau$ is $\St_2(\sigma)$, $\Sp_2(\sigma)$ or $\sigma\times\widetilde\sigma$ with $\dim\sigma>1$;
		\item $\tau=\chi\St_2$ with $\chi\ne\one_1$.
	\end{itemize}
\end{remark}

\subsection{The case \texorpdfstring{$n=3$}{n=3}: statement}\label{ss:n3-statement}

We now treat the case $n=3$, which serves as the base of the induction in the general case. Our aim is to prove the following theorem.

\begin{theorem}\label{thm:n3}
	Let $s\in\R$. An irreducible smooth representation of $G_3$ is $(H_{1,2},\chi_s)$-distinguished if and only if it is listed for $(3,s)$.
\end{theorem}

We first make the notation of \S\ref{ss:dps} explicit for $n=3$. In this case, we have $A_s=[-s]$, $L_s=-s-1$, $R_s=1-s$, $A_s^-=[-s-1,-s]$ and $A_s^+=[-s,1-s]$, and therefore
$$\xione_s=\Ss([-s-1,-s])\times\Ss([s+1]),\qquad\xitwo_s=\Ss([s-1])\times\Ss([-s,1-s]).$$
Twisting Definition~\ref{def:Qm} and Lemma~\ref{lem:Qm-Langlands} by $\nu^{ec}$, with $c\in\R$, we obtain
\begin{equation}\label{eq:Q3}
	\begin{aligned}
		\nu^{ec}Q_3&=\Ss(\{[c,c+1],[c+2]\})=\Qq(\{[c],[c+1,c+2]\}),\\
		\nu^{ec}\Qt_3&=\Ss(\{[c-2],[c-1,c]\})=\Qq(\{[c-2,c-1],[c]\}).
	\end{aligned}
\end{equation}

For $n=3$, Theorem~\ref{thm:orbit-bound} and Corollaries~\ref{cor:dIII-k1} and~\ref{cor:dIII-kn1} take the following explicit form. For $\rho,\rho'\in\Irr(G_1)$, we put $\mathrm m(\rho,\rho')=1$ if $\rho\simeq\rho'$, and $\mathrm m(\rho,\rho')=0$ otherwise. We begin with induction from $P_{(1,2)}$.

\begin{lemma}\label{lem:n3-bound-k1}
	Let $\pi_1\in\Irr(G_1)$, $\pi_2\in\Irr(G_2)$ and $s\in\R$. Then $\dim\Hom_{H_{1,2}}(\pi_1\times\pi_2,\chi_s)\le\dI+\dII+\dIII$, where
	\begin{align*}
		\dI&=\mathrm m(\pi_1,\nu^{-es})\,\dim\Hom_{H_{1,1}}(\pi_2,\chi_{s+\frac12}),\\
		\dII&=\mathrm m(\pi_1,\nu^{e(s-1)})\,\dim\Hom_{G_2}(\pi_2,\nu^{\frac e2-es}),\\
		\dIII&=\dim\Hom_{G_2}(\pi_2,\widetilde\pi_1\times\nu^{-es}).
	\end{align*}
	Moreover, if $\dI\ne0$ or $\dII\ne0$, then $\pi_1\times\pi_2$ is $(H_{1,2},\chi_s)$-distinguished.
\end{lemma}

Similarly, for induction from $P_{(2,1)}$, we have the following.

\begin{lemma}\label{lem:n3-bound-k2}
	Let $\pi_1\in\Irr(G_2)$, $\pi_2\in\Irr(G_1)$ and $s\in\R$. Then $\dim\Hom_{H_{1,2}}(\pi_1\times\pi_2,\chi_s)\le\dI+\dII+\dIII$, where
	\begin{align*}
		\dI&=\dim\Hom_{G_2}(\pi_1,\nu^{-\frac e2-es})\,\mathrm m(\pi_2,\nu^{e(s+1)}),\\
		\dII&=\dim\Hom_{H_{1,1}}(\pi_1,\chi_{s-\frac12})\,\mathrm m(\pi_2,\nu^{-es}),\\
		\dIII&=\dim\Hom_{G_2}(\pi_1,\nu^{-es}\times\widetilde\pi_2).
	\end{align*}
	Moreover, if $\dI\ne0$ or $\dII\ne0$, then $\pi_1\times\pi_2$ is $(H_{1,2},\chi_s)$-distinguished.
\end{lemma}

\begin{proof}[Proof of Lemmas~\ref{lem:n3-bound-k1} and~\ref{lem:n3-bound-k2}]
	These are Theorem~\ref{thm:orbit-bound} together with Corollary~\ref{cor:dIII-k1} (respectively, Corollary~\ref{cor:dIII-kn1}), for $n=3$ and $k=1$ (respectively, $k=2$). Here, we use the convention $\chi_c=\nu^{ec}$ on $H_{1,0}=G_1$, and the fact that, for $\rho\in\Irr(G_1)$ and $c\in\R$, we have $\dim\Hom_{G_1}(\rho,\nu^{ec})=\mathrm m(\rho,\nu^{ec})$.
\end{proof}

The terms $\dIII$ are computed as follows. For $\sigma_1,\sigma_2\in\Irr(G_1)$, the representation $\sigma_1\times\sigma_2$ has an irreducible socle: it is either irreducible or, by (T2), of length two with a unique irreducible subrepresentation. Hence, by Lemma~\ref{lem:socle-hom}, each term $\dIII$ in Lemmas~\ref{lem:n3-bound-k1} and~\ref{lem:n3-bound-k2} is at most $1$, with equality if and only if $\pi_2$ (respectively, $\pi_1$) is isomorphic to this socle. In the cases we need, the relevant representation is of the form $\nu^{ea}\times\nu^{-ea}$. By Lemma~\ref{lem:St2}, its unique irreducible subrepresentation is $\tau_a$ if $a\ne\pm\frac12$, the character $\one_2$ if $a=-\frac12$, and $\St_2$ if $a=\frac12$.

\subsection{Exceptional representations of \texorpdfstring{$G_3$}{G3}}\label{ss:n3-exceptional}

In this subsection, we show that certain twists of $Q_3$ and $\Qt_3$ are not distinguished. They occur as quotients of degenerate principal series at the exceptional values of $s$, and must therefore be excluded separately. In each case, we realize the representation as a quotient of a suitable induced representation, and we show, using Lemmas~\ref{lem:n3-bound-k1} and~\ref{lem:n3-bound-k2}, that the induced representation is not distinguished. We begin with the case $s=0$.

\begin{lemma}\label{lem:n3-Q3-s0}
	The representation $\nu^{-e}Q_3$ is not $(H_{1,2},\chi_0)$-distinguished.
\end{lemma}

\begin{proof}
	By \eqref{eq:Q3}, we have $\nu^{-e}Q_3=\Qq(\{[-1],[0,1]\})$. Since $[-1]\prec[0,1]$, the ordering $([0,1],[-1])$ is admissible. Hence, by Theorem~\ref{thm:lang} and Lemma~\ref{lem:St2}, $\nu^{-e}Q_3$ is a quotient of $\nu^{e/2}\St_2\times\nu^{-e}$. We apply Lemma~\ref{lem:n3-bound-k2} with $s=0$, $\pi_1=\nu^{e/2}\St_2$ and $\pi_2=\nu^{-e}$. We have $\dI=0$, because $\nu^{e/2}\St_2$ is not a character, and $\dII=0$, because $\nu^{-e}\not\simeq\one_1$. Moreover,
	$$
	\dIII=\dim\Hom_{G_2}\big(\nu^{e/2}\St_2,\one_1\times\nu^{e}\big)=0,$$
	because the unique irreducible subrepresentation of $\one_1\times\nu^e$ is the character $\nu^{e/2}$ (Lemmas~\ref{lem:St2} and~\ref{lem:socle-hom}). Hence, $\nu^{e/2}\St_2\times\nu^{-e}$ is not $(H_{1,2},\chi_0)$-distinguished. Therefore, by Lemma~\ref{lem:dist-quotient}, neither is its quotient $\nu^{-e}Q_3$.
\end{proof}

The same argument applies to the corresponding twist of $\Qt_3$.

\begin{lemma}\label{lem:n3-Qt3-s0}
	The representation $\nu^{e}\Qt_3$ is not $(H_{1,2},\chi_0)$-distinguished.
\end{lemma}

\begin{proof}
	By \eqref{eq:Q3}, we have $\nu^e\Qt_3=\Qq(\{[-1,0],[1]\})$. Since $[-1,0]\prec[1]$, the ordering $([1],[-1,0])$ is admissible. Hence, by Theorem~\ref{thm:lang} and Lemma~\ref{lem:St2}, $\nu^e\Qt_3$ is a quotient of $\nu^e\times\nu^{-e/2}\St_2$. We apply Lemma~\ref{lem:n3-bound-k1} with $s=0$, $\pi_1=\nu^e$ and $\pi_2=\nu^{-e/2}\St_2$. We have $\dI=\dII=0$, since $\nu^e\notin\{\one_1,\nu^{-e}\}$. Moreover,
	$$
	\dIII=\dim\Hom_{G_2}\big(\nu^{-e/2}\St_2,\nu^{-e}\times\one_1\big)=0,$$
	since the unique irreducible subrepresentation of $\nu^{-e}\times\one_1$ is the character $\nu^{-e/2}$. Hence, $\nu^e\times\nu^{-e/2}\St_2$ is not $(H_{1,2},\chi_0)$-distinguished, and therefore, by Lemma~\ref{lem:dist-quotient}, neither is $\nu^e\Qt_3$.
\end{proof}

Finally, we treat a twist of $Q_3$ at the exceptional value $s=-\frac32$.

\begin{lemma}\label{lem:n3-Q3-s-32}
	The representation $\nu^{-e/2}Q_3$ is not $(H_{1,2},\chi_{-3/2})$-distinguished.
\end{lemma}

\begin{proof}
	By \eqref{eq:Q3}, we have $\nu^{-e/2}Q_3=\Ss(\{[-\frac12,\frac12],[\frac32]\})$. Since $[-\frac12,\frac12]\prec[\frac32]$, the ordering $([\frac32],[-\frac12,\frac12])$ is admissible. Hence, by Lemma~\ref{lem:S-quotient}, $\nu^{-e/2}Q_3$ is a quotient of $\one_2\times\nu^{3e/2}$. We apply Lemma~\ref{lem:n3-bound-k2} with $s=-\frac32$, $\pi_1=\one_2$ and $\pi_2=\nu^{3e/2}$. We have $\dI=0$, since $\one_2\not\simeq\nu^{-\frac e2-es}=\nu^{e}$. Next, $\dII=0$ by Lemma~\ref{lem:char-dist}, since $\one_2$ is not $(H_{1,1},\chi_{s-\frac12})$-distinguished when $s-\frac12=-2\ne0$. Finally,
	$$
	\dIII=\dim\Hom_{G_2}\big(\one_2,\nu^{3e/2}\times\nu^{-3e/2}\big)=0,$$
	since $\nu^{3e/2}\times\nu^{-3e/2}=\tau_{3/2}$ is irreducible and not a character. Hence, $\one_2\times\nu^{3e/2}$ is not $(H_{1,2},\chi_{-3/2})$-distinguished, and therefore, by Lemma~\ref{lem:dist-quotient}, neither is $\nu^{-e/2}Q_3$.
\end{proof}

\subsection{The key lemma for \texorpdfstring{$n=3$}{n=3}}\label{ss:n3-key}

The following lemma is the main step in the proof of Theorem~\ref{thm:n3}. It treats the distinguished representations of $G_3$ that are quotients of a representation induced from $\nu^{-es}$ and an irreducible representation of $G_2$.

\begin{lemma}\label{lem:n3-key}
	Let $s\in\R$ and $\tau'\in\Irr(G_2)$, and let $\pi\in\Irr(G_3)$ be $(H_{1,2},\chi_s)$-distinguished. Suppose that $\pi$ is a quotient of $\nu^{-es}\times\tau'$ or of $\tau'\times\nu^{-es}$. Then $\pi$ is listed for $(3,s)$.
\end{lemma}

\begin{proof}
	Case $1$. \emph{Either $\tau'\in\T$, or $\tau'=\one_2$ and $s=\pm\frac12$.} By Lemma~\ref{lem:dist-quotient}, the induced representation of which $\pi$ is a quotient is distinguished.
	
	Suppose first that $\pi$ is a quotient of $\nu^{-es}\times\tau'$, and apply Lemma~\ref{lem:n3-bound-k1} with $\pi_1=\nu^{-es}$ and $\pi_2=\tau'$. If $\dI\ne0$, then $\tau'$ is $(H_{1,1},\chi_{s+\frac12})$-distinguished. Hence, by Corollary~\ref{cor:T-chi-t} and Lemma~\ref{lem:char-dist} (for $n=2$), either $\tau'\in\T$, or $\tau'=\one_2$ and $s=-\frac12$. If $\dII\ne0$, then $\nu^{-es}=\nu^{e(s-1)}$ and $\tau'=\nu^{\frac e2-es}$, that is, $s=\frac12$ and $\tau'=\one_2$. If $\dIII\ne0$, then, by Lemma~\ref{lem:socle-hom}, $\tau'$ is the unique irreducible subrepresentation of $\nu^{es}\times\nu^{-es}$. By Lemma~\ref{lem:St2} and Corollary~\ref{cor:T-members}, this is $\tau_s\in\T$ if $s\ne\pm\frac12$, $\St_2\in\T$ if $s=\frac12$, and $\one_2$ if $s=-\frac12$.
	
	Suppose now that $\pi$ is a quotient of $\tau'\times\nu^{-es}$, and apply Lemma~\ref{lem:n3-bound-k2} with $\pi_1=\tau'$ and $\pi_2=\nu^{-es}$. If $\dI\ne0$, then $\tau'=\nu^{-\frac e2-es}$ and $\nu^{-es}=\nu^{e(s+1)}$, so that $s=-\frac12$ and $\tau'=\one_2$. If $\dII\ne0$, then $\tau'$ is $(H_{1,1},\chi_{s-\frac12})$-distinguished, so either $\tau'\in\T$, or $\tau'=\one_2$ and $s=\frac12$. If $\dIII\ne0$, then $\tau'$ is the unique irreducible subrepresentation of $\nu^{-es}\times\nu^{es}$, which is $\tau_s\in\T$ if $s\ne\pm\frac12$, $\one_2$ if $s=\frac12$, and $\St_2\in\T$ if $s=-\frac12$.
	
	Case $2$. \emph{The case $\tau'=\one_2$ and $s=\pm\frac12$.} In this case, the segments $[-s]$ and $[-\frac12,\frac12]$ are not linked, since $[-s]\subseteq[-\frac12,\frac12]$. Hence, by Lemma~\ref{lem:unlinked-product}, both $\nu^{-es}\times\one_2$ and $\one_2\times\nu^{-es}$ are irreducible, and isomorphic to $\pi$. For $s=-\frac12$, this representation is $\one_2\times\nu^{e/2}=\xione_{-1/2}$, and for $s=\frac12$, it is $\nu^{-e/2}\times\one_2=\xitwo_{1/2}$. Both are listed, by Lemmas~\ref{lem:dps-irred} and~\ref{lem:dps2-irred}.
	
	Case $3$.  \emph{The case $\tau'\in\T$.} If $\nu^{-es}\times\tau'$ is irreducible, then $\pi\simeq\nu^{-es}\times\tau'\simeq\tau'\times\nu^{-es}$, which is listed in (c). Otherwise, Lemma~\ref{lem:reducible-nu-tau} applies. Moreover, since $\nu^{-es}\times\tau'$ and $\tau'\times\nu^{-es}$ have the same constituents, $\pi$ is a constituent of $\nu^{-es}\times\tau'$. We go through the cases of Lemma~\ref{lem:reducible-nu-tau}.
	
	(3a) \emph{$\tau'=\St_2$ and $s=\pm\frac32$.} By Lemma~\ref{lem:two-segments} and \eqref{eq:Q3}, the constituents of $\Qq([-\frac32])\times\Qq([-\frac12,\frac12])$ are $\Qq(\{[-\frac32],[-\frac12,\frac12]\})=\nu^{-3e/2}Q_3$ and the discrete series representation $\Qq([-\frac32,\frac12])$. Similarly, the constituents of $\Qq([\frac32])\times\Qq([-\frac12,\frac12])$ are $\Qq(\{[-\frac12,\frac12],[\frac32]\})=\nu^{3e/2}\Qt_3$ and $\Qq([-\frac12,\frac32])$. The discrete series representations are not distinguished, by Theorem~\ref{thm:ds-not-dist}. On the other hand, $\nu^{-3e/2}Q_3$ (for $s=\frac32$) and $\nu^{3e/2}\Qt_3$ (for $s=-\frac32$) are listed.
	
	In the remaining cases, $\tau'=\tau_x$ for some $x\in\{L_s,R_s\}=\{-s-1,1-s\}$ with $x\ne\pm\frac12$. By Lemma~\ref{lem:constituents}, the constituents of $\nu^{-es}\times\tau_x=\Ss([-s])\times\Ss([x])\times\Ss([-x])$ are among the representations $\Ss(\bb)$ with $\bb\le\mm=\{[-s],[x],[-x]\}$.
	
	(3b) \emph{$x=1-s$ and $s\ne0$.} Here, the condition $x\ne\pm\frac12$ means that $s\notin\{\frac12,\frac32\}$. The segments $[-s]$ and $[1-s]$ are linked. The segment $[s-1]$ is linked neither to $[1-s]$ nor to $[-s,1-s]$, and it is linked to $[-s]$ only if $s=1$. In the latter case, $[-1]$ and $[0]$ merge to $[-1,0]=[-s,1-s]$ as well. Hence, the constituents are $\Ss(\mm)$ and $\Ss(\{[-s,1-s],[s-1]\})\simeq\xitwo_s$, which is irreducible and listed by Lemma~\ref{lem:dps2-irred}.
	
	Suppose that $\Ss(\mm)$ is distinguished. The ordering $([s-1],[1-s],[-s])$ of $\mm$ is admissible. Hence, by Lemmas~\ref{lem:S-quotient} and~\ref{lem:quotient-constituents}, $\Ss(\mm)$ is a quotient of $X\times\nu^{e(s-1)}$ for some constituent $X$ of $\nu^{-es}\times\nu^{e(1-s)}$. By Lemma~\ref{lem:St2}, we have $X=\nu^{e(\frac12-s)}$ or $X=\nu^{e(\frac12-s)}\St_2$. In the first case, $X\times\nu^{e(s-1)}=\xi^{(2)\prime}_s$ is irreducible by Lemma~\ref{lem:dps2-irred}. Hence, $\Ss(\mm)\simeq\Ss(\{[-s,1-s],[s-1]\})$, which contradicts Theorem~\ref{thm:zel}. In the second case, we apply Lemma~\ref{lem:n3-bound-k2} with $\pi_1=X$ and $\pi_2=\nu^{e(s-1)}$. We have $\dI=0$, since $X$ is not a character, and $\dII=0$, since $\nu^{e(s-1)}\ne\nu^{-es}$ as $s\ne\frac12$. Moreover, $\dIII=\dim\Hom_{G_2}(X,\nu^{-es}\times\nu^{e(1-s)})=0$, since the unique irreducible subrepresentation of $\nu^{-es}\times\nu^{e(1-s)}$ is the character $\nu^{e(\frac12-s)}$. Hence, $X\times\nu^{e(s-1)}$ is not distinguished, which contradicts Lemma~\ref{lem:dist-quotient}. Therefore, $\Ss(\mm)$ is not distinguished.
	
	(3c) \emph{$x=-s-1$ and $s\ne0$.} Here, the condition $x\ne\pm\frac12$ means that $s\notin\{-\frac12,-\frac32\}$. The segments $[-s-1]$ and $[-s]$ are linked. The segment $[s+1]$ is linked neither to $[-s-1]$ nor to $[-s-1,-s]$, and it is linked to $[-s]$ only if $s=-1$. In the latter case, $[0]$ and $[1]$ merge to $[0,1]=[-s-1,-s]$ as well. Hence, the constituents are $\Ss(\mm)$ and $\Ss(\{[-s-1,-s],[s+1]\})\simeq\xione_s$, which is irreducible and listed by Lemma~\ref{lem:dps-irred}.
	
	Suppose that $\Ss(\mm)$ is distinguished. The ordering $([-s],[-s-1],[s+1])$ of $\mm$ is admissible; for $s=-1$, note that $[0]\prec[1]=[-s]$. Hence, as in (3b), $\Ss(\mm)$ is a quotient of $\nu^{e(s+1)}\times Y$ for some constituent $Y$ of $\nu^{-e(s+1)}\times\nu^{-es}$, that is, $Y\in\{\nu^{-e(s+\frac12)},\nu^{-e(s+\frac12)}\St_2\}$. If $Y$ is the character, then $\nu^{e(s+1)}\times Y=\xi^{(1)\prime}_s$ is irreducible by Lemma~\ref{lem:dps-irred}, and we obtain a contradiction as in (3b). If $Y=\nu^{-e(s+\frac12)}\St_2$, we apply Lemma~\ref{lem:n3-bound-k1} with $\pi_1=\nu^{e(s+1)}$ and $\pi_2=Y$. We have $\dI=0$, since $\nu^{e(s+1)}\ne\nu^{-es}$ as $s\ne-\frac12$, and $\dII=0$, since $Y$ is not a character. Moreover, $\dIII=\dim\Hom_{G_2}(Y,\nu^{-e(s+1)}\times\nu^{-es})=0$, since the unique irreducible subrepresentation of $\nu^{-e(s+1)}\times\nu^{-es}$ is the character $\nu^{-e(s+\frac12)}$. Hence, $\nu^{e(s+1)}\times Y$ is not distinguished, which again contradicts Lemma~\ref{lem:dist-quotient}. Therefore, $\Ss(\mm)$ is not distinguished.
	
	(3d) \emph{$s=0$ and $\tau'=\tau_1$.} In this case, $\mm=\{[-1],[0],[1]\}$, and the multisegments $\bb\le\mm$ are $\mm$, $\{[-1,0],[1]\}$, $\{[-1],[0,1]\}$ and $\{[-1,1]\}$. By \eqref{eq:Q3}, the second and the third give $\nu^{-e}Q_3$ and $\nu^e\Qt_3$, respectively, which are not distinguished by Lemmas~\ref{lem:n3-Q3-s0} and~\ref{lem:n3-Qt3-s0}. The last one gives $\one_3$, which is listed. Finally, by Lemma~\ref{lem:S-quotient}, $\Ss(\mm)$ is the unique irreducible quotient of $\nu^{-e}\times\one_1\times\nu^e$. By \S\ref{ss:segments}, the unique irreducible quotient of this representation is $\Qq([-1,1])=\St_3$. Thus, $\Ss(\mm)=\St_3$, which is not distinguished by Theorem~\ref{thm:ds-not-dist}.
	
	It remains to check that the cases (3a)--(3d) are exhaustive. By Lemma~\ref{lem:reducible-nu-tau}, either we are in case (3a), or $\tau'=\tau_x$ with $x\in\{-s-1,1-s\}$. Since $\tau_x\simeq\tau_{-x}$, the representations $\tau_{-s-1}$ and $\tau_{1-s}$ coincide only if $-s-1=s-1$, that is, $s=0$. In that case, $\tau'=\tau_1$, which is case (3d). For $s\ne0$, we are in case (3b) or (3c).
\end{proof}

\subsection{Proof of Theorem~\ref{thm:n3}}\label{ss:n3-proof}

We now combine the results of this section to prove Theorem~\ref{thm:n3}. By Proposition~\ref{prop:listed-dist}, every listed representation is distinguished. Conversely, let $\pi\in\Irr(G_3)$ be $(H_{1,2},\chi_s)$-distinguished. By Lemma~\ref{lem:quotient-pi-sigma}, $\pi$ is a quotient of $\pi'\times\sigma$ for some $\pi'\in\Irr(G_2)$ and $\sigma\in\Irr(G_1)$. Hence, by Proposition~\ref{prop:three-cases}, one of the cases (I), (II) and (III) of that proposition holds. For $n=3$, case (III) reads as follows: there is $\rho\in\Irr(G_1)$ with $\sigma\simeq\widetilde\rho$ and
\begin{itemize}
	\item $\pi'=\nu^{e(\frac12-s)}$ if $\rho\simeq\nu^{e(1-s)}$;
	\item $\pi'=\nu^{e(\frac12-s)}\Qt_2=\nu^{-e(s+\frac12)}\St_2$ if $\rho\simeq\nu^{-e(s+1)}$, by \eqref{eq:Q2};
	\item $\pi'=\nu^{-es}\times\rho$, which is irreducible, otherwise.
\end{itemize}
We treat the cases in turn.

\emph{Case \textup{(I)}.} Here $\pi'\times\sigma=\xione_s$. If $s\notin\{-\frac32,0\}$, then $\xione_s$ is irreducible by Lemma~\ref{lem:dps-irred}, and therefore $\pi\simeq\xione_s$, which is listed. If $s=0$, then $\pi$ is the unique irreducible quotient $\nu^{-e}Q_3$ of $\xione_0$ by Lemma~\ref{lem:dps1-exceptional}. This is impossible by Lemma~\ref{lem:n3-Q3-s0}. If $s=-\frac32$, then $\pi$ is the unique irreducible quotient of $\xione_{-3/2}$, namely the character $\nu^{e/2}$ (Lemma~\ref{lem:dps1-exceptional}). This is impossible by Lemma~\ref{lem:char-dist}, since $s\ne0$.

\emph{Case \textup{(II)}.} Here $\sigma\simeq\nu^{e(L_s+1)}=\nu^{-es}$, so $\pi$ is a quotient of $\pi'\times\nu^{-es}$. Hence, Lemma~\ref{lem:n3-key} applies with $\tau'=\pi'$, and $\pi$ is listed.

\emph{Case \textup{(III)} with $\rho\simeq\nu^{e(1-s)}$.} In this case, $\pi'\times\sigma=\Ss([-s,1-s])\times\Ss([s-1])=\xi^{(2)\prime}_s$. By Lemmas~\ref{lem:dps2-irred} and~\ref{lem:dps2-exceptional}, we have $\pi\simeq\xitwo_s$ if $s\notin\{0,\frac32\}$, $\pi\simeq\one_3$ if $s=0$, and $\pi\simeq\nu^{-3e/2}Q_3$ if $s=\frac32$. All of these representations are listed.

\emph{Case \textup{(III)} with $\rho\simeq\nu^{-e(s+1)}$.} In this case, $\pi'=\Qq([-s-1,-s])$ and $\sigma=\Qq([s+1])$, so $\pi$ is a quotient of $\Qq([-s-1,-s])\times\Qq([s+1])$. The two segments are linked if and only if $s\in\{0,-\frac32\}$.
\begin{itemize}
	\item If $s=0$, then $[-1,0]\prec[1]$. Hence, by Lemma~\ref{lem:two-segments}, the unique irreducible quotient of $\Qq([-1,0])\times\Qq([1])$ is $\Qq([-1,1])=\St_3$, which is not distinguished by Theorem~\ref{thm:ds-not-dist}.
	\item If $s=-\frac32$, then $[-\frac12]\prec[\frac12,\frac32]$. Hence, by Lemma~\ref{lem:two-segments} and \eqref{eq:Q3}, the unique irreducible quotient of $\Qq([\frac12,\frac32])\times\Qq([-\frac12])$ is $\Qq(\{[-\frac12],[\frac12,\frac32]\})=\nu^{-e/2}Q_3$, which is not distinguished by Lemma~\ref{lem:n3-Q3-s-32}.
	\item If $s\notin\{0,-\frac32\}$, then the two segments are not linked. Hence, by Lemma~\ref{lem:two-segments}, $\pi\simeq\nu^{-e(s+\frac12)}\St_2\times\nu^{e(s+1)}\simeq\nu^{e(s+1)}\times\nu^{-e(s+\frac12)}\St_2$. We apply Lemma~\ref{lem:n3-bound-k1} to the latter, with $\pi_1=\nu^{e(s+1)}$ and $\pi_2=\nu^{-e(s+\frac12)}\St_2$. We have $\dII=0$, since $\pi_2$ is not a character. Moreover,
	$$\dIII=\dim\Hom_{G_2}\big(\nu^{-e(s+\frac12)}\St_2,\nu^{-e(s+1)}\times\nu^{-es}\big)=0,
	$$
	since the unique irreducible subrepresentation of $\nu^{-e(s+1)}\times\nu^{-es}$ is the character $\nu^{-e(s+\frac12)}$. Since $\pi$ is distinguished, it follows that $\dI\ne0$. This forces $\nu^{e(s+1)}=\nu^{-es}$, that is, $s=-\frac12$. In this case, $\pi\simeq\nu^{e/2}\times\St_2$. This representation is listed in (c), since $\St_2\in\T$ by Corollary~\ref{cor:T-members}, and $\nu^{e/2}\times\St_2$ is irreducible.
\end{itemize}

\emph{Case \textup{(III)} with $\rho\notin\{\nu^{e(1-s)},\nu^{-e(s+1)}\}$.} In this case, $\pi'=\nu^{-es}\times\rho$, so $\pi$ is a quotient of $\nu^{-es}\times(\rho\times\widetilde\rho)$, by transitivity of induction. Hence, by Lemma~\ref{lem:quotient-constituents}, $\pi$ is a quotient of $\nu^{-es}\times\tau'$ for some constituent $\tau'\in\Irr(G_2)$ of $\rho\times\widetilde\rho$. Therefore, Lemma~\ref{lem:n3-key} applies, and $\pi$ is listed.

This completes the proof of Theorem~\ref{thm:n3}. \qed

\section{Auxiliary results for the inductive step}\label{sec:aux}

Throughout this section, we assume that $n\ge3$, we put $H=H_{1,n-1}$, and we keep the notation of \S\ref{ss:dps}. We collect here the consequences of Proposition~\ref{prop:three-alternatives}, together with the non-distinction results for exceptional representations, that will be used in the inductive step of Section~\ref{sec:induction}.

\subsection{Quotients of representations induced from a character of \texorpdfstring{$G_1$}{G1}}\label{ss:quotient-char-G1}

We first determine the possible shapes of the inducing data when a distinguished representation is a quotient of a representation induced from an irreducible representation of $G_{n-1}$ and a character of $G_1$. We begin with the case where the character of $G_1$ is the second factor.

\begin{lemma}\label{lem:quot-X-nu}
	Let $n\ge3$ and $s,y\in\R$, and let $\pi\in\Irr(G_n)$ be $(H,\chi_s)$-distinguished. Suppose that $\pi$ is a quotient of $X\times\nu^{ey}$ for some $X\in\Irr(G_{n-1})$. Then one of the following holds:
	\begin{enumerate}[label=\textup{(\alph*)}]
		\item $y=s+\frac{n-1}2$ and $X\simeq\nu^{-e(s+\frac12)}$;
		\item $y=-s-\frac{n-3}2$, and $X$ is $(H_{1,n-2},\chi_{s-\frac12})$-distinguished;
		\item $y=s-\frac{n-1}2$ and $X\simeq\nu^{-e(s-\frac12)}$;
		\item $y=s+\frac{n-1}2$ and $X\simeq\nu^{-e(s-\frac12)}\Qt_{n-1}$;
		\item $y\notin\{s\pm\frac{n-1}2\}$, and $X\simeq\nu^{-es}\times\nu^{-ey}$, where the latter representation is induced from $G_{n-2}\times G_1$ and is irreducible.
	\end{enumerate}
\end{lemma}

\begin{proof}
	We apply Proposition~\ref{prop:three-cases} with $\pi'=X$ and $\sigma=\nu^{ey}$. Case (I) of that proposition is exactly (a). Case (II) gives (b), since $\nu^{\frac{e(3-n)}2-es}=\nu^{ey}$ means that $y=-s-\frac{n-3}2$. In case (III), we have $\rho\simeq\widetilde\sigma=\nu^{-ey}$. Moreover, the conditions $\rho\simeq\nu^{eR_s}$ and $\rho\simeq\nu^{eL_s}$ mean that $y=s-\frac{n-1}2$ and $y=s+\frac{n-1}2$, respectively. Since $\nu^{\frac e2-es}=\nu^{-e(s-\frac12)}$, the three possibilities in case (III) give (c), (d) and (e), respectively.
\end{proof}

Next, we treat the case where the character of $G_1$ is the first factor.

\begin{lemma}\label{lem:quot-nu-X}
	Let $n\ge3$ and $s,y\in\R$, and let $\pi\in\Irr(G_n)$ be $(H,\chi_s)$-distinguished. Suppose that $\pi$ is a quotient of $\nu^{ey}\times X$ for some $X\in\Irr(G_{n-1})$. Then one of the following holds:
	\begin{enumerate}[label=\textup{(\alph*)}]
		\item $y=-s+\frac{n-3}2$, and $X$ is $(H_{1,n-2},\chi_{s+\frac12})$-distinguished;
		\item $y=s-\frac{n-1}2$ and $X\simeq\nu^{-e(s-\frac12)}$;
		\item $y=s+\frac{n-1}2$ and $X\simeq\nu^{-e(s+\frac12)}$;
		\item $y=s-\frac{n-1}2$ and $X\simeq\nu^{-e(s+\frac12)}Q_{n-1}$;
		\item $y\notin\{s\pm\frac{n-1}2\}$, and $X\simeq\nu^{-ey}\times\nu^{-es}$, where the latter representation is induced from $G_1\times G_{n-2}$ and is irreducible.
	\end{enumerate}
\end{lemma}

\begin{proof}
	By Lemma~\ref{lem:dist-quotient}, the representation $\nu^{ey}\times X$ is distinguished. We apply Proposition~\ref{prop:three-alternatives} with $k=1$, $\pi_1=\nu^{ey}$ and $\pi_2=X$.
	\begin{itemize}
		\item Condition (I) says that $\nu^{ey}\simeq\nu^{\frac{e(n-3)}2-es}$ and that $X$ is $(H_{1,n-2},\chi_{s+\frac12})$-distinguished. This is (a).
		\item Condition (II) says that $\nu^{ey}\simeq\chi_{s-\frac{n-1}2}=\nu^{e(s-\frac{n-1}2)}$ on $H_{1,0}=G_1$, and that $X\simeq\nu^{\frac e2-es}=\nu^{-e(s-\frac12)}$. This is (b).
		\item Condition (III) says that there is $\sigma\in\Irr(G_1)$ such that $\nu^{ey}\in\mathcal O_1(s,1;\sigma)=\{\sigma\}$ and $X\in\mathcal O_2(s,n-1;\sigma)$; thus, $\sigma\simeq\nu^{ey}$. We apply \eqref{eq:O2} with $k=1$, so that $y_s=-es$ and $\nu^{-\frac{ek}2-es}=\nu^{-e(s+\frac12)}$. Since $es+\frac{e(n-1)}2=e(s+\frac{n-1}2)$ and $es-\frac{e(n-2k+1)}2=e(s-\frac{n-1}2)$, we obtain (c), (d) and (e), according as $y=s+\frac{n-1}2$, $y=s-\frac{n-1}2$, or neither.\qedhere
	\end{itemize}
\end{proof}

\subsection{Quotients of representations induced from a character of \texorpdfstring{$G_{n-2}$}{Gn-2}}\label{ss:quotient-char-Gn2}

We now consider distinguished representations that are quotients of a representation induced from the character $\nu^{-es}$ of $G_{n-2}$ and an irreducible representation of $G_2$. In this situation, only one of the three orbits can contribute. We begin with the case where $\nu^{-es}$ is the first factor.

\begin{lemma}\label{lem:quot-nu-tau}
	Let $n\ge4$, $s\in\R$ and $\tau\in\Irr(G_2)$, and let $\pi\in\Irr(G_n)$ be $(H,\chi_s)$-distinguished. If $\pi$ is a quotient of the representation $\nu^{-es}\times\tau$ induced from $G_{n-2}\times G_2$, then $\tau$ is $(H_{1,1},\chi_{s+\frac{n-2}2})$-distinguished.
\end{lemma}

\begin{proof}
	By Lemma~\ref{lem:dist-quotient}, the representation $\nu^{-es}\times\tau$ is distinguished. We apply Proposition~\ref{prop:three-alternatives} with $k=n-2\ge2$, $\pi_1=\nu^{-es}$ and $\pi_2=\tau$.
	\begin{itemize}
		\item Condition (II) would imply that the character $\nu^{-es}$ of $G_{n-2}$ is $(H_{1,n-3},\chi_{s-1})$-distinguished. Since $n-2\ge2$, Lemma~\ref{lem:char-dist} would then force $s-1=0$ and $\nu^{-es}=\one_{n-2}$, that is, $s=1$ and $s=0$, which is absurd.
		\item Condition (III) would give $\nu^{-es}\in\mathcal O_1(s,n-2;\sigma)$ for some $\sigma\in\Irr(G_1)$. Since $k=n-2\ge2$, Remark~\ref{rem:O-characters} shows that the only character that can lie in such a set is $\nu^{\frac{e(n-k)}2-es}=\nu^{e-es}$, which is different from $\nu^{-es}$.
	\end{itemize}
	Hence, condition (I) holds. Since $\nu^{\frac{e(n-k)-2e}2-es}=\nu^{-es}$ for $k=n-2$, condition (I) says exactly that $\tau$ is $(H_{1,1},\chi_{s+\frac{n-2}2})$-distinguished.
\end{proof}

The case where $\nu^{-es}$ is the second factor is analogous.

\begin{lemma}\label{lem:quot-tau-nu}
	Let $n\ge4$, $s\in\R$ and $\tau\in\Irr(G_2)$, and let $\pi\in\Irr(G_n)$ be $(H,\chi_s)$-distinguished. If $\pi$ is a quotient of the representation $\tau\times\nu^{-es}$ induced from $G_2\times G_{n-2}$, then $\tau$ is $(H_{1,1},\chi_{s-\frac{n-2}2})$-distinguished.
\end{lemma}

\begin{proof}
	By Lemma~\ref{lem:dist-quotient}, the representation $\tau\times\nu^{-es}$ is distinguished. We apply Proposition~\ref{prop:three-alternatives} with $k=2$, $\pi_1=\tau$ and $\pi_2=\nu^{-es}$.
	\begin{itemize}
		\item Condition (I) would imply that the character $\nu^{-es}$ of $G_{n-2}$ is $(H_{1,n-3},\chi_{s+1})$-distinguished. As in the proof of Lemma~\ref{lem:quot-nu-tau}, Lemma~\ref{lem:char-dist} would force $s=-1$ and $s=0$, which is absurd.
		\item Condition (III) would give $\nu^{-es}\in\mathcal O_2(s,n-2;\sigma)$ for some $\sigma\in\Irr(G_1)$. Since $n-k=n-2\ge2$, Remark~\ref{rem:O-characters} shows that the only character that can lie in such a set is $\nu^{-\frac{ek}2-es}=\nu^{-e-es}$, which is different from $\nu^{-es}$.
	\end{itemize}
	Hence, condition (II) holds. Since $\nu^{\frac{e(2-k)}2-es}=\nu^{-es}$ for $k=2$, condition (II) says exactly that $\tau$ is $(H_{1,1},\chi_{s-\frac{n-2}2})$-distinguished.
\end{proof}

Combining these two lemmas with the classification for $n=2$, we obtain the following consequences.

\begin{lemma}\label{lem:tau-T-left}
	Let $n\ge4$, $s\in\R$ and $\tau\in\Irr(G_2)$, and let $\pi\in\Irr(G_n)$ be $(H,\chi_s)$-distinguished. If $\pi$ is a quotient of $\nu^{-es}\times\tau$, then either $\tau\in\T$, or $\tau=\one_2$ and $s=-\frac{n-2}2$.
\end{lemma}

\begin{proof}
	By Lemma~\ref{lem:quot-nu-tau}, $\tau$ is $(H_{1,1},\chi_{s+\frac{n-2}2})$-distinguished. If $\tau$ is infinite-dimensional, then $\tau\in\T$ by Corollary~\ref{cor:T-chi-t}. Otherwise, $\tau$ is a character by Lemma~\ref{lem:G2-finite-dim}, and Lemma~\ref{lem:char-dist} for $G_2$ gives $\tau=\one_2$ and $s+\frac{n-2}2=0$.
\end{proof}

\begin{lemma}\label{lem:tau-T-right}
	Let $n\ge4$, $s\in\R$ and $\tau\in\Irr(G_2)$, and let $\pi\in\Irr(G_n)$ be $(H,\chi_s)$-distinguished. If $\pi$ is a quotient of $\tau\times\nu^{-es}$, then either $\tau\in\T$, or $\tau=\one_2$ and $s=\frac{n-2}2$.
\end{lemma}

\begin{proof}
	By Lemma~\ref{lem:quot-tau-nu}, $\tau$ is $(H_{1,1},\chi_{s-\frac{n-2}2})$-distinguished. The rest of the argument is the same as in the proof of Lemma~\ref{lem:tau-T-left}, and gives either $\tau\in\T$, or $\tau=\one_2$ and $s-\frac{n-2}2=0$.
\end{proof}

\subsection{Twists of \texorpdfstring{$Q_n$}{Qn} and \texorpdfstring{$\widetilde Q_n$}{Q~n}}\label{ss:Qn-twists}

The following lemmas strengthen Lemmas~\ref{lem:Q-dist} and~\ref{lem:Qtilde-dist}. In those lemmas the twist is tied to the parameter $s$; here we allow arbitrary twists. We first treat twists of $Q_n$ for $n\ge4$.

\begin{lemma}\label{lem:Qn-twist-dist}
	Let $n\ge4$ and $s,t\in\R$. Then $\nu^{-et}Q_n$ is $(H,\chi_s)$-distinguished if and only if $s=t=\frac n2$.
\end{lemma}

\begin{proof}
	The ``if'' part is Lemma~\ref{lem:Q-dist}. Conversely, suppose that $\nu^{-et}Q_n$ is distinguished. By Lemma~\ref{lem:Qm-quotient}, $\nu^{-et}Q_n$ is a quotient of $\pi_1\times\pi_2$, where $\pi_1=\nu^{e(\frac n2-t)}\St_2\in\Irr(G_2)$ and $\pi_2=\nu^{-et}\in\Irr(G_{n-2})$. Hence, $\pi_1\times\pi_2$ is distinguished by Lemma~\ref{lem:dist-quotient}, and we may apply Proposition~\ref{prop:three-alternatives} with $k=2$.
	
	Condition (I) fails, since $\pi_1$ is not a character. Suppose that condition (II) holds. Then $\pi_2\simeq\nu^{-es}$, so $t=s$, and $\pi_1$ is $(H_{1,1},\chi_{s-\frac{n-2}2})$-distinguished. By Lemma~\ref{lem:central}, the central character $z\mapsto|z|_{\mathrm{F}}^{2e^2(\frac n2-t)}$ of $\pi_1$ is trivial, so $t=\frac n2$. Hence, $s=t=\frac n2$.
	
	Suppose now that condition (III) holds, so that $\pi_1\in\mathcal O_1(s,2;\sigma)$ for some $\sigma\in\Irr(G_1)$. By \eqref{eq:O1} and \eqref{eq:Q2}, the elements of the sets $\mathcal O_1(s,2;\sigma)$ are characters, the twist $\nu^{\frac{e(n-2)}2-es}\Qt_2=\nu^{\frac{e(n-4)}2-es}\St_2$, and irreducible representations induced from $G_1\times G_1$. By Remark~\ref{rem:O-characters}, the latter are not twists of $\St_2$. Hence, $\pi_1=\nu^{\frac{e(n-4)}2-es}\St_2$, that is, $t=s+2$, and $\sigma\simeq\nu^{\frac{e(n-5)}2-es}$. Then $\pi_2=\nu^{-e(s+2)}\in\mathcal O_2(s,n-2;\sigma)$. Since $n-2\ge2$, Remark~\ref{rem:O-characters} shows that the only character that can belong to $\mathcal O_2(s,n-2;\sigma)$ is $\nu^{-e-es}$, which is different from $\nu^{-e(s+2)}$. This is a contradiction.
\end{proof}

The analogous statement holds for twists of $\Qt_n$.

\begin{lemma}\label{lem:Qtn-twist-dist}
	Let $n\ge4$ and $s,t\in\R$. Then $\nu^{-et}\Qt_n$ is $(H,\chi_s)$-distinguished if and only if $s=t=-\frac n2$.
\end{lemma}

\begin{proof}
	The ``if'' part is Lemma~\ref{lem:Qtilde-dist}. Conversely, suppose that $\nu^{-et}\Qt_n$ is distinguished. By Lemma~\ref{lem:Qtm-quotient}, $\nu^{-et}\Qt_n$ is a quotient of $\pi_1\times\pi_2$, where $\pi_1=\nu^{-et}\in\Irr(G_{n-2})$ and $\pi_2=\nu^{-e(t+\frac n2)}\St_2\in\Irr(G_2)$. Hence, $\pi_1\times\pi_2$ is distinguished by Lemma~\ref{lem:dist-quotient}, and we may apply Proposition~\ref{prop:three-alternatives} with $k=n-2\ge2$.
	
	Suppose that condition (I) holds. Then $\pi_1\simeq\nu^{-es}$, so $t=s$, and $\pi_2$ is $(H_{1,1},\chi_{s+\frac{n-2}2})$-distinguished. Hence, the central character of $\pi_2$ is trivial by Lemma~\ref{lem:central}, which gives $t=-\frac n2$. Thus, $s=t=-\frac n2$. Condition (II) fails, since $\pi_2$ is not a character.
	
	Suppose now that condition (III) holds. Then the character $\pi_1$ of $G_{n-2}$ lies in some set $\mathcal O_1(s,n-2;\sigma)$. By Remark~\ref{rem:O-characters} and \eqref{eq:O1}, this forces $\pi_1=\nu^{e-es}$, that is, $t=s-1$, and $\sigma\simeq\nu^{\frac{e(n-1)}2-es}$. Moreover, $\pi_2\in\mathcal O_2(s,2;\sigma)$. By \eqref{eq:O2} with $k=n-2$ and by \eqref{eq:Q2}, the only twist of $\St_2$ that can lie in such a set is $\nu^{-\frac{e(n-2)}2-es}Q_2=\nu^{-\frac{e(n-4)}2-es}\St_2$, and it occurs only if $\sigma\simeq\nu^{es+\frac{e(n-5)}2}$. Since $\pi_2$ is a twist of $\St_2$, we obtain $\nu^{\frac{e(n-1)}2-es}\simeq\nu^{es+\frac{e(n-5)}2}$, that is, $s=1$, and hence $t=0$. But then $\pi_2=\nu^{-\frac{en}2}\St_2$ would be isomorphic to $\nu^{-\frac{e(n-4)}2-e}\St_2=\nu^{-\frac{e(n-2)}2}\St_2$, which is absurd.
\end{proof}

For $n=3$, we deduce the corresponding statements from Theorem~\ref{thm:n3}.

\begin{lemma}\label{lem:Q3-twist-dist}
	Let $s,t\in\R$. Then $\nu^{-et}Q_3$ is $(H_{1,2},\chi_s)$-distinguished if and only if $s=t=\frac32$.
\end{lemma}

\begin{proof}
	By Theorem~\ref{thm:n3}, it suffices to determine when a twist of $Q_3$ is listed for $(3,s)$. We treat twists of $\Qt_3$ at the same time, for use in Lemma~\ref{lem:Qt3-twist-dist}. By \eqref{eq:Q3}, the Zelevinsky multisegment of such a twist consists of two linked segments on the line of $\one_1$, one of which is a singleton. We compare this with the Zelevinsky multisegments of the listed representations.
	\begin{itemize}
		\item When they are irreducible, the representations $\xione_s$ and $\xitwo_s$ have Zelevinsky multisegments consisting of two unlinked segments (Lemmas~\ref{lem:sum-multisegments} and~\ref{lem:two-segments}).
		\item The Zelevinsky multisegment of $\one_3$ consists of one segment.
		\item Let $\nu^{-es}\times\tau$ be irreducible, with $\tau\in\T$, and write $\tau=\Ss(\bb)$. By Lemma~\ref{lem:sum-multisegments}, its Zelevinsky multisegment is $\{[-s]\}+\bb$. If the supercuspidal support of $\tau$ lies on the line of $\one_1$, then, by the proof of Lemma~\ref{lem:reducible-nu-tau}, either $\tau=\St_2=\Ss(\{[-\frac12],[\frac12]\})$ or $\tau=\tau_x=\Ss(\{[x],[-x]\})$. In both cases, $\{[-s]\}+\bb$ consists of three segments. Otherwise, $\bb$ contains a segment that is not on the line of $\one_1$.
	\end{itemize}
	In none of these cases is the Zelevinsky multisegment of the form described above. Hence, by Theorem~\ref{thm:zel} and Lemma~\ref{lem:Q-vs-Qt}, a twist of $Q_3$ is listed for $(3,s)$ only if it is $\nu^{-3e/2}Q_3$ and $s=\frac32$. Similarly, a twist of $\Qt_3$ is listed for $(3,s)$ only if it is $\nu^{3e/2}\Qt_3$ and $s=-\frac32$. Conversely, $\nu^{-3e/2}Q_3$ is listed for $(3,\frac32)$. Finally, by Theorem~\ref{thm:zel}, we have $\nu^{-et}Q_3\simeq\nu^{-3e/2}Q_3$ if and only if $t=\frac32$.
\end{proof}

\begin{lemma}\label{lem:Qt3-twist-dist}
	Let $s,t\in\R$. Then $\nu^{-et}\Qt_3$ is $(H_{1,2},\chi_s)$-distinguished if and only if $s=t=-\frac32$.
\end{lemma}

\begin{proof}
	By the proof of Lemma~\ref{lem:Q3-twist-dist}, a twist of $\Qt_3$ is listed for $(3,s)$ if and only if it is $\nu^{3e/2}\Qt_3$ and $s=-\frac32$. Moreover, $\nu^{-et}\Qt_3\simeq\nu^{3e/2}\Qt_3$ if and only if $t=-\frac32$, by Theorem~\ref{thm:zel}. Hence, Theorem~\ref{thm:n3} gives the assertion.
\end{proof}

Combining Lemmas~\ref{lem:Qn-twist-dist}, \ref{lem:Qtn-twist-dist}, \ref{lem:Q3-twist-dist} and~\ref{lem:Qt3-twist-dist}, we obtain, for every $m\ge3$ and all $s,t\in\R$,
\begin{equation}\label{eq:Q-twist-all}
	\begin{aligned}
		\nu^{-et}Q_m\ &\text{is }(H_{1,m-1},\chi_s)\text{-distinguished}\iff s=t=\tfrac m2,\\
		\nu^{-et}\Qt_m\ &\text{is }(H_{1,m-1},\chi_s)\text{-distinguished}\iff s=t=-\tfrac m2.
	\end{aligned}
\end{equation}

\subsection{Twists of \texorpdfstring{$S_n$}{Sn} and \texorpdfstring{$\widetilde S_n$}{S~n}}\label{ss:Sn-twists}

Recall the representations $S_n=\Ss(\sss_n)$ and $\Stil_n=\Ss(\widetilde\sss_n)$ of \S\ref{ss:Sn}. We now show that no twist of these representations is distinguished. In both cases, we realize the twist as a quotient of a representation induced from a character of $G_1$ and a representation of $G_{n-1}$, and we apply Lemmas~\ref{lem:quot-X-nu} and~\ref{lem:quot-nu-X}.

\begin{lemma}\label{lem:Sn-not-dist}
	Let $n\ge4$ and $s,t\in\R$. Then $\nu^{-et}S_n$ is not $(H,\chi_s)$-distinguished.
\end{lemma}

\begin{proof}
	The Zelevinsky multisegment of $\pi=\nu^{-et}S_n$ is $\{\Delta_{n-2}(-t),[\frac{n-1}2-t],[\frac{n+1}2-t]\}$, and the ordering $([\frac{n+1}2-t],[\frac{n-1}2-t],\Delta_{n-2}(-t))$ of it is admissible. Hence, by Lemmas~\ref{lem:S-quotient} and~\ref{lem:quotient-constituents}, $\pi$ is a quotient of $X\times\nu^{ey}$, where $y=\frac{n+1}2-t$ and $X$ is a constituent of $\Ss(\Delta_{n-2}(-t))\times\Ss([\frac{n-1}2-t])$. Since $\Delta_{n-2}(-t)\prec[\frac{n-1}2-t]$, Lemma~\ref{lem:two-segments} shows that either $X=\Ss(\Delta_{n-1}(\frac12-t))=\nu^{e(\frac12-t)}$, which corresponds to the union of the two segments, or
	$$X=\Ss(\{\Delta_{n-2}(-t),[\tfrac{n-1}2-t]\})=\nu^{-e(t+\frac12)}Q_{n-1}.$$
	Suppose that $\pi$ is distinguished. We go through the cases of Lemma~\ref{lem:quot-X-nu}.
	\begin{enumerate}[label=(\alph*)]
		\item Suppose that $y=s+\frac{n-1}2$ and $X\simeq\nu^{-e(s+\frac12)}$. Then $X$ is the character $\nu^{e(\frac12-t)}$, so $s=t-1$. Moreover, $y=\frac{n+1}2-t$ gives $s=1-t$. Hence, $t=1$ and $s=0$, and $\pi$ is a quotient of $\nu^{-e/2}\times\nu^{e\frac{n-1}2}=\xione_0$. By Lemma~\ref{lem:dps1-exceptional}, the unique irreducible quotient of $\xione_0$ is $\nu^{-e}Q_n$. However, $\pi=\nu^{-e}S_n\not\simeq\nu^{-e}Q_n$ by Theorem~\ref{thm:zel}, since their Zelevinsky multisegments have three and two segments, respectively.
		\item Suppose that $y=-s-\frac{n-3}2$ and $X$ is $(H_{1,n-2},\chi_{s-\frac12})$-distinguished. If $X$ is the character $\nu^{e(\frac12-t)}$ of $G_{n-1}$, then Lemma~\ref{lem:char-dist} gives $s=t=\frac12$. In that case, $y=\frac n2$, whereas $-s-\frac{n-3}2=-\frac{n-2}2$, a contradiction. If $X=\nu^{-e(t+\frac12)}Q_{n-1}$, then \eqref{eq:Q-twist-all} in rank $n-1\ge3$ gives $s-\frac12=t+\frac12=\frac{n-1}2$. In that case, $y=\frac32$, whereas $-s-\frac{n-3}2=-\frac{2n-3}2$, again a contradiction.
		\item Suppose that $y=s-\frac{n-1}2$ and $X\simeq\nu^{-e(s-\frac12)}$. Then $X$ is the character $\nu^{e(\frac12-t)}$, so $s=t$. Moreover, $\frac{n+1}2-t=t-\frac{n-1}2$, so $s=t=\frac n2$. Hence, $\pi$ is a quotient of $\nu^{e\frac{1-n}2}\times\nu^{e/2}=\xi^{(2)\prime}_{n/2}$, whose unique irreducible quotient is $\nu^{-en/2}Q_n$ by Lemma~\ref{lem:dps2-exceptional}. However, $\pi=\nu^{-en/2}S_n\not\simeq\nu^{-en/2}Q_n$ by Theorem~\ref{thm:zel}.
		\item[\textup{(d), (e)}] These cases do not occur. Indeed, the Zelevinsky multisegment of $X$ is either a single segment or consists of two linked segments, the singleton following the other segment. Hence, $X$ is neither a twist of $\Qt_{n-1}$ (Lemmas~\ref{lem:Qm-not-char} and~\ref{lem:Q-vs-Qt}) nor an irreducible representation induced from $G_{n-2}\times G_1$ (Remark~\ref{rem:O-characters}).
	\end{enumerate}
	In all cases, we reach a contradiction. Therefore, $\pi$ is not distinguished.
\end{proof}

The argument for $\Stil_n$ is symmetric, with the character of $G_1$ now in the first factor.

\begin{lemma}\label{lem:Stn-not-dist}
	Let $n\ge4$ and $s,t\in\R$. Then $\nu^{-et}\Stil_n$ is not $(H,\chi_s)$-distinguished.
\end{lemma}

\begin{proof}
	The Zelevinsky multisegment of $\pi=\nu^{-et}\Stil_n$ is $\{\Delta_{n-2}(-t),[-\frac{n-1}2-t],[-\frac{n+1}2-t]\}$, and the ordering $(\Delta_{n-2}(-t),[-\frac{n-1}2-t],[-\frac{n+1}2-t])$ of it is admissible. Hence, by Lemmas~\ref{lem:S-quotient} and~\ref{lem:quotient-constituents}, $\pi$ is a quotient of $\nu^{ey}\times X$, where $y=-\frac{n+1}2-t$ and $X$ is a constituent of $\Ss([-\frac{n-1}2-t])\times\Ss(\Delta_{n-2}(-t))$. Since $[-\frac{n-1}2-t]\prec\Delta_{n-2}(-t)$, Lemma~\ref{lem:two-segments} shows that either $X=\Ss(\Delta_{n-1}(-t-\frac12))=\nu^{-e(t+\frac12)}$, or
	$$X=\Ss(\{[-\tfrac{n-1}2-t],\Delta_{n-2}(-t)\})=\nu^{e(\frac12-t)}\Qt_{n-1}.$$
	Suppose that $\pi$ is distinguished. We go through the cases of Lemma~\ref{lem:quot-nu-X}.
	\begin{enumerate}[label=(\alph*)]
		\item Suppose that $y=-s+\frac{n-3}2$ and $X$ is $(H_{1,n-2},\chi_{s+\frac12})$-distinguished. If $X$ is the character $\nu^{-e(t+\frac12)}$, then $s=t=-\frac12$ by Lemma~\ref{lem:char-dist}. In that case, $y=-\frac n2$, whereas $-s+\frac{n-3}2=\frac{n-2}2$, a contradiction. If $X=\nu^{e(\frac12-t)}\Qt_{n-1}$, then \eqref{eq:Q-twist-all} in rank $n-1\ge3$ gives $s+\frac12=t-\frac12=-\frac{n-1}2$. In that case, $y=-\frac32$, whereas $-s+\frac{n-3}2=\frac{2n-3}2$, again a contradiction.
		\item Suppose that $y=s-\frac{n-1}2$ and $X\simeq\nu^{-e(s-\frac12)}$. Then $X$ is the character $\nu^{-e(t+\frac12)}$, so $-s+\frac12=-t-\frac12$. Moreover, $-\frac{n+1}2-t=s-\frac{n-1}2$. Hence, $t=-1$ and $s=0$, and $\pi$ is a quotient of $\nu^{-e\frac{n-1}2}\times\nu^{e/2}=\xitwo_0$. By Lemma~\ref{lem:dps2-exceptional}, the unique irreducible quotient of $\xitwo_0$ is $\nu^e\Qt_n$. However, $\pi=\nu^e\Stil_n\not\simeq\nu^e\Qt_n$ by Theorem~\ref{thm:zel}.
		\item Suppose that $y=s+\frac{n-1}2$ and $X\simeq\nu^{-e(s+\frac12)}$. Then $X$ is the character $\nu^{-e(t+\frac12)}$, so $s=t$. Moreover, $-\frac{n+1}2-t=t+\frac{n-1}2$, so $s=t=-\frac n2$. Hence, $\pi$ is a quotient of $\nu^{-e/2}\times\nu^{e\frac{n-1}2}=\xi^{(1)\prime}_{-n/2}$, whose unique irreducible quotient is $\nu^{en/2}\Qt_n$ by Lemma~\ref{lem:dps1-exceptional}. However, $\pi=\nu^{en/2}\Stil_n\not\simeq\nu^{en/2}\Qt_n$ by Theorem~\ref{thm:zel}.
		\item[\textup{(d), (e)}] These cases do not occur. Indeed, the Zelevinsky multisegment of $X$ is either a single segment or consists of two linked segments, the singleton preceding the other segment. Hence, $X$ is neither a twist of $Q_{n-1}$ (Lemmas~\ref{lem:Qm-not-char} and~\ref{lem:Q-vs-Qt}) nor an irreducible representation induced from $G_1\times G_{n-2}$ (Remark~\ref{rem:O-characters}).
	\end{enumerate}
	In all cases, we reach a contradiction. Therefore, $\pi$ is not distinguished.
\end{proof}

\subsection{Constituents of reducible representations \texorpdfstring{$\nu^{-es}\times\tau$}{nu times tau}}\label{ss:constituents-reducible}

Throughout this subsection, we assume that $n\ge4$. By Lemma~\ref{lem:reducible-nu-tau}, if $\tau\in\T$ and $\nu^{-es}\times\tau$ is reducible, then either $\tau=\St_2$ and $s=\pm\frac n2$, or $\tau=\tau_x$ for some $x\in\{L_s,R_s\}$ with $x\ne\pm\frac12$. We treat these cases one by one, beginning with $\tau=\St_2$ and $s=\frac n2$.

\begin{lemma}\label{lem:constituents-St2-plus}
	Let $n\ge4$ and $s=\frac n2$. If a constituent $\pi$ of $\nu^{-es}\times\St_2$ is $(H,\chi_s)$-distinguished, then $\pi\simeq\nu^{-en/2}Q_n$.
\end{lemma}

\begin{proof}
	We have $\nu^{-es}\times\St_2=\nu^{-en/2}(\one_{n-2}\times\nu^{en/2}\St_2)$, which has the same constituents as $\nu^{-en/2}(\nu^{en/2}\St_2\times\one_{n-2})$. Hence, by Lemma~\ref{lem:Sn-embedding} and \eqref{eq:b-le-sn}, $\pi$ is one of the representations $\nu^{-en/2}S_n$, $\nu^{-en/2}Q_n$, the character $\nu^{e(1-\frac n2)}$, and
	$$Z=\nu^{-en/2}\Ss\big(\{\Delta_{n-2},[\tfrac{n-1}2,\tfrac{n+1}2]\}\big).$$
	The representation $\nu^{-en/2}S_n$ is not distinguished by Lemma~\ref{lem:Sn-not-dist}, and the character is not distinguished by Lemma~\ref{lem:char-dist}, since $s\ne0$. It remains to show that $Z$ is not distinguished.
	
	We have $Z=\Ss(\{A_s,[-\frac12,\frac12]\})$, where $A_s=[-\frac{2n-3}2,-\frac32]\prec[-\frac12,\frac12]$. Hence, by Lemma~\ref{lem:S-quotient}, $Z$ is a quotient of $\Ss(A_s)\times\Ss([-\frac12,\frac12])=\nu^{-es}\times\one_2$, induced from $G_{n-2}\times G_2$. Suppose that $Z$ is distinguished. Then $\nu^{-es}\times\one_2$ is distinguished by Lemma~\ref{lem:dist-quotient}, and we apply Proposition~\ref{prop:three-alternatives} with $k=n-2$.
	\begin{itemize}
		\item Condition (I) would imply that $\one_2$ is $(H_{1,1},\chi_{s+\frac{n-2}2})=(H_{1,1},\chi_{n-1})$-distinguished, which contradicts Lemma~\ref{lem:char-dist}.
		\item Condition (II) would imply that the character $\nu^{-es}$ of $G_{n-2}$ is $(H_{1,n-3},\chi_{s-1})$-distinguished. By Lemma~\ref{lem:char-dist}, this forces $s=1$, which is not the case.
		\item Condition (III) would imply that $\nu^{-es}$ lies in some set $\mathcal O_1(s,n-2;\sigma)$. By Remark~\ref{rem:O-characters}, the only character in such a set is $\nu^{e-es}\ne\nu^{-es}$.
	\end{itemize}
	This is a contradiction. Hence, $Z$ is not distinguished, and therefore $\pi\simeq\nu^{-en/2}Q_n$.
\end{proof}

The case $s=-\frac n2$ is symmetric.

\begin{lemma}\label{lem:constituents-St2-minus}
	Let $n\ge4$ and $s=-\frac n2$. If a constituent $\pi$ of $\nu^{-es}\times\St_2$ is $(H,\chi_s)$-distinguished, then $\pi\simeq\nu^{en/2}\Qt_n$.
\end{lemma}

\begin{proof}
	We have $\nu^{-es}\times\St_2=\nu^{en/2}(\one_{n-2}\times\nu^{-en/2}\St_2)$. Hence, by Lemma~\ref{lem:Stn-embedding} and the list following \eqref{eq:b-le-sn}, $\pi$ is one of the representations $\nu^{en/2}\Stil_n$, $\nu^{en/2}\Qt_n$, the character $\nu^{e\frac{n-2}2}$, and
	$$Z=\nu^{en/2}\Ss\big(\{\Delta_{n-2},[-\tfrac{n+1}2,-\tfrac{n-1}2]\}\big).$$
	The representation $\nu^{en/2}\Stil_n$ is not distinguished by Lemma~\ref{lem:Stn-not-dist}, and the character is not distinguished by Lemma~\ref{lem:char-dist}, since $s\ne0$. It remains to show that $Z$ is not distinguished.
	
	We have $Z=\Ss(\{[-\frac12,\frac12],A_s\})$, where $[-\frac12,\frac12]\prec A_s=[\frac32,\frac{2n-3}2]$. Hence, by Lemma~\ref{lem:S-quotient}, $Z$ is a quotient of $\one_2\times\nu^{-es}$, induced from $G_2\times G_{n-2}$. Suppose that $Z$ is distinguished. Then $\one_2\times\nu^{-es}$ is distinguished by Lemma~\ref{lem:dist-quotient}, and we apply Proposition~\ref{prop:three-alternatives} with $k=2$.
	\begin{itemize}
		\item Condition (I) would give $\one_2\simeq\nu^{\frac{e(n-4)}2-es}=\nu^{e(n-2)}$, which is false since $n\ge4$.
		\item Condition (II) would imply that $\one_2$ is $(H_{1,1},\chi_{s-\frac{n-2}2})=(H_{1,1},\chi_{-(n-1)})$-distinguished, which contradicts Lemma~\ref{lem:char-dist}.
		\item Condition (III) would imply that $\nu^{-es}$ lies in some set $\mathcal O_2(s,n-2;\sigma)$. By Remark~\ref{rem:O-characters}, the only character in such a set is $\nu^{-e-es}\ne\nu^{-es}$.
	\end{itemize}
	This is a contradiction. Hence, $Z$ is not distinguished, and therefore $\pi\simeq\nu^{en/2}\Qt_n$.
\end{proof}

We now turn to the case $\tau=\tau_x$ with $x=R_s$.

\begin{lemma}\label{lem:constituents-R}
	Let $n\ge4$ and $s\in\R$, and put $x=R_s$. Suppose that $x\ne\pm\frac12$ and $-x\ne L_s$. If a constituent $\pi$ of $\nu^{-es}\times\tau_x$ is $(H,\chi_s)$-distinguished, then $\pi\simeq\xitwo_s$, and $\xitwo_s$ is irreducible.
\end{lemma}

\begin{proof}
	By Lemma~\ref{lem:constituents}, the constituents of $\nu^{-es}\times\tau_x$ are among the representations $\Ss(\bb)$ with $\bb\le\mm=\{A_s,[R_s],[-R_s]\}$. The only elementary operation on $\mm$ replaces $A_s$ and $[R_s]$ by $A_s^+$. Indeed, $[-R_s]$ is not linked to $[R_s]$, since $x\ne\pm\frac12$. If $R_s\ne0$, then $[-R_s]$ is linked neither to $A_s$ nor to $A_s^+$, since $-R_s\ne L_s$ and $R_s\ne-\frac12$. If $R_s=0$, then $[-R_s]=[R_s]$, and the two copies of $[R_s]$ give the same elementary operation. Thus, the only new multisegment is $\{A_s^+,[-R_s]\}$, and hence $\pi=\Ss(\mm)$ or $\pi=\Ss(\{A_s^+,[-R_s]\})=\xitwo_s$. Moreover, $\xitwo_s$ is irreducible by Lemma~\ref{lem:dps2-irred}, since $s=0$ would give $-R_s=L_s$, and $s=\frac n2$ would give $x=-\frac12$.
	
	Suppose that $\Ss(\mm)$ is distinguished. The ordering $([-R_s],[R_s],A_s)$ of $\mm$ is admissible. Hence, by Lemmas~\ref{lem:S-quotient} and~\ref{lem:quotient-constituents}, $\Ss(\mm)$ is a quotient of $X\times\nu^{-eR_s}$, where $X$ is a constituent of $\Ss(A_s)\times\Ss([R_s])$, that is, $X\in\{\Ss(\{A_s,[R_s]\}),\Ss(A_s^+)\}$. Here, $\Ss(A_s^+)=\nu^{-e(s-\frac12)}$, and, twisting Definition~\ref{def:Qm}, $\Ss(\{A_s,[R_s]\})=\nu^{-e(s+\frac12)}Q_{n-1}$. We apply Lemma~\ref{lem:quot-X-nu} with $y=-R_s=s-\frac{n-1}2$. Since $y\ne s+\frac{n-1}2$, cases (a), (d) and (e) of that lemma do not occur, and case (b) occurs only if $s=\frac12$. In case (c), and in case (b) with $X=\Ss(A_s^+)$, the representation $\Ss(\mm)$ is a quotient of the irreducible representation $\Ss(A_s^+)\times\Ss([-R_s])=\xi^{(2)\prime}_s$, and hence isomorphic to it. This contradicts Theorem~\ref{thm:zel}. In case (b) with $X=\nu^{-e}Q_{n-1}$ (and $s=\frac12$), \eqref{eq:Q-twist-all} in rank $n-1$ shows that $X$ is not $(H_{1,n-2},\chi_0)$-distinguished, which is again a contradiction.
\end{proof}

The case $x=L_s$ is symmetric.

\begin{lemma}\label{lem:constituents-L}
	Let $n\ge4$ and $s\in\R$, and put $x=L_s$. Suppose that $x\ne\pm\frac12$ and $-x\ne R_s$. If a constituent $\pi$ of $\nu^{-es}\times\tau_x$ is $(H,\chi_s)$-distinguished, then $\pi\simeq\xione_s$, and $\xione_s$ is irreducible.
\end{lemma}

\begin{proof}
	As in the proof of Lemma~\ref{lem:constituents-R}, the constituents are among $\Ss(\mm)$, with $\mm=\{A_s,[L_s],[-L_s]\}$, and $\Ss(\{A_s^-,[-L_s]\})=\xione_s$. Moreover, $\xione_s$ is irreducible by Lemma~\ref{lem:dps-irred}, since $s=0$ would give $-L_s=R_s$, and $s=-\frac n2$ would give $L_s=\frac12$.
	
	Suppose that $\Ss(\mm)$ is distinguished. The ordering $(A_s,[L_s],[-L_s])$ of $\mm$ is admissible. Hence, by Lemmas~\ref{lem:S-quotient} and~\ref{lem:quotient-constituents}, $\Ss(\mm)$ is a quotient of $\nu^{-eL_s}\times X$, where
	$$X\in\{\Ss(\{[L_s],A_s\}),\Ss(A_s^-)\}=\{\nu^{e(\frac12-s)}\Qt_{n-1},\nu^{-e(s+\frac12)}\}.$$
	We apply Lemma~\ref{lem:quot-nu-X} with $y=-L_s=s+\frac{n-1}2$. Since $y\ne s-\frac{n-1}2$, cases (b), (d) and (e) of that lemma do not occur, and case (a) occurs only if $s=-\frac12$. In case (c), and in case (a) with $X=\Ss(A_s^-)$, the representation $\Ss(\mm)$ is a quotient of the irreducible representation $\Ss([-L_s])\times\Ss(A_s^-)=\xi^{(1)\prime}_s$, and hence isomorphic to it. This contradicts Theorem~\ref{thm:zel}. In case (a) with $X=\nu^{e}\Qt_{n-1}$ (and $s=-\frac12$), \eqref{eq:Q-twist-all} in rank $n-1$ shows that $X$ is not $(H_{1,n-2},\chi_0)$-distinguished, which is again a contradiction.
\end{proof}

It remains to treat the case $s=0$, in which $\tau_{L_0}\simeq\tau_{R_0}$.

\begin{lemma}\label{lem:constituents-s0}
	Let $n\ge4$, $s=0$ and $x=R_0=\frac{n-1}2$. If a constituent $\pi$ of $\one_{n-2}\times\tau_{R_0}$ is $(H,\chi_0)$-distinguished, then $\pi\simeq\one_n$.
\end{lemma}

\begin{proof}
	Put $a=\frac{n-3}2\ge\frac12$, so that $A_0=[-a,a]$, $L_0=-a-1$ and $R_0=a+1$, and let $\mm=\{A_0,[L_0],[R_0]\}$. The segments $[L_0]$ and $[R_0]$ are not linked. Hence, the multisegments $\bb\le\mm$ are $\mm$, $\{A_0^-,[R_0]\}$, $\{[L_0],A_0^+\}$ and $\{\Delta_n\}$, and the corresponding representations are $\Ss(\mm)$, $\nu^{-e}Q_n$, $\nu^{e}\Qt_n$ and $\one_n$, respectively. The second and the third are not $(H,\chi_0)$-distinguished by \eqref{eq:Q-twist-all}, and $\one_n$ is listed. It remains to show that $\Ss(\mm)$ is not distinguished.
	
	The ordering $([R_0],A_0,[L_0])$ of $\mm$ is admissible. Hence, by Lemma~\ref{lem:S-quotient}, $\Ss(\mm)$ is the unique irreducible quotient of $\Ss([L_0])\times\Ss(A_0)\times\Ss([R_0])$. Since $[-a,a-1]\prec[a]$, Lemma~\ref{lem:two-segments} shows that $\Ss(A_0)$ is a quotient of $\Ss([a])\times\Ss([-a,a-1])$. Hence, $\Ss(\mm)$ is a quotient of
	$$
	\Ss([-a-1])\times\Ss([a])\times\Ss([-a,a-1])\times\Ss([a+1])\simeq\big(\Ss([a])\times\Ss([a+1])\big)\times\big(\Ss([-a-1])\times\Ss([-a,a-1])\big).$$
	Here, we have interchanged unlinked factors: since $a\ge\frac12$, none of the pairs $\{[-a-1],[a]\}$, $\{[-a,a-1],[a+1]\}$ and $\{[-a-1],[a+1]\}$ is linked. By Lemma~\ref{lem:quotient-constituents}, $\Ss(\mm)$ is a quotient, and hence a constituent, of $W_1\times W_2$, where $W_1$ is a constituent of $\Ss([a])\times\Ss([a+1])$ and $W_2$ is a constituent of $\Ss([-a-1])\times\Ss([-a,a-1])$. By Lemma~\ref{lem:two-segments},
	$$
	W_1\in\big\{\nu^{e\frac{n-2}2}\St_2,\ \nu^{e\frac{n-2}2}\big\},\qquad W_2\in\big\{\Ss(\{[-a-1],[-a,a-1]\}),\ \Ss([-a-1,a-1])=\nu^{-e}\big\},$$
	where $\nu^{-e}$ is regarded as a character of $G_{n-2}$. If $W_1$ and $W_2$ are both characters, then $W_1\times W_2=\Ss([a,a+1])\times\Ss([-a-1,a-1])$. Its constituents are among the representations $\Ss(\bb)$ with $\bb=\{[a,a+1],[-a-1,a-1]\}$ or $\bb=\{\Delta_n\}$ (Lemma~\ref{lem:constituents}), and neither of these multisegments is $\mm$. Hence, at least one of $W_1$ and $W_2$ is not a character.
	
	Suppose that $\Ss(\mm)$ is distinguished. By Lemma~\ref{lem:dist-quotient}, $W_1\times W_2$ is distinguished, and we apply Proposition~\ref{prop:three-alternatives} with $k=2$ and $s=0$.
	\begin{itemize}
		\item Condition (I) requires $W_1\simeq\nu^{e\frac{n-4}2}$, which is false.
		\item Condition (II) requires $W_2\simeq\one_{n-2}$, which is false.
		\item Condition (III) requires $W_1\in\mathcal O_1(0,2;\sigma)$ and $W_2\in\mathcal O_2(0,n-2;\sigma)$ for some $\sigma\in\Irr(G_1)$. By \eqref{eq:O1} and \eqref{eq:Q2}, the elements of the sets $\mathcal O_1(0,2;\sigma)$ are the character $\nu^{e\frac{n-2}2}$ (for $\sigma\simeq\nu^{e\frac{n-1}2}$), the twist $\nu^{e\frac{n-2}2}\Qt_2=\nu^{e\frac{n-4}2}\St_2$, and irreducible representations induced from $G_1\times G_1$, which are not twists of $\St_2$. Hence, $W_1=\nu^{e\frac{n-2}2}\St_2$ is impossible, and $W_1=\nu^{e\frac{n-2}2}$ forces $\sigma\simeq\nu^{e\frac{n-1}2}$. In that case, \eqref{eq:O2} with $k=2$ and $n-2\ge2$ gives $\mathcal O_2(0,n-2;\sigma)=\{\nu^{-e}\}$, so $W_2=\nu^{-e}$. But then both $W_1$ and $W_2$ are characters, which we have excluded.
	\end{itemize}
	This is a contradiction. Hence, $\Ss(\mm)$ is not distinguished, and therefore $\pi\simeq\one_n$.
\end{proof}

Combining Lemmas~\ref{lem:constituents-St2-plus}--\ref{lem:constituents-s0}, we obtain the following result.

\begin{proposition}\label{prop:constituents-reducible}
	Let $n\ge4$, $s\in\R$ and $\tau\in\T$, and suppose that $\nu^{-es}\times\tau$ is reducible. If a constituent $\pi$ of $\nu^{-es}\times\tau$ is $(H,\chi_s)$-distinguished, then $\pi$ is listed for $(n,s)$ in \textup{(a)} or \textup{(b)}.
\end{proposition}

\begin{proof}
	By Lemma~\ref{lem:reducible-nu-tau}, either $\tau=\St_2$ and $s=\pm\frac n2$, or $\tau=\tau_x$ for some $x\in\{L_s,R_s\}$ with $x\ne\pm\frac12$. The first case is treated in Lemmas~\ref{lem:constituents-St2-plus} and~\ref{lem:constituents-St2-minus}. In the second case, if $s\ne0$, then $-R_s\ne L_s$ and $-L_s\ne R_s$, so Lemmas~\ref{lem:constituents-R} and~\ref{lem:constituents-L} apply. If $s=0$, then $\tau_{L_0}=\tau_{-R_0}\simeq\tau_{R_0}$, and Lemma~\ref{lem:constituents-s0} applies. In every case, $\pi$ is one of $\xione_s$, $\xitwo_s$, $\one_n$, $\nu^{-en/2}Q_n$ and $\nu^{en/2}\Qt_n$, at the relevant values of $s$. All of these are listed in (a) or (b).
\end{proof}

As a consequence, we can treat all quotients of representations induced from $\nu^{-es}$ and an irreducible representation of $G_2$.

\begin{corollary}\label{cor:quotient-nu-tau-listed}
	Let $n\ge4$ and $s\in\R$. Let $\pi\in\Irr(G_n)$ be $(H,\chi_s)$-distinguished, and suppose that $\pi$ is a quotient of $\nu^{-es}\times\tau$ or of $\tau\times\nu^{-es}$ for some $\tau\in\Irr(G_2)$. Then $\pi$ is listed for $(n,s)$.
\end{corollary}

\begin{proof}
	By Lemmas~\ref{lem:tau-T-left} and~\ref{lem:tau-T-right}, either $\tau\in\T$, or $\tau=\one_2$ and $s=\mp\frac{n-2}2$. Suppose first that $\tau=\one_2$. Then, by Lemma~\ref{lem:two-segments}, the representations $\nu^{-es}\times\one_2$ (for $s=-\frac{n-2}2$) and $\one_2\times\nu^{-es}$ (for $s=\frac{n-2}2$) have the unique irreducible quotients
	$$
	\Ss([-\tfrac12,\tfrac{2n-5}2])\times\Ss([\tfrac12])=\xione_s\qquad\text{and}\qquad\Ss([-\tfrac{2n-5}2,\tfrac12])\times\Ss([-\tfrac12])\simeq\xitwo_s,$$
	respectively. These representations are irreducible and listed, since $s\notin\{0,\pm\frac n2\}$. Suppose now that $\tau\in\T$. If $\nu^{-es}\times\tau$ is irreducible, then $\pi\simeq\nu^{-es}\times\tau\simeq\tau\times\nu^{-es}$, which is listed in (c). Otherwise, $\pi$ is a constituent of $\nu^{-es}\times\tau$, and Proposition~\ref{prop:constituents-reducible} applies.
\end{proof}

\section{The classification for \texorpdfstring{$n\ge4$}{n>=4}}\label{sec:induction}

In this section, we prove part (i) of Theorem~\ref{thm:A} by induction on $n$, starting from the case $n=3$ (Theorem~\ref{thm:n3}). We keep the notation of \S\ref{ss:dps}, and we put $H=H_{1,n-1}$.

\subsection{The inductive input}\label{ss:inductive-input}

Let $n\ge4$, and suppose that the classification holds in rank $n-1$. That is, for every $t\in\R$, every $(H_{1,n-2},\chi_t)$-distinguished irreducible representation of $G_{n-1}$ is listed for $(n-1,t)$. Let $s\in\R$, and let $\pi'\in\Irr(G_{n-1})$ be $(H_{1,n-2},\chi_{s-\frac12})$-distinguished. We now express the representations listed for $(n-1,s-\frac12)$ in terms of the segments attached to $(n,s)$.

Put $s'=s-\frac12$, and denote by $L'_{s'}$, $R'_{s'}$, $A'_{s'}$, $\xi'^{(1)}_{s'}$ and $\xi'^{(2)}_{s'}$ the objects of \S\ref{ss:dps} in rank $n-1$. Then
$$
L'_{s'}=-s'-\tfrac{n-2}2=L_s+1,\qquad R'_{s'}=-s'+\tfrac{n-2}2=R_s,\qquad A'_{s'}=\Delta_{n-3}(\tfrac12-s)=[L_s+2,R_s-1],$$
and
$$
\xi'^{(1)}_{s'}=\nu^{-es}\times\nu^{-e(L_s+1)},\qquad\xi'^{(2)}_{s'}=\nu^{-eR_s}\times\nu^{e(1-s)},$$
induced from $G_{n-2}\times G_1$ and from $G_1\times G_{n-2}$, respectively. Moreover, the exceptional values $s'=0$, $s'=-\frac{n-1}2$ and $s'=\frac{n-1}2$ correspond to $s=\frac12$, $s=1-\frac n2$ and $s=\frac n2$, respectively. Hence, by the induction hypothesis, $\pi'$ is isomorphic to one of the following representations:
\begin{enumerate}[label=\textup{(\alph*)}]
	\item $\one_{n-1}$, if $s=\frac12$;
	\item $\nu^{e\frac{n-1}2}\Qt_{n-1}$, if $s=1-\frac n2$;
	\item $\nu^{-e\frac{n-1}2}Q_{n-1}$, if $s=\frac n2$;
	\item $\nu^{-es}\times\nu^{-e(L_s+1)}$, induced from $G_{n-2}\times G_1$, if $s\notin\{1-\frac n2,\frac12\}$;
	\item $\nu^{-eR_s}\times\nu^{e(1-s)}$, induced from $G_1\times G_{n-2}$, if $s\notin\{\frac12,\frac n2\}$;
	\item $\nu^{e(\frac12-s)}\times\tau$, induced from $G_{n-3}\times G_2$, where $\tau\in\T$ and $\nu^{e(\frac12-s)}\times\tau$ is irreducible.
\end{enumerate}
Here, (a)--(c) are the exceptional members of the two degenerate principal series in rank $n-1$, and (d)--(f) are the three families of Definition~\ref{def:listed}, with $n$ replaced by $n-1$ and $s$ by $s-\frac12$.

\subsection{A lemma for \texorpdfstring{$n=4$}{n=4}}\label{ss:n4-lemma}

In the inductive step, case (f) requires an additional argument when $n=4$, and we give this argument first. We use the following description of the Langlands multisegments of the representations listed for $n=3$.

\begin{lemma}\label{lem:n3-Langlands-listed}
	Let $t\in\R$, and let $\pi$ be listed for $(3,t)$. Then the Langlands multisegment of $\pi$ satisfies one of the following:
	\begin{enumerate}
		\item it consists of three singletons;
		\item it is $\{[-\frac12,\frac12],[\frac32]\}$, and then $t=-\frac32$;
		\item it is $\{[-\frac32],[-\frac12,\frac12]\}$, and then $t=\frac32$;
		\item it contains the singleton $[-t]$.
	\end{enumerate}
\end{lemma}

\begin{proof}
	Suppose that $\pi=\xione_t$ is irreducible. Then $\pi=\Ss(A_t^-)\times\Ss([-L_t])$, so by Lemmas~\ref{lem:char-Langlands} and~\ref{lem:sum-multisegments}, its Langlands multisegment consists of the three singletons $[L_t]$, $[L_t+1]$ and $[-L_t]$. The same argument applies to $\xitwo_t$ and to $\one_3$. By \eqref{eq:Q3}, we have $\nu^{3e/2}\Qt_3=\Qq(\{[-\frac12,\frac12],[\frac32]\})$ and $\nu^{-3e/2}Q_3=\Qq(\{[-\frac32],[-\frac12,\frac12]\})$; these are listed only for $t=-\frac32$ and $t=\frac32$, respectively. Finally, suppose that $\pi=\nu^{-et}\times\tau$ is irreducible, with $\tau=\Qq(\bb)\in\T$. Then $\pi=\Qq(\{[-t]\}+\bb)$ by Lemma~\ref{lem:sum-multisegments}.
\end{proof}

We can now prove the additional lemma needed for $n=4$.

\begin{lemma}\label{lem:n4-special}
	Let $s\in\R$ and $\tau\in\T$. Then no irreducible quotient of the representation $\tau\times\nu^{-es}\St_2$ induced from $G_2\times G_2$ is $(H_{1,3},\chi_s)$-distinguished.
\end{lemma}

\begin{proof}
	Let $\pi$ be an irreducible quotient of $\tau\times\nu^{-es}\St_2$, and suppose that $\pi$ is distinguished. By Lemma~\ref{lem:dist-quotient}, $\tau\times\nu^{-es}\St_2$ is distinguished, and we may apply Proposition~\ref{prop:three-alternatives} with $n=4$, $k=2$, $\pi_1=\tau$ and $\pi_2=\nu^{-es}\St_2$. Conditions (I) and (II) fail, because neither $\tau$ nor $\nu^{-es}\St_2$ is a character. Hence, condition (III) holds, that is, there is $\sigma\in\Irr(G_1)$ such that $\tau\in\mathcal O_1(s,2;\sigma)$ and $\nu^{-es}\St_2\in\mathcal O_2(s,2;\sigma)$.
	
	By \eqref{eq:O2} with $n=4$ and $k=2$, together with \eqref{eq:Q2}, the only twist of $\St_2$ that can lie in a set $\mathcal O_2(s,2;\sigma)$ is $\nu^{-e-es}Q_2=\nu^{-es}\St_2$, and it occurs only for $\sigma\simeq\nu^{es-\frac e2}$. We now apply \eqref{eq:O1} with $n=4$ and $k=2$, so that $x_s=\frac e2-es$.
	\begin{itemize}
		\item If $s=1$, then $\sigma\simeq\nu^{\frac{3e}2-es}$, and $\tau$ is a character, which is impossible.
		\item If $s=0$, then $\sigma\simeq\nu^{-\frac e2-es}$, and $\tau=\nu^{e-es}\Qt_2=\St_2$.
		\item Otherwise, $\tau=\nu^{e(\frac12-s)}\times\nu^{e(s-\frac12)}=\tau_x$ with $x=\frac12-s\ne\pm\frac12$.
	\end{itemize}
	We show that each of the remaining possibilities leads to a contradiction.
	
	\emph{The case $s=0$ and $\tau=\St_2$.} The representation $\St_2\times\St_2=\Qq([-\frac12,\frac12])\times\Qq([-\frac12,\frac12])$ is irreducible by Lemma~\ref{lem:two-segments}, so $\pi\simeq\St_2\times\St_2$. By Lemma~\ref{lem:St2}, $\St_2$ is a quotient of $\nu^{-e/2}\times\nu^{e/2}$. Hence, $\pi$ is a quotient of $\nu^{-e/2}\times X$, where $X=\nu^{e/2}\times\St_2\in\Rep(G_3)$. Since $[\frac12]\subseteq[-\frac12,\frac12]$, the representation $X$ is irreducible, and $X=\Qq(\{[\frac12],[-\frac12,\frac12]\})$ by Lemma~\ref{lem:sum-multisegments}. We apply Lemma~\ref{lem:quot-nu-X} with $n=4$, $s=0$ and $y=-\frac12$. Since $y\notin\{\frac12,\pm\frac32\}$, only case (e) can occur, and it gives $X\simeq\nu^{e/2}\times\one_2$, with this representation irreducible. But then $X=\Qq(\{[\frac12],[\frac12],[-\frac12]\})$ by Lemmas~\ref{lem:char-Langlands} and~\ref{lem:sum-multisegments}, which contradicts Theorem~\ref{thm:lang}.
	
	\emph{The case $s\notin\{0,1\}$ and $\tau=\tau_x$, with $x=\frac12-s$.} Since $\tau_x=\nu^{ex}\times\nu^{-ex}$, the representation $\pi$ is a quotient of $\nu^{ex}\times Y$, where
	$$
	Y=\nu^{-ex}\times\nu^{-es}\St_2=\Qq([s-\tfrac12])\times\Qq([-s-\tfrac12,-s+\tfrac12])\in\Rep(G_3).$$
	These two segments are linked if and only if $s-\frac12\in\{-s-\frac32,-s+\frac32\}$, that is, $s\in\{-\frac12,1\}$.
	
	Suppose first that $s\ne-\frac12$. Then $Y$ is irreducible, and $Y=\Qq(\{[s-\frac12],[-s-\frac12,-s+\frac12]\})$. We apply Lemma~\ref{lem:quot-nu-X} with $n=4$ and $y=\frac12-s$.
	\begin{itemize}
		\item Since $y=-s+\frac{n-3}2$, case (a) may occur. In that case, $Y$ is $(H_{1,2},\chi_{s+\frac12})$-distinguished, and hence listed for $(3,s+\frac12)$ by Theorem~\ref{thm:n3}. However, the Langlands multisegment of $Y$ consists of one singleton $[s-\frac12]$ and one segment of length two, and $[s-\frac12]=[-(s+\frac12)]$ only if $s=0$. Hence, by Lemma~\ref{lem:n3-Langlands-listed}, either $Y=\nu^{-3e/2}Q_3$ and $s+\frac12=\frac32$, which is excluded since $s\ne1$, or $Y=\nu^{3e/2}\Qt_3$ and $s+\frac12=-\frac32$. In the latter case, $s=-2$, and $Y=\Qq(\{[-\frac52],[\frac32,\frac52]\})\ne\Qq(\{[-\frac12,\frac12],[\frac32]\})$. Thus, case (a) is impossible.
		\item Cases (b) and (d) require $y=s-\frac32$, that is, $s=1$, and case (c) requires $y=s+\frac32$, that is, $s=-\frac12$. Both are excluded.
		\item In case (e), $Y\simeq\nu^{e(s-\frac12)}\times\nu^{-es}$, where $\nu^{-es}$ is a character of $G_2$, and this representation is irreducible. Hence, $Y=\Qq(\{[s-\frac12],[-s-\frac12],[-s+\frac12]\})$ by Lemmas~\ref{lem:char-Langlands} and~\ref{lem:sum-multisegments}, which contradicts Theorem~\ref{thm:lang}.
	\end{itemize}
	
	Suppose now that $s=-\frac12$, so that $x=1$. Then $Y=\Qq([-1])\times\Qq([0,1])$. By Lemma~\ref{lem:two-segments} and \eqref{eq:Q3}, the constituents of $Y$ are $\Qq([-1,1])=\St_3$ and $\Qq(\{[-1],[0,1]\})=\nu^{-e}Q_3$. By Lemma~\ref{lem:quotient-constituents}, $\pi$ is a quotient of $\nu^e\times X$ for some $X\in\{\St_3,\nu^{-e}Q_3\}$. We apply Lemma~\ref{lem:quot-nu-X} with $n=4$, $s=-\frac12$ and $y=1$.
	\begin{itemize}
		\item Case (a) would imply that $X$ is $(H_{1,2},\chi_0)$-distinguished. This is false for $\St_3$ by Theorem~\ref{thm:ds-not-dist}, and for $\nu^{-e}Q_3$ by Lemma~\ref{lem:n3-Q3-s0}.
		\item Case (c) would give $X\simeq\nu^{-e(s+\frac12)}=\one_3$, which is false.
		\item Cases (b) and (d) require $y=s-\frac32=-2$, and case (e) requires $y\notin\{s\pm\frac32\}=\{1,-2\}$. Neither condition holds.
	\end{itemize}
	In all cases, we reach a contradiction, which completes the proof.
\end{proof}

\subsection{The inductive step}\label{ss:inductive-step}

We can now carry out the inductive step. Combined with Theorem~\ref{thm:n3} and Proposition~\ref{prop:listed-dist}, the following theorem proves part (i) of Theorem~\ref{thm:A}.

\begin{theorem}\label{thm:inductive-step}
	Let $n\ge4$, and suppose that, for every $t\in\R$, every $(H_{1,n-2},\chi_t)$-distinguished irreducible representation of $G_{n-1}$ is listed for $(n-1,t)$. Then, for every $s\in\R$, every $(H_{1,n-1},\chi_s)$-distinguished irreducible representation of $G_n$ is listed for $(n,s)$.
\end{theorem}

\begin{proof}
	Let $\pi\in\Irr(G_n)$ be $(H,\chi_s)$-distinguished. By Lemma~\ref{lem:quotient-pi-sigma}, $\pi$ is a quotient of $\pi'\times\sigma$ for some $\pi'\in\Irr(G_{n-1})$ and $\sigma\in\Irr(G_1)$. By Proposition~\ref{prop:three-cases}, one of the cases (I), (II) and (III) of that proposition occurs. In each case, we show either that $\pi$ is listed, or that the case cannot occur.
	
	\emph{Case \textup{(I)}.} Here $\pi'\times\sigma=\xione_s$. If $s\notin\{-\frac n2,0\}$, then $\xione_s$ is irreducible by Lemma~\ref{lem:dps-irred}, and therefore $\pi\simeq\xione_s$, which is listed. If $s=0$, then $\pi\simeq\nu^{-e}Q_n$ by Lemma~\ref{lem:dps1-exceptional}, and this representation is not $(H,\chi_0)$-distinguished by Lemma~\ref{lem:Qn-twist-dist}. If $s=-\frac n2$, then $\pi$ is the character $\nu^{e\frac{n-2}2}$ by Lemma~\ref{lem:dps1-exceptional}, which is impossible by Lemma~\ref{lem:char-dist}, since $s\ne0$.
	
	\emph{Case \textup{(III)}.} Let $\rho\in\Irr(G_1)$, with $\sigma\simeq\widetilde\rho$, be as in Proposition~\ref{prop:three-cases}. We distinguish three subcases.
	
	(III.1) \emph{$\rho\simeq\nu^{eR_s}$.} Then $\pi'\times\sigma=\Ss(A_s^+)\times\Ss([-R_s])=\xi^{(2)\prime}_s$. Hence, by Lemmas~\ref{lem:dps2-irred} and~\ref{lem:dps2-exceptional}, we have $\pi\simeq\xitwo_s$ if $s\notin\{0,\frac n2\}$, $\pi\simeq\one_n$ if $s=0$, and $\pi\simeq\nu^{-en/2}Q_n$ if $s=\frac n2$. All of these representations are listed.
	
	(III.2) \emph{$\rho\simeq\nu^{eL_s}$.} Then $\sigma=\nu^{-eL_s}$ and $\pi'=\Ss(\{[L_s],A_s\})$. Since $[L_s]\prec A_s$, Lemma~\ref{lem:S-quotient} shows that $\pi'$ is the unique irreducible quotient of $\Ss([L_s])\times\Ss(A_s)$. Hence, by Lemma~\ref{lem:quotient-constituents}, $\pi$ is a quotient of $\nu^{eL_s}\times X$ for some constituent $X$ of $\Ss(A_s)\times\Ss([-L_s])$. We apply Lemma~\ref{lem:quot-nu-X} with $y=L_s=-s-\frac{n-1}2$. Since $y\ne-s+\frac{n-3}2$, case (a) of that lemma does not occur. Moreover, cases (b) and (d) require $y=s-\frac{n-1}2$, that is, $s=0$, and case (c) requires $y=s+\frac{n-1}2$, that is, $s=-\frac{n-1}2$.
	\begin{itemize}
		\item Suppose that $s\notin\{0,-\frac{n-1}2\}$. Then only case (e) can occur, and it gives $X\simeq\nu^{-eL_s}\times\nu^{-es}$. Hence, $\pi$ is a quotient of $(\nu^{eL_s}\times\nu^{-eL_s})\times\nu^{-es}$. By Lemma~\ref{lem:quotient-constituents}, $\pi$ is therefore a quotient of $\tau'\times\nu^{-es}$ for some constituent $\tau'$ of $\nu^{eL_s}\times\nu^{-eL_s}$, and Corollary~\ref{cor:quotient-nu-tau-listed} applies.
		\item Suppose that $s=-\frac{n-1}2$, so that $L_s=0$. Then only case (c) can occur, and it gives $X=\Ss(A_s^-)$. Hence, $\pi$ is a quotient of $\Ss([L_s])\times\Ss(A_s^-)=\Ss([-L_s])\times\Ss(A_s^-)=\xi^{(1)\prime}_s$. This representation is irreducible by Lemma~\ref{lem:dps-irred}, since $s\notin\{-\frac n2,0\}$. Therefore, $\pi\simeq\xione_s$, which is listed.
		\item Suppose that $s=0$, so that cases (b) and (d) may occur. In case (b), we have $X=\Ss(A_0^+)$, and $\pi$ is a quotient of $\Ss([L_0])\times\Ss(A_0^+)=\xitwo_0$, since $L_0=-R_0$. Hence, $\pi\simeq\nu^e\Qt_n$ by Lemma~\ref{lem:dps2-exceptional}, which is impossible by Lemma~\ref{lem:Qtn-twist-dist}. In case (d), we have $X=\nu^{-e/2}Q_{n-1}=\Ss(\{A_0,[R_0]\})$. Hence, $\pi$ is a constituent of $\Ss([L_0])\times\Ss(A_0)\times\Ss([R_0])$, which has the same constituents as $\one_{n-2}\times\tau_{R_0}$. By Lemma~\ref{lem:constituents-s0}, we obtain $\pi\simeq\one_n$, which is listed.
	\end{itemize}
	
	(III.3) \emph{$\rho\notin\{\nu^{eR_s},\nu^{eL_s}\}$.} Then $\pi'=\nu^{-es}\times\rho$ is irreducible, and $\pi$ is a quotient of $\nu^{-es}\times(\rho\times\widetilde\rho)$. By Lemma~\ref{lem:quotient-constituents}, $\pi$ is a quotient of $\nu^{-es}\times\tau'$ for some constituent $\tau'$ of $\rho\times\widetilde\rho$, and Corollary~\ref{cor:quotient-nu-tau-listed} applies.
	
	\emph{Case \textup{(II)}.} Here $\sigma=\nu^{e\beta}$, with $\beta=L_s+1=-s-\frac{n-3}2$, and $\pi'$ is $(H_{1,n-2},\chi_{s-\frac12})$-distinguished. By the induction hypothesis, $\pi'$ is one of the representations (a)--(f) of \S\ref{ss:inductive-input}. We treat these in turn.
	
	(II.a) \emph{$s=\frac12$ and $\pi'=\one_{n-1}=\Ss(A_s^+)$.} Then $\beta=-R_s$, so $\pi'\times\sigma=\xi^{(2)\prime}_s$. This representation is irreducible by Lemma~\ref{lem:dps2-irred}, and therefore $\pi\simeq\xitwo_{1/2}$, which is listed.
	
	(II.b) \emph{$s=1-\frac n2$ and $\pi'=\nu^{e\frac{n-1}2}\Qt_{n-1}=\Ss(\{[L_s],A_s\})$.} Here $L_s=-\frac12$ and $\beta=\frac12$. As in (III.2), $\pi$ is a quotient of $\Ss([-\frac12])\times\Ss(A_s)\times\Ss([\frac12])$. Since $[\frac12]$ is the left end of $A_s=[\frac12,\frac{2n-5}2]$, the last two factors are not linked, and hence commute. Thus, $\pi$ is a quotient of $(\nu^{-e/2}\times\nu^{e/2})\times\nu^{-es}$. By Lemma~\ref{lem:quotient-constituents} and Corollary~\ref{cor:quotient-nu-tau-listed}, $\pi$ is listed.
	
	(II.c) \emph{$s=\frac n2$ and $\pi'=\nu^{-e\frac{n-1}2}Q_{n-1}$.} We have $\pi'=\Ss(\{E,[\frac12]\})$, where $E=\Delta_{n-2}(1-s)=[L_s+2,R_s]=[-\frac{2n-5}2,-\frac12]$. Since $E\prec[\frac12]$, the representation $\pi'$ is the unique irreducible quotient of $\Ss(E)\times\Ss([\frac12])$. Here $\beta=-\frac{2n-3}2$, so $[\beta]\prec E$ and $[\beta]\cup E=A_s^+$. The singletons $[\frac12]$ and $[\beta]$ are not linked. Hence, $\pi$ is a quotient of $(\Ss(E)\times\Ss([\beta]))\times\nu^{e/2}$, and therefore of $X\times\nu^{e/2}$ for some constituent $X\in\{\Ss(\{[\beta],E\}),\Ss(A_s^+)\}$ of $\Ss(E)\times\Ss([\beta])$. We apply Lemma~\ref{lem:quot-X-nu} with $y=\frac12=s-\frac{n-1}2$. Since $y$ differs from $s+\frac{n-1}2$ and from $-s-\frac{n-3}2$, only case (c) can occur. Then $X=\Ss(A_s^+)$, and $\pi$ is a quotient of $\Ss(A_s^+)\times\Ss([-R_s])=\xi^{(2)\prime}_{n/2}$. Hence, $\pi\simeq\nu^{-en/2}Q_n$ by Lemma~\ref{lem:dps2-exceptional}, which is listed.
	
	(II.d) \emph{$s\notin\{1-\frac n2,\frac12\}$ and $\pi'=\nu^{-es}\times\nu^{-e\beta}$.} Then $\pi$ is a quotient of $\nu^{-es}\times(\nu^{-e\beta}\times\nu^{e\beta})$. By Lemma~\ref{lem:quotient-constituents}, $\pi$ is a quotient of $\nu^{-es}\times\tau'$ for some constituent $\tau'$ of $\nu^{-e\beta}\times\nu^{e\beta}$, and Corollary~\ref{cor:quotient-nu-tau-listed} applies.
	
	(II.e) \emph{$s\notin\{\frac12,\frac n2\}$ and $\pi'=\nu^{-eR_s}\times\nu^{e(1-s)}$.} We have $\pi'=\Ss([-R_s])\times\Ss(E)$, where $E=\Delta_{n-2}(1-s)=[L_s+2,R_s]$. Since $[\beta]\prec E$ and $[\beta]\cup E=A_s^+$, the representation $\pi$ is a quotient of $\nu^{-eR_s}\times X$ for some $X\in\{\Ss(\{[\beta],E\}),\Ss(A_s^+)\}$. Twisting Definition~\ref{def:Qm}, we have $\Ss(\{[\beta],E\})=\nu^{e(\frac32-s)}\Qt_{n-1}$. We apply Lemma~\ref{lem:quot-nu-X} with $y=-R_s=s-\frac{n-1}2$. Then case (c) does not occur, case (a) occurs only for $s=\frac{n-2}2$, and cases (b) and (d) may occur.
	\begin{itemize}
		\item In case (b), we have $X=\Ss(A_s^+)$, and $\pi$ is a quotient of $\xitwo_s$. If $s\ne0$, then $\xitwo_s$ is irreducible, and $\pi\simeq\xitwo_s$ is listed. If $s=0$, then $\pi\simeq\nu^e\Qt_n$ by Lemma~\ref{lem:dps2-exceptional}, which is impossible by Lemma~\ref{lem:Qtn-twist-dist}.
		\item Case (d) does not occur, since neither possibility for $X$ is a twist of $Q_{n-1}$ (Lemmas~\ref{lem:Qm-not-char} and~\ref{lem:Q-vs-Qt}).
		\item In case (a), $X$ is $(H_{1,n-2},\chi_{\frac{n-1}2})$-distinguished. This is impossible for the character $\Ss(A_s^+)$, by Lemma~\ref{lem:char-dist}. It is also impossible for $\nu^{e(\frac32-s)}\Qt_{n-1}=\nu^{-e\frac{n-5}2}\Qt_{n-1}$, by \eqref{eq:Q-twist-all} in rank $n-1$, since $\frac{n-1}2\ne-\frac{n-1}2$.
	\end{itemize}
	
	(II.f) \emph{$\pi'=\nu^{e(\frac12-s)}\times\tau$, induced from $G_{n-3}\times G_2$ and irreducible, with $\tau\in\T$.} Put $K=\Delta_{n-3}(\frac12-s)=[L_s+2,R_s-1]$, so that $\nu^{e(\frac12-s)}=\Ss(K)$, $[\beta]\prec K$ and $[\beta]\cup K=A_s$. Since $\pi'$ is irreducible, we have $\pi'\simeq\tau\times\Ss(K)$. Hence, $\pi$ is a quotient of $\tau\times(\Ss(K)\times\Ss([\beta]))$, and therefore, by Lemma~\ref{lem:quotient-constituents}, of $\tau\times Z$ for some $Z\in\{\nu^{-es},\Ss(\{[\beta],K\})\}$. If $Z=\nu^{-es}$, then Corollary~\ref{cor:quotient-nu-tau-listed} applies. It remains to treat the case $Z=\Ss(\{[\beta],K\})$, which, twisting Definition~\ref{def:Qm}, is $\nu^{e(1-s)}\Qt_{n-2}$.
	
	If $n=4$, then $K=[\frac12-s]$, and $Z=\Ss(\{[-s-\frac12],[-s+\frac12]\})=\nu^{-es}\St_2$ by Lemma~\ref{lem:St2}. Hence, Lemma~\ref{lem:n4-special} shows that $\pi$ is not distinguished, which is a contradiction.
	
	If $n\ge5$, then $\tau\times Z$ is distinguished by Lemma~\ref{lem:dist-quotient}, and we apply Proposition~\ref{prop:three-alternatives} with $k=2$, $\pi_1=\tau$ and $\pi_2=Z$. Condition (I) fails because $\tau$ is not a character, and condition (II) fails because $Z$ is not a character (Lemma~\ref{lem:Qm-not-char}). By \eqref{eq:O2}, the elements of the sets $\mathcal O_2(s,n-2;\sigma')$ are characters, twists of $Q_{n-2}$, and irreducible representations induced from $G_1\times G_{n-3}$. Since $n-2\ge3$, the twist $Z$ of $\Qt_{n-2}$ is of none of these forms, by Lemmas~\ref{lem:Qm-not-char} and~\ref{lem:Q-vs-Qt} and Remark~\ref{rem:O-characters}. Hence, condition (III) fails as well, which is a contradiction.
	
	This completes the proof.
\end{proof}

\subsection{Proof of Theorem~\ref{thm:A}(i)}\label{ss:proof-Ai}

We now combine the case $n=3$ with the inductive step. Since the listed representations are exactly those of Theorem~\ref{thm:A}(i), the following theorem proves part (i) of Theorem~\ref{thm:A}.

\begin{theorem}\label{thm:classification-all}
	Let $n\ge3$ and $s\in\R$. An irreducible smooth representation of $G_n$ is $(H_{1,n-1},\chi_s)$-distinguished if and only if it is listed for $(n,s)$.
\end{theorem}

\begin{proof}
	For $n=3$, this is Theorem~\ref{thm:n3}. Let $n\ge4$. By Proposition~\ref{prop:listed-dist}, every representation listed for $(n,s)$ is $(H_{1,n-1},\chi_s)$-distinguished. Conversely, we argue by induction on $n$, the case $n=3$ being Theorem~\ref{thm:n3}. Assume that the assertion holds in rank $n-1$ for every value of the parameter. Then the hypothesis of Theorem~\ref{thm:inductive-step} is satisfied, and that theorem shows that every $(H_{1,n-1},\chi_s)$-distinguished irreducible representation of $G_n$ is listed for $(n,s)$.
\end{proof}

Specializing Theorem~\ref{thm:classification-all} to $s=0$, we obtain the classification of $H_{1,n-1}$-distinguished representations, that is, Corollary~\ref{cor:C}.

\begin{corollary}\label{cor:linear-periods}
	Let $n\ge3$. An irreducible smooth representation $\pi$ of $G_n$ is $H_{1,n-1}$-distinguished if and only if either $\pi\simeq\one_n$, or $\pi\simeq\one_{n-2}\times\tau$, where $\tau\in\T$ and $\one_{n-2}\times\tau$ is irreducible.
\end{corollary}

\begin{proof}
	Since $\chi_0$ is trivial, $\pi$ is $H_{1,n-1}$-distinguished if and only if it is $(H_{1,n-1},\chi_0)$-distinguished. By Theorem~\ref{thm:classification-all}, this holds if and only if $\pi$ is listed for $(n,0)$. For $s=0$, the families (a) and (b) of Definition~\ref{def:listed} contribute only $\one_n$, since $0$ is an exceptional value for both. Moreover, family (c) consists of the irreducible representations $\nu^{0}\times\tau=\one_{n-2}\times\tau$ with $\tau\in\T$. This proves the corollary.
\end{proof}

\section{Multiplicity one}\label{sec:multiplicity}

In this section, we prove part (ii) of Theorem~\ref{thm:A}. Let $n\ge3$, and put $H=H_{1,n-1}$. By Theorem~\ref{thm:classification-all}, we have $\Hom_{H}(\pi,\chi_s)=0$ unless $\pi$ is listed for $(n,s)$. The listed representations are of the following six kinds:
\begin{gather*}
	\xione_s\ (s\notin\{-\tfrac n2,0\}),\qquad\xitwo_s\ (s\notin\{0,\tfrac n2\}),\qquad\one_n\ (s=0),\\
	\nu^{\frac{en}2}\Qt_n\ (s=-\tfrac n2),\qquad\nu^{-\frac{en}2}Q_n\ (s=\tfrac n2),\qquad\nu^{-es}\times\tau\ (\tau\in\T, \nu^{-es}\times\tau\ \text{irreducible}).
\end{gather*}
We prove multiplicity one for each kind separately. Apart from two elementary cases, we realize $\pi$ as a quotient of a representation $\pi_1\times\pi_2$ induced from a maximal parabolic subgroup, for which the three orbit contributions of Theorem~\ref{thm:orbit-bound} add up to at most one. Then Lemma~\ref{lem:dist-quotient} gives $\dim\Hom_H(\pi,\chi_s)\le1$. The two elementary cases are the trivial representation, and $\nu^{-es}\times\tau$ for $n=3$ and $s=0$, where we use Theorem~\ref{thm:ALMSY-mult}. Since $\pi$ is distinguished by Proposition~\ref{prop:listed-dist}, the dimension is exactly one in all cases.

We begin with the irreducible degenerate principal series.

\begin{proposition}\label{prop:mult-xi1}
	Let $n\ge3$ and $s\notin\{-\frac n2,0\}$. Then $\dim\Hom_H(\xione_s,\chi_s)=1$.
\end{proposition}

\begin{proof}
	Take $k=n-1$, $\pi_1=\nu^{-es-\frac e2}\in\Irr(G_{n-1})$ and $\pi_2=\nu^{e(s+\frac{n-1}2)}\in\Irr(G_1)$, so that $\pi_1\times\pi_2=\xione_s$. Then we have the following.
	\begin{itemize}
		\item $\dI=1$. Indeed, $\pi_1=\nu^{\frac{e(n-k-2)}2-es}$, and $\pi_2=\chi_{s+\frac{n-1}2}$ on $H_{1,0}=G_1$.
		\item $\dII=0$. Indeed, since $n-1\ge2$, Lemma~\ref{lem:char-dist} shows that the character $\pi_1$ of $G_{n-1}$ is $(H_{1,n-2},\chi_{s-\frac12})$-distinguished only if $s=\frac12$ and $\pi_1=\one_{n-1}$. The latter condition forces $s=-\frac12$, which is a contradiction.
		\item $\dIII=0$. Indeed, by Corollary~\ref{cor:dIII-kn1}, we have $\dIII=\dim\Hom_{G_{n-1}}(\pi_1,\nu^{-es}\times\nu^{-e(s+\frac{n-1}2)})$. By Lemma~\ref{lem:char-sigma} with $m=n-1$, $y=-es$ and $\sigma=\nu^{y-e\frac{n-1}2}$, the unique irreducible subrepresentation of $\nu^{-es}\times\nu^{-e(s+\frac{n-1}2)}$ is $\nu^{-es+\frac e2}\Qt_{n-1}$, which is not a character by Lemma~\ref{lem:Qm-not-char}. Hence, $\dIII=0$ by Lemma~\ref{lem:socle-hom}.
	\end{itemize}
	Therefore, $\dim\Hom_H(\xione_s,\chi_s)\le1$ by Theorem~\ref{thm:orbit-bound}.
\end{proof}

The argument for $\xitwo_s$ is symmetric.

\begin{proposition}\label{prop:mult-xi2}
	Let $n\ge3$ and $s\notin\{0,\frac n2\}$. Then $\dim\Hom_H(\xitwo_s,\chi_s)=1$.
\end{proposition}

\begin{proof}
	Take $k=1$, $\pi_1=\nu^{e(s-\frac{n-1}2)}\in\Irr(G_1)$ and $\pi_2=\nu^{-es+\frac e2}\in\Irr(G_{n-1})$, so that $\pi_1\times\pi_2=\xitwo_s$. Then we have the following.
	\begin{itemize}
		\item $\dII=1$. Indeed, $\pi_1=\chi_{s-\frac{n-1}2}$ on $H_{1,0}=G_1$, and $\pi_2=\nu^{\frac{e(2-k)}2-es}$.
		\item $\dI=0$. Indeed, by Lemma~\ref{lem:char-dist}, the character $\pi_2$ of $G_{n-1}$ is $(H_{1,n-2},\chi_{s+\frac12})$-distinguished only if $s=-\frac12$ and $\pi_2=\one_{n-1}$. The latter condition forces $s=\frac12$, which is a contradiction.
		\item $\dIII=0$. Indeed, by Corollary~\ref{cor:dIII-k1}, we have $\dIII=\dim\Hom_{G_{n-1}}(\pi_2,\nu^{-es+e\frac{n-1}2}\times\nu^{-es})$. By Lemma~\ref{lem:sigma-char} with $m=n-1$, $y=-es$ and $\sigma=\nu^{y+e\frac{n-1}2}$, the unique irreducible subrepresentation of this representation is $\nu^{-es-\frac e2}Q_{n-1}$, which is not a character. Hence, $\dIII=0$ by Lemma~\ref{lem:socle-hom}.
	\end{itemize}
	Therefore, Theorem~\ref{thm:orbit-bound} gives the assertion.
\end{proof}

The trivial representation is immediate.

\begin{proposition}\label{prop:mult-trivial}
	Let $n\ge3$. Then $\dim\Hom_H(\one_n,\chi_0)=1$.
\end{proposition}

\begin{proof}
	This is clear, since $\one_n$ is one-dimensional and $\chi_0$ is trivial.
\end{proof}

Next, we treat the exceptional representations at $s=\mp\frac n2$.

\begin{proposition}\label{prop:mult-Qtn}
	Let $n\ge3$ and $s=-\frac n2$. Then $\dim\Hom_H(\nu^{en/2}\Qt_n,\chi_s)=1$.
\end{proposition}

\begin{proof}
	By Lemma~\ref{lem:dps1-exceptional}, $\nu^{en/2}\Qt_n$ is a quotient of $\xi^{(1)\prime}_{-n/2}=\pi_1\times\pi_2$, where $k=1$, $\pi_1=\nu^{-e/2}\in\Irr(G_1)$ and $\pi_2=\nu^{e\frac{n-1}2}\in\Irr(G_{n-1})$. With $s=-\frac n2$, we have the following.
	\begin{itemize}
		\item $\dI=0$, since $\pi_1\not\simeq\nu^{\frac{e(n-3)}2-es}=\nu^{e\frac{2n-3}2}$.
		\item $\dII=0$, since $\pi_1\not\simeq\chi_{s-\frac{n-1}2}=\nu^{-e\frac{2n-1}2}$.
		\item $\dIII=1$. Indeed, by Corollary~\ref{cor:dIII-k1}, we have $\dIII=\dim\Hom_{G_{n-1}}(\pi_2,\nu^{e/2}\times\nu^{en/2})$. By Lemma~\ref{lem:sigma-char} with $m=n-1$, $y=\frac{en}2$ and $\sigma=\nu^{e/2}=\nu^{y-e\frac{n-1}2}$, the unique irreducible subrepresentation of $\nu^{e/2}\times\nu^{en/2}$ is $\nu^{y-\frac e2}=\nu^{e\frac{n-1}2}=\pi_2$. Hence, $\dIII=1$ by Lemma~\ref{lem:socle-hom}.
	\end{itemize}
	Therefore, by Theorem~\ref{thm:orbit-bound} and Lemma~\ref{lem:dist-quotient}, we obtain $\dim\Hom_H(\nu^{en/2}\Qt_n,\chi_s)\le1$.
\end{proof}

\begin{proposition}\label{prop:mult-Qn}
	Let $n\ge3$ and $s=\frac n2$. Then $\dim\Hom_H(\nu^{-en/2}Q_n,\chi_s)=1$.
\end{proposition}

\begin{proof}
	By Lemma~\ref{lem:dps2-exceptional}, $\nu^{-en/2}Q_n$ is a quotient of $\xi^{(2)\prime}_{n/2}=\pi_1\times\pi_2$, where $k=n-1$, $\pi_1=\nu^{-e\frac{n-1}2}\in\Irr(G_{n-1})$ and $\pi_2=\nu^{e/2}\in\Irr(G_1)$. With $s=\frac n2$, we have the following.
	\begin{itemize}
		\item $\dI=0$, since $\pi_1\not\simeq\nu^{-\frac e2-es}=\nu^{-e\frac{n+1}2}$.
		\item $\dII=0$, by Lemma~\ref{lem:char-dist} applied to the character $\pi_1$ of $G_{n-1}$, since $s-\frac12\ne0$.
		\item $\dIII=1$. Indeed, by Corollary~\ref{cor:dIII-kn1}, we have $\dIII=\dim\Hom_{G_{n-1}}(\pi_1,\nu^{-en/2}\times\nu^{-e/2})$. By Lemma~\ref{lem:char-sigma} with $m=n-1$, $y=-\frac{en}2$ and $\sigma=\nu^{-e/2}=\nu^{y+e\frac{n-1}2}$, the unique irreducible subrepresentation of $\nu^{-en/2}\times\nu^{-e/2}$ is $\nu^{y+\frac e2}=\nu^{-e\frac{n-1}2}=\pi_1$. Hence, $\dIII=1$ by Lemma~\ref{lem:socle-hom}.
	\end{itemize}
	Therefore, Theorem~\ref{thm:orbit-bound} and Lemma~\ref{lem:dist-quotient} give the assertion.
\end{proof}

Finally, we treat the representations $\nu^{-es}\times\tau$, first for $n\ge4$ and then for $n=3$.

\begin{proposition}\label{prop:mult-c-n4}
	Let $n\ge4$, $s\in\R$ and $\tau\in\T$, and suppose that $\nu^{-es}\times\tau$ is irreducible. Then $\dim\Hom_H(\nu^{-es}\times\tau,\chi_s)=1$.
\end{proposition}

\begin{proof}
	Take $k=n-2\ge2$, $\pi_1=\nu^{-es}\in\Irr(G_{n-2})$ and $\pi_2=\tau$. Then we have the following.
	\begin{itemize}
		\item $\dI\le1$. Indeed, the first factor of $\dI$ is $\dim\Hom_{G_{n-2}}(\nu^{-es},\nu^{-es})=1$, and the second is $\dim\Hom_{H_{1,1}}(\tau,\chi_{s+\frac{n-2}2})\le1$ by Proposition~\ref{prop:rank2-mult}.
		\item $\dII=0$, since $\tau$ is irreducible and infinite-dimensional, and hence has no character as a quotient.
		\item $\dIII=0$. Indeed, in the notation of \S\ref{ss:open-orbit}, we have $A=G_{n-3}$ with $n-3\ge1$, and $J_1=r_{(n-3,1)}(\nu^{-es})=\nu^{-es-\frac e2}\otimes\nu^{-es+e\frac{n-3}2}$ by Lemma~\ref{lem:jacquet-char}. Since $A$ acts on $J_1$ by the character $\nu^{-es-\frac e2}$, which is different from $\nu^{x_s}=\nu^{\frac e2-es}$, the space $U$ of Proposition~\ref{prop:U-T} is zero. Hence, $\dIII=\dim\Hom_B(U,\widetilde T)=0$ by Proposition~\ref{prop:open-condition}.
	\end{itemize}
	Therefore, Theorem~\ref{thm:orbit-bound} gives the assertion.
\end{proof}

\begin{proposition}\label{prop:mult-c-n3}
	Let $s\in\R$ and $\tau\in\T$, and suppose that the representation $\nu^{-es}\times\tau$ induced from $G_1\times G_2$ is irreducible. Then $\dim\Hom_{H_{1,2}}(\nu^{-es}\times\tau,\chi_s)=1$.
\end{proposition}

\begin{proof}
	Put $\pi=\nu^{-es}\times\tau$. If $s=0$, then $\chi_0$ is trivial, and $\dim\Hom_{H_{1,2}}(\pi,\one)\le1$ by Theorem~\ref{thm:ALMSY-mult} with $(p,q)=(1,2)$.
	
	Let $s\ne0$. We first apply Lemma~\ref{lem:n3-bound-k1} with $\pi_1=\nu^{-es}$ and $\pi_2=\tau$. We have $\dI\le1$ by Proposition~\ref{prop:rank2-mult}, and $\dII=0$, since $\tau$ is not a character. Moreover, $\dIII=\dim\Hom_{G_2}(\tau,\nu^{es}\times\nu^{-es})$, which, by Lemma~\ref{lem:socle-hom}, vanishes unless $\tau$ is the unique irreducible subrepresentation of $\nu^{es}\times\nu^{-es}$. Since $\tau\in\T$, by Lemma~\ref{lem:St2} this happens only if either $\tau=\tau_s$ with $s\ne\pm\frac12$, or $\tau=\St_2$ with $s=\frac12$. Outside these two cases, we obtain $\dim\Hom_{H_{1,2}}(\pi,\chi_s)\le\dI\le1$.
	
	Suppose that $\tau=\St_2$ and $s=\frac12$. Then $\pi\simeq\St_2\times\nu^{-e/2}$, and we apply Lemma~\ref{lem:n3-bound-k2} with $\pi_1=\St_2$ and $\pi_2=\nu^{-e/2}$. We have $\dI=0$, since $\St_2$ is not a character, and $\dII\le1$ by Proposition~\ref{prop:rank2-mult}. Moreover, $\dIII=\dim\Hom_{G_2}(\St_2,\nu^{-e/2}\times\nu^{e/2})=0$, since the unique irreducible subrepresentation of $\nu^{-e/2}\times\nu^{e/2}$ is $\one_2$.
	
	Suppose that $\tau=\tau_s$ with $s\notin\{0,\pm\frac12\}$. Then $\pi=\nu^{-es}\times\nu^{es}\times\nu^{-es}$, and the segments $[-s]$, $[s]$, $[-s]$ are pairwise unlinked. Hence, by Lemma~\ref{lem:unlinked-product}, $\pi\simeq\nu^{es}\times\tau''$, where $\tau''=\nu^{-es}\times\nu^{-es}$ is irreducible. We apply Lemma~\ref{lem:n3-bound-k1} with $\pi_1=\nu^{es}$ and $\pi_2=\tau''$. We have $\dI=0$, since $\nu^{es}\ne\nu^{-es}$, and $\dII=0$, since $\nu^{es}\ne\nu^{e(s-1)}$. Moreover, $\dIII=\dim\Hom_{G_2}(\tau'',\nu^{-es}\times\nu^{-es})=1$.
	
	In all cases, $\dim\Hom_{H_{1,2}}(\pi,\chi_s)\le1$, and equality holds by Proposition~\ref{prop:listed-dist}.
\end{proof}

We can now state the multiplicity one theorem in all ranks.

\begin{theorem}\label{thm:multiplicity}
	Let $n\ge1$, $s\in\R$ and $\pi\in\Irr(G_n)$. Then $\dim\Hom_{H_{1,n-1}}(\pi,\chi_s)\le1$. Moreover, for $n\ge3$,
	$$
	\dim\Hom_{H_{1,n-1}}(\pi,\chi_s)=\begin{cases}1&\text{if }\pi\text{ is listed for }(n,s),\\0&\text{otherwise.}\end{cases}$$
\end{theorem}

\begin{proof}
	For $n=1$, this is Proposition~\ref{prop:n1}, and for $n=2$, it is Theorem~\ref{thm:main-n2}(2). For $n\ge3$, Theorem~\ref{thm:classification-all} shows that the space vanishes unless $\pi$ is listed for $(n,s)$. The listed representations are treated in Propositions~\ref{prop:mult-xi1}--\ref{prop:mult-c-n3}.
\end{proof}

\section{Langlands parameters}\label{sec:langlands}

In this section, we prove part (iii) of Theorem~\ref{thm:A}. Throughout, we assume that $n\ge3$. We use the Langlands parameters of Definition~\ref{def:L-param} and the notation of \S\ref{ss:dps}. Let $\aaa_s$ denote the multisegment consisting of the singletons $[y]$ with $y\in A_s$.

We first compute the Langlands parameter of the character $\nu^{-es}$ of $G_{n-2}$.

\begin{lemma}\label{lem:L-nu-es}
	Let $n\ge3$ and $s\in\R$. As a character of $G_{n-2}$, we have $\nu^{-es}=\Qq(\aaa_s)$ and
	$$\varphi_{\nu^{-es}}\simeq\bigoplus_{j=0}^{n-3}|\cdot|^{e(-s+\frac{n-3}2-j)}\otimes\Sp_e,
	$$
	which is a representation of dimension $e(n-2)$.
\end{lemma}

\begin{proof}
	Since $A_s=\Delta_{n-2}(-s)$, this follows from Lemmas~\ref{lem:char-Langlands}, \ref{lem:L-additive} and~\ref{lem:L-segment}.
\end{proof}

Next, we identify the parameters attached to a pair of opposite singletons.

\begin{lemma}\label{lem:L-two-singletons}
	Let $x\in\R$. The parameter $|\cdot|^{ex}\otimes\Sp_e\oplus|\cdot|^{-ex}\otimes\Sp_e$ is equal to $\varphi_{\tau_x}$ if $x\ne\pm\frac12$, and to $\varphi_{\one_2}$ if $x=\pm\frac12$. Consequently, a representation of $G_n$ whose Langlands multisegment is $\aaa_s+\{[x],[-x]\}$ has Langlands parameter $\varphi_{\nu^{-es}}\oplus\varphi_{\tau_x}$ if $x\ne\pm\frac12$, and $\varphi_{\nu^{-es}}\oplus\varphi_{\one_2}$ if $x=\pm\frac12$.
\end{lemma}

\begin{proof}
	For $x\ne\pm\frac12$, the representation $\tau_x=\Qq([x])\times\Qq([-x])$ is irreducible, so $\tau_x=\Qq(\{[x],[-x]\})$ by Lemma~\ref{lem:sum-multisegments}. Moreover, $\one_2=\Qq(\{[\frac12],[-\frac12]\})$ by Lemma~\ref{lem:char-Langlands}. The assertions now follow from Lemmas~\ref{lem:L-additive} and~\ref{lem:L-segment}.
\end{proof}

We now determine the Langlands multisegments of all listed representations.

\begin{lemma}\label{lem:Langlands-listed}
	Let $n\ge3$ and $s\in\R$. The Langlands multisegments of the listed representations are as follows.
	\begin{enumerate}
		\item If $s\notin\{-\frac n2,0\}$, then $\xione_s=\Qq(\aaa_s+\{[L_s],[-L_s]\})$.
		\item If $s\notin\{0,\frac n2\}$, then $\xitwo_s=\Qq(\aaa_s+\{[R_s],[-R_s]\})$.
		\item $\one_n=\Qq(\aaa_0+\{[L_0],[R_0]\})$, where $R_0=-L_0=\frac{n-1}2$.
		\item $\nu^{en/2}\Qt_n=\Qq(\aaa_{-n/2}+\{[-\frac12,\frac12]\})$.
		\item $\nu^{-en/2}Q_n=\Qq(\aaa_{n/2}+\{[-\frac12,\frac12]\})$.
		\item If $\tau=\Qq(\bb)\in\T$ and $\nu^{-es}\times\tau$ is irreducible, then $\nu^{-es}\times\tau=\Qq(\aaa_s+\bb)$.
	\end{enumerate}
\end{lemma}

\begin{proof}
	(1) Since $A_s^-=[L_s]\cup A_s$, the character $\nu^{-es-\frac e2}=\Ss(A_s^-)$ of $G_{n-1}$ is equal to $\Qq(\aaa_s+\{[L_s]\})$ by Lemma~\ref{lem:char-Langlands}. Moreover, $\nu^{-eL_s}=\Qq([-L_s])$. Since $\xione_s$ is irreducible, the assertion follows from Lemma~\ref{lem:sum-multisegments}.
	
	(2) The proof is the same as for (1), using $A_s^+=A_s\cup[R_s]$.
	
	(3) Since $\Delta_n=[L_0]\cup A_0\cup[R_0]$, this follows from Lemma~\ref{lem:char-Langlands}.
	
	(4) and (5) These follow by twisting Lemma~\ref{lem:Qm-Langlands} by $\nu^{\pm en/2}$. For instance, twisting by $\nu^{en/2}$, the singletons $[j]$ with $|j|\le\frac{n-3}2$ in the Langlands multisegment of $\Qt_n$ become the singletons of $A_{-n/2}=[\frac32,\frac{2n-3}2]$, and the segment $[-\frac{n+1}2,-\frac{n-1}2]$ becomes $[-\frac12,\frac12]$. The case of $Q_n$ is similar, with $A_{n/2}=[-\frac{2n-3}2,-\frac32]$.
	
	(6) This follows from Lemmas~\ref{lem:L-nu-es} and~\ref{lem:sum-multisegments}.
\end{proof}

We can now prove part (iii) of Theorem~\ref{thm:A}.

\begin{theorem}\label{thm:langlands}
	Let $n\ge3$ and $s\in\R$. An irreducible smooth representation $\pi$ of $G_n$ is $(H_{1,n-1},\chi_s)$-distinguished if and only if
	$\varphi_{\pi}\simeq\varphi_{\nu^{-es}}\oplus\varphi_{\tau_0}$
	for some $\tau_0\in\Irr(G_2)$ satisfying one of the following conditions:
	\begin{itemize}
		\item $\tau_0=\one_2$ and $s=\pm\frac{n-2}2$;
		\item $\tau_0\in\T$.
	\end{itemize}
	Here, $\nu^{-es}$ is regarded as a character of $G_{n-2}$, so that $\varphi_{\nu^{-es}}$ has dimension $e(n-2)$ and $\varphi_{\tau_0}$ has dimension $2e$.
\end{theorem}

\begin{proof}
	Since the representations of $W_{\mathrm{F}}\times\mathrm{SL}_2(\C)$ under consideration are semisimple, the condition of the theorem means that
	\begin{equation}\label{eq:L-condition}
		\varphi_{\pi}\simeq\varphi_{\nu^{-es}}\oplus\varphi_{\tau_0},\qquad\text{where either }\tau_0\in\T,\ \text{or }\tau_0=\one_2\text{ and }s=\pm\tfrac{n-2}2.
	\end{equation}
	
	\emph{Distinguished representations satisfy \eqref{eq:L-condition}.} Let $\pi$ be distinguished. By Theorem~\ref{thm:classification-all}, $\pi$ is listed for $(n,s)$, and we use Lemma~\ref{lem:Langlands-listed}.
	\begin{itemize}
		\item If $\pi=\nu^{-es}\times\tau$ is as in (c), then $\varphi_{\pi}=\varphi_{\nu^{-es}}\oplus\varphi_{\tau}$ by Lemma~\ref{lem:Langlands-listed}(6) and Lemma~\ref{lem:L-additive}.
		\item Suppose that $\pi=\xione_s$ with $s\notin\{-\frac n2,0\}$. Then the Langlands multisegment of $\pi$ is $\aaa_s+\{[L_s],[-L_s]\}$. Note that $L_s=\frac12$ would mean $s=-\frac n2$, which is excluded, and that $L_s=-\frac12$ means $s=-\frac{n-2}2$. Hence, by Lemma~\ref{lem:L-two-singletons}, the quotient $\varphi_{\pi}/\varphi_{\nu^{-es}}$ is $\varphi_{\one_2}$ if $s=-\frac{n-2}2$, and $\varphi_{\tau_{L_s}}$ with $\tau_{L_s}\in\T$ (Corollary~\ref{cor:T-members}) otherwise.
		\item The case $\pi=\xitwo_s$, with $s\notin\{0,\frac n2\}$, is the same, with $R_s$ in place of $L_s$. In this case, the quotient is $\varphi_{\one_2}$ exactly when $s=\frac{n-2}2$.
		\item If $s=0$ and $\pi=\one_n$, then the quotient is $\varphi_{\tau_{(n-1)/2}}$, and $\tau_{(n-1)/2}\in\T$ since $\frac{n-1}2\ge1$.
		\item If $s=-\frac n2$ and $\pi=\nu^{en/2}\Qt_n$, or $s=\frac n2$ and $\pi=\nu^{-en/2}Q_n$, then the quotient is $\varphi_{\Qq([-\frac12,\frac12])}=\varphi_{\St_2}$, and $\St_2\in\T$.
	\end{itemize}
	
	\emph{Representations satisfying \eqref{eq:L-condition} are distinguished.} Write $\tau_0=\Qq(\bb)$. By Lemmas~\ref{lem:L-nu-es}, \ref{lem:L-additive} and~\ref{lem:L-injective}, we have $\pi\simeq\Qq(\aaa_s+\bb)$. We show that $\pi$ is listed for $(n,s)$; it is then distinguished by Proposition~\ref{prop:listed-dist}.
	\begin{itemize}
		\item Suppose that $\tau_0=\one_2$ and $s=-\frac{n-2}2$. Then $L_s=-\frac12$, and $\aaa_s+\bb=\aaa_s+\{[L_s],[-L_s]\}$ is the Langlands multisegment of $\xione_s$, which is irreducible since $s\notin\{-\frac n2,0\}$. Hence, $\pi\simeq\xione_s$. Similarly, if $s=\frac{n-2}2$, then $R_s=\frac12$ and $\pi\simeq\xitwo_s$.
		\item Suppose that $\tau_0\in\T$ and $\nu^{-es}\times\tau_0$ is irreducible. Then $\pi\simeq\nu^{-es}\times\tau_0$ by Lemma~\ref{lem:Langlands-listed}(6), which is listed in (c).
		\item Suppose that $\tau_0\in\T$ and $\nu^{-es}\times\tau_0$ is reducible. Then Lemma~\ref{lem:reducible-nu-tau} applies. If $\tau_0=\St_2$ and $s=\pm\frac n2$, then $\bb=\{[-\frac12,\frac12]\}$, and $\pi$ is $\nu^{-en/2}Q_n$ or $\nu^{en/2}\Qt_n$, respectively, by Lemma~\ref{lem:Langlands-listed}(4) and~(5). If $\tau_0=\tau_x$ with $x\in\{L_s,R_s\}$ and $x\ne\pm\frac12$, then $\bb=\{[x],[-x]\}$. For $x=L_s$, the multisegment $\aaa_s+\bb$ is that of $\one_n$ if $s=0$, and otherwise that of the irreducible representation $\xione_s$; note that $s\ne-\frac n2$, since $L_s\ne\frac12$. For $x=R_s$, it is that of $\one_n$ if $s=0$, and otherwise that of the irreducible representation $\xitwo_s$; note that $s\ne\frac n2$, since $R_s\ne-\frac12$.
	\end{itemize}
	In all cases, $\pi$ is listed for $(n,s)$, which completes the proof.
\end{proof}

\subsection{Proof of Theorem~\ref{thm:A}}\label{ss:proof-A}

We can now complete the proof of our main theorem. Let $n\ge3$. Part (i) of Theorem~\ref{thm:A} is Theorem~\ref{thm:classification-all}, together with the description of the listed representations following Definition~\ref{def:listed}. Part (ii) is Theorem~\ref{thm:multiplicity}, and part (iii) is Theorem~\ref{thm:langlands}. \qed


\begin{thebibliography}{ALMSY24}
	
	\bibitem[AGRS10]{AGRS10}
	A. Aizenbud, D. Gourevitch, S. Rallis, and G. Schiffmann, \emph{Multiplicity one theorems}, Ann. of Math. (2) \textbf{172} (2010), no.~2, 1407--1434.
	
	\bibitem[ALMSY24]{ALMSY24}
	U.~K.~Anandavardhanan, H.~Lu, N.~Matringe, V.~S\'echerre, and C.~Yang,
	\emph{The sign of linear periods},
	preprint (2024), \href{https://arxiv.org/abs/2402.12106}{arXiv:2402.12106}.
	
	\bibitem[BR07]{BR07}
	A. I. Badulescu and D. Renard, \emph{Zelevinsky involution and Moeglin--Waldspurger algorithm for $\mathrm{GL}_n(D)$}, Functional analysis IX, Various Publ. Ser. (Aarhus), vol.~48, Univ. Aarhus, Aarhus, 2007, pp.~9--15.
	
	\bibitem[BZ76]{BZ76}
	I. N. Bernstein and A. V. Zelevinsky, \emph{Representations of the group $\mathrm{GL}(n,F)$, where $F$ is a local non-Archimedean field}, Uspehi Mat. Nauk \textbf{31} (1976), no.~3(189), 5--70; English transl., Russian Math. Surveys \textbf{31} (1976), no.~3, 1--68.
	
	\bibitem[BZ77]{BZ77}
	I. N. Bernstein and A. V. Zelevinsky, \emph{Induced representations of reductive $\mathfrak p$-adic groups. I}, Ann. Sci. \'Ecole Norm. Sup. (4) \textbf{10} (1977), no.~4, 441--472.
	
	\bibitem[CH26]{CH26}
	Q.~Chen and J.~Huang,
	\emph{Classification of $\mathrm{GL}_1(\mathrm{D})\times\mathrm{GL}_{n-1}(\mathrm{D})$-distinguished representations},
	preprint (2026), \href{https://arxiv.org/abs/2611.19149}{arXiv:2611.19149}.
	
	\bibitem[CS20]{CS20}
	F. Chen and B. Sun, \emph{Uniqueness of twisted linear periods and twisted Shalika periods}, Sci. China Math. \textbf{63} (2020), no.~1, 1--22.
	
\bibitem[DS26]{DS26}
P.~Dagar and H.~Sharma,
\emph{On representations of $\mathrm{GL}_n(\mathrm{D})$ admitting a generalized linear period},
preprint (2026), \href{https://arxiv.org/abs/2609.34307}{arXiv:2609.34307}.
	
	\bibitem[DSV27]{DSV27}
	P. Dagar, H. Sharma, and M. K. Verma, \emph{On representations of $\mathrm{GL}_n$ distinguished by $\mathrm{GL}_1\times\mathrm{GL}_{n-1}$ over a quaternion division algebra}, J. Algebra \textbf{711} (2027), 405--420.
	
	\bibitem[DKV84]{DKV84}
	P. Deligne, D. Kazhdan, and M.-F. Vign\'eras, \emph{Repr\'esentations des alg\`ebres centrales simples $p$-adiques}, Repr\'esentations des groupes r\'eductifs sur un corps local, Travaux en Cours, Hermann, Paris, 1984, pp.~33--117.
	
	\bibitem[Dra83]{Dra83}
	P. K. Draxl, \emph{Skew fields}, London Mathematical Society Lecture Note Series, vol.~81, Cambridge University Press, Cambridge, 1983.
	
	\bibitem[FJ93]{FJ93}
	S. Friedberg and H. Jacquet, \emph{Linear periods}, J. Reine Angew. Math. \textbf{443} (1993), 91--139.
	
	\bibitem[GGP23]{GGP23}
	W. T. Gan, B. H. Gross, and D. Prasad, \emph{Twisted GGP problems and conjectures}, Compositio Math. \textbf{159} (2023), no.~9, 1916--1973.
	
	\bibitem[GT10]{GT10}
	W. T. Gan and S. Takeda, \emph{On Shalika periods and a theorem of Jacquet--Martin}, Amer. J. Math. \textbf{132} (2010), no.~2, 475--528.
	
	\bibitem[HT01]{HT01}
	M. Harris and R. Taylor, \emph{The geometry and cohomology of some simple Shimura varieties}, Annals of Mathematics Studies, vol.~151, Princeton University Press, Princeton, NJ, 2001.
	
	\bibitem[Hen00]{Hen00}
	G. Henniart, \emph{Une preuve simple des conjectures de Langlands pour $\mathrm{GL}(n)$ sur un corps $p$-adique}, Invent. Math. \textbf{139} (2000), no.~2, 439--455.
	
	\bibitem[JR96]{JR96}
	H. Jacquet and S. Rallis, \emph{Uniqueness of linear periods}, Compositio Math. \textbf{102} (1996), no.~1, 65--123.
	
	\bibitem[LM16]{LM16}
	E. Lapid and A. M\'inguez, \emph{On parabolic induction on inner forms of the general linear group over a non-archimedean local field}, Selecta Math. (N.S.) \textbf{22} (2016), no.~4, 2347--2400.
	
	\bibitem[LLX26]{LLX26}
	Y.~Li, H.~Lu, and J.~Xiao,
	\emph{Uniqueness of twisted linear periods},
	preprint (2026), \href{https://arxiv.org/abs/2610.26983}{arXiv:2610.26983}.
	
	\bibitem[Lu22]{Lu22}
	H. Lu, \emph{The generalized linear periods}, J. Number Theory \textbf{233} (2022), 261--284.
	
	\bibitem[Lu26]{Lu26}
	H. Lu, \emph{The generalized linear period and the Shalika period over a division algebra}, J. Number Theory \textbf{279} (2026), 238--255.
	
	\bibitem[MW86]{MW86}
	C. M\oe glin and J.-L. Waldspurger, \emph{Sur l'involution de Zelevinski}, J. Reine Angew. Math. \textbf{372} (1986), 136--177.
	
	\bibitem[NS25]{NS25}
	S. Nadimpalli and M. Sheth, \emph{Twisted Jacquet modules: a conjecture of D. Prasad}, Proc. Amer. Math. Soc. \textbf{153} (2025), no.~3, 1349--1361.
	
	\bibitem[Off17]{Off17}
	O. Offen, \emph{On parabolic induction associated with a $p$-adic symmetric space}, J. Number Theory \textbf{170} (2017), 211--227.
	
	\bibitem[PR94]{PR94}
	V. Platonov and A. Rapinchuk, \emph{Algebraic groups and number theory}, Pure and Applied Mathematics, vol.~139, Academic Press, Boston, MA, 1994.
	
	\bibitem[Pra93]{Pra93}
	D. Prasad, \emph{On the decomposition of a representation of $\mathrm{GL}(3)$ restricted to $\mathrm{GL}(2)$ over a $p$-adic field}, Duke Math. J. \textbf{69} (1993), no.~1, 167--177.
	
	\bibitem[Pra99]{Pra99}
	D. Prasad, \emph{Some remarks on representations of a division algebra and of the Galois group of a local field}, J. Number Theory \textbf{74} (1999), no.~1, 73--97.
	
	\bibitem[Pra01]{Pra01}
	D. Prasad, \emph{The space of degenerate Whittaker models for $\mathrm{GL}(4)$ over $p$-adic fields}, Cohomology of arithmetic groups, $L$-functions and automorphic forms (Mumbai, 1998/1999), Tata Inst. Fund. Res. Stud. Math., vol.~15, Tata Inst. Fund. Res., Bombay, 2001, pp.~103--115.
	
	\bibitem[PR00]{PR00}
	D. Prasad and A. Raghuram, \emph{Kirillov theory for $\mathrm{GL}_2(\mathcal{D})$ where $\mathcal{D}$ is a division algebra over a non-Archimedean local field}, Duke Math. J. \textbf{104} (2000), no.~1, 19--44.
	
	\bibitem[Rag99]{Rag99}
	A. Raghuram, \emph{Some topics in algebraic groups: Representation theory of $\mathrm{GL}_2(\mathcal{D})$ where $\mathcal{D}$ is a division algebra over a non-Archimedean local field}, Ph.D. thesis, Tata Institute of Fundamental Research, University of Mumbai, 1999.
	
	\bibitem[Rag07]{Rag07}
	A. Raghuram, \emph{On the restriction to $\mathcal{D}^*\times\mathcal{D}^*$ of representations of $p$-adic $\mathrm{GL}_2(\mathcal{D})$}, Canad. J. Math. \textbf{59} (2007), no.~5, 1050--1068.
	
	\bibitem[Sch13]{Sch13}
	P. Scholze, \emph{The local Langlands correspondence for $\mathrm{GL}_n$ over $p$-adic fields}, Invent. Math. \textbf{192} (2013), no.~3, 663--715.
	
	\bibitem[Smi17]{Smi17}
	J. M. Smith, \emph{Construction of relative discrete series representations for $p$-adic $\mathrm{GL}_n$}, Ph.D. thesis, University of Toronto, 2017.
	
	\bibitem[Tad90]{Tad90}
	M. Tadi\'c, \emph{Induced representations of $\mathrm{GL}(n,A)$ for $p$-adic division algebras $A$}, J. Reine Angew. Math. \textbf{405} (1990), 48--77.
	
	\bibitem[Ven13]{Ven13}
	C. G. Venketasubramanian, \emph{On representations of $\mathrm{GL}(n)$ distinguished by $\mathrm{GL}(n-1)$ over a $p$-adic field}, Israel J. Math. \textbf{194} (2013), no.~1, 1--44.
	
	\bibitem[Wei74]{Wei74}
	A. Weil, \emph{Basic number theory}, 3rd ed., Grundlehren der Mathematischen Wissenschaften, vol.~144, Springer-Verlag, New York, 1974.
	
	\bibitem[Zel80]{Zel80}
	A. V. Zelevinsky, \emph{Induced representations of reductive $\mathfrak p$-adic groups. II. On irreducible representations of $\mathrm{GL}(n)$}, Ann. Sci. \'Ecole Norm. Sup. (4) \textbf{13} (1980), no.~2, 165--210.

%
%
%
%
%
%
%
%
%
%
%
%
%
%
%
%
%
%
%
%
%
%
%
%
%
%
%
%
%
%
%
%
%
%
%
%
%
%
%
%

\end{thebibliography}
\end{document}